\documentstyle[twoside,11pt]{report}
\def\qed{\relax\ifmmode\hskip2em \Box\else\unskip\nobreak\hskip1em $\Box$\fi}

\newtheorem{thm}{Theorem}
\newtheorem{cor}[thm]{Corollary}

\newtheorem{lem}[thm]{Lemma}
\newtheorem{prop}[thm]{Proposition}
\newtheorem{defn}[thm]{Definition}
\newtheorem{notat}[thm]{Notation}
\newtheorem{exmp}[thm]{Example}
\newtheorem{rem}[thm]{Remark}
\input amssym.def
\input amssym.tex

\newcommand{\prodi}{\frak Y}

\newcommand{\proof}{\noindent{\bf Proof.}}
\newcommand{\nsem}{$(2-\epsilon)$}
\newcommand{\BB}{\Bbb B}
\newcommand{\beq}{\begin{equation}}
\newcommand{\boldC}{{\bf C}}
\newcommand{\boldf}{{\bf f}}
\newcommand{\boldff}{{\bf F}}
\newcommand{\boldg}{{\bf g}}
\newcommand{\boldh}{{\bf h}}
\newcommand{\boldy}{{\bf y}}
\newcommand{\cadlag}{c\`adl\`ag}

\newcommand{\CAL}{{\cal L}}
\newcommand{\CLW}{{\cal W}}

\newcommand{\dual}{\mbox{dual}}
\newcommand{\eeq}{\end{equation}}

\newcommand{\esssup}{\mbox{esssup}}

\newcommand{\FF}{\Bbb F}
\newcommand{\fgotik}{\frak{F}}

\newcommand{\g}{\phi}
\newcommand{\h}{\psi}
\newcommand{\inv}{\mbox{\tiny{inv}}}
\newcommand{\ind}{\gamma}
\newcommand{\loc}{\lbrack\!\lbrack}
\newcommand{\lei}{\loc}

\newcommand{\lmuv}{\lambda_m\Cap\lei u,v\rei}
\newcommand{\lsi}{\circleddash}
\newcommand{\kuv}{\kappa\Cap\lei u,v\rei}
\newcommand{\MM}{\Bbb M}

\newcommand{\NN}{\Bbb N}
\newcommand{\osc}{\mbox{Osc}}

\newcommand{\prtn}{\frak{P}}
\newcommand{\pindex}{\upsilon}
\newcommand{\pv}{A}
\newcommand{\PP}{\frak{T}}
\newcommand{\RR}{\Bbb R}
\newcommand{\roc}{\rbrack\!\rbrack}
\newcommand{\rei}{\roc}

\newcommand{\SYS}{SY\!S}
\newcommand{\ugotik}{\frak{U}}
\newcommand{\unc}{\xi}
\newcommand{\unit}{1\!\!\Bbb I}
\newcommand{\RYS}{RY\!S}
\newcommand{\ZZ}{\Bbb Z}
\newcommand{\WW}{{\cal W}}
\newcommand{\YS}{Y\!S}

\title{QUADRATIC VARIATION, $p$-VARIATION AND INTEGRATION \\ 
WITH APPLICATIONS TO STOCK PRICE MODELLING}

\author{Rimas Norvai\v sa\thanks{This research was partially supported 
by Lithuanian State Science and Studies Foundation Grant K014.}\\
Institute of Mathematics and Informatics, Vilnius, Lithuania}

\date{ \today}

\begin{document}
\maketitle

\newpage
\tableofcontents

\chapter{Introduction and results}\label{introduction}
\pagenumbering{arabic}
\setcounter{thm}{0}

This long paper deals with several aspects of calculus for real-valued
functions having a quadratic variation.
Here the notion of a quadratic variation is a property of a 
``deterministic" function rather than the well-known property of a 
Brownian motion discovered by L\'evy \cite{PL40}.
Let $f$ be a \emph{regulated} real-valued function on a closed interval 
$[a,b]$, that is, for $a\leq s<t\leq b$, there exist the limits 
$f(t-):=\lim_{u\uparrow t}f(u)$ and $f(s+):=\lim_{u\downarrow s}f(u)$.
Let $\lambda=\{\lambda_m\colon\,m\geq 1\}$ be a nested
sequence of partitions $\lambda_m=\{t_i^m\colon\,i=0,\dots,n(m)\}$
of $[a,b]$ such that $\cup_m\lambda_m$ is dense in $[a,b]$;
the class of all such $\lambda$ is denoted by $\Lambda [a,b]$.
We say that $f$ has the \emph{quadratic $\lambda$-variation} on $[a,b]$,
if there is a regulated function $V$ on $[a,b]$ such that $V(a)=0$ and
for each $a\leq s < t\leq b$, letting
$x_i^m:=(t_i^m\wedge t)\vee s$ for $i=0,\dots,n(m)$,
\beq\label{qv-var1}
V(t)-V(s)=\lim_{m\to\infty}\sum_{i=1}^{n(m)}\big [f(x_i^m)-f(x_{i-1}^m)
\big ]^2,
\eeq
\beq\label{qv-var2}
V(t)-V(t-)=[f(t)-f(t-)]^2\quad\mbox{and}\quad
V(s+)-V(s)=[f(s+)-f(s)]^2.
\eeq
The function $[f]_{\lambda}:=V$, if exists, is nondecreasing and is
called the \emph{bracket function} of $f$. 
Actually, existence of $V$ with the stated properties is equivalent to 
the definition of the quadratic $\lambda$-variation (see Definition 
\ref{function}.\ref{qv} and Proposition \ref{function}.\ref{alpha}).
The class of all functions having the quadratic $\lambda$-variation on
$[a,b]$ is denoted by $Q_{\lambda}[a,b]$.
F\"olmer \cite{HFa} introduced a quadratic variation for a regulated and
right-continuous function on $[0,\infty)$ to be a Radon measure on 
$[0,\infty)$ if exists as a limit of sums of point masses with respect 
to a sequence of partitions of a half line $[0,\infty)$ with vanishing 
mesh.
The bracket function $[f]_{\lambda}$ being a function rather than a 
measure, allows us to use a Stieltjes type integration for all
classes of functions appearing in this paper.
We notice that the quadratic variation of a function defined by
Wiener \cite{NW24}, if exists, does not depend on partitions, and is
closely related to the local $2$-variation of a function defined
by Love and Young \cite{LY38}.
To recall its definition let $f$ be a function on $[a,b]$, and let $0<p<\infty$.
For a partition $\kappa=\{t_i\colon\,i=0,\dots,n\}$ of $[a,b]$, let 
\beq\label{s_p}
s_p(f;\kappa):=\sum_{i=1}^n|f(t_i)-f(t_{i-1})|^p.
\eeq
We say that $f$ has the \emph{local $p$-variation},
$1<p<\infty$, if the limit
\beq\label{local-p-var}
\lim_{\kappa,\prtn}s_p(f;\kappa)
\eeq
exists in the sense of refinements of partitions
(see Appendix A for details).
Wiener's quadratic variation of $f$ means that (\ref{local-p-var}) holds
for $p=2$ and the limit in the sense of refinements is replaced by
the limit as the mesh of partitions tends to zero.
Notice that the limit (\ref{qv-var1}) is taken along a fixed sequence
of partitions. 
Both Wiener's quadratic variation and the local $2$-variation
do \emph{not} exist for almost all sample functions of a Brownian motion.

Again let $f$ be a function on $[a,b]$.
For $0<p<\infty$, the \emph{$p$-variation} of $f$  is defined by
\beq\label{-2intr}
v_p(f;[a,b]):=\sup\big\{s_p(f;\kappa)\colon\,\kappa\in\Xi [a,b]\},
\eeq
where $\Xi [a,b]$ denotes the class of all partitions of $[a,b]$.                                                                 
We say that $f$ has \emph{bounded $p$-variation} if $v_p(f;[a,b])<\infty$,
and denote the class of all such functions by ${\cal W}_p[a,b]$.
Usefulness of the $p$-variation property hinges on its relation to the
Stieltjes type integrability theory established around the late of thirties 
by L.\ C.\ Young in connection to applications in the theory of Fourier series.
Let $g$ and $f$ be a pair of functions on $[a,b]$ such that
\beq\label{-1intr}
f\in {\cal W}_p[a,b],\quad g\in {\cal W}_q[a,b]\quad\mbox{and}\quad 
p^{-1}+q^{-1}>1\quad\mbox{for}\quad p,q >0.
\eeq
Then by the L.\ C.\ Young Theorem on Stieltjes integrability  \cite{LCY36}
the integral $\smallint_a^bg\,df$
exists $(a)$ in the sense of Riemann-Stieltjes if $f$ and $g$ 
have no common discontinuities, $(b)$ in the sense of refinement 
Riemann-Stieltjes if $f$ and $g$ have no common discontinuities on the 
right and no common discontinuities on the left, $(c)$ always in the sense
defined by L.\ C.\ Young.
The functional calculus developed by many authors later on
has a full-fledged application to the class of functions having bounded 
$p$-variation for some $0<p<2$ (see \cite{DNc}).
In particular, the theory applies path by path to L\'evy processes without a 
Gaussian component, and to a fractional Brownian motion having the Hurst 
exponent  $H\in (1/2,1)$ because their almost all sample functions have
bounded $p$-variation for some $0<p<2$.
Let $B=\{B(t)\colon\,t\geq 0\}$ be a standard Brownian motion,
and let $X=\{X(t)\colon\,t\geq 0\}$ be a stochastic process
both defined on a complete probability space $(\Omega,{\cal F},\Pr)$.
Then by L.\ C.\ Young's Theorem on Stieltjes integrability,
for almost all $\omega\in\Omega$, the Riemann-Stieltjes integral
\beq\label{7intr}
(RS)\int_0^TX(t,\omega)\,dB(t,\omega)
\eeq
exists provided for some $0<p<2$, $X$ has almost all sample functions of 
bounded $p$-variation on $[0,T]$.
This fact was used by P.\ L\'evy to solve suitable Riemann-Stieltjes 
integral equations with respect to a Brownian motion (see e.g. \cite{PL53}).
In (\ref{7intr}), $X$ may be an $\alpha$-stable L\'evy process with
$0<\alpha<2$, or $X$ may be a fractional Brownian motion with
the Hurst index $1/2<H<1$.
However, one cannot replace $X$ in (\ref{7intr}) by $B$ (see Section 
\ref{nonex} for details).

A function $f$ on $[a,b]$ having bounded $p$-variation for some $1\leq p<2$ 
also has the quadratic $\lambda$-variation for each sequence of partitions
$\lambda=\{\lambda_m\colon\,m\geq 1\}\in\Lambda [a,b]$ such that  
\beq\label{two-sided-access}
\{t\in (a,b)\colon\,[f(t+)-f(t)][f(t)-f(t-)]\not =0\}\subset\cup_m\lambda_m.
\eeq
In this case, the bracket function $[f]_{\lambda}$ at $t\in [a,b]$
is a sum of squared jumps over $[a,t]$ 
(see Corollary \ref{function}.\ref{p-var-qv-var} below).
The converse is not true: given $\lambda\in\Lambda [0,1]$, almost all sample 
functions of a Brownian motion have the quadratic $\lambda$-variation on $[0,1]$, 
and have the $p$-variation on $[0,1]$ unbounded for each $1\leq p\leq 2$.
Another related example is provided by a subclass of self-affine functions 
introduced by K\^ono \cite{NK86} (see Example \ref{function}.\ref{NK}
below).
An extension of the $p$-variation calculus to Brownian motion path like functions
require a more refined constructions than the Riemann-Stieltjes integral.
The quadratic $\lambda$-variation and related integral constructions introduced
below make the core of our attempts towards building a desirable calculus.

\paragraph*{Evolution representation problem.}
The results in this paper are motivated by and applied to an evolution 
representation problem.
To sketch the problem, let $\BB=(\BB,\|\cdot\|)$ be a Banach algebra with 
unit $\unit$.
A family $U=\{U(t,s)\colon\,a\leq s\leq t\leq b\}\subset\BB$ is
called an {\em evolution} in $\BB$ if it is multiplicative, that is,
\beq\label{B-evolution}
\left\{\begin{array}{ll}
U(r,t)U(t,s)=U(r,s) &\mbox{for $a\leq s\leq t\leq r\leq b$},\\
U(t,t)=\unit &\mbox{for $a\leq t\leq b$}.
\end{array}\right.
\eeq
The notion of an evolution generalizes the concept of a one-parameter
semigroup of bounded linear operators on a Banach space.
The classical Hille-Yosida theorem describes any strongly continuous,
contractive semigroup in terms of its generator.
Analogous pairing results have been established for evolutions
under various conditions on the function $U_a$ defined by
$U_a(t):=U(t,a)$ for $t\in [a,b]$.
For example, $U_a$ is analytic, or almost differentiable, or
of bounded variation is one set of conditions.
Other type of conditions imply that the generators will be
unbounded operators in some cases, bounded in other.

Without additional assumptions a representation of an evolution may not
be unique.
We are interested in an evolution $U$ which satisfies the condition: 
there is a function $A$ on the simplex 
$\{(s,t)\colon\,a\leq s\leq t\leq b\}$ such that 
\beq\label{Frechet}
A(t,s)=\lim_{\kappa,\prtn}S(U;[s,t],\kappa),
\eeq
where $S(U;[s,t],\kappa):=\sum_{i=1}^n[U(t_i,t_{i-1})-\unit ]$
for a partition $\kappa=\{t_i\colon\,i=0,\dots,n\}$ of $[s,t]$.
Actually in this paper we are more concerned with the case when the limit under
refinement of partitions in (\ref{Frechet}) is replaced by a limit 
along a sequence of partitions.
In Chapter \ref{modelling} below, an evolution with such a condition 
and $\BB$ being the set of real numbers, is used to build a continuous 
time stock price model.
According to Dobrushin \cite{RLD53}, condition (\ref{Frechet}) appeared in 
earlier works of M.\ Fr\'echet and E.\ B.\ Dynkin, in connection 
of extending the Chapman-Kolmogorov equations to more general
transition probabilities of a Markov process.
The function $A$ defined by (\ref{Frechet}) is additive, and so
letting $h:=A(\cdot,a)$, we have $A(t,s)=h(t)-h(s)$ for each
$a\leq s\leq t\leq b$.
In some cases an evolution $U$ satisfying (\ref{Frechet}) can
be represented by a product integral with respect to the
function $h=A(\cdot,a)$.
For a $\BB$-valued function $h$ on $[a,b]$, and for a partition 
$\kappa=\{t_i\colon\,i=0,\dots,n\}$ of $[a,b]$, let
\beq\label{PandS}
P(h;\kappa)=P(h;[a,b],\kappa)
:=\prod_{i=1}^n\big [\unit +h(t_i)-h(t_{i-1})\big ].
\eeq
The \emph{product integral} on $[a,b]$ with respect to $h$ is defined to be 
the limit of products $P(h;[a,b],\kappa)$, 
if it exists in the sense of refinements of partitions $\kappa$, that 
is, using the notion of a limit of a directed function (see Appendix A),
\beq\label{prod-int}
\prodi_a^b(\unit +dh):=\lim_{\kappa,\prtn}P(h;[a,b],\kappa).
\eeq
Then the \emph{product integral representation} of an evolution $U$ means
the relation
$$
U(t,s)=\prodi_s^t(\unit +dh),\qquad a\leq s\leq t\leq b,
$$
where $h=A(\cdot,a)$ and $A$ is defined by (\ref{Frechet}).
Apparently the product integral representation of an evolution $U$
was known for a long time, but in all such cases
an evolution $U$ satisfies the condition:
\beq\label{classM}
\sup\Big\{\sum_{i=1}^n\|U(t_i,t_{i-1})-\unit\|\colon\,
a=t_0 < t_1 <\cdots <t_n=b,\,\,\, n\geq 1\Big\}<\infty.
\eeq
More often a product integral representation of an evolution was
used under various additional smoothness conditions
(see Masani \cite{PRM79} for a survey).

A most general product integral representation of an evolution  $U$
satisfying (\ref{classM}) was proved by Mac Nerney \cite{JSMN63} in the form
of a duality between multiplicative and additive functions.
Let ${\cal U}_1={\cal U}_1([a,b];\BB)$ be the class of all evolutions $U$ 
such that (\ref{classM}) holds,
and let ${\cal W}_1={\cal W}_1([a,b];\BB)$ be the class of all
$\BB$-valued functions on $[a,b]$ with bounded variation, that is, 
with bounded $1$-variation.
By Theorem 3.3 of Mac Nerney \cite{JSMN63}, there is a reversible 
mapping ${\cal P}$, from ${\cal W}_1$ onto ${\cal U}_1$, such that for 
$U\in {\cal U}_1$ and $h\in {\cal W}_1$, the following three statements 
are equivalent: 
\begin{enumerate}
\item[$(a)$] $U={\cal P}(h)$;
\item[$(b)$] $U(t,s)=\lim_{\kappa,\prtn}P(h;[s,t],\kappa)$
for each $a\leq s\leq t\leq b$;
\item[$(c)$] $h(t)-h(s)=\lim_{\kappa,\prtn}S(U;[s,t],\kappa)$
for each $a\leq s\leq t\leq b$.
\end{enumerate}
Thus $(b)$ gives the product integral representation of the evolution $U$
provided (\ref{classM}) and (\ref{Frechet}) hold.

We illustrate the preceding representation in the case when an evolution
$U$ on $\BB=\MM^d$, the Banach algebra of $d\times d$ matrices for some 
integer $d\geq 1$, is continuous and translation-invariant.
That is, suppose that the function $U_a\equiv U(\cdot,a)$, is 
continuous and
\beq\label{transl-invar}
U(t+u,s+u)=U(t,s)\qquad\mbox{for each $a\leq s\leq t\leq t+u\leq b$.} 
\eeq
Then letting $T(t-s):=U(t,s)$ for $a\leq s\leq t\leq b$, by 
(\ref{B-evolution}), it follows that the family
$T=\{T(u)\colon\,u\geq 0\}$ is a continuous one-parameter
semigroup:
$$
\left\{\begin{array}{ll}
T(u)T(v)=T(u+v) &\mbox{for each $u, v\geq 0$},\\
T(0)=\unit. &
\end{array}\right.
$$
In this case it is well-known that there exists the right-derivative
$T'(0)$ and $T(u)=\exp\{uT'(0)\}$ for $u\geq 0$, where the exponential 
$\exp {x}:=\sum_{k=0}^{\infty}x^k/k!$ for $x\in\MM^d$
(see e.g.\ Theorem I.2.9 and Proposition I.2.8 in Engel and Nagel \cite{EandN}).  
For any partition $\kappa=\{t_i\colon\,i=0,\dots,n\}$ of $[a,b]$,
we have
$$
\sum_{i=1}^n\big\|U(t_i,t_{i-1})-\unit\big\|
=\sum_{i=1}^n\big\|T(t_i-t_{i-1})-\unit\big\|
\leq (b-a)\sup_{u>0}u^{-1}\|T(u)-\unit\|< \infty,
$$
and so (\ref{classM}) holds for $U$.
Therefore the above stated Mac Nerney \cite{JSMN63} evolution representation
when applied to $U$ gives the representation $U={\cal P}(h)$ with
the function $h$ such that $h(t)-h(s)=(t-s)T'(0)$ and
$$
\prodi_s^t(\unit +dh)=U(t,s)=T(t-s)=\exp\{(t-s)T'(0)\}=\exp\{
h(t)-h(s)\}
$$
for all $a\leq s\leq t\leq b$.

As we have seen, an evolution $U$ in a Banach algebra $\BB$ has the 
product integral representation provided (\ref{classM}) holds.
Also in this case, an evolution gives a solution of a linear integral equation
as follows.
A $\BB$-valued function $f$ on $[a,b]$ is invertible if $f(t)$
is invertible for each $t\in [a,b]$, and $f^{\inv}$ is the reciprocal
function of $f$ defined by $f^{\inv}(t):=[f(t)]^{-1}$ for $t\in [a,b]$.
If the function $U_a\equiv U(\cdot,a)$ has bounded reciprocal 
$U_a^{\inv}$, then the reverse ${\cal P}^{-1}$ of Mac Nerney's 
mapping ${\cal P}$, has the Left Cauchy-Stieltjes integral representation.
Let $f$, $g$ be $\BB$-valued functions defined on $[a,b]$, and
for any partition $\kappa=\{t_i\colon\,i=0,\dots,n\}$ of $[a,b]$, let
$$
S_{LC}(f,g;[a,b],\kappa)
:=\sum_{i=1}^n\big [f(t_i)-f(t_{i-1})\big ]g(t_{i-1}).
$$
The \emph{Left Cauchy-Stieltjes} integral is defined by
$$
(LCS)\int_a^bdf\,g:=\lim_{\kappa,\prtn}S_{LC}(f,g;[a,b],\kappa)
$$
provided the limit exists in the sense of refinement of partitions.
Let the reciprocal $U_a^{\inv}$ of the function $U_a$ be defined
and bounded on $[a,b]$.
By (\ref{B-evolution}), for any partition 
$\kappa=\{t_i\colon\,i=0,\dots,n\}$ of $[a,b]$, we have
$$
S(U;[a,b],\kappa)=\sum_{i=1}^n\big [U(t_i,t_{i-1})-\unit\big ]
=\sum_{i=1}^n\big [U_a(t_i)-U_a(t_{i-1})\big ]U_a^{\inv}(t_{i-1})
=S_{LC}(U_a,U_{a}^{\inv};[a,b],\kappa).
$$
Also (\ref{classM}) holds if $U_a$ has bounded variation.
Therefore the integral $(LCS)\smallint_a^bdU_a\,U_a^{\inv}$ exists,
and the representation $U={\cal P}(h)$ holds with the function $h$
defined by
$$
h(t)=(LCS)\int_a^tdU_a\,U_a^{\inv},\qquad a\leq t\leq b.
$$
A representation of a solution of a linear Left Cauchy-Stieltjes integral 
equation in terms of a flow in $\BB$ is another well-known characterization 
of the evolution representation.
Suppose that $U={\cal P}(h)$ for some $U\in {\cal U}_1([a,b];\BB)$ and
$h\in {\cal W}_1([a,b];\BB)$.
Then by Theorem 4.1 of Mac Nerney \cite{JSMN63}, in the class ${\cal W}_1$, 
the linear Left Cauchy-Stieltjes integral equation
$$
f(t)=f(a)+(LCS)\int_a^tdh\,f,\qquad a\leq t\leq b,
$$
has the unique solution $f(t)=U(t,a)f(a)$, $t\in [a,b]$.

Replacing boundedness of the variation of the function $U_a$, which is
condition (\ref{classM}) above,
by a weaker condition provides an interesting challenge to the
evolution representation problem.
Suppose that the reciprocal function $U_a^{\inv}$ of the
function $U_a\equiv U(\cdot,a)$ exists and is bounded.
In this case, boundedness of the $p$-variation of $U_a$ for some
$1\leq p<\infty$ is one such weaker condition.
The main advantage of this candition is a well established
integration theory for a class of functions with bounded
$p$-variation.
Freedman \cite{MAF}  applied L.\ C.\ Young's results and ideas to extend
the product integral representation of an evolution $U$ provided
the function $U_a$ is continuous and has bounded $p$-variation
for some $1\leq p<2$.
Using the Banach fixed point theorem Freedman \cite[Theorem 5.1]{MAF}
proved the existence of the product integral as a unique solution
of the linear Riemann-Stieltjes integral equation.
A non-linear Riemann-Stieltjes integral equation with respect
to a continuous function of bounded $p$-variation, $1\leq p<2$,
was solved by Lyons \cite{TL94}.

Dudley and Norvai\v sa \cite[Part II]{DNa}  started their investigation
of the product integral (\ref{prod-int})
with respect to a Banach algebra valued function $h$ from scratch.
By their Theorem 4.4, for a real valued function $f$,
the product integral with respect to $f$ exists and is non-zero
if and only if $f$ is regulated function satisfying $(a)$ and $(b)$:
\begin{enumerate} 
\item[$(a)$] $f$ has local $2$-variation, that is, 
(\ref{local-p-var}) holds with $p=2$;
\item[$(b)$] for any $a\leq s<t\leq b$,
\beq\label{dfnot=-1}
\Delta^{+}f(s):=f(s+)-f(s)\not =-1\not =f(t)-f(t-)=:\Delta^{-}f(t).
\eeq
\end{enumerate}
Moreover, if the product integral exists then for any $a\leq s<t\leq b$,
$$
\prodi_s^t(1+df)=e^{f(t)-f(s)}\prod_{[s,t]}\big [(1+\Delta^{-}f)
(1+\Delta^{+}f)\big ]e^{-\Delta^{-}f-\Delta^{+}f},
$$
where $\Delta^{-}f(s):=0$ and $\Delta^{+}f(t):=0$.
Assuming that $f$ is a regulated function, a new proof of this fact is given 
below by Theorem \ref{function}.\ref{prodint}.
For a function $h$ with values in a Banach algebra $\BB$, let
${\cal P}_a(h)=\{{\cal P}_a(h;t)\colon\,a\leq t\leq b\}$ be the
indefinite product integral defined by
$$
{\cal P}_a(h;t):=\prodi_a^t(\unit +dh),\qquad a\leq t\leq b,
$$
and let ${\cal W}_p([a,b];\BB)$ be the Banach space of 
all $\BB$-valued functions with bounded $p$-variation on $[a,b]$.
The main result in \cite[Part II]{DNa}
is that the mapping ${\cal P}_a$ from ${\cal W}_p([a,b];\BB)$ into itself, 
is Fr\'echet differentiable and ${\cal P}_a$ is analytic when 
restricted to right- or left-continuous functions $h$. 
To formulate the evolution representation in this case,
for a $\BB$-valued function $f$, let ${\cal L}_a(f)
=\{{\cal L}_a(f;t)\colon\,a\leq t\leq b\}$ be the indefinite
Left Young integral on $[a,b]$ defined by
$$
{\cal L}_a(f;t):= (LY)\int_a^tdf\,f^{\inv},\qquad a\leq t\leq b,
$$
provided the Left Young integral $(LY)\smallint_a^bdf\,f^{\inv}$
exists (see Definition \ref{variation}.\ref{LYandRY} below).
The following is Theorem 6.14 in \cite[Part II]{DNa}.

\begin{thm}
Let $0<p<2$, let $U$ be an evolution in $\BB$ such that $U_a\equiv U(\cdot,a)
\in {\cal W}_p([a,b];\BB)$, and let the reciprocal function
$U_a^{\inv}$ exists and is bounded.
Then $(a)$, $(b)$ and $(c)$ hold{\rm :}
\begin{enumerate}
\item[$(a)$] the indefinite Left Young integral ${\cal L}_a(U_a)$ on $[a,b]$
is defined and is in ${\cal W}_p([a,b];\BB)${\rm ;} 
\item[$(b)$] the product integral on $[a,b]$
with respect to ${\cal L}_a(U_a)$ is defined{\rm ;} 
\item[$(c)$] the function ${\cal L}_a(U_a)$ in ${\cal W}_p([a,b];\BB)$
is unique up to an additive function such that
$$
U(t,s)=\prodi_s^t (\unit +d{\cal L}_a(U_a)),\qquad a\leq s\leq t\leq b.
$$ 
\end{enumerate}
\end{thm}

The main result of this paper is a solution of a representation 
problem for a class of evolutions $U$ in $\RR$ such that
$U_a\equiv U(\cdot,a)$ has the quadratic $\lambda$-variation
for some $\lambda\in\Lambda [a,b]$. 
Before formulating this result we summarize other results.

\paragraph*{$(\lambda,p)$-decomposable function and 
$(2-\epsilon)$-semimartingale.}
Let $\lambda\in\Lambda [a,b]$ and $1\leq p<2$.
As it was mentioned earlier, a function $f$ having bounded $p$-variation
also has the quadratic $\lambda$-variation provided a two sided jump of $f$
occure only at a point of a partition of $\lambda$.
By the preceding remarks such functions in the class $Q_{\lambda}[a,b]$ 
of all functions having the quadratic $\lambda$-variation may be considered 
as special: for such functions we have a calculus quite similar to 
the classical calculus for smooth functions (see \cite{DNc}).
Therefore when proving something for a function having the quadratic
$\lambda$-variation it is desirable to make the new result concerted with
the calculus applicable to functions having bounded $p$-variation.
To this aim we consider sums $g+h$ such that the $p$-variation
is bounded for $h$ and may be unbounded for $g$, but $g$ still has the
quadratic $\lambda$-variation.
We also found useful to assume when $p>1$ that $g$ has bounded $q$-variation 
for some $q>2$ such that $p^{-1}+q^{-1}>1$.
This condition is satisfied by our typical examples mentioned earlier.  
Under an additional mild condition on jumps of $g$ and $h$,
the sum $g+h$ also has the quadratic $\lambda$-variation
(see Corollary \ref{function}.\ref{qvclassQ} below). 
We also treat stochastic processes with a similar decomposability property
except that the first component of such a decomposition has to be
a local martingale.

Recall that  ${\cal W}_p[a,b]$, $0<p<\infty$, denotes the class of all 
real-valued functions $f$ on a closed bounded interval $[a,b]$
with bounded $p$-variation defined by  (\ref{-2intr}).
Also, a real-valued  function $f$ on $[a,b]$ is regulated on $[a,b]$
if for $a\leq s<t\leq b$, there exist the limits 
$f(t-):=\lim_{u\uparrow t}f(u)$ and $f(s+):=\lim_{u\downarrow s}f(u)$.
If for some  $p<\infty$, $f$ has bounded $p$-variation on $[a,b]$ 
then $f$ is regulated on $[a,b]$. 
A point $t\in (a,b]$ is \emph{left discontinuity point of $f$} if
$\Delta^{-}f(t):=f(t)-f(t-)\not =0$, and a point $s\in [a,b)$ is \emph{right 
discontinuity point of $f$} if $\Delta^{+}f(s):=f(s+)-f(s)\not =0$.
The class of all regulated functions on $[a,b]$ is denoted by ${\cal R}[a,b]$.
For  $1\leq p<\infty$, let
\beq\label{dual-Wp}
\dual ({\cal W}_p)[a,b]:=\left\{ \begin{array}{ll} 
\cup\big\{{\cal W}_q[a,b]\colon\,\frac{1}{p}+\frac{1}{q}>1\big\}
&\mbox{if $1<p<\infty$,}\\
{\cal R}[a,b]&\mbox{if $p=1$.}
\end{array}\right.
\eeq
Thus $f\in \dual ({\cal W}_p)[a,b]$ with $1<p<\infty$ if $f$
has bounded $q$-variation for some $q<p/(p-1)$, and
$f\in \dual ({\cal W}_1)[a,b]$ if $f$ is a regulated function. 
Now we are ready to define the class of $(\lambda,p)$-decomposable 
considered in Chapter \ref{function}.

\begin{defn}\label{decomposable}
{\rm Let $f$ be a regulated function on $[a,b]$, let $\lambda\in \Lambda [a,b]$, 
and let $1\leq p<2$.
We say that $f$ is \emph{{\rm ($\lambda,p)$}-decomposable} on $[a,b]$, or
$f\in {\cal D}_{\lambda,p}[a,b]$,
if there exists a pair $(g,h)$ of regulated functions on $[a,b]$ 
such that $g\in Q_{\lambda}[a,b]\cap \dual ({\cal W}_p)[a,b]$, 
$h\in {\cal W}_p[a,b]$, $g(a)=h(a)=0$ and $f\equiv C+g+h$ for some
constant $C$.
Any such pair $(g,h)$ is called a $(\lambda,p)$-dual pair 
$(\lambda,p)$-decomposing $f$.
The set of all $(\lambda,p)$-dual pairs which
$(\lambda,p)$-decompose $f$ is denoted by $D_{\lambda,p}(f)$.}
\end{defn}

Among all $(\lambda,p)$-dual pairs we select those pairs which have
left or right discontinuities at points of partitions of $\lambda$. 
Given $\lambda=\{\lambda_m\colon\,m\geq 1\}\in\Lambda [a,b]$, we often
assume to hold
\beq\label{right-jumps-access}
\mbox{either}\quad N_{(a,b)}(\Delta^{+}f):=
\{t\in (a,b)\colon\,(\Delta^{+}f)(t)\not =0\}\subset\cup\lambda
\eeq
\beq\label{left-jumps-access}
\mbox{or}\qquad N_{(a,b)}(\Delta^{-}f):=
\{t\in (a,b)\colon\,(\Delta^{-}f)(t)\not =0\}\subset\cup\lambda \,,
\eeq
where $\cup\lambda :=\cup_m\lambda_m$.
We say that the \emph{right discontinuity points of $f$ are accessible by
$\lambda$} if $N_{(a,b)}(\Delta^{+}f)\subset\cup\lambda$.
Similarly, we say that the \emph{left discontinuity points of $f$ are accessible 
by $\lambda$} if $N_{(a,b)}(\Delta^{-}f)\subset\cup\lambda$.
Clearly, the first condition holds for each $\lambda\in\Lambda [a,b]$
if $f$ is right-continuous, and the second one holds if $f$ is left-continuous.
In some cases we use a weaker assumption (see (\ref{two-sided-access}) above): 
we say that the \emph{two-sided
discontinuity points of $f$ are accessible by $\lambda$} if 
\beq\label{jumps-access}
N_{(a,b)}(\Delta^{-}f)\cap N_{(a,b)}(\Delta^{+}f)=\{t\in (a,b)\colon\,
\big (\Delta^{-}f\Delta^{+}f\big )(t)\not =0\}\subset\cup\lambda\,.
\eeq
Now let $f\in {\cal D}_{\lambda,p}[a,b]$ for some $\lambda\in\Lambda [a,b]$ 
and $1\leq p<2$.
We denote by $D_{\lambda,p}^{+}(f)$ the set of all pairs 
$(g,h)\in D_{\lambda,p}(f)$ such that
\beq\label{+accessible}
N_{(a,b)}(\Delta^{+}g)\cup N_{(a,b)}(\Delta^{+}h)\subset\cup\lambda.
\eeq
Similarly,  $D_{\lambda,p}^{-}(f)$ denotes the set of all
pairs $(g,h)\in D_{\lambda,p}(f)$ such that (\ref{+accessible})
holds with $\Delta^{+}$ replaced by $\Delta^{-}$.

Suppose that $f\equiv C+g+h$ for some $(g,h)\in D_{\lambda,p}^{+}(f)
\cup D_{\lambda,p}^{-}(f)$. 
By Corollary \ref{function}.\ref{qvclassQ},  the function $f$ has 
the quadratic $\lambda$-variation, and its bracket function is given by
\beq\label{4qvclassQ}
[f]_{\lambda}(t)=[g]_{\lambda}^c(t)+\sum_{(a,t]}\big (\Delta^{-}
f\big )^2+\sum_{[a,t)}\big (\Delta^{+}f\big )^2\qquad a\leq t\leq b.
\eeq
Thus the representation $f\equiv C+g+h$ is not unique because
$(f-f(a),0)\in D_{\lambda,p}^{+}(f)\cup D_{\lambda,p}^{-}(f)$ 
is another $(\lambda,p)$-decomposition of $f$.
Nevertheless several constructions defined in this paper for functions
$f$ of the class ${\cal D}_{\lambda,p}[a,b]$, do not depend on
a representation of $f$ by a $(\lambda,p)$-dual pair in
$D_{\lambda,p}^{+}(f)\cup D_{\lambda,p}^{-}(f)$.
For example, if $D_{\lambda,p}^{+}(f)\cup D_{\lambda,p}^{-}(f)\not =\emptyset$
then the continuous part of the bracket function $[f]_{\lambda}$ does not 
depend on a $(\lambda,p)$-decomposition.
This follows from formula (\ref{4qvclassQ}).
Therefore we define $[f]_{\lambda}^c:=[g]_{\lambda}^c$ with a $g$ being such 
that $(g,h)\in D_{\lambda,p}^{+}(f)\cup D_{\lambda,p}^{-}(f)$.
Below we define two integrals, the Left and Right $\lambda$-integrals, 
with respect to a function $f$ of the class ${\cal D}_{\lambda,p}[a,b]$ 
which under mild additional assumptions also do not depend on nonuniqueness 
of $(\lambda,p)$-decomposability of $f$.

Another nice feature of the class of $(\lambda,p)$-decomposable functions 
is that it is closed under taking  a composition with the $C^2$ class functions.

\begin{thm}\label{invar}
For $\lambda\in\Lambda [a,b]$ and $1\leq p<2$, let $f$ be a 
$(\lambda,p)$-decomposable function on $[a,b]$, and let $\phi$ be a 
$C^2$ class function.
If $D_{\lambda,p}^{+}(f)\cup D_{\lambda,p}^{-}(f)\not =\emptyset$ then 
the composition  $\phi{\circ}f$ is a $(\lambda,p)$-decomposable function.
\end{thm}

The proof of this statement given at the end of Section \ref{l-integrals},
shows that the composition $\phi{\circ}f$ has at least two 
different $(\lambda,p)$-decompositions if
$D_{\lambda,p}^{+}(f)\cap D_{\lambda,p}^{-}(f)\not =\emptyset$.

Let $1\leq p<2$.
We say that a stochastic process $X$ is a \emph{$p$-semimartingale}
if $X-X(0)=M+A$, where $M(0)=A(0)=0$, $M$ is a local martingale
and $A$ is a stochastic process having almost all sample functions of
bounded $p$-variation on each bounded time interval.
In the case $p=1$, $X$ is usually called a semimartingale.
We also say that $X$ is a \emph{$(2-\epsilon)$-semimartingale} if it is
a $p$-semimartingale for some $1\leq p<2$.
This class of proceesses may be considered as a stochastic analog
of the class of $(\lambda,p)$-decompasable functions.
In Chapter \ref{process} we prove several elementary properties
for such processes.
Given a $(2-\epsilon)$-semimartingale $X$ we define the stochastic
$\lsi$-integral with respect to $X$ to be a sum of two integrals with respect
to each component of a decomposition and show that values of such sums
do not depend on a decomposition in most interesting cases.  
Analogous uniqueness reult is proved with respect to the continuous
part of a quadratic variation of a $(2-\epsilon)$-semimartingale.
Finally in Chapter \ref{process} we prove that in the class of all
$(2-\epsilon)$-semimartingales, an extended stochastic Dol\'eans 
exponential is a unique solution of the linear stochastic $\lsi$-integral
equation.

\paragraph*{ Integrals of functions having the quadratic $\lambda$-variation.}
As we mentiened earlier, the L.\ C.\ Young Theorem on Stieltjes integrability
does not hold in general when the integrand and integrator both
have the quadratic $\lambda$-variation.
For such functions we use restricted Stieltjes type integrals defined to be
limits, if exist, along the sequence of partitions $\lambda$.
An existence of such an integral with a special form of the integrand can
be proved by replicating a proof of It\^o's formula (see F\"olmer \cite{HFa}).
In this paper we follow the same route except tha we
do not assume functions to be right or left continuous,
and use a Riemann-Stieltjes type integral to integrate with respect to a bracket 
function.

Let $\lambda=\{\lambda_m\colon\,m\geq 1\}\in\Lambda [a,b]$
with $\lambda_m=\{t_i^m\colon\,i=0,\dots,n(m)\}$ for $m\geq 1$.
For regulated functions $g$ and $f$ on $[a,b]$,
the \emph{Left Cauchy $\lambda$-integral} $(LC)\smallint g\,d_{\lambda}f$
is defined on $[a,b]$ if there exists a regulated function $\Phi$ 
on $[a,b]$ such that $\Phi (a)=0$ and  for each $a\leq s<t\leq b$,
letting $x_i^m=(t_i^m\wedge t)\vee s$ for $i=0,\dots,n(m)$,
$$
\Phi (t)-\Phi (s)=\lim_{m\to\infty}\sum_{i=1}^{n(m)}g(x_{i-1}^m)
\big [f(x_i^m)-f(x_{i-1}^m)\big ],
$$
$$
\Delta^{-}\Phi (t)=\big (g_{-}\Delta^{-}f\big )(t)
\qquad\mbox{and}\qquad
\Delta^{+}\Phi (s)=\big (g\Delta^{+}f\big )(s).
$$
If such $\Phi$ exists then we also say that $g$ is Left Cauchy 
$\lambda$-integrable on $[a,b]$ with respect to $f$, and  for $a\leq s<t\leq b$, 
let $(LC)\int_s^tg\,d_{\lambda}f:=\Phi (t)-\Phi (s)$.
In fact we define slightly more general integral when $g$ may not be regulated
(Definition \ref{function}.\ref{LCint}). 
Such extension is used to prove the Black-Scholes formula in Section
\ref{option}.
Similarly the Right Cauchy $\lambda$-integral is defined
(Definition \ref{function}.\ref{RCint}).

The following statement is a corollary of Theorem \ref{function}.\ref{chainrule}.
We omit its proof because it is symmetric to the proof of
Proposition \ref{function}.\ref{1chqv} concerning the Right Cauchy 
$\lambda$-integral.

\begin{prop}\label{chqv}
Let $f$ be a regulated function on $[a,b]$, and let $\lambda\in\Lambda [a,b]$ 
be such that the right discontinuity points of $f$ are accessible by $\lambda$.
The following statements about $f$ are equivalent{\rm :}
\begin{enumerate}
\item[$(a)$] $f$ has the quadratic $\lambda$-variation on $[a,b]${\rm ;}
\item[$(b)$] for a $C^1$ class function $\psi$, the composition $\psi{\circ}f$
is Left Cauchy $\lambda$-integrable on $[a,b]$ with respect to $f$, 
and for any $a\leq s<t\leq b$,
\begin{eqnarray*}
(\Psi{\circ}f)(t)&=&(\Psi{\circ}f)(s)
+(LC)\int_s^t(\psi{\circ}f)\,d_{\lambda}f
+\frac{1}{2}(RS)\int_s^t(\psi'{\circ}f)\,d[f]_{\lambda}^c\\[2mm]
& &+\sum_{(s,t]}\Big\{\Delta^{-}(\Psi{\circ}f)-(\psi{\circ}f)_{-}
\Delta^{-}f\Big\}+\sum_{[s,t)}\Big\{\Delta^{+}(\Psi{\circ}f)-
(\psi{\circ}f)\Delta^{+}f\Big\},
\end{eqnarray*}
where the two sums converge unconditionally
and $\Psi (u)=\Psi (0)+\smallint_0^u\psi (x)\,dx$, $u\in\RR${\rm ;}
\item[$(c)$] $f$ is Left Cauchy $\lambda$-integrable on $[a,b]$ with respect to
itself.
\end{enumerate}
If any of the three statements holds, then for $a\leq t\leq b$,
$$
(LC)\int_a^tf\,d_{\lambda}f=\frac{1}{2}\left\{f^2(t)-f^2(a)
-[f]_{\lambda}(t)\right\}.
$$
\end{prop}

Now we define an integral with respect to a $(\lambda,p)$-decomposable
function.
Let $\lambda\in\Lambda [a,b]$, let $1\leq p<2$ and let $f$ be
a $(\lambda,p)$-decomposable function.
For a regulated function $F$ and a set $A\subset D_{\lambda,p}(f)$,
we say that the \emph{Left $\lambda$-integral $(L)\smallint F\,d_{\lambda}f$ 
is defined on $[a,b]$ with respect to $A$} if for any $(g,h)\in A$, on the right side
$$
(L)\int F\,d_{\lambda}f:=
(LC)\int F\,d_{\lambda}g+(LY)\int F\,dh
$$
the Left Cauchy $\lambda$-integral and the Left Young integral are both 
defind on $[a,b]$, and sums of values at $a\leq s<t\leq b$
do not depend on $(g,h)\in A$.
In Proposition \ref{function}.\ref{property3} we show that the Left
$\lambda$-integral $(L)\smallint F\,d_{\lambda}f$ is defined  on $[a,b]$
with respect to $D_{\lambda,p}^{+}(f)\not =\emptyset$ provided
$F\in\dual ({\cal W}_p)[a,b]$ and $(LC)\smallint F\,d_{\lambda}g$
is defined on $[a,b]$ for some $(g,h)\in D_{\lambda,p}^{+}(f)$.
Analogous statement holds for the Right $\lambda$-integral
defined symmetrically.

Next we summarize results concerning an extension of the product integral 
(\ref{prod-int}).
Let $\lambda=\{\lambda_m\colon\,m\geq 1\}\in\Lambda [a,b]$
with $\lambda_m=\{t_i^m\colon\,i=0,\dots,n(m)\}$ for $m\geq 1$.
Let $f$ be a regulated function on $[a,b]$.
By Proposition \ref{function}.\ref{prod-int2}, the product $\lambda$-integral
$\prodi (1+d_{\lambda}f)$ is defined on $[a,b]$ if and only if there exists
a regulated function $H$ on $[a,b]$ such that $H(a)=1$, $|H|\gg 0$
and for $a\leq s<t\leq b$, letting $x_i^m=(t_i^m\wedge t)\vee s$ for 
$i=0,\dots,n(m)$,
$$
H(t)/H(s)=\lim_{m\to\infty}\prod_{i=1}^{n(m)}
\big [1+f(x_i^m)-f(x_{i-1}^m)\big ],
$$
$$
H(t)/H(t-)=1+(\Delta^{-}f)(t)\qquad\mbox{and}\qquad
H(s+)/H(s)=1+(\Delta^{+}f)(s). 
$$
If such $H$ exists then for $a\leq s<t\leq b$, let 
$\prodi_s^t (1+d_{\lambda}f):=H(t)/H(s)$.
By Proposition \ref{function}.\ref{prod->l-prod}, if the product integral
$\prodi (1+df)$ exists and is non-zero on $[a,b]$, and if the
two-sided discontinuity points of $f$ are accessible by $\lambda$, then
the product $\lambda$-integral $\prodi (1+d_{\lambda}f)$
is defined on $[a,b]$, and the two integrals have the same values.

Let $f$ be a regulated function on $[a,b]$ having the quadratic 
$\lambda$-variation for some $\lambda\in\Lambda [a,b]$.
The Dol\'eans exponential ${\cal E}_{\lambda,a}(f)$ on $[a,b]$ is defined by
$$
{\cal E}_{\lambda,a}(f;t):=\left\{ \begin{array}{ll}
\exp\big\{f(t)-f(a)-\frac{1}{2}[f]_{\lambda}^c(t)\big\}
\prod_{[a,t]}(1+\Delta f)e^{-\Delta f}&
                          \mbox{if $t\in (a,b]$}\\
1 &\mbox{if $t=a$,}
\end{array}\right. 
$$
where the product on the right side is defined by (\ref{1bvofV}).
By Theorem \ref{function}.\ref{ExistLprod}, if 
the two-sided discontinuity points of $f$ are accessible by $\lambda$,
and if (\ref{dfnot=-1}) holds, then the product $\lambda$-integral
$\prodi (1+d_{\lambda}f)$ is defined on $[a,b]$ and for $a\leq t\leq b$,
$$
\prodi_a^t(1+d_{\lambda}f)={\cal E}_{\lambda,a}(f;t).
$$

For $\lambda\in\Lambda [a,b]$ and $1\leq p<2$, let $f$ be a function
on $[a,b]$ which is $(\lambda,p)$-decomposable and $D_{\lambda,p}^{+}(f)$
is nonempty.
Thus $f$ has the quadratic $\lambda$-variation.
Consider the forward linear Left $\lambda$-integral equation
\beq\label{6intr}
F(t)=1+(L)\int_a^tF\,d_{\lambda}f,\qquad a\leq t\leq b,
\eeq
provided the Left $\lambda$-integral is defined on $[a,b]$ with respect
to $D_{\lambda,p}^{+}(f)$.
Using a chain rule it is easy to check that the Dol\'eans exponential 
${\cal E}_{\lambda,a}(f)$ is a solution of (\ref{6intr}).
To show its uniqueness in the class $\dual ({\cal W}_p)[a,b]$ we augment 
equation (\ref{6intr}) as formulated in Definition \ref{function}.\ref{solution}.
The augmented equation is solved by Theorem \ref{LIE} in Section 
\ref{uniqueness}.
Analogous result for the backward linear Right $\lambda$-integral 
equation is proved by Theorem \ref{BLIE} in the same Section
\ref{uniqueness}.

\paragraph*{Evolution and the quadratic $\lambda$-variation.}
Now we are prepared to formulate the main result.
For a function $f$ on $[a,b]$, $f\gg 0$ means that there is a constant $C>0$ 
such that $f(t)\geq C$ for each $t\in [a,b]$.
By Proposition \ref{function}.\ref{regul-evol2}, an
\emph{evolution $U$ in $\RR$ is regulated} if and only if $U_a\equiv U(\cdot,a)$
is a regulated function on $[a,b]$ and $|U_a|\gg 0$.
Let $\lambda=\{\lambda_m\colon\,m\geq 1\}\in\Lambda [a,b]$
with $\lambda_m=\{t_i^m\colon\,i=0,\dots,n(m)\}$ for $m\geq 1$.
By Proposition \ref{function}.\ref{generator}, a regulated evolution $U$ in $\RR$
has the \emph{$\lambda$-generator} if and only if there exists a regulated 
function $G$ on $[a,b]$ such that $G(a)=0$ and for each $a\leq s<t\leq b$,
letting $x_i^m=(t_i^m\wedge t)\vee s$ for $i=0,\dots,n(m)$,
$$
G(t)-G(s)=\lim_{m\to\infty}\sum_{i=1}^{n(m)}\big [U(x_{i-1}^m,x_i^m)-1\big ],
$$
$$
\Delta^{-}G(t)=U_a(t)/U_a(t-)-1\quad\mbox{and}\quad
\Delta^{+}G(s)=U_a(s+)/U_a(s)-1.
$$
The function $G$ is called the $\lambda$-generator of $U$.
This function generalizes the notion of the generator of a one-parameter 
semigroup as follows.
Let $U$ be a continuous translation-invariant evolution on $[a,b]$
as defined by (\ref{transl-invar}).
Then letting $T(y-z):=U(z,y)$ it follows that the family
$T=\{T(u)\colon\,u\geq 0\}$ is a continuous one-parameter semigroup.
By the results of Mac Nerney \cite{JSMN63}, we have that for each
$\lambda\in\Lambda [a,b]$, $U$ has the $\lambda$-generator $G$
such that $G(t)-G(s)=(t-s)T'(0)$ for $a\leq s<t\leq b$.

To define a mapping analogous to Mac Nerney's duality mapping ${\cal P}$ 
between  ${\cal W}_1$ and ${\cal U}_1$, we use two sets of
functions having the quadratic $\lambda$-variation.
For $\lambda\in\Lambda [a,b]$, let
\beq\label{class-L}
L_{\lambda}[a,b]:=\big\{f\in Q_{\lambda}[a,b]\colon\,
f(a)=1,\quad f\gg 0\quad\mbox{and}\quad N_{(a,b)}(\Delta^{+}f)
\subset\cup\lambda\,\big\}\quad\mbox{and}
\eeq
\beq\label{class-E}
E_{\lambda}[a,b]:=\big\{g\in Q_{\lambda}[a,b]\colon\,
g(a)=0,\quad 1+\Delta^{-}g\wedge\Delta^{+}g\gg 0\quad\mbox{and}\quad
N_{(a,b)}(\Delta^{+}g)\subset\cup\lambda\,\big\}.
\eeq
We consider a regulated evolution $U$ such that $U_a\in L_{\lambda}[a,b]$
for some $\lambda\in\Lambda [a,b]$. 
To show that $U$ has the $\lambda$-generator and to
give its representation, we use the Left Cauchy $\lambda$-integral
and the product $\lambda$-integral.
We say that the function ${\cal L}_{\lambda}f$ is defined on $[a,b]$ if the 
Left Cauchy $\lambda$-integral $(LC)\smallint f^{-1}\,d_{\lambda}f$ is defined 
on $[a,b]$, and then let 
${\cal L}_{\lambda}f(t):=(LC)\smallint_a^tf^{-1}\,d_{\lambda}f$ for $t\in [a,b]$.

\begin{thm}\label{evolution}
Let  $\lambda\in\Lambda [a,b]$, and let $U$ be a regulated evolution 
in $\RR$ such that $U_a\equiv U(\cdot,a)\in L_{\lambda}[a,b]$.
Then  statements $(a)$, $(b)$ and $(c)$ hold, where{\rm :}
\begin{enumerate}
\item[$(a)$] the function ${\cal L}_{\lambda}U_a$ is defined on $[a,b]$, 
it is in $E_{\lambda}[a,b]$ and it is the $\lambda$-generator of $U${\rm ;}
\item[$(b)$] the product $\lambda$-integral
$\prodi (1+d_{\lambda}({\cal L}_{\lambda}U_a))$ is defined on $[a,b]${\rm ;}
\item[$(c)$] the function ${\cal L}_{\lambda}U_a$ in $E_{\lambda}[a,b]$ 
is unique up to an additive constant such that
\beq\label{1evolution}
 U(t,s)=\prodi_s^t\big (1+d_{\lambda}\big ({\cal L}_{\lambda}U_a\big )
\big )\qquad a\leq s\leq t\leq b.
\eeq
\end{enumerate}
\end{thm}

\begin{proof}
Statement $(a)$ holds by Proposition \ref{function}.\ref{Elogarithm}
and Theorem \ref{function}.\ref{generator1}.
Statement $(b)$ holds by Theorem \ref{function}.\ref{ratio}.
By (\ref{1ratio}) of Theorem \ref{function}.\ref{ratio}, we have that for all 
$a\leq s\leq t\leq b$,
$$
\prodi_s^t\big (1+d_{\lambda}\big ({\cal L}_{\lambda}U_a\big )\big )
=U_a(t)/U_a(s)=U(t,s),
$$
and so (\ref{1evolution}) holds.
Suppose that $g\in E_{\lambda}[a,b]$ is such that 
$\prodi_a^t (1+d_{\lambda}g)=U(t,a)$ for $t\in [a,b]$.
Then ${\cal L}_{\lambda}U_a=g+C$ for some constant $C$
by Theorem \ref{function}.\ref{Ulogarithm}, 
proving uniqueness of the representation, and hence statement $(c)$.
The proof of Theorem \ref{evolution} is complete.
\qed\end{proof}

In Chapter \ref{modelling}, the solution of the evolution representation
problem is used to build up an asset pricing model.
More specifically, stock price changes during a time period $[0,T]$ are 
considered as an evolution $U$ such that 
$U_0\equiv U(\cdot,0)\in L_{\lambda}[0,T]$.
By Theorem \ref{evolution}, the evolution $U$ has the $\lambda$-generator 
which is used to model a return of a stock in our asset pricing model.
In this chapter there are two new results: an almost sure approximation of
a continuous time Black-Scholes model by a dicrete time binomial model,
and the option pricing formula when a stock price is modelled by
a function having the quadratic $\lambda$-variation.
The second result clarifies and solves the problem posed in Bick
and Willinger \cite{BandW94}.

\chapter{The $p$-variation and integrals of Stieltjes type}\label{variation}
\setcounter{thm}{0}

\vspace*{0.2truein}
\begin{quotation}{\footnotesize
One thing at least we may say: no theory of integration is complete
which does not enable us to define, by some kind of integration process,
the integral of a continuous function with respect to itself, or,
more generally, the integration of an ordinary function, such as
the square, of a continuous function with respect to that
continuous function, when the latter has not bounded variation.

W.\ H.\ Young. ``The progress of mathematical analysis in the 
twentieth century'' (Presidential Address). Proc. London Math. Soc.
Ser. 2, Vol. 24 (1925), 421-434. }
\end{quotation}
\vspace*{0.2truein}

The problem of integration of functions having unbounded variation 
was approached in several different ways.
Although only very few of these approaches have reached a development
to a full fledged calculus.
Integration theories dealing with suitable families of functions,
such as stochastic processes or supports of probability distributions,
are the most advanced ones.
That is not surprising because of Stochastic Analysis ability
to model uncertainty in natural and social sciences.
The classical approach to integration of a single function having
unbounded variation received less attention partly because there were 
no immediate applications, at least in other areas of analysis.
We believe that further developments of the classical approach
to integration have potentials to shed a new light on modelling
problems of a real life.

\section{Additive and multiplicative functions}\label{add&mult}

It is well-known that the Lebesgue-Stieltjes integral and the 
Riemann-Stieltjes integral are both additive, but in a different
sense.
Namely, let $\nu$ be a function defined on the class of
all subintervals, closed or open at either end, of an interval $[a,b]$.
The function $\nu$ is additive over disjoint subintervals if 
$\nu (A\cup B)=\nu (A)+\nu (B)$ for any disjoint subintervals 
$A, B$ of $[a,b]$ such that $A\cup B$ is a subinterval.
In particular, for any $a\leq s < t < r\leq b$, we have
$$
\nu ([s,r])=\nu ([s,t))+\nu ([t,r])=\nu ([s,t])+\nu ((t,r]).
$$
This is the additivity property enjoyed by the Lebesgue-Stieltjes
integral.
On the other hand, the Riemann-Stieltjes integral may be considered
as a function $\mu$ defined on the simplex
$$
S[a,b]:=\big\{(s,t)\colon\,s,t\in [a,b],\quad s\leq t\big\}
$$ 
of an interval $[a,b]$, and with this notation the integral possess the 
additivity property:  for any $a\leq s\leq t\leq r\leq b$,
$$
\mu (s,r)=\mu (s,t)+\mu (t,r).
$$
As compared to the additivity over disjoint subintervals, the latter 
property amounts to the additivity over adjacent closed subintervals.
The latter property also holds for all extended Stieltjes type integrals
discussed in this paper, and so the resulting integration theories differ 
from the theory of integration with respect to a measure.
Notice that an evolution $U$ defined by (\ref{B-evolution}),
is just a multiplicative function when it is viewed as a function 
defined on the simplex $S[a,b]$.

In this paper we consider functions defined on the simplex of a lineary 
ordered set, which are additive or multiplicative.
This approach provides a concise formulation of various properties
related to a sequence $\lambda$ of partitions of an interval $[a,b]$.
For example, conditions (\ref{qv-var1}) and (\ref{qv-var2}) defining 
the quadratic $\lambda$-variation will become equivalent to a single 
property of an additive upper continuous function defined on the simplex 
of an extended interval to be defined next.

\paragraph*{Regulated functions and extended intervals.}
Let $J$ be an interval of a real line, which may be bounded or
unbounded, and open or closed at either end.
A real-valued function $f$ defined on $J$ is called {\em regulated on $J$}
\index{regulated}
if the right limit $f(x+):=\lim_{y\downarrow x}f(y)$ exists and
is finite for $x\in\bar{J}$ not equal to the right endpoint of $J$,
and if the left limit $f(x-):=\lim_{y\uparrow x}f(y)$ exists and
is finite for $x\in\bar{J}$ not equal to the left endpoint of $J$.
See Section 2 in \cite[Part III]{DNa} for more information
about regulated functions.
The set of all regulated function on $J$ will be denoted by 
${\cal R}(J)$.

For a regulated function $f$ on $J$, define a function $\Delta^{+}f$
on $J$ by $(\Delta^{+}f)(x):=f(x+)-f(x)$ for all $x\in J$ except at the
right endpoint, where it is defined to be zero.
Similarly, define a function $\Delta^{-}f$ on $J$ by $(\Delta^{-}f)(x)
:=f(x)-f(x-)$ for all $x\in J$ except at the left endpoint,
where it is defined to be zero.
Also, let $\Delta^{\pm}f:=\Delta^{+}f+\Delta^{-}f$ on $J$.

With each point $x$ of a real line one can adjoin two symbols
$x-$ and $x+$.
Then a set formed by an interval with adjoined symbols extends
naturally the domain of a regulated function.
More formally, given a closed subinterval $[a,b]$ of an interval 
$J\subset \RR$,  define the {\em extended interval} $\lei a,b\rei$ by
\beq\label{ext-int}
\lei a,b\rei :=
\big (\{a\}\times\{0,+\}\big )\cup\big ((a,b)\times\{-,0,+\}\big )
\cup\big (\{b\}\times \{-,0\}\big )
\eeq
if $a<b$, and $\lei a,b\rei :=\emptyset$ if $a=b$.
For $x\in [a,b]$, we identify $(x,-)$ with $x-$,  
$(x,0)$ with $x$ and $(x,+)$ with $x+$.
With this identification, we have that an ordinary interval $[a,b]$
is a subset of the extended interval $\lei a,b\rei$.
Letting $z+<x- <x <x+ <y-$ for all $z<x<y$ with $z, y\in [a,b]$, 
gives a linear ordering of points of $\lei a,b\rei$.
If $f$ is a regulated function on $[a,b]$ then its natural extension on
$\lei a,b\rei$ again will be denoted by $f$.
On $\lei a,b\rei$ we have the topology induced by the linear ordering
with a base given by the class of all sets $\{u\colon\,v<u<w\}$,
$\{u\colon\,a\leq u<w\}$ and $\{u\colon\,v<u\leq b\}$
for $v, w\in \lei a,b\rei$.
The extended interval $\lei a,b\rei$ is compact, not metrizable and is 
totally disconnected.
The set ${\cal R}[a,b]$ can be identified with the set of all continuous
real-valued functions on $\lei a,b\rei$.
This identification also has a Banach algebra aspect
(see Theorem 5 of Berberian \cite{SKB78}).

Let $[a,b]$ be a closed interval with endpoints $a<b$.
The {\em simplex of the extended interval} $\lei a,b\rei$ is defined by
$$
S\lei a,b\rei:=\big\{(u,v)\colon\, u,v\in\lei a,b\rei,\quad
u\leq v\big\}.
$$ 
In view of the above mentioned identification of $x\in [a,b]$ with
$(x,0)$ in (\ref{ext-int}), we can and do assume that the simplex 
$S[a,b]$ of the ordinary interval $[a,b]$ is a subset of the simplex 
$S\lei a,b\rei $ of the extended interval $\lei a,b\rei$.
For $(u,v)\in S\lei a,b\rei $, let
$$
\lei u,v\rei :=\left\{ \begin{array}{ll}
\big\{w\in\lei a,b\rei\colon\,u\leq w\leq v\big\} &\mbox{if $u\not =v$,}\\
\emptyset &\mbox{if $u=v$.}
\end{array}\right.
$$
This definition agree with (\ref{ext-int}) if $u,v\in [a,b]$ are
ordinary points.
In that case, that is, for ordinary points $(s, t)\in S[a,b]$, 
$\lei s,t\rei$ may be called the extended closed interval.
Also, for ordinary points $(s, t)\in S[a,b]$, $s\not =t$, the extended 
open interval is defined by $\lei s+,t-\rei$, the extended left-open 
interval is defined by $\lei s+,t\rei$, and the extended right-open 
interval is defined by $\lei s,t-\rei$.

\paragraph*{Functions on $S\lei a,b\rei$.}
\index{function on $S(\lei a,b\rei$}
Let $[a,b]$ be a closed interval with endpoints $a<b$.
Any function $\mu$ defined on the simplex $S\lei a,b\rei $ of the 
extended interval $\lei a,b\rei$ is called a function on $S\lei a,b\rei $.
Next we define additive, multiplicative and upper continuous
functions on $S\lei a,b\rei $.

\begin{defn}
{\rm A function $\mu$ on $S\lei a,b\rei $ is called \emph{additive} if 
$$
\left\{ \begin{array}{ll}
\mu (u,v)=\mu (u,w)+\mu (w,v) &\mbox{for each $u, w, v\in \lei a,b\rei$
such that $u\leq w\leq v$,}\\
\mu (u,u) =0 &\mbox{for each $u\in \lei a,b\rei$.} 
\end{array}\right.
$$}
\index{additive}
\end{defn}

 Here the second condition follows from the first one.
To give an example of an additive function on $S\lei a,b\rei $,
let $f$ be a regulated function on $[a,b]$.
Define the function $\mu (f)$ on $S\lei a,b\rei $ by 
\beq\label{muf}
\mu (f;u,v):=f(v)-f(u),\qquad\mbox{for}\quad (u,v) \in S\lei a,b\rei.
\eeq
Clearly,  $\mu (f)$ is additive function on $S\lei a,b\rei$
such that $\mu (f;t-,t)=\Delta^{-}f(t)$ for $t\in (a,b]$
and $\mu (f;s,s+ )=\Delta^{+}f(s)$ for $s\in [a,b)$.
Other examples of additive functions on $S\lei a,b\rei $
will be given later on.
As we already mentioned, the quadratic $\lambda$-variation is defined by
means of a suitable additive function on $S\lei a,b\rei $
(see Definition \ref{qv} below).

\begin{defn}
{\rm A function $\mu$ on $S\lei a,b\rei $ is called \emph{multiplicative} if 
$$
\left\{ \begin{array}{ll}
\mu (u,v)=\mu (u,w)\mu (w,v) &\mbox{for each $u, w, v\in \lei a,b\rei$
such that $u\leq w\leq v$,}\\
\mu (u,u) =1 &\mbox{for each $u\in \lei a,b\rei$.} 
\end{array}\right.
$$}
\index{multiplicative}
\end{defn}

Here the second condition follows from the first one if 
$\mu (a,u)\not =0$ for each $u\in \lei a,b\rei$.
To give an example of a multiplicative function on $S\lei a,b\rei $,
let $f$ be a regulated function on $[a,b]$ such that $1/f$ is
bounded on $[a,b]$.
Define the function $\pi (f)$ on $S\lei a,b\rei $ by 
\beq\label{pif}
\pi (f;u,v):=f(v)/f(u),\qquad\mbox{for}\quad (u,v)\in S\lei a,b\rei.
\eeq
Clearly, $\pi (f)$ is multiplicative function on $S\lei a,b\rei $.

For extended intervals as for other sets, $\lei u_k,v_k\rei\uparrow \lei 
u,v\rei$ means that $u_1\geq u_2\geq\cdots$, $v_1\leq v_2\leq\cdots$ and 
$\cup_{k=1}^{\infty}\lei u_k,v_k\rei=\lei u,v\rei$, while 
$\lei u_k,v_k\rei\downarrow \lei u,v\rei$ means that $u_1\leq u_2\leq\cdots$, 
$v_1\geq v_2\geq\cdots$ and 
$\cap_{k=1}^{\infty}\lei u_k,v_k\rei=\lei u,v\rei$.

\begin{defn}
{\rm A function $\mu$ on $S\lei a,b\rei $ is called {\it upper 
continuous} if $\mu (u_k,v_k)\to\mu (u,v)$ for any sequence
$(u,v), (u_1,v_1), (u_2,v_2 ),\dots\in S\lei a,b\rei $ such that 
$\lei u_k,v_k\rei\downarrow \lei u,v\rei$.
If the same is true whenever $u=v$ then
a function $\mu$ on $S\lei a,b\rei $ is called {\it upper continuous
at $\emptyset$} (empty set).}
\index{upper continuous}
\end{defn}

For a function $\mu$ on $S\lei a,b\rei$, define the right
distribution function $R_{\mu}$ and the left distribution function
$L_{\mu}$ on $\lei a,b\rei$ respectively by
\beq\label{RandLdf}
R_{\mu}(u):=\mu (a,u),\quad u\in\lei a,b\rei
\quad\mbox{and}\quad
L_{\mu}(u):=\mu (u,b),\quad u\in\lei a,b\rei.
\eeq
By definition of the interval topology, a real-valued function
$f$ on $\lei a,b\rei$ is {\em continuous at} $u\in\lei a,b\rei$, which 
is not an endpoint, if for each $\epsilon >0$ there are $u_1, 
u_2\in\lei a,b\rei$ such that $u_1< u<u_2$ and 
$|f(u)-f(v)|<\epsilon$ for each $v\in\lei a,b\rei$, $u_1< v<u_2$.
Continuity at an endpoint means the one-sided variant of the continuity
just formulated.
A function $f$ is {\em continuous on} $\lei a,b\rei$ if it is
continuous at each $u\in\lei a,b\rei$.
Any function $f$ defined on $\lei a,b\rei$ is
continuous at an ordinary point $t\in [a,b]\subset\lei a,b\rei$ 
because always one can take $u_1:=t-$ if $t>a$ and $u_2:=t+$ if $t<b$.
Thus the continuity on $\lei a,b\rei$ is a restriction of a function
on its behaviour on the left side of each point $t-$ and on the
right side of each point $t+$.  

Let $\mu$ be a function on an extended interval $\lei a,b\rei$.
We say that $\mu$ is \emph{nondegenerate} if $\mu (a,b)\not =0$, 
and $\mu$ is {\em bounded} if $\|\mu\|_{\infty}:=\sup\{|\mu (u,v)|\colon\,
(u,v)\in S\lei a,b\rei\}<\infty$.

\begin{thm}\label{interv}
Let a function $\mu$ on $S\lei a,b\rei $ be either additive or
multiplicative, nondegenerate and bounded.
Then the following statements are equivalent{\rm :}
\begin{enumerate}
\item[$(i)$] $\mu$ is upper continuous{\rm ;}
\item[$(ii)$] $\mu$ is upper continuous at $\emptyset${\rm ;}
\item[$(iii)$] $R_{\mu}$ is continuous on $\lei a,b\rei${\rm ;}
\item[$(iv)$] when restricted to $[a,b]$, $R_{\mu}$ is regulated on 
$[a,b]$, $\lim_{s\uparrow t}R_{\mu}(s)=\mu (a,t-)$ for $t\in (a,b]$ and
$\lim_{s\downarrow t}R_{\mu}(s)=\mu (a,t+)$ for $t\in [a,b)${\rm ;}
\item[$(v)$] $L_{\mu}$ is continuous on $\lei a,b\rei${\rm ;}
\item[$(vi)$] when restricted to $[a,b]$, $L_{\mu}$ is regulated on 
$[a,b]$, $\lim_{s\uparrow t}L_{\mu}(s)=\mu (t-,b)$ for $t\in (a,b]$ and
$\lim_{s\downarrow t}L_{\mu}(s)=\mu (t+,b)$ for $t\in [a,b)${\rm ;}
\item[$(vii)$] $\mu (s_k+,t_k-)\to\mu (u,u)$ whenever extended open 
intervals $\lei s_k+,t_k-\rei\downarrow\lei u,u\rei$ for some
$u\in\lei a,b\rei$, and 
\beq\label{singletons}
\forall\,\epsilon >0,\qquad\left\{\begin{array}{ll}
\mbox{card $\{t\in (a,b]$}\colon\,|\mu (t-,t)-\mu (t,t)|
>\epsilon\} <\infty, \\
\mbox{card $\{t\in [a,b)$}\colon\,|\mu (t,t+)-
\mu (t,t)|>\epsilon\} <\infty.
\end{array}
\right.
\eeq
\end{enumerate}
\end{thm}

\begin{proof}
$(i)\Leftrightarrow (ii)$.
Clearly, $(i)$ implies $(ii)$.
For the converse implication, let extended intervals 
$A_k:=\lei u_k,v_k\rei\downarrow \lei u,v\rei =:A$.
Then $A_k=B_k\cup A\cup C_k$ for extended intervals $B_k$, $C_k$ 
such that $B_k\downarrow\lei u,u\rei$ and
$C_k\downarrow\lei v,v\rei$.
Thus $\mu (u_k,v_k)\to\mu (u,v)$ by additivity or multiplicativity.

$(ii)\Rightarrow (iii)$.
Let $u\in\lei a,b\rei$.
Due to above remarks following the definition of continuity on extended
interval and by symmetry, one can assume that $u=t-$ for some $t\in (a,b]$.
Let $\epsilon >0$.
Since $\lei v,t-\rei\downarrow\lei t-,t-\rei$ as $v\uparrow t-$,
by the assumption there exists $a\leq u_1<t-$ such that 
$|\mu (v,t-)-\mu (t-,t-)|<\epsilon$ for each $u_1< v<t-$.
Otherwise we could extract a sequence 
$\lei v_k,t-\rei\downarrow\lei t-,t-\rei$
such that $\mu (v_k,t-)\not\to\mu (t-,t-)$.
Let $u_2:=t$.
If $\mu$ is additive then
$
|R_{\mu}(t-)-R_{\mu}(v)|=|\mu (v,t-)|<\epsilon
$
for each $u_1< v<u_2$.
If $\mu$ is multiplicative then
$$
|R_{\mu}(t-)-R_{\mu}(v)|=|\mu (a,v)|
|\mu (v,t-)-1|<\epsilon\|\mu\|_{\infty}
$$
for each $u_1<v<u_2$.
Since $\epsilon >0$ is arbitrary, $R_{\mu}$ is continuous at $t-$,
and so $R_{\mu}$ is continuous on $\lei a,b\rei$.

$(iii)\Rightarrow (iv)$.
For $t\in (a,b]$, since $R_{\mu}$ is continuous at $t-$ and
$R_{\mu}(t-)=\mu (a,t-)$, it follows that $\lim_{s\uparrow t}
R_{\mu}(s)=\mu (a,t-)$.
Similarly it follows that $\lim_{s\downarrow t}
R_{\mu}(s)=\mu (a,t+)$ for $t\in [a,b)$.
Thus $R_{\mu}$ is regulated when restricted to $[a,b]$.

$(iv)\Rightarrow (vii)$.
Let $\bar R_{\mu}$ be the restriction of $R_{\mu}$ to $[a,b]$.
Since $\bar R_{\mu}$ is regulated, $\bar R_{\mu}(t-)$ is defined and equals
to $\mu (a,t-)$ for each $t\in (a,b]$ by the assumption.
Similarly, $\bar R_{\mu}(t+)=\mu(a,t+)$ for each $t\in [a,b)$.
First suppose that $\mu$ is additive.
Since $\mu (t-,t)=\Delta^{-}\bar R_{\mu}(t)$, $\mu (t,t+)=
\Delta^{+}\bar R_{\mu}(t)$ and $\mu (t,t)=0$,  (\ref{singletons}) holds 
by the analogous property for regulated functions.
Let open extended intervals $\lei s_k+,t_k-\rei\downarrow\emptyset$ as
$k\to\infty$, where each $(s_k,t_k)\in S[a,b]$.
Then either for some $t\in (a,b]$ and all sufficiently large $k$,
$t_k=t$ and $s_k\uparrow t$, or  for some $s\in [a,b)$ and all sufficiently 
large $k$, $s_k=s$ and $t_k\downarrow s$.
Since $\mu$ is additive and $\bar R_{\mu}$ is regulated, 
$\mu (s_k+,t_k-)=\bar R_{\mu}(t-)-\bar R_{\mu}(s_k+)\to 0$ as 
$k\to\infty$ in the first case,
and $\mu (s_k,t_k-)=\bar R_{\mu}(t_k-)-\bar R_{\mu}(s+)\to 0$ as 
$k\to\infty$ in the second case.
Thus the statement $(vii)$ holds when $\mu$ is additive.
Now suppose that $\mu$ is multiplicative, and hence nondegenerate
and bounded.
For $t\in (a,b]$, we have 
\beq\label{1singletons}
\bar R_{\mu}(t)-\bar R_{\mu}(t-)=\bar R_{\mu}(t-)\big\{\mu (t-,t)-1\big\}.
\eeq
We claim that $C:=\inf_t|\bar R_{\mu}(t-)|>0$.
Suppose not.
Then there exists an infinite sequence $s_j\to s\in [a,b]$ such that
either $\bar R_{\mu}(s_j-)\to \bar R_{\mu}(s-)=0$, or
$\bar R_{\mu}(s_j-)\to \bar R_{\mu}(s+)=0$.
In either case due to multiplicativity of $\mu$, it follows
that $\mu (a,b)=0$.
The contradiction proves the claim.
Then by (\ref{1singletons}), $|\Delta^{-}\bar R_{\mu}(t)|\geq
C|\mu (t-,t)-1|$ for $t\in (a,b]$.
Similarly, one can conclude that
$|\Delta^{+}\bar R_{\mu}(t)|\geq C|\mu (t,t+)-1|$ for $t\in [a,b)$
and some $C>0$.
Thus (\ref{singletons}) holds by the
analogous property for regulated functions.
For $(s_k,t_k)\in S[a,b]$, since 
$$
\bar R_{\mu}(t_k-)-\bar R_{\mu}(s_k+)=\bar R_{\mu}(s_k+)\big\{
\mu (s_k+,t_k-)-1\big\}, 
$$
one can concluded as in the additive case
that $\mu (s_k+,t_k-)\to 1$ if $(s_k,t_k)\downarrow\emptyset$.
Thus the statement $(vii)$ holds when $\mu$ is multiplicative.

$(vii)\Rightarrow (ii)$.
Let extended intervals $\lei u_k,v_k\rei\downarrow\lei w,w\rei$
for some $w\in\lei a,b\rei$.
Then either for some $s\in (a,b]$, $v_k=s-=w$ for all sufficiently 
large $k$, or for some $t\in [a,b)$, $u_k=t+=w$ for all sufficiently
large $k$.
Using additivity if $\mu$ is additive, or continuity of the
multiplication if $\mu$ is multiplicative, in each of the two cases
$\mu (u_k,v_k)\to\mu (w,w)$ follows by the statement $(vii)$.

The proof of implications $(ii)\Rightarrow (v)\Rightarrow (vi)\Rightarrow 
(vii)\Rightarrow (ii)$ is based on symmetric arguments which we omit.
The proof of Theorem \ref{interv} is complete.
\qed\end{proof}

Let $f$ be a regulated function on $[a,b]$.
Since the right distribution function $R_{\mu (f)}=f-f(a)$ 
of the additive function $\mu (f)$ on $\lei a,b\rei$ defined by
(\ref{muf}), satisfies statement $(iv)$ of the preceding theorem, 
the function $\mu (f)$ is upper continuous.
Similarly, the multiplicative function $\pi (f)$ defined by (\ref{pif}) is
upper continuous, if in addition, $f$ is bounded away from zero.
In the sense of the following two statements, these two examples are
canonical.
Let ${\cal AC}\lei a,b\rei $ be the set of all additive upper 
continuous functions on $S\lei a,b\rei $, and let
$R_{\mu,[a,b]}$ denotes the right distribution function $R_{\mu}$
restricted to $[a,b]$.

\begin{cor}\label{projad}
Let $a<b$.
Then one-to-one linear operators between the vector spaces
${\cal AC}\lei a,b\rei $ and $\{f\in {\cal R}[a,b]\colon\,
f(a)=0\}$ are given by $f=R_{\mu,[a,b]}$ and $\mu=\mu (f)$.
\end{cor} 

\begin{proof}
For $\mu\in {\cal AC}\lei a,b\rei $, let $f:=R_{\mu,[a,b]}$.
By statement $(iv)$ of Theorem \ref{interv}, $f\in {\cal R}[a,b]$
and $f(a)=\mu (a,a)=0$.
Conversely, for $f\in {\cal R}[a,b]$ such that $f(a)=0$, let
$\mu :=\mu (f)$ be defined by (\ref{muf}).
Since $R_{\mu (f),[a,b]}=f$, $\mu$ is upper continuous by the
implication $(iv)\Rightarrow (i)$ of Theorem \ref{interv},
and so $\mu\in {\cal AC}\lei a,b\rei $.
Clearly, the two mappings are linear and one-to-ne, proving the 
conclusion.
\qed\end{proof}

For a function $f$ on $[a,b]$, we write $|f|\gg 0$ if there is
a constant $C>0$ such that $|f(t)|\geq C$ for all $t\in [a,b]$.
Let ${\cal MC}\lei a,b\rei $ be the set of all nondegenerate,
bounded, multiplicative and upper continuous functions defined
on $S\lei a,b\rei $.

\begin{cor}\label{projmult}
Let $a<b$.
Then one-to-one mappings between the sets ${\cal MC}\lei a,b\rei $
and $\{f\in {\cal R}[a,b]\colon\, f(a)=1$ {\mbox and} $|f|\gg 0\}$ are
given by $f=R_{\mu,[a,b]}$ and $\mu =\pi (f)$.
\end{cor}

\begin{proof}
For $\mu\in {\cal MC}\lei a,b\rei $, let $f:=R_{\mu,[a,b]}$.
By statement $(iv)$ of Theorem \ref{interv}, $f\in {\cal R}[a,b]$
and $f(a)=\mu (a,a)=1$.
Letting $C:=\inf\{|\mu (a,t)|\colon\,t\in [a,b]\}$, we have that
$C>0$.
If not, due to compactness of $[a,b]$, there is $u\in\lei a,b\rei$
such that $\mu (a,u)=0$.
Then $\mu (a,b)=\mu (a,u)\mu (u,b)=0$, a contradiction proving
that $|f(t)|=|\mu (a,t)|\geq C>0$ for each $t\in [a,b]$.
Conversely, for $f\in {\cal R}[a,b]$ such that $f(a)=1$
and $|f|\gg 0$, let $\mu :=\pi (f)$ be defined by (\ref{pif}).
Then $\mu$ is nondegenerate, bounded and multiplicative function
defined on $S\lei a,b\rei $.
Since $R_{\mu (f),[a,b]}=f$, $\mu$ is upper continuous by the
implication $(iv)\Rightarrow (i)$ of Theorem \ref{interv},
and so $\mu\in {\cal MC}\lei a,b\rei $, proving the conclusion.
\qed\end{proof}

\section{The Wiener class}\label{Wiener}

Recall that $\Xi [a,b]$ denotes the set of all partitions of a closed
interval of real numbers $[a,b]$. 
For $\kappa,\lambda\in \Xi[a,b]$, the relation $\kappa\subset\lambda$
means that the partition $\lambda$ is a refinement of the partition $\kappa$.
Let $\prtn$ be the family of all sets $P(\kappa):=\{\lambda\in
\Xi [a,b]\colon\,\lambda\supset\kappa\}$, $\kappa\in \Xi [a,b]$.
Then $\prtn$ is the direction in the sense of Definition
\ref{convergence}.\ref{direction},
which we call the {\em direction of partitions} of $[a,b]$.
Recall also that ${\cal W}_p={\cal W}_p[a,b]$, $0<p<\infty$,
denotes the set of all real-valued functions $f$ on $[a,b]$ having bounded
$p$-variation, that is, (\ref{-2intr}) holds.
For a function $f$ on $[a,b]$ and for $0<p<\infty$, let
$$
v_p^{\ast}(f)=v_p^{\ast}(f;[a,b]):=\limsup_{\kappa,\prtn}s_p(f;\kappa)
:=\inf_{\kappa}\sup\{s_p(f;\lambda)\colon\,\lambda\in \Xi [a,b],\quad
\lambda\supset\kappa\},
$$
$$
\sigma_p^{\ast}(f)=\sigma_p^{\ast}(f;[a,b])
:=\liminf_{\kappa,\prtn}s_p(f;\kappa)
:=\sup_{\kappa}\inf\{s_p(f;\lambda)\colon\,\lambda\in \Xi [a,b],\quad
\lambda\supset\kappa\}.
$$
Notation $s_p(f;\kappa)$ is defined by (\ref{s_p}).
For a regulated function $f$ on $[a,b]$ and for $0<p<\infty$, let
\beq\label{wiener0}
\sigma_p(f)=\sigma_p(f;[a,b]):=\sum_{[a,b]}\big |\Delta f\big |^p:=
\lim_{\sigma,\fgotik}\sum_{\sigma}\big |\Delta f\big |^p,
\eeq
where the limit is defined in Definition \ref{convergence}.\ref{unconditional},
and it always exists finite or infinite.
Here and below, for a finite set $\sigma\subset [a,b]$,
$$
\sum_{\sigma}\big |\Delta f\big |^p:=
\sum_{x\in\sigma}\Big\{|\Delta_a^{-}f(x)|^p+
|\Delta_b^{+}f(x)|^p\Big\}, 
$$
where $\Delta_a^{-}f(x)=f(x)-f(x-)$ if $x\in (a,b]$,
$\Delta_a^{-}f(a):=0$, and $\Delta_b^{+}f(x)=f(x+)-f(x)$ if $x\in [a,b)$,
$\Delta_b^{+}f(b):=0$.
Then the following relations hold for the four quantities just defined.

\begin{lem}\label{wiener}
Let $f$ be a regulated function on $[a,b]$, and let $0<p<\infty$.
Then
\beq\label{wiener3}
\sigma_p(f)\leq \sigma_p^{\ast}(f)\leq v_p^{\ast}(f)\leq v_p(f).
\eeq
\end{lem}

{\proof}
Suppose that $\sigma_p(f)<+\infty$.
Then for each finite $\sigma\subset [a,b]$ and each $\epsilon >0$,
there exists a partition $\kappa$ of $[a,b]$ such that
$\sigma\subset\kappa$ and for each refinement $\lambda$ of $\kappa$,
$s_p(f;\lambda)\geq \sum_{\sigma}|\Delta f|^p-\epsilon$,
and so $\sigma_p(f)\leq\sigma_p^{\ast}(f)$.
If $\sigma_p(f)=+\infty$, then similarly it follows that
$\sigma_p^{\ast}(f)=+\infty$.
Therefore the first inequality in (\ref{wiener3}) holds.
To prove the second inequality, one can assume that $v_p^{\ast}(f)$
is finite.
For each $\kappa\in \Xi [a,b]$, let 
$$
I(\kappa):=\inf\{s_p(f;\lambda)\colon\,\lambda\in \Xi [a,b],\,
\lambda\supset\kappa\}
\quad\mbox{and}\quad
S(\kappa):=\sup\{s_p(f;\lambda)\colon\,\lambda\in \Xi [a,b],\,
\lambda\supset\kappa\}.
$$
Suppose that $v_p^{\ast}(f)=\inf_{\kappa}S(\kappa) <
\sigma_p^{\ast}(f)=\sup_{\kappa}I(\kappa)$.
Then there exist $\kappa_1, \kappa_2\in \Xi [a,b]$ such that
$S(\kappa_1)< I(\kappa_2)$.
Let $\kappa_3:=\kappa_1\cup\kappa_2$.
Since $\kappa_3\supset \kappa_1$ and $\kappa_3\supset\kappa_2$,
it follows that 
$$
S(\kappa_3)\leq S(\kappa_1)< I(\kappa_2)\leq I(\kappa_3).
$$
This contradiction proves the second inequality in (\ref{wiener3}).
Since the third inequality is obvious, the proof of Lemma \ref{wiener} 
is complete.
\qed

By Proposition 2.12 in \cite[Part II]{DNa},
for a regulated function $f$ on $[a,b]$ and for $0<p<1$, we have
$$
\sigma_p(f)= \sigma_p^{\ast}(f)= v_p^{\ast}(f)=v_p(f).
$$
In fact, each function $f$ which has the $p$-variation on $[a,b]$ bounded 
for some $0<p<1$, is a pure jump function on $[a,b]$
(see Theorem 2.11 in \cite[Part II]{DNa}).
Recall that a regulated function $f$ on $[a,b]$ is called a {\it pure jump} 
function, if for each $a<x\leq b$, the directed function
$\sum_{\mu\cap (a,x)}\{\Delta^{-}f+\Delta^{+}f\}$ has a limit
$\sum_{(a,x)}\{\Delta^{-}f+\Delta^{+}f\}$, which satisfies the relation
$$
f(x)=f(a)+\Delta^{+}f(a)+\sum_{(a,x)}\{\Delta^{-}f+\Delta^{+}f\} 
+\Delta^{-}f(x).
$$
In the case $p=1$, by Proposition 2.13 in
\cite[Part II]{DNa}, $v_1^{\ast}(f)=v_1(f)$ for any function $f$ on $[a,b]$,
while $\sigma_1(f)=v_1^{\ast}(f)$ if and only if $f$ is a pure jump 
function.
The situation again becomes different when $1<p<\infty$.
By Theorem \ref{liminf} to be proved next, for a regulated function $f$,
$\sigma_p^{\ast}(f)=\sigma_p(f)$ provided $1<p<\infty$.

It appears that functions $f\in {\cal W}_p$, $1<p<\infty$, such that
$\sigma_p^{\ast}(f)= v_p^{\ast}(f)$, have especially nice properties.
Define the \emph{Wiener class} ${\cal W}_p^{\ast}$ by
$$
{\cal W}_p^{\ast}:=
{\cal W}_p^{\ast}[a,b]:=\big\{f\in{\cal W}_p[a,b]\colon\,
\sigma_p^{\ast}(f)= v_p^{\ast}(f)\big\}.
$$
By Theorem 11.4 of McShane and Botts \cite[p.\ 55]{MB},
$f\in {\cal W}_p^{\ast}$ if and only if the directed function
$\big (s_p(f;\cdot ),\prtn\big )$ has a limit $\lim_{\kappa,\prtn}
s_p(f;\kappa)$.
Thus recalling (\ref{local-p-var}), functions having the local $p$-variation
constitute the Wiener class ${\cal W}_p^{\ast}$.
This in conjunction with Theorem \ref{liminf}, yields 
that for a regulated function $f$ and $1<p<\infty$, 
$f\in {\cal W}_p^{\ast}$ if and only if
the  directed function $\big (s_p(f;\cdot ),\prtn\big )$ converges and
$$
\lim_{\kappa,\prtn}s_p(f;\kappa)=\sigma_p(f).
$$
Love and Young \cite{LY38} defined the Wiener class ${\cal W}_p^{\ast}$
to be the set of all functions $f\in {\cal W}_p$ for which
$v_p^{\ast}(f)=\sigma_p(f)$.
Due to the following Theorem \ref{liminf}, the original and the present 
definitions of ${\cal W}_p^{\ast}$ give the same class of functions.
Love and Young \cite[\S 12, p.\ 32]{LY38} showed that, for $1<p<\infty$,
$$
\cup\big\{{\cal W}_q\colon\,1\leq q<p\big\}\subset
{\cal W}_p^{\ast}\subset {\cal W}_p.
$$

\begin{thm}\label{liminf}
Let $f$ be a regulated function on $[a,b]$, and let $1<p<\infty$.
Then $\sigma_p^{\ast}(f)=\sigma_p(f)$.
\end{thm}

\begin{proof}
If $\sigma_p(f)=+\infty$ then the conclusion holds by Lemma \ref{wiener}.
Suppose that $\sigma_p(f)<+\infty$.
Let $\epsilon \in (0,1)$.
It is enough to show that for a given partition $\kappa$, there exists
a refinement $\lambda$ of $\kappa$ such that
\beq\label{1liminf}
s_p(f;\lambda)\leq \sigma_p(f)+\epsilon.
\eeq
Let $\kappa$ be a partition of $[a,b]$, and let $K:=(4^{p}/2)(1+\osc (f))$.
Choose a partition $\mu=\{y_j\colon\,j=0,\dots,m\}$ of $[a,b]$
such that $\kappa\subset\mu$ and
\beq\label{2liminf}
\sum_{\mu}|\Delta f|^p\geq \sigma_p(f)-\epsilon/(3K).
\eeq
Then choose a set $\{u_{j-1},v_j\colon\,j=1,\dots,m\}\subset [a,b]$ 
such that $y_{j-1}<u_{j-1}<v_j<y_j$ for $j=1,\dots,m$ and
\beq\label{3liminf}
\Big |\sum_{\mu}|\Delta f|^p-\sum_{j=1}^m\Big\{|f(u_{j-1})-f(y_{j-1})|^p
+|f(y_j)-f(v_j)|^p\Big\}\Big |<\epsilon /3.
\eeq
Also, without loss of generality we can assume that $f$ is continuous
at each of the points  $\{u_{j-1},v_j\colon\,j=1,\dots,m\}$.
We claim that for each $j=1,\dots,m$ there exists a partition 
$\lambda_j$ of $[u_{j-1},v_j]$ such that
\beq\label{5liminf}
s_p(f;\lambda_j)\leq K\sigma_p(f;[u_{j-1},v_j])+\epsilon/(3m).
\eeq
Assuming that this is true, and
letting $\lambda :=\mu\cup \big (\cup_{j=1}^m\lambda_j\big )$, 
it then follows that (\ref{1liminf})
holds by (\ref{2liminf}), (\ref{3liminf}) and (\ref{5liminf}),
as desired.

To prove the claim, let $j\in\{1,\dots,m\}$, and suppose that
$\sigma_j:=\sigma_p(f;[u_{j-1},v_j])>0$.
If $|f(v_j)-f(u_{j-1})|^p\leq K\sigma_j$ then taking $\lambda_j=
\{u_{j-1},v_j\}$, (\ref{5liminf}) holds.
Otherwise, either $f(v_j)-f(u_{j-1}) >(K\sigma_j)^{1/p}$ or
$f(v_j)-f(u_{j-1})<-(K\sigma_j)^{1/p}$.
We construct $\lambda_j$ only in the first case because the second case
is symmetric.
It is enough to construct a partition $\lambda_j^{+}$ of $[u_{j-1},v_j]$
having partition points $t+$ with $t\in (u_{j-1},v_j)$ and such that
\beq\label{claim}
s_p(f;\lambda_j^{+})\leq K\sigma_j=K\sigma_p(f;[u_{j-1},v_j]).
\eeq
Approximating $f(t+)$ by values of $f$ at ordinary points, we get 
a partition $\lambda_j$ such that $s_p(f;\lambda_j)$ is arbitrary
close to $s_p(f;\lambda_j^{+})$, and hence (\ref{5liminf}) follows.
To begin the construction let $d_j:=f(v_j)-f(u_{j-1})$.
Since $d_j\leq \osc (f)<K< d_j^p/\sigma_j$,
we have $0<\sigma_j^{1/(p-1)}<d_j$.
Choose the least integer $M\geq 2$ such that 
\beq\label{1claim}
d_j/M\leq \sigma_j^{1/(p-1)}<d_j.
\eeq
If $f(t)-f(u_{j-1})\leq d_j/M$ for each $t\in (u_{j-1},v_j)$, then
by continuity of $f$ at $v_j$, $d_j\leq d_j/M$,
a contradiction.
Thus, there is a $t\in (u_{j-1},v_j)$ such that $f(t)-f(u_{j-1})>d_j/M$.
Let
$$
t_1:=\inf\big\{t\in (u_{j-1},v_j)\colon\,f(t)-f(u_{j-1})>d_j/M\big\}.
$$
Since $f$ is continuous at $u_{j-1}$, $t_1>u_{j-1}$.
Then we have $f(t_1+)-f(u_{j-1})\geq d_j/M$ and 
$f(t_1-)-f(u_{j-1})\leq d_j/M$.
Thus
\beq\label{2claim}
d_j/M\leq f(t_1+)-f(u_{j-1})\leq\Delta^{\pm}f(t_1)+d_j/M.
\eeq
Since for each $t\in (u_{j-1},v_j)$,
$$
|\Delta^{\pm}f(t)|\leq 2^{p-1}\big [|\Delta^{+}f(t)|^p
+|\Delta^{-}f(t)|^p\big ]^{1/p}\leq (K\sigma_j)^{1/p}<d_j.
$$
we have $0\leq\Delta^{\pm}f (t_1)\leq d_j$.
Suppose that $d_j-d_j/M\leq f(t_1+)-f(u_{j-1})\leq d_j+d_j/M$, 
that is, 
$$
-d_j/M\leq f(v_j)-f(t_1+)=d_j-[f(t_1+)-f(u_{j-1})]\leq d_j/M.
$$
Letting $\lambda_j^{+}:=\{u_{j-1},t_1+,v_j\}$, by (\ref{2claim})
and (\ref{1claim}), we then have
\begin{eqnarray*}
s_p(f;\lambda_j^{+})&=&\big |f(t_1+)-f(u_{j-1})\big |^p+\big |f(v_j)-f(t_1+)
\big |^p\leq 2^{p-1}\big\{|\Delta^{\pm}f(t_1)|^p+2(d_j/M)^p\big\}\\[2mm]
&\leq &4^{p-1}\big\{\sigma_j+\osc (f)(d_j/M)^{p-1}\big\}
\leq K\sigma_p(\mu;[u_{j-1},v_j]),
\end{eqnarray*}
and so (\ref{claim}) holds.
If $M=2$ then the construction stops, proving the claim.
Otherwise, we can assume that $M>2$ and $f(t_1+)-f(u_{j-1})<d_j-d_j/M$.
Suppose that for some integer $2\leq k<M$,
we have constructed $u_{j-1}=:t_0<t_1<\cdots <t_{k-1}<v_j$
such that $f(t_{k-1}+)-f(u_{j-1})<d_j-d_j/M$ and for $i=1,\dots,k-1$,
\beq\label{3claim}
d_j/M\leq f(t_i+)-f(t_{i-1}+)\leq \Delta^{\pm} f(t_i)+d_j/M.
\eeq
Thus $f(v_j)-f(t_{k-1}+)=d_j-[f(t_{k-1}+)-f(u_{j-1})]>d_j/M$.
Let
$$
t_k:=\inf\big\{t\in (t_{k-1},v_j)\colon\,f(t)-f(t_{k-1}+)>d_j/M\},
$$
which exists since $f$ is continuous at $v_j$, and $t_k>t_{k-1}$.
As in (\ref{2claim}), we have that (\ref{3claim})
holds with $i=k$, and so,
$$
f(t_k+)-f(u_{j-1})=\big [f(t_{k-1}+)-f(u_{j-1})\big ]
+\big [f(t_k+)-f(t_{k-1}+)\big ]<d_j+\Delta^{\pm}f(t_k).
$$
Suppose that $f(t_k+)-f(u_{j-1})\geq d_j-d_j/M$, and so, 
$$
-\Delta^{\pm}f(t_k)\leq f(v_j)-f(t_k+)=d_j-\big [f(t_k+)-f(u_{j-1})
\big ]\leq d_j/M.
$$
Letting $\lambda_j^{+}=\{t_i+\colon\,i=0,\dots,k+1\}$ with $t_{k+1}:=v_j$
(recall that $f$ is continuous at $u_{j-1}, v_j$), 
by (\ref{3claim}) and (\ref{5liminf}), we then have
\begin{eqnarray*}
s_p(f;\lambda_j^{+})&=&\sum_{i=1}^{k+1}\big |f(t_i+)-f(t_{i-1}+)\big 
|^p\\[2mm]
&\leq& 2^{p-1}\Big\{\sum_{i=1}^k|\Delta^{\pm}f(t_i)|^p
+|\Delta^{\pm}f(t_k)|^p+(k+1)(d_j/M)^p\Big\}\\[2mm]
&\leq &(4^{p}/2)\big\{\sigma_j+\osc (f)(d_j/M)^{p-1}\big\}
\leq K\sigma_p(f;[u_{j-1},v_j]),
\end{eqnarray*}
so (\ref{claim}) holds.
Otherwise, $f(t_k+)-f(u_{j-1})<d_j-d_j/M$.
On the other hand, by (\ref{3claim}) with $i=1,\dots,k$,
$$
f(t_k+)-f(u_{j-1})=\sum_{i=1}^k[f(t_i+)-f(t_{i-1}+)]
\geq k(d_j/M),
$$
which yields that $k\leq M-1$.
We continue construction which  stops after at most $M$ steps, 
proving the claim for all $j$ such that $\sigma_j>0$.

Let $j\in\{1,\dots,m\}$ be such that $\sigma_j=0$, if such $j$ exists.
If $f(v_j)=f(u_{j-1})$ then (\ref{5liminf}) holds with 
$\lambda_j=\{u_{j-1},v_j\}$.
Otherwise, either $f(v_j)>f(u_{j-1})$ or $f(v_j)<f(u_{j-1})$.
Again, we construct $\lambda_j$ only in the first case because the second 
case is symmetric.
Let $M\geq 2$ be the least integer such that 
$(1/M)^{p-1}\leq \epsilon/(3m\osc (f)^p)$.
Let $t_0:=u_{j-1}$, and for $i=1,\dots,M$, let
$$
t_i:=\inf\{t\in [t_{i-1},v_j]\colon\,f(t)=f(u_{j-1})+(i/M)[f(v_j)
-f(u_{j-1})]\}.
$$ 
Then taking $\lambda_j=\{t_i\colon\,i=0,\dots,M\}$, we have
$$
s_p(f;\lambda_j)=\sum_{i=1}^M\big |f(t_i)-f(t_{i-1})\big |^p
=M^{1-p}|f(v_j)-f(u_{j-1})|^p<\epsilon /(3m),
$$
so (\ref{5liminf}) holds in the case $\sigma_j=0$.
The proof of Theorem \ref{liminf} is now complete.
\qed\end{proof}

The following characterization of the Wiener class is
often used in what follows.
We extend definition (\ref{-2intr}) of the $p$-variation over 
a closed interval $J=[a,b]$ to the case when $J$ is open or closed
at either end, by 
$$
v_p(f;J):=\sup\big\{s_p(f;\kappa)\colon\,\,\kappa=\{
x_i\colon\,i=0,\dots,n\},\,\,x_0, x_n\in J\big\}.
$$

\begin{lem}\label{lemma1}
Let $f$ be a regulated function on $[a,b]$, and let $1<p<\infty$.
Then the following statements are equivalent{\rm :}
$(a)$ $f\in{\cal W}_p^{\ast}${\rm ;}
\begin{enumerate}
\item[$(b)$] for every $\epsilon >0$, there is a partition 
$\{z_j\colon\,j=0,\dots,m\}$ of $[a,b]$  such that
\beq\label{wiener1}
\sum_{j=1}^mv_p(f;(z_{j-1},z_j))<\epsilon;
\eeq
\item[$(c)$] $\sigma_p(f)<\infty$, and for every $\epsilon >0$, there
is a partition $\lambda$ of $[a,b]$ such that
\beq\label{wiener2}
\sum_{i=1}^n|f(x_i-)-f(x_{i-1}+)|^p<\epsilon
\eeq
for each refinement $\{x_i\colon\,i=0,\dots,n\}$ of $\lambda$.
\end{enumerate}
\end{lem}

\begin{proof}
$(a)\Rightarrow (b)$.
Let $f\in{\cal W}_p^{\ast}$ and $\epsilon >0$.
By the definitions of $v_p^{\ast}(f)$ and $\sigma_p(f)$, there exists 
a partition $\kappa=\{z_j\colon\,j=0,\dots,m\}$ of $[a,b]$ such that
$$
\sum_{j=1}^mv_p(f;[z_{j-1},z_j])=\sup\big\{s_p(f;\lambda)\colon\,
\lambda\in \Xi [a,b],\quad \lambda\supset\kappa\big\}
<v_p^{\ast}(f)+\epsilon/2
$$
$$
\mbox{and}\qquad
\sum_{j=1}^m\big\{|\Delta^{+}f(z_{j-1})|^p+|\Delta^{-}f(z_j)|^p
\big\}>\sigma_p(f)-\epsilon/2.
$$
Let $\{u_{j-1},v_j\colon\,j=1,\dots,m\}$ be a set of points in
$(a,b)$ such that $z_{j-1}<u_{j-1}<v_j<z_j$ for $j=1,\dots,m$.
Then we have, for $j=1,\dots,m$,
$$
v_p(f;[u_{j-1},v_j])\leq v_p(f;[z_{j-1},z_j])
-\Big\{v_p(f;[z_{j-1},u_{j-1}])+v_p(f;[v_j,z_j])\Big\}.
$$
For each $j=1,\dots,m$, letting $u_{j-1}\downarrow z_{j-1}$ and
$v_j\uparrow z_j$, by Lemma 2.19 of Dudley and Norvai\v sa 
\cite[Part II]{DNa}, it  follows that
\begin{eqnarray*}
\sum_{j=1}^mv_p(f;(z_{j-1},z_j))\!&\leq&\!\sum_{j=1}^m\Big\{
v_p(f;[z_{j-1},z_j])
-|\Delta^{+}f(z_{j-1})|^p-|\Delta^{-}f(z_j)|^p\Big\}\\
&<&v_p^{\ast}(f)+\epsilon/2-\sigma_p(f)+\epsilon/2=\epsilon.
\end{eqnarray*}
The last equality holds by the definition of the Wiener class
${\cal W}_p^{\ast}$ and by Theorem \ref{liminf}, proving statement $(b)$.

$(b)\Rightarrow (a)$.
Let statement $(b)$ hold.
Then it is easy to see that $f\in{\cal W}_p$.
Assume however that $f\not\in{\cal W}_p^{\ast}$.
Therefore, since (\ref{wiener3}) always holds,
$v_p^{\ast}(f)-\sigma_p^{\ast} (f)\geq C$ for some positive 
constant $C$.
Let $\kappa=\{z_j\colon\,j=0,\dots,m\}\in\Xi [a,b]$.
Then we have
$$
\sum_{j=1}^m\Big\{|\Delta^{+}f(z_{j-1})|^p+v_p(f;(z_{j-1},z_j))
+|\Delta^{-}f(z_j)|^p\Big\}\geq v_p^{\ast}(f).
$$
By Theorem \ref{liminf}, it then follows that
\begin{eqnarray*}
\sum_{j=1}^mv_p(f;(z_{j-1},z_j))&\geq& v_p^{\ast}(f)-\sum_{j=1}^m
\Big \{\Delta^{+}f(z_{j-1})|^p+|\Delta^{-}f(z_j)|^p\Big \}\\
&\geq& v_p^{\ast}(f)-\sigma_p (f)
=v_p^{\ast}(f)-\sigma_p^{\ast} (f)\geq C>0.
\end{eqnarray*}
Since $\kappa$ is arbitrary, (\ref{wiener1}) can't hold for each 
$\epsilon >0$.
This contradiction implies that $f\in{\cal W}_p^{\ast}$.

$(b)\Rightarrow (c)$.
Let statement $(b)$ hold.
Then it is easy to see that $f\in{\cal W}_p$, and hence
$\sigma_p (f)<\infty$ by (\ref{wiener3}).
Let $\epsilon >0$ and let $\lambda :=\{z_j\colon\,j=0,\dots,m\}\in\Xi [a,b]$
be such that (\ref{wiener1}) holds.
It is clear that (\ref{wiener1}) also holds when the partition
$\lambda$ is replaced by any its refinement.
Let $\{x_i\colon\,i=0,\dots,n\}\in\Xi [a,b]$ be a refinement of $\lambda$.
Then we have
\begin{eqnarray*}
\sum_{i=1}^n|f(x_i-)-f(x_{i-1}+)|^p&=&\sum_{i=1}^n\lim_{y_{i-1}\downarrow
x_{i-1},\,y_i\uparrow x_i}|f(y_i)-f(y_{i-1})|^p\\
&\leq&\sum_{i=1}^nv_p(f;(x_{i-1},x_i))<\epsilon.
\end{eqnarray*}
Therefore statement $(c)$ holds.

$(c)\Rightarrow (b)$.
Let $\epsilon >0$ and let $\lambda\in\Xi [a,b]$ be such that 
(\ref{wiener2}) holds for all its refinements.
Since $\sigma_p(f)<\infty$,
there exists a finite set $\nu\subset [a,b]$ such that
$$
\sum_{x\in\mu}\Big\{|\Delta^{-}f(x)|^p+|\Delta^{+}f(x)|^p
\Big \}<\epsilon
$$
whenever a finite set $\mu\subset (a,b)$ and $\mu\cap\nu
=\emptyset$.
Let $\{z_j\colon\,j=0,\dots,m\}:=\lambda\cup\nu\in\Xi [a,b]$
and, for each $j=1,\dots,m$, let $\kappa (j):=\{x_i^j\colon\,
i=0,\dots,n(j)\}$ be a set of points in $(a,b)$ such that
$z_{j-1}<x_0^j<\cdots <x_{n(j)}^j<z_j$.
Then we have
\begin{eqnarray*}
\sum_{j=1}^ms_p(f;\kappa (j))&\leq& 4^{p-1}\sum_{j=1}^m\sum_{i=1}^{n(j)}
\Big \{|\Delta^{+}f(x_{i-1}^j)|^p+|\Delta^{-}f(x_i^j)|^p\Big \}\\
& &+4^{p-1}\sum_{j=1}^m\sum_{i=1}^{n(j)}|f(x_i^j-)-f(x_{i-1}^j+)|^p
<4^p\epsilon /2.
\end{eqnarray*}
Since all partitions $\kappa (j)$ of $(z_{j-1},z_j)$ are arbitrary,
it follows that
$$
\sum_{j=1}^mv_p(f;(z_{j-1},z_j))\leq 4^p\epsilon /2.
$$
This proves statement $(b)$.
The proof of Lemma \ref{lemma1} is complete.
\qed\end{proof}

\section{Extended Riemann-Stieltjes integrals}\label{extended}

The refinement Riemann-Stieltjes and the Henstock-Kurzweil
integrals are the best known extensions of the Riemann-Stieltjes integral.
We refer to Graves \cite{LMG}, Hildebrandt \cite{THH63}
and McLeod \cite{RMM} for detail expositions of the properties of the
three integrals.
Less well known are extensions of the Riemann-Stieltjes
integral suggested by several members of the family of Young.
This section is just a quick glance at relations between
several different extensions of the Riemann-Stieltjes integral
in the light of the $p$-variation conditions originated in the Theorem
on Stieltjes integrability of L.\ C.\ Young \cite[p.\ 264]{LCY36}. 

The Riemann-Stieltjes integral, or the $RS$ integral, of a function
$f$ with respect to a function $g$ is defined to exist and equal the limit
of Riemann-Stieltjes sums as the mesh of partitions tends to zero.
The $RS$ integral exists first, if $f\in {\cal W}_p$, $g\in {\cal W}_q$
with $p^{-1}+q^{-1}>1$, and second if $f$ and $g$ have no a discontinuity
at the same point.
The first condition cannot be replaced in general by the
condition $p^{-1}+q^{-1}=1$, and the second condition is
necessary for the existence of the $RS$ integral.
In this section we recall several extensions of the $RS$
integral which exist when either one of the two conditions,
or both of them, fail to hold for a pair of functions $(f,g)$.
We will say that an integral $I_2$ is an $({\cal A},{\cal B})$
extension of an integral $I_1$, or shortly
\beq\label{extend}
I_1\stackrel{({\cal A},{\cal B})}{\longrightarrow} I_2,
\eeq
if the following two statements hold:
\begin{enumerate}
\item[$(a)$] for a class ${\cal A}$ of pairs of functions, if the $I_1$
integral exists then the $I_2$ integral so does;
\item[$(b)$] for a class ${\cal B}\subset {\cal A}$, the converse
implication holds, and the two integrals have the same value.
\end {enumerate}
Typically, ${\cal A}$ is a class of pairs of regulated 
functions, and ${\cal B}$ is a subclass of the class ${\cal D}_p$ 
of all pairs  $(f,g)$ such that $f\in {\cal W}_p$, $g\in{\cal W}_q^{\ast}$ 
with $p^{-1}+q^{-1}=1$, or vice versa.
A summary of relations (\ref{extend}) between extended Riemann-Stieltjes 
integrals is given at the end of this section.

\begin{notat}\label{partition}
{\rm Let $[a,b]$ be a closed interval of real numbers.
A set $\kappa=\{x_i\colon\,i=0,\dots,n\}$ of points in $[a,b]$ is
called a {\it partition of} $[a,b]$ if $a=x_0<x_1<\dots <x_n=b$.
For $i=1,\dots,n$, a point $y_i\in [x_{i-1},x_i]$ 
attached to a subinterval $[x_{i-1},x_i]$ is called a {\it tag}, 
and the set $\tau =\tau (\kappa)=\{([x_{i-1},x_i],y_i)\colon\,
i=1,\dots,n\}$ is called a {\it tagged partition} of $[a,b]$
associated to $\kappa$.
A tagged partition $\tau=\tau (\kappa)$ is a refinement of a partition
$\lambda$ if $\kappa\supset\lambda$.
Let $\Xi [a,b]$ and $\Theta [a,b]$ be the sets of all partitions and
all tagged partitions of $[a,b]$, respectively.
Given a partition $\lambda$ of $[a,b]$, let $R(\lambda)$ be 
the set of all tagged partitions which are
refinements of $\lambda$.
Let $\PP$ be the family of all sets $R=R(\lambda)$, $\lambda\in\Xi [a,b]$.  
Then $\PP$ is a direction in $\Theta [a,b]$ in the sense of Definition
\ref{convergence}.\ref{direction}.}
\end{notat}

\paragraph*{The refinement Riemann-Stieltjes integral.}
Let $f$ and $g$ be real-valued functions defined on $[a,b]$.
The \emph{Riemann-Stieltjes sum} is a function $S_{RS}=S_{RS}(f,g)$ 
defined on $\Theta [a,b]$ by
$$
S_{RS}(\tau)=S_{RS}(f,g;\tau):=
\sum_{i=1}^nf(y_i)[g(x_i)-g(x_{i-1})]
$$
for a tagged partition $\tau=\{([x_{i-1},x_i],y_i)\colon\,i=1,\dots,n\}
\in\Theta [a,b]$.
Then $(S_{RS},\PP)$ is the directed function.
The {\it refinement Riemann-Stieltjes} integral 
$(RRS)\smallint_a^bf\,dg$ is defined to exist and equal the limit of the
directed function $(S_{RS},\PP)$, that is,
\beq\label{RRSint}
(RRS)\int_a^bf\,dg:=\lim_{\PP,\tau}S_{RS}(f,g;\tau)
\eeq
provided the limit exists.
Therefore $(RRS)\smallint_a^bf\,dg$ is defined to equal $A$
if given $\epsilon >0$ there exists a partition $\lambda$ of $[a,b]$ 
such that
$
\big |S_{RS}(f,g;\tau)-A|<\epsilon
$
for each tagged partition $\tau$ which is a refinement 
of $\lambda$.
It is easy to ascertain that if $(RRS)\smallint_a^bf\,dg$
exists then for each $a\leq x<y\leq b$,
$$
\lim_{z\downarrow x}\big [f(z)-f(x)\big ]\big [g(z)-g(x)\big ]
=\lim_{z\uparrow y}\big [f(z)-f(y)\big ]\big [g(z)-g(y)\big ]=0.
$$
Therefore $f$ and $g$ cannot have common discontinuities on the same
side at the same point provided their refinement Riemann-Stieltjes
integral exists.
This restriction does not allow to integrate an indicator function with respect 
to itself.

\begin{thm}\label{meshRS}
Let $g$ and $f$ be bounded functions on $[a,b]$.
The integral $\smallint_a^bg\,df$ exists as the Riemann-Stieltjes
integral if and only if it exists as the refinement Riemann-Stieltjes
integral and the two functions $g$, $f$ have no common discontinuities 
on $[a,b]$.
\end{thm}

The proof of this theorem is given in \cite[p.\ 264]{LMG}
and in \cite[p.\ 51]{THH63}.
We use extensively the preceding theorem when dealing with chain rules
in the next chapter for functions having the quadratic $\lambda$-variation.

\paragraph*{The refinement Young-Stieltjes integral.}
One way to resolve the problem of common discontinuities for the
refinement Riemann-Stieltjes integral is to modify the Riemann-Stieltjes
sum as was suggested by W.H.\ Young (1914)\footnote{In this work
W.H.\ Young developed the integral equivalent to the Lebesgue-Stieltjes
integral in a way which was latter extended by Daniell (1918)}.
Let $g$ be a regulated function on $[a,b]$.
The \emph{refinement Young-Stieltjes integral}
$(\RYS)\smallint_a^bf\,dg$ is defined to be the limit, if it exists, as in
(\ref{RRSint}) except that the Riemann-Stieltjes sum
$S_{RS}(f,g;\tau)$ is replaced by the {\em Young-Stieltjes sum}
\beq\label{YS}
S_{YS}(f,g;\tau):=
\sum_{i=1}^n\Big\{[f\Delta^{+}g](x_{i-1})+f(y_i)[g(x_i-)-g(x_{i-1}+)]
+[f\Delta^{-}g](x_i)\Big\},
\eeq
and the direction is restricted to all tagged partitions
$\tau =\{([x_{i-1},x_i],y_i)\colon\,i=1,\dots,n\}$ of $[a,b]$ such that
$x_{i-1}<y_i<x_i$ for $i=1,\dots,n$.
Clearly, the refinement Young-Stieltjes integral exists if either $f$, or 
$g$ (or both) is a step function.
Also, the refinement Young-Stieltjes integral exists if the refinement
Riemann-Stieltjes integral so does, and the two integrals have the 
same value.
This is so because any Young-Stieltjes sum (\ref{YS}) can be approximated 
arbitrarily closely by Riemann-Stieltjes sums based on refinements of 
$\{x_i\colon\,i=1,\dots,n\}\in \Xi [a,b]$.
More information about the refinement Young-Stieltjes integral one can find
in \cite[pp.\ 88-95]{THH63}.

\paragraph{The Central Young integral.}
L.\ C.\ Young \cite{LCY36} suggested  another way to bypass the problem of 
common discontinuities for the refinement Riemann-Stieltjes integral.
His integral is defined when both the integrand and integrator
are regulated functions, and for such functions extends
the refinement Young-Stieltjes integral (see Proposition 3.17 in \cite[Part II]{DNa}).
For a regulated function $f$ on $[a,b]$, define its left continuous
modification $f_{-}^{(a)}$ on $[a,b]$ by
$$
f_{-}^{(a)}(x):=f_{-}(x):=f(x-)\quad\mbox{for $a<x\leq b$, and}
\quad f_{-}^{(a)}(a):=f(a).
$$	
Similarly, define the right continuous modification
$f_{+}^{(b)}$ on $[a,b]$ by
$$
f_{+}^{(b)}(x):=f_{+}(y):=f(y+)\quad\mbox{for $a\leq y< b$, and}
\quad f_{-}^{(b)}(b):=f(b).
$$	
Let $g$ and $f$ be regulated functions on $[a,b]$.
Define the $Y_1$ integral on $[a,b]$ by
\beq\label{Y1}
(Y_1)\int_a^bf\,dg :=(RRS)\int_a^bg_{+}^{(b)}\,df_{-}^{(a)}
-\sum_{[a,b)}\Delta^{+}f[g_{+}-g_{-}^{(c)}]+[f\Delta^{-}g](b)
\eeq
provided the $RRS$ integral exists, and the sum converges unconditionally.
Also, define the $Y_2$ integral on $[a,b]$ by
$$
(Y_2)\int_a^bf\,dg:=(RRS)\int_a^bg_{-}^{(a)}\,df_{+}^{(b)}
+[f\Delta^{+}g](a)+\sum_{(a,b]}\Delta^{-}f[g_{+}^{(b)}-g_{-}]
$$
provided the $RRS$ integral exists, and the sum converges unconditionally.
The $Y_1$ and $Y_2$ integrals satisfy standard properties of integrals, 
such as linearity and additivity over adjacent intervals.
Also, the $Y_1$ integral exists if and only if the $Y_2$ integral
so does, and the integrals have the same value.
Therefore for regulated functions $g$ and $f$, 
the \emph{Central Young} integral on $[a,b]$ is defined by
$$
(CY)\int_a^bf\,dg :=(Y_1)\int_a^b f\,dg =(Y_2)\int_a^b f\,dg
$$
provided either the $Y_1$ or the $Y_2$ integral exists on $[a,b]$.
All these facts are proved in \cite[Section 3 of Part II]{DNa} 
(see also \cite{DNb} for simpler proofs).
The integrals $Y_1$ and $Y_2$ defined Dudley \cite{RMD92} to clarify
the original definition of L.\ C.\ Young \cite{LCY36}.
The word ``Central'' in the name of the Central Young integral
is used because there are its two other modifications: 
the Left Young and Right Young integrals discussed in Section
\ref{LYRYintegrals} below.

\paragraph{The Henstock-Kurzweil integral.}
The above two extensions of the refinement Riemann-Stieltjes integral
are defined for regulated functions.
This restriction is not a problem if one wishes to apply these integrals
to stochastic processes, which usually have regulated, and hence
bounded sample functions.
However there is another extension of the refinement Riemann-Stieltjes
integral which may exist when neither the integrand nor the integrator
are regulated, or bounded functions.
Ward \cite{AJW36} working with the Perron-Stieltjes integral,
modified its definition so that the new integral extends
the refinement Riemann-Stieltjes integral and is defined for
functions which need not be regulated.
An equivalent definition of Ward's extension was given indepepndently
by Kurzweil \cite{JK57} and Henstock \cite{RH63}.
The resulting integral is often called the Henstock-Kurzweil integral.
Namely, for real-valued functions $g$ and $f$ on $[a,b]$,
the Henstock-Kurzweil integral $(HK)\smallint_a^bf\,dg$
is defined to exist and equal $A$ if given $\epsilon >0$ there 
exists a positive function $\delta (\cdot)$ on $[a,b]$ such that 
Riemann-Stieltjes sums $S_{RS}(g,f;\tau)$ differ from $A$ 
at most by $\epsilon$ for all tagged partitions $\tau =\{([x_{i-1},x_i],
y_i)\colon\,i=1,\dots,n\}$ such that $y_i-\delta (y_i)\leq x_{i-1}
\leq y_i\leq x_i\leq y_i+\delta (y_i)$ for $i=1,\dots,n$.
Several comparison results between the Hestock-Kurzweil and the 
refinement Young-Stieltjes integral are given in
\cite[Appendix F in Part I]{DNa}.
In particular, the Henstock-Kurzweil integral extends strictly the 
Lebesgue-Stieltjes integral.
Most likely this fact made the Henstock-Kurzweil integral so popular
among those currently working  in Real Analysis.
In this paper the Henstock-Kurzweil integral is used to improve earlier results
related to the Black-Scholes formula in Financial Mathematics (see Theorem
\ref{modelling}.\ref{BandS}).

\paragraph{The L.\ C.\ Young Stieltjes integrability theorem.}
As it is mentioned in the introduction, in terms of the $p$-variation,
the first and best possible existence conditions for the extended 
Riemann-Stieltjes integrals are due to L.\ C.\ Young \cite{LCY36}
in 1936.
Two years latter he published in \cite{LCY38b} a different proof of a more
general result using the $\phi$-variation property of a function, 
an extention of the $p$-variation.
In the first paper L.\ C.\ Young used, what we now call the $Y_1$ integral,
while in the second one he used the refinement Young-Stieltjes integral.
By Theorem F.2 in \cite[Part I]{DNa}, L.\ C.\ Young's Stieltjes integrability 
theorem also applies to the Henstock-Kurzweil integral.

We formulate the L.\ C.\ Young Stieltjes integrability theorem in a form suitable 
for further references in the present notes.
Recall that ${\cal W}_{\infty}[a,b]$ is the class ${\cal R}[a,b]$ of all regulated
functions on $[a,b]$.

\begin{thm}\label{LCY}
Let $1\leq p<\infty$, and $q$ be such that $p^{-1}+q^{-1}>1$ if $p>1$
and $q=\infty$ if $p=1$.
For $f\in {\cal W}_p[a,b]$ and $g\in {\cal W}_q[a,b]$,
the integral $\smallint_a^b f\,dg$ exists
\begin{enumerate}
\item[$(a)$] as the Riemann-Stieltjes integral if $f$ and $g$ have no
discontinuities at the same point{\rm ;}
\item[$(b)$] as the refinement Riemann-Stieltjes integral 
if $f$ and $g$ have no common discontinuities on the same side 
at the same point{\rm ;}
\item[$(c)$] always as the refinement Young-Stieltjes integral, 
as  the Central Young integral and as the Henstock-Kurzweil
integral.
\end{enumerate}
In whichever of the five senses the integral exists, 
there exists a constant $K$ if $p>1$ which depends only on $p$ and $q$ 
such that the inequality
\beq\label{1LCY}
\Big |\int_a^b\big\{g-g(y)\big\}\,df\Big |
\leq \left\{\begin{array}{ll}
\osc\, (g;[a,b])\,v_1(f;[a,b]) &\mbox{if $p=1$,}\\
Kv_q(g;[a,b])\,v_p(f;[a,b]) &\mbox{if $1<p<\infty$,}
\end{array}\right.
\eeq
holds for any $y\in [a,b]$.
\end{thm}

In the case $p=1$ the conclusion is well-known and easy to prove.
In the case $p>1$ the stated theorem with the Central Young integral in 
statement $(c)$, follows from the theorem on Stieltjes integrability of 
L.\ C.\ Young \cite{LCY36}.
Statement $(c)$ for the other two integrals then follows from Corollary 3.20 
of \cite[Part II]{DNa} and from Theorem F.2 of \cite[Part I]{DNa}.
There are other proofs of the same result with a one of the three integrals
in statement $(c)$ and with a slightly different constant $K$ in (\ref{1LCY}), 
in addition to already mentioned work of L.\ C.\ Young \cite{LCY38b}.
Necessary and sufficient Stieltjes integrability conditions in terms of the
$\phi$-variation have been established by D'ya\v ckov \cite{AMD88}. 
These results also can be found \cite{DNb}.

Suppose that a function $g$ is integrable with respect to a function $f$
in the sense of any of the five integrals defined so far,
and suppose that inequality (\ref{1LCY}) holds for some $1\leq p<\infty$.
Let $\Phi (x):=C+\smallint_a^xf\,dg$, $x\in [a,b]$, 
be the corresponding indefinite integral.
Since each of the five integrals is additive over adjacent intervals,
then in the case $p> 1$, we have the bound
$$
s_p(\Phi;\kappa)\leq 2^{p-1}\Big\{K^pv_q(g)^{p/q}
+\|g\|_{\sup}^p\Big\}v_p(f)
$$
for any partition $\kappa$ of $[a,b]$, and so the $p$-variation of $\Phi$
is bounded in this case.
A similar bound holds in the case $p=1$, and so we get:

\begin{cor}\label{indefinite}
Under the assumptions of the preceding theorem, each of the five
indefinite integrals $\Phi$ has bounded $p$-variation.
\end{cor}

\paragraph*{The symmetric Young-Stieltjes integral.}
L.\ C.\ Young \cite{LCY36} showed that the $p$-variation
conditions (\ref{-1intr}) in his Theorem on Stieltjes integrability
cannot be replaced in general by the same condition with
$p^{-1}+q^{-1}=1$.
However, if one restricts to integrals of the form 
\beq\label{special}
\int_a^b\big (\phi{\circ}f\big )\,df 
\eeq
where the integrand is a composition of the integrator
with a smooth function $\phi$, then the conditions
(\ref{-1intr}) can be weakened.
For example, by Proposition 4.4 of \cite{RN99a}, the integral
$(\RYS)\smallint_a^bf\,df$ exists and satisfies the relation
$$
(\RYS)\int_a^bf\,df=\frac{1}{2}\Big\{f(b)^2-f(a)^2+\sum_{(a,b]}
\big\{\Delta^{-}f\big\}^2-\sum_{[a,b)}\big\{\Delta^{+}f\big\}^2
\Big\}
$$
if and only if $f\in {\cal W}_2^{\ast}[a,b]$.
It is interesting that one can modify Young-Stieltjes sums further so that 
the resulting integral, called the symmetric Young-Stieltjes
integral, of the form (\ref{special}) exists for some
functions with infinite $2$-variation.
In particular,  the symmetric Young-Stieltjes integral exists almost surely
for a sample function of Brownian motion with respect to itself.
At the same time the new integral agrees with the refinement 
Young-Stieltjes integral if the conditions (\ref{-1intr})
of L.\ C.\ Young's theorem are satisfied for functions with suitable
values at their discontinuity points.
 
Let $g$ and $f$ be regulated functions on $[a,b]$.
If all tags of a tagged parition $\tau=\{([x_{i-1},x_i],y_i)\colon$
$i=1,\dots,n\}$ of $[a,b]$  satisfy the condition $x_{i-1}<y_i<x_i$ for
$i=1,\dots,n$, then we call $\tau$ a \emph{Young tagged partition}. 
For a Young tagged partition $\tau$, let
$$
C_{\YS}(g,f;\tau):=\sum_{i=1}^n\Big\{\big [\Delta^{+}g\Delta^{+}f\big ]
(x_{i-1})-\big [\Delta g\Delta f\big ](x_{i-1},y_i]
+\big [\Delta g\Delta f\big ][y_i,x_i)
-\big [\Delta^{-}g\Delta^{-}f\big ](x_i)\Big\},
$$
and let
\begin{eqnarray}\label{SYSsum}
&&T_{\YS}(g,f;\tau):=S_{\YS}(g,f;\tau)+\frac{1}{2}C_{\YS}(g,f;\tau)
=\frac{1}{2}\sum_{i=1}^n\Big\{\big [(g_{+}+g)\Delta^{+}f\big ](x_{i-1})\\
&&\qquad +\big [g_{+}(x_{i-1})+g(y_i)\big ]
\Delta f(x_{i-1},y_i]+\big [g(y_i)+g_{-}(x_i)\big ]
\Delta f[y_i,x_i)+\big [(g+g_{-})\Delta^{-}f\big ](x_i)\Big\},\nonumber
\end{eqnarray}
where $\Delta f(u,v]:=f(v)-f(u+)$, $\Delta f[u,v):=f(v-)-f(u)$, 
and $[\Delta g\Delta f](J):=\Delta g(J)\Delta f(J)$.
Let $R(\lambda)$ be the set of all Young tagged partitions
which are refinements of a partition $\lambda$, and let
$\frak R$ be the family of all sets $R(\lambda)$, $\lambda\in\Xi [a,b]$.
Then $\frak R$ is the direction in the sense of Definition
\ref{convergence}.\ref{direction}.

\begin{defn}\label{SYSint}
{\rm Let $g$, $f$ be regulated functions on $[a,b]$.
Define the \emph{symmetric Young-Stieltjes} integral $(\SYS)\smallint_a^b
g\,df$ to exist and equal the limit of the directed function
$(T_{\YS}(g,f;\cdot),\frak R)$.
Also, define the \emph{Young-Stieltjes quadratic covariation}
$C_{\YS}(g,f)$ on $[a,b]$ to exist and equal the limit of
the directed function $(C_{\YS}(g,f;\cdot ),\frak R)$.}
\end{defn}

Sufficient existence conditions for the symmetric Young-Stieltjes
integral are formulated in the next section.
Here we present a formal relation between symmetric
and refinement Young-Stieltjes integrals. 

\begin{thm}\label{SYSandRYS}
Let $g$, $f$ be regulated functions on $[a,b]$ satisfying the relations
\beq\label{normalization}
\left\{ \begin{array}{ll}
\mbox{$f=(f_{+}+f_{-})/2$ on $(a,b)$}\\
\mbox{$g=(g_{+}+g_{-})/2$ on $(a,b)$, $g(a+)=g(a)$ and $g(b-)=g(b)$.}
\end{array}\right.
\eeq
Then the following holds:
\begin{enumerate}
\item[$(a)$] The integral $(\RYS)\smallint_a^b f\,dg$ exists if and only if
both the integral $(\SYS)\smallint_a^bf\,dg$ exists and the
quadratic covariation $C_{\YS}(f,g)=0$.
\item[$(b)$] Assume in addition that for $p^{-1}+q^{-1}=1$, 
one of the two functions $f, g$
is in ${\cal W}_p^{\ast}[a,b]$ and the other one is in 
${\cal W}_q[a,b]$.
Then the integrals $(\RYS)\smallint_a^bf\,dg$ and $(\SYS)\smallint_a^b
f\,dg$ both exist or not simultaneously, and if they exist then both have
the same value.
\end{enumerate}
\end{thm}

Statements $(a)$ and $(b)$ are Theorem 5.5 and Proposition 5.6
proved in \cite{RN99a}.

\paragraph{Summary.}
The relations (\ref{extend}) between extended
Riemann-Stieltjes integrals can be summarized as follows:
$$
RS\stackrel{({\cal A}_1,{\cal B}_1)}{\longrightarrow} RRS
\stackrel{({\cal A}_2,{\cal B}_2)}{\longrightarrow} \RYS
\left\{\begin{array}{lll}
\stackrel{({\cal A}_3,{\cal B}_3)}{\longrightarrow} CY\\
\stackrel{({\cal A}_4,{\cal B}_4)}{\longrightarrow} \SYS\\
\stackrel{{\cal A}_5}{\longrightarrow} HK
\end{array}
\right.
$$
{\it Legend}: Let ${\cal D}_p=\{(f,g)\colon\,f\in {\cal W}_p,
g\in {\cal W}_{p'}\}$ with $1/p+1/p'=1$. Then we have:
\begin{enumerate}
\item 
$
\left\{ \begin{array}{ll}
{\cal A}_1=\{(f,g)\colon\,f,g\in\RR^{[a,b]}\}\\
{\cal B}_1=\{(f,g)\colon\,\mbox{$\forall x\in [a,b]$, either
$f$ or $g$ is continuous at $x$}\};
\end{array}\right.
$
\item 
$
\left\{ \begin{array}{ll}
{\cal A}_2=\{(f,g)\colon\,f\in\RR^{[a,b]},\,g\in {\cal R}[a,b]\}\\
{\cal B}_2=\{(f,g)\colon\, \mbox{either $\Delta^{+}f=\Delta^{-}g
=0$ or $\Delta^{-}f=\Delta^{+}g=0$}\};
\end{array}\right.
$
\item ${\cal A}_3=\{(f,g)\colon\,f,g\in {\cal R}[a,b]\}$ 
and ${\cal B}_3={\cal D}_p$;
\item ${\cal A}_4={\cal A}_3$ and
${\cal B}_4={\cal D}_p\cap\{(f,g)\colon\,\mbox{$f=(f^{+}+f^{-})/2$
and $g=(g^{+}+g^{-})/2$}\}$;
\item ${\cal A}_5={\cal D}_p$.
\end{enumerate}

\paragraph*{The refinement Riemann-Stieltjes integral with monotone  
integrator.}
The refinement Riemann-Stieltjes and the Riemann-Stieltjes integrals with 
respect to a monotone function appear in chain rules for a function having the
quadratic $\lambda$-variation in the next chapter.
Therefore we include here few remarks concerning these integrals.
In text-books a Stieltjes type integral is often defined by restricting to
monotone integrators, or to integrators having bounded variation. 
In that case the refinement Riemann-Stieltjes integral can be defined 
without explicit use of the limit under refinements of partitions,
so that there may not be seen obvious conection between such integrals.
Usually in text-books,  an extension of the Riemann integral
developed by Darboux, Peano and others is used.
A Stieltjes type  extension of such an integral was due to Pollard \cite{SP23}.
We state it next as the necessary and 
sufficient condition for the existence of the refinement Riemann-Stieltjes 
integral in the case when $h$ is nondecreasing and $f$ is bounded.

Let $f$ be a bounded function on $[a,b]$ and let $h$ be a nondecreasing
function on $[a,b]$.
For a partition $\kappa=\{x_i\colon\,i=0,\dots,n\}$ of $[a,b]$, define the 
{\it lower and upper Riemann-Stieltjes sums}, respectively, by
$$
\left\{ \begin{array}{ll}
L_{RS}(\kappa):=
L_{RS}(f,h;\kappa):=\sum_{i=1}^nm(f;[x_{i-1},x_i])
[h(x_i)-h(x_{i-1})]\\
U_{RS}(\kappa):=
U_{RS}(f,h;\kappa):=\sum_{i=1}^nM(f;[x_{i-1},x_i])
[h(x_i)-h(x_{i-1})],
\end{array}\right.
$$
where, for any set $J$ in the domain of $f$, 
\beq\label{mandM}
m(f;J):=\inf\{f(x)\colon\,x\in J\}
\quad\mbox{and}\quad
M(f;J):=\sup\{f(x)\colon\,x\in J\}.
\eeq
Since $h$ is nondecreasing, it follows that 
$L_{RS}(\lambda)\leq U_{RS}(\kappa)$
and $U_{RS}(\lambda)\geq U_{RS}(\kappa)$ whenever $\kappa\supset\lambda$
for $\kappa,\lambda\in\Xi [a,b]$.
Also, letting
$$
\left\{ \begin{array}{ll}
L_{RS}:=L_{RS}(f,h;[a,b])
:=\sup\{L_{RS}(f,h;\kappa)\colon\,\kappa\in \Xi [a,b]\}\\
U_{RS}:=U_{RS}(f,h;[a,b])
:=\inf\{U_{RS}(f,h;\kappa)\colon\,\kappa\in \Xi [a,b]\}
\end{array}\right.
$$
we have the relation $L_{RS}\leq U_{RS}$.
Then one says that the \emph{Darboux integral} $(D)\smallint_a^b
f\,dh$ is defined and equals $b\in\RR$ if $L_{RS}=U_{RS}=b$.

Let $\PP$ be the direction in the set of all tagged partitions
of $[a,b]$ (see Notation \ref{partition}),
and let $S_{RS}(f,h;\tau)$ be the Riemann-Stieltjes sum
based on a tagged partition $\tau\in \Theta [a,b]$.
Since $h$ is non decreasing it is easy to see that
$$
L_{RS}(f,h;[a,b])=\liminf_{\PP}S_{RS}(f,h)
\quad\mbox{and}\quad
U_{RS}(f,h;[a,b])=\limsup_{\PP}S_{RS}(f,h),
$$
where the right sides are defined by (\ref{1orderconv}).
By Theorem \ref{convergence}.\ref{orderconv},
the integral $\smallint_a^bf\,dh$ exists as the
Darboux integral if and only if it exists as the refinement
Riemann-Stieltjes integrals, and then both integrals have the same
value.
It seems that the definition of an integral in terms
of the limit in the sense of refinements of partitions is preferable 
because it is not restricted to integrands having bounded variation. 

\paragraph*{The refinement Young-Stieltjes integral with monotone integrator.}
As in the case of the refinement Riemann-Stieltjes integral, 
the refinement Young-Stieltjes integration with
respect to monotone functions afford certain simplifications.
Let $f$ be a bounded function on $[a,b]$ and let $h$ be a nondecreasing
function on $[a,b]$. 
For $\kappa=\{x_i\colon\,i=0,\dots,n\}\in\Xi [a,b]$, define the lower
and upper Young sums, respectively, by
\beq\label{lower}
L_Y(f,h;\kappa):=\frak S(f,h;\kappa)+\sum_{i=1}^nm(f;(x_{i-1},x_i))
[h(x_i-)-h(x_{i-1}+)]
\eeq
\beq\label{upper}
U_Y(f,h;\kappa):=\frak S(f,h;\kappa)+\sum_{i=1}^nM(f;(x_{i-1},x_i))
[h(x_i-)-h(x_{i-1}+)],
\eeq
where $m$, $M$ are defined by (\ref{mandM}), and 
$$
\frak S(f,h;\kappa):=\sum_{i=0}^nf(x_i)[h_{+}^{(b)}(x_i)-h_{-}^{(a)}
(x_{i-1})]=[f\Delta^{+}f](a)+\sum_{i=1}^{n-1}[f\Delta^{\pm}h](x_i)
+[f\Delta^{-}h](b).
$$
Again as in the preceding paragraph, let $\PP$ be the direction in the 
set of all tagged partitions of $[a,b]$,
and let $S_{\YS}(f,h;\tau)$, $\tau\in \Theta [a,b]$, be the Young-Stieltjes sum
defined by (\ref{YS}).
Since $h$ is nondecreasing function,we have that
$$
\left\{ \begin{array}{ll}
L_{Y}:=L_{Y}(f,h;[a,b])
:=\sup\{L_{Y}(f,h;\kappa)\colon\,\kappa\in Q([a,b])\}
=\liminf_{\PP}S_{YS}(f,h)\\
U_{Y}:=U_{Y}(f,h;[a,b])
:=\inf\{U_{Y}(f,h;\kappa)\colon\,\kappa\in Q([a,b])\}
=\limsup_{\PP}S_{YS}(f,h)
\end{array}\right.
$$
Then by Theorem \ref{convergence}.\ref{orderconv}, the following holds:

\begin{thm}\label{DarbouxYoung}
Let $f$ be a bounded function on $[a,b]$ and let $h$ be a nondecreasing 
function on $[a,b]$.
Then the integral $(\RYS)\smallint_a^bf\,dh$ exists if and only if 
\beq\label{3Darboux}
A:=L_{Y}(f,h;[a,b])=U_{Y}(f,h;[a,b]),
\eeq 
and then $A$ is the value of the $\RYS$ integral.
\end{thm}

Ross \cite[\S 35]{KARb}, \cite{KARa} defined the integral
$\smallint_a^bf\,dh:=A$
if (\ref{3Darboux}) holds and proved several its properties based on this 
definition.
However, he did not relate the integral so defined
to the refinement Young-Stieltjes integral.
Love \cite[Section 5]{ERL93} does relate the integral of Ross
with the refinement Young-Stieltjes integral but makes no a reference to 
its origin.
In fact, as earlier as 1936, essentially the same integral have been 
developed by Glivenko \cite[Section V]{VIG}.

\section{Chain rules for extended Riemann-Stieltjes integrals}

Here by a chain rule we mean an integral representation of the composition
$\phi{\circ}g$ of two functions $g$ and $\phi$ such that $\phi$ is smooth
and $g$ has possibly unbounded variation.
The best known example of a chain rule is It\^o's formula giving a stochastic 
integral representation of a smooth function with a semimartingale.

\paragraph*{The refinement Young-Stieltjes integral representation.}
Let $B=\{B(t)\colon\,t\geq 0\}$ be a standard Brownian motion
on a complete probability space $(\Omega,{\cal F},\Pr)$.
For any $p>2$, we have that $v_p(B(\cdot,\omega);[0,1])<\infty$ for
almost all $\omega\in\Omega$, that is the $p$-variation index
for almost all sample functions of a Brownian motion is $2$.
Therefore it follows from the L.\ C.\ Young Stieltjes integrability theorem
(Theorem \ref{LCY}) that the Riemann-Stieltjes integral
\beq\label{1BM}
(RS)\int_0^1f\,dB(\cdot,\omega)
\eeq
is defined for almost all $\omega\in\Omega$ and for any
function $f$ having bounded $p$-variation for some
$p<2$.
Since $v_2(B(\cdot,\omega);[0,1])=+\infty$ for almost all
$\omega\in\Omega$, we cannot take $f$ in (\ref{1BM}) from
the class ${\cal W}_2$, still less can we take
$f$ to be a sample function of a Brownian motion (see Section \ref{nonex}
below).
The condition $p^{-1}+q^{-1}>1$ in Theorem \ref{LCY}
is best possible in general.
It can be improved to the condition $p^{-1}+q^{-1}=1$
provided one restricts the class of integrands to functions
of a special form.
That is, the existence part of the preceding theorem can be 
extended when the integrand $f$ is the composition $\phi{\circ}g$ 
for some smooth function $\phi$.
Instead of the inequality (\ref{1LCY}), a chain rule formula 
then holds.

Next is a chain rule for a composition representation by the refinement 
Young-Stieltjes integral; its proof is in \cite[Theorem 1.1]{RN99a}.

\begin{thm}\label{YSchrule}
For $0<\alpha \leq 1$, let $g\in{\cal W}_{1+\alpha}^{\ast}[a,b]$,
and let $\phi$ be differentiable with a derivative $\phi'$ satisfying
a Lipschitz condition of order $\alpha$.
The refinement Young-Stieltjes integral 
$(\RYS)\smallint_a^b (\phi'{\circ}g)\,dg$
is defined, and its value is determined by the relation
\beq\label{1YSchrule}
(\phi{\circ}g)(b)-(\phi{\circ}g)(a)=(\RYS)\int_a^b(\phi'{\circ}g)\,dg
\eeq
$$
+\sum_{(a,b]}\Big\{\Delta^{-}(\phi{\circ}g)-(\phi'{\circ}g)
\Delta^{-}g\Big\}+\sum_{[a,b)}\Big\{\Delta^{+}(\phi{\circ}g)
-(\phi'{\circ}g)\Delta^{+}g\Big\},
$$
where the two sums are unconditional.
\end{thm} 

The statement of the preceding theorem also holds for the Central Young 
and Henstock-Kurzweil integrals in place of the refinement Young-Stieltjes
integral.
In Theorem \ref{YSchrule}, $g$ cannot be taken to be a sample function
of a Brownian motion. 
To extend the preceding chain rule formula to functions $g\in {\cal W}_p$
with $p<3$ one can use the symmetric Young-Stieltjes integral.
The proof of the next theorem is in \cite[Theorem 1.4]{RN99a}.

\begin{thm}
For $0<\alpha\leq 1$, let $g\in {\cal W}_{2+\alpha}^{\ast}([a,b])$,
and let $\phi$ be twice differentiable with a second derivative
satisfying a Lipschitz condition of order $\alpha$.
The symmetric Riemann-Stieltjes integral $(SRS)\smallint_a^b
(\phi'{\circ}g)\,dg$ is defined, and its value is determined by the relation
$$
(\phi{\circ}g)(b)-(\phi{\circ}g)(a)=(SRS)\int_a^b(\phi'{\circ}g)\,dg
$$
$$
+\sum_{(a,b]}\Big\{\Delta^{-}(\phi{\circ}g)-(\phi'{\circ}g)\Delta^{-}g+
\frac{\Delta^{-}(\phi'{\circ}g)}{2}\Delta^{-}g\Big\}
+\sum_{[a,b)}\Big\{\Delta^{+}(\phi{\circ}g)-(\phi'{\circ}g)\Delta^{+}g-
\frac{\Delta^{+}(\phi'{\circ}g)}{2}\Delta^{+}g\Big\}.
$$
\end{thm}

\paragraph*{Integration by parts formula.}
Using the multivariate variant of the chain rule for the refinement
Young-Stieltjes integral we get an integration by parts formula:
for $g, f\in {\cal W}_2^{\ast}[a,b]$,
if $(\RYS)\smallint_a^bg\,df$ exists then so does $(\RYS)\smallint_a^b
f\,dg$ and
\beq\label{2RYSbyparts}
(\RYS)\int_a^bg\,df +(\RYS)\int_a^bf\,dg =f(b)g(b)-f(a)g(a)
+\sum_{(a,b]}\Delta^{-}f\Delta^{-}g-\sum_{[a,b)}\Delta^{+}f\Delta^{+}g,
\eeq
where the two sums converge unconditionally (see \cite[Corollary 4.3]{RN99a}).
This formula under similar conditions was known already
to L.\ C.\ Young, who used it in the proof of other results 
(see p.\ 611 in \cite{LCY38b}). 
According to Hildebrandt \cite[p.\ 276]{THH38},
formula (\ref{2RYSbyparts}) was proved by de Finetti and Jacob
\cite{deFJ} under the assumption that $g$ and $f$ have bounded variation
on $[a,b]$.
We do not know at this writing whether the assumption $g, f\in 
{\cal W}_2^{\ast}$ could be replaced just by the assumption that 
the two sums in (\ref{2RYSbyparts}) converge unconditionally.
However, a simple integration by parts formula holds for the symmetric
Young-Stieltjes integral.

\begin{thm}\label{SYSbyparts}
Let $g$, $f$ be regulated functions on $[a,b]$.
For $\sharp =\SYS$, if $(\sharp)\smallint_a^bg\,df$ exists then so does
$(\sharp)\smallint_a^bf\,dg$ and
\beq\label{1SYSbyparts}
(\sharp)\int_a^bg\,df +(\sharp)\int_a^bf\,dg =f(b)g(b)-f(a)g(a).
\eeq
\end{thm}

\begin{proof}
Let $\tau=\{([x_{i-1},x_i],y_i)\colon\,i=0,\dots,n\}$ be a Young
tagged partition of $[a,b]$.
By definition (\ref{SYSsum}) of $\SYS$ sums, we have
\begin{eqnarray*}
\lefteqn{T_{\YS}(g,f;\tau)+T_{\YS}(f,g;\tau)}\\[2mm]
&&=\sum_{i=1}^n\Big\{\Delta^{+}
(gf)(x_{i-1})+\Delta (gf)(x_{i-1},y_i]+\Delta (gf)[y_i,x_i)
+\Delta^{-}(gf)(x_i)\Big\}=\Delta (gf)[a,b].
\end{eqnarray*}
Thus the directed function $(T_{\YS}(g,f;\cdot ),\frak R)$ has a limit
if and only if the directed function $(T_{\YS}(f,g;\cdot ),\frak R)$ does,
and (\ref{1SYSbyparts}) holds provided both limits exist.
The proof is complete.
\qed\end{proof}

Under the conditions of Theorem \ref{SYSandRYS}, relating the
symmetric and refinement Young-Stieltjes integral,
we get the same integration by parts formula for the
refinement Young-Stieltjes integral.

\begin{cor}\label{RYSbyparts}
Let $g$, $f$ be regulated functions on $[a,b]$ satisfying
{\rm (\ref{normalization})}.
For $\sharp =\RYS$, if $(\sharp)\smallint_a^bg\,df$ exists then so does
$(\sharp)\smallint_a^bf\,dg$ and {\rm (\ref{1SYSbyparts})} holds.
\end{cor}

\begin{proof}
Suppose that $(\RYS)\smallint_a^bg\,df$ exists.
Then by statement $(a)$ of Theorem \ref{SYSandRYS}, 
$(\SYS)\smallint_a^bg\,df$ exists, and 
\beq\label{1RYSbyparts}
0=C_{\YS}(g,f)=C_{\YS}(f,g).
\eeq
By the form of the $\SYS$ sums, it also follows that
the two integrals have the same value.
Hence by the preceding theorem, $(\SYS)\smallint_a^bf\,dg$ exists
and the relations
$$
(\RYS)\int_a^bg\,df=(\SYS)\int_a^bg\,df=\Delta (gf)[a,b]
-(\SYS)\int_a^bf\,dg.
$$
hold.
This in conjunction with (\ref{1RYSbyparts}) allow us to apply
statement $(a)$ of Theorem \ref{SYSandRYS} in the opposite direction
to conclude that $(\RYS)\smallint_a^bf\,dg$ exists and (\ref{1SYSbyparts}) 
holds with $\sharp=\RYS$.
The proof of Corollary \ref{RYSbyparts} is complete.
\qed\end{proof}

It is easy to see that if regulated functions $g$, $f$ have
discontinuities satisfying condition (\ref{normalization}), then
the two sums in formula (\ref{2RYSbyparts}) disappear and we get
the integration by parts formula (\ref{1SYSbyparts}) with
$\sharp=\RYS$ without the assumption $g, f\in {\cal W}_2^{\ast}$.

If instead of condition (\ref{normalization}), the regulated
functions $g$, $f$ on $[a,b]$ satisfy the same conditions except 
that jumps of $g$ at the endpoints are allowed, then
statement $(a)$ of Theorem \ref{SYSandRYS} holds with 
the value 
$$
C_{\YS}(g,f)=\Delta^{+}g(a)\Delta^{+}f(a)-\Delta^{-}g(b)\Delta^{-}f(b)
:=A.
$$
Following then the arguments of the proof of the preceding
corollary, we get that the integration by parts formula
$$
(\RYS)\int_a^bg\,df +(\RYS)\int_a^bf\,dg =f(b)g(b)-f(a)g(a)-A
$$
holds for the refinement Young-Stieltjes integral.
This is a result recently proved by Love \cite[Theorem 25]{ERL98}.
We notice that the $R^3S$ integral in Love's work is the
refinement Young-Stieltjes integral.

\section{The Left and Right Young integrals}\label{LYRYintegrals}

Next we define two different versions of the Central Young integral,
that is the Left and Right Young integrals.
These two integrals are better suited to solving linear integral equations
(see Section 5.4 in \cite[Part II]{DNa}).
For a regulated function $f$ on $[a,b]$, the left-continuous
function $f_{-}^{(a)}$ and the right-continuous function $f_{+}^{(b)}$
are defined by
\beq\label{lmodification}
f_{-}^{(a)}(x):=f_{-}(x):=f(x-)\quad\mbox{for}\quad a<x\leq b,
\quad f_{-}^{(a)}(a):=f(a)
\eeq
and
\beq\label{rmodification}
f_{+}^{(b)}(x):=f_{+}(x):=f(x+)\quad\mbox{for}\quad
a\leq x<b,
\quad f_{+}^{(b)}(b):=f(b),
\eeq
respectively.

\begin{defn}\label{LYandRY}
{\rm Let $f$ and $g$ be regulated functions on $[a,b]$.
Define the {\it Left Young integral}, or the $LY$ integral,
by
$$
(LY)\int_a^bg\,df:=(RRS)\int_a^bg_{-}^{(a)}\,df_{+}^{(b)}
+\big [g\Delta^{+}f\big ](a)+\sum_{(a,b)}\Delta^{-}g\Delta^{+}f
$$
provided $a<b$, the $RRS$ integral exists and the sum  converges
unconditionally.
We say that $g$ is $LY$ integrable with respect to $f$ on $[a,b]$.
Define the {\it Right Young integral}, or the $RY$ integral, by
$$
(RY)\int_a^bg\,df:=(RRS)\int_a^bg_{+}^{(b)}\,df_{-}^{(a)}
-\sum_{(a,b)}\Delta^{+}g\Delta^{-}f+\big [g\Delta^{-}f\big ](b)
$$
provided $a<b$, the $RRS$ integral exists and the sum  converges 
unconditionally.
We say that $g$ is $RY$ integrable with respect to $f$ on $[a,b]$.
If $a=b$ then both integrals are defined as $0$.}
\end{defn}

The Left Young and Right Young integrals are defined in \cite{DNa}
(see relation (3.44) there) for Banach algebra valued functions.
The present variants of the two integrals are special cases of the integrals
with two integrands used to extend the Duhamel formula for product integrals.
It is easy to see that if $f$ is right-continuous then 
\beq\label{LYMPS}
(LY)\int_a^bg\,df=(RRS)\int_a^bg_{-}^{(a)}\,df
\eeq
for each regulated function $g$.
Similarly, if $f$ is left-continuous then
$$
(RY)\int_a^bg\,df=(RRS)\int_a^bg_{+}^{(b)}\,df.
$$
The Left Young and Right Young integrals can be approximated by 
Left and Right Cauchy sums, respectively. 
For $\kappa=\{x_i\colon\,i=0,\dots,n\}\in\Xi [a,b]$, define
the Left Cauchy sum
$$
S_{LC}(g,f;\kappa):= 
\sum_{i=1}^ng(x_{i-1})[f(x_i)-f(x_{i-1})],
$$
and the Right Cauchy sum 
$$
S_{RC}(g,f;\kappa)
:=\sum_{i=1}^ng(x_{i})[f(x_i)-f(x_{i-1})].
$$
These sums define respectively the Left Cauchy-Stieltjes integral
$(LCS)\smallint_a^bg\,df$ and the Right Cauchy-Stieltjes integral
$(RCS)\smallint_a^bg\,df$ in the sense of refinements of partitions.
Let $\Lambda [a,b]$ be the set of all nested sequences $\lambda=
\{\lambda_m\colon\,m\geq 1\}$ of partitions of $[a,b]$ such that
$\cup_m\lambda_m$ is dense in $[a,b]$.

\begin{thm}\label{approximation}
For $1\leq p<\infty$, let $f\in {\cal W}_p[a,b]$ and
{\rm $g\in \dual\, ({\cal W}_p)[a,b]$}.
Then $g$ is Left Young integrable with respect to $f$ on $[a,b]$
and the relation
\beq\label{1approximation}
\lim_{m\to\infty}S_{LC}(g,f;\lambda_m)=(LY)\int_a^bg\,df
\eeq
holds for each $\{\lambda_m\colon\,m\geq 1\}\in\Lambda [a,b]$ 
such that 
\beq\label{0approximation}
\{x\in (a,b)\colon\,[\Delta^{-}g\Delta^{+}f](x)\not =0\}
\subset\cup_m\lambda_m.
\eeq
Also,  $g$ is Right Young integrable with respect to $f$ on $[a,b]$
and the relation
$$
\lim_{m\to\infty}S_{RC}(g,f;\lambda_m)=(RY)\int_a^bg\,df
$$
holds for each $\{\lambda_m\colon\,m\geq 1\}\in\Lambda [a,b]$ 
such that 
$$
\{x\in (a,b)\colon\,[\Delta^{+}g\Delta^{-}f](x)\not =0\}
\subset\cup_m\lambda_m.
$$
\end{thm}

\begin{proof}
We prove the theorem only for the Left Young integral and in the case $p=1$
because a proof for the Right Young integral is similar and 
a proof for the case $p>1$ is given in \cite[Theorem 3]{RN99b}.
Since $f_{+}^{(b)}$ has bounded variation and $g_{-}^{(a)}$ is regulated,
the integral $(RRS)\smallint_a^bg_{-}^{(a)}\,df_{+}^{(b)}$ exists
by Theorem \ref{LCY}.$(b)$.
Clearly the sum $\sum_{(a,b)}|\Delta^{-}g\Delta^{+}f|$ converges 
absolutely, and hence unconditionaly.
Thus the integral $(LY)\smallint_a^bg\,df$ is defined.
To prove relation (\ref{1approximation}) consider a sequence 
of partitions $\lambda_m=\{x_i^m\colon\,i=0,\dots,n(m)\}$ 
satisfying the stated assumptions.
For each $m\geq 1$, let $g_m(x):=g(x_{i-1}^m)$ if
$x\in [x_{i-1}^m,x_i^m)$ for some $i\in\{1,\dots,n(m)\}$ and let 
$g_m(b):=g(b)$.
Using notation (\ref{lmodification}) and (\ref{rmodification})
for $[a,b]=[x_{i-1}^m,x_i^m]$, for each $i=1,\dots,
n(m)$, we have
$$
(RRS)\int_{x_{i-1}^m}^{x_i^m}(g_m)_{-}^{(x_{i-1}^m)}\,df_{+}^{(x_i^m)}
=g(x_{i-1}^m)[f(x_i^m)-f(x_{i-1}^m+)].
$$
Using additivity of the $LY$ integral proved in \cite[Theorem 5]{RN99b},
and since  $\Delta^{-}(g_m)=0$ on $(x_{i-1}^m,x_i^m)$, it then follows that
\begin{eqnarray}
(LY)\int_a^bg_m\,df&=&\sum_{i=1}^{n(m)}\Big\{(RRS)\int_{x_{i-1}^m}^{x_i^m}
(g_m)_{-}^{(x_{i-1}^m)}\,df_{+}^{(x_i^m)}+[g_m\Delta^{+}f](x_{i-1}^m)
\Big\}\nonumber\\
&=&S_{LCS}(g,f;\lambda_m).
\label{2approximation}
\end{eqnarray}
On the other hand, since $g_m$ has only finitely many jumps, we have
the representation
\begin{eqnarray}
(LY)\int_a^bg_m\,df\!&=&\!(RRS)\int_a^b(g_m)_{-}^{(a)}\,df_{+}^{(b)}
+\sum_{i=1}^{n(m)-1}[\Delta^{-}g\Delta^{+}f](x_i^m)\label{3approximation}\\
& &\!+\!\sum_{i=1}^{n(m)-1}[g(x_i^m-)-g(x_{i-1}^m)]\Delta^{+}f(x_i^m)
+[g\Delta^{+}f](a).\nonumber
\end{eqnarray}
To prove
\beq\label{4approximation}
\lim_{m\to\infty}(RRS)\int_a^b(g_m)_{-}^{(a)}\,df_{+}^{(b)}
=(RRS)\int_a^bg_{-}^{(a)}\,df_{+}^{(b)}.
\eeq
we use the Osgood theorem on term by term integration (see e.g.
\cite[Theorem II.15.6]{THH63}).
It is clear that $g_m$ are uniformly bounded.
To show that $(g_m)_{-}^{(a)}$ converge to $g_{-}^{(a)}$, let $x\in (a,b]$.
For each $m\geq 1$, there is $i\leq n(m)$ such that $x_{i-1}^m<x\leq x_i^m$.
Then $g_m(x-)=g(x_{i-1}^m)\to g(x-)$ as $m\to\infty$.
Thus all the hypotheses of the Osgood theorem are satisfied, and hence
(\ref{4approximation}) holds.
To prove
\beq\label{5approximation}
\lim_{m\to\infty}\sum_{i=1}^{n(m)-1}[\Delta^{-}g\Delta^{+}f]
(x_i^m)=\sum_{(a,b)}[\Delta^{-}g\Delta^{+}f]
\eeq
notice that each finite sum of terms $[\Delta^{-}g\Delta^{+}f](x)$ is included
into approximating sum for sufficiently large $m$ by (\ref{0approximation}).
To prove
\beq\label{6approximation}
\lim_{m\to\infty}\sum_{i=1}^{n(m)-1}[g(x_i^m-)-g(x_{i-1}^m)]\Delta^{+}f(x_i^m)
=0
\eeq
notice that each term tends to zero as $m\to\infty$ and any finite sum of them is
bounded by $\sum |\Delta^{+}f|$.
Now (\ref{1approximation}) follows from (\ref{2approximation}) -- 
(\ref{6approximation}).
The proof is complete.
\qed\end{proof}

\begin{cor}
Under the hypotheses of Theorem {\rm \ref{approximation}},
the $LCS$ and $RCS$ integrals 
exist and equal respectively to the $LY$ and $RY$ integrals.
\end{cor}

\begin{proof}
Suppose that the conclusion concerning the integral 
$(LCS)\smallint_a^bg\,df$ does not hold.
Then there exists $\epsilon >0$ such that for any $\lambda\in\Xi [a,b]$
the least upper bound of $|S_{LCS}(g,f;\kappa)-(LY)\smallint_a^b
g\,df|$ over all refinements $\kappa$ of $\lambda$ is bigger
than $\epsilon$.
Then one can construct recursively a sequence $\{\lambda_m\colon\,
m\geq 1\}\in\Lambda [a,b]$ satisfying (\ref{0approximation}) but 
not (\ref{1approximation}).
A contradiction proves the first part of the corollary.
The proof of the second one is similar. 
\end{proof}

Next is the Love-Young inequality for the Left Young and Right Young integrals.
Recall that ${\cal W}_{\infty}[a,b]$ is the class ${\cal R}[a,b]$ of all regulated 
functions on $[a,b]$.

\begin{thm}\label{LYineq}
Let $1\leq p<\infty$, and $q$ be such that $p^{-1}+q^{-1}>1$ if $p>1$
and $q=\infty$ if $p=1$.
For $f\in {\cal W}_p[a,b]$ and $g\in {\cal W}_q[a,b]$, the integrals 
$(LY)\smallint_a^b g\,df$ and $(RY)\smallint_a^b g\,df$ are defined.
Moreover,
there exists a constant $K$ if $p>1$ which depends only on $p$ and $q$ 
such that the inequality
$$
\Big |(\sharp)\int_a^b\big [g-g(y)\big ]\,df\Big |
\leq \left\{\begin{array}{ll}
\osc\, (g;[a,b])\,v_1(f;[a,b]) &\mbox{if $p=1$,}\\
Kv_q(g;[a,b])\,v_p(f;[a,b]) &\mbox{if $1<p<\infty$,}
\end{array}\right.
$$ 
holds for $\sharp =LY$, $\sharp =RY$ and for any $y\in [a,b]$.
\end{thm}

The proof follows easily from Theorem \ref{LCY}, H\"older's inequality
and Theorem \ref{approximation}.

We finish with a formulation of several useful properties of
the Left Young and Right Young integrals used in this paper.

\begin{thm}\label{additivity}
Let $g$, $f$ be regulated functions on $[a,b]$, and let 
$a\leq c\leq b$.
For $\sharp =LY$ or $RY$,
$(\sharp)\smallint_a^bg\,df$ exists if and only if both $(sharp)\smallint_a^c
g\,df$ and $(sharp)\smallint_c^bg\,df$ exist, and then
$$
(\sharp)\int_a^bg\,df=(\sharp)\int_a^cg\,df+(\sharp)\int_c^bg\,df.
$$
\end{thm}

The direct proof of the preceding additivity property is given in
\cite[Theorem 5]{RN99b}.

Next is an integration by parts formula for the Left Young and Right Young
integrals.

\begin{thm}\label{byparts}
Let $g$, $f$ be regulated functions on $[a,b]$.
If either of the two integrals $(LY)\smallint_a^bg\,df$
and $(RY)\smallint_a^bf\,dg$ exists then both exist, and
$$
(LY)\int_a^bg\,df+(RY)\int_a^bf\,dg=f(b)g(b)-f(a)g(a).
$$
\end{thm}

The proof of the theorem is given in \cite[Theorem 7]{RN99b}.

Let $\Phi$ and $\Psi$ be \emph{indefinite} Left Young and Right Young 
integrals defined respectively by
$$
\Phi (x):=const +(LY)\int_a^xg\,dh\quad\mbox{and}\quad
\Psi (y):=const +(RY)\int_y^bg\,dh
$$
for $x,y\in [a,b]$.
By the following statement $\Phi$ and $\Psi$, are regulated functions.

\begin{prop}\label{LYjumps}
For regulated functions $g$ and $h$ on $[a,b]$ the following hold{\rm :}
\begin{enumerate}
\item[$(a)$] If $g$ is Left Young integrable with respect to $h$ on $[a,b]$
then $\Phi$ is a regulated function on $[a,b]$, and for
$a\leq y<x\leq b$,
$$
(\Delta^{-}\Phi)(x)=[g_{-}\Delta^{-}h](x)\quad\mbox{ and }\quad
(\Delta^{+}\Phi)(y)=[g\Delta^{+}h](y).
$$
\item[$(b)$] If $g$ is Right Young integrable with respect to $h$ on $[a,b]$
then $\Psi$ is a regulated function on $[a,b]$, and for
$a\leq y<x\leq b$,
$$
(\Delta^{-}\Psi)(x)=-[g\Delta^{-}h](x)\quad\mbox{ and }\quad
(\Delta^{+}\Psi)(y)=-[g_{+}\Delta^{+}h](y).
$$
\end{enumerate}
\end{prop}

The proof of the proposition is given in \cite[Proposition 8]{RN99b}.

Next is the \emph{substitution rule} for the refinement Riemann-Stieltjes
integrals.
It is used in Chapter \ref{process} and to derive the corresponding results 
for the Left Young and Right Young integrals.                    

\begin{prop}\label{MPSsrule}
Let $h\in{\cal W}_p[a,b]$ and $f,g\in{\cal W}_q[a,b]$ for some 
$p,q >0$ with $p^{-1}+q^{-1}>1$. 
Suppose that the pairs $h, g$ and $h, f$ have no common 
discontinuities on the left and on the right at the same point.
Then  $g$ and $fg$ are $RRS$ integrable with respect to $h$,
$f$ is $RRS$ integrable with respect to the indefinite $RRS$ integral
$\Phi (y):=const+(RRS)\smallint_a^yg\,dh$, $y\in[a,b]$,
and 
$$             
(RRS)\int_a^bf\,d\Phi =(RRS)\int_a^b fg\,dh\,.
$$
\end{prop}

The proof of this proposition is given in \cite[Proposition 9]{RN99b}.

Next are substitution rules for the Left Young and Right Young integrals.

\begin{thm}\label{srule}
Let $h\in\WW_p[a,b]$ and $f,g\in\WW_q[a,b]$ for some 
$p,q >0$ with $p^{-1}+q^{-1}>1$. 
Then $g$ and $fg$ are $LY$ integrable with respect to $h$, 
$f$ is $LY$ integrable with respect to the indefinite 
$LY$ integral $\Psi$ and 
$$
(LY)\int_a^bf\,d\Psi =(LY)\int_a^bfg\,dh\,.
$$
Similarly, $g$ and $fg$ are $RY$ integrable with respect to $h$, 
$f$ is $RY$ integrable with respect to the indefinite 
$RY$ integral $\Phi$ and 
$$
(RY)\int_a^bf\,d\Phi =(RY)\int_a^bfg\,dh\,.
$$
\end{thm}

The proof of this theorem is given in \cite[Theorem  10]{RN99b}.

The final statement provides chain rules for the Left Young
and Right Young integral representations of a composition.
Its proof is given in \cite[Theorems 11 and 14]{RN99b}.

\begin{thm}\label{main}
For $\alpha\in [0,1]$, let  $f\in{\cal W}_{1+\alpha}^{\ast}[a,b]$ 
let $\g$ be differentiable with the first derivative satisfying  a
H\"older condition of order $\alpha$,
and let $\h$ be a regulated function on $[a,b]$.
Then the two equalities
$$
(LY)\int_a^b\h\,d(\g{\circ}f)= (LY)\int_a^b\h(\g'{\circ}f)\,df
$$
$$
+\sum_{(a,b]}h_{-}\big [\Delta^{-}(\g{\circ} f)-
(\g'{\circ}f)_{-}\Delta^{-}f\big ]
+\sum_{[a,b)}h\big [\Delta^{+}(\g{\circ} f)-
(\g_l'{\circ}f)\Delta^{+}f\big ]
$$
and
$$
(RY)\int_a^b\h\,d(\g{\circ}f)= (RY)\int_a^b\h(\g'{\circ}f)\,df
$$
$$
+\sum_{(a,b]}h\big [\Delta^{-}(\g{\circ} f)-
(\g'{\circ}f)\Delta^{-}f\big ]
+\sum_{[a,b)}h_{+}\big [\Delta^{+}(\g{\circ} f)-
(\g_l'{\circ}f)_{+}\Delta^{+}f\big ]
$$
hold both meaning that the two integrals exist provided any
one of the two exists, and the two sums converge unconditionally.
\end{thm}

\newpage
\chapter{Quadratic variation and restricted integrals of Stieltjes 
type}\label{function}
\setcounter{thm}{0}

\vspace*{0.2truein}
\begin{quotation}
{\footnotesize
In mathematics, fruitful points of view are more powerful tools
than so called key theorems.
To say more precisely, they are eyes of a researcher, ardently trying
to understand the essence of mathematical objects.

If the reality we seek to understand is rich and complicated
then to grasp its depth and exquisitiness, we need as many such
"eyes" as possible.

A.\  Grothendieck. "Recoltes et Semailles".}
\end{quotation}
\vspace*{0.2truein}

Recall that a regulated function $f$ on $[a,b]$ has the local $2$-variation,
that is belongs to the Wiener class ${\cal W}_2^{\ast}[a,b]$,
 if and only if the limit
\beq\label{0exmp}
\lim_{\kappa,\prtn}s_2(f;\kappa)\quad 
\Big (=\sum_{(a,b]}(\Delta^{-}f)^2
+\sum_{[a,b)}(\Delta^{+}f)^2\Big )
\eeq
exists in the sense of refinements of partitions of 
the interval $[a,b]$ (see (\ref{local-p-var}) and Section \ref{Wiener}).
In this chapter we consider a class of functions such that instead of
(\ref{0exmp}), the limit exists along a nested sequence $\lambda$ of 
partitions densely filling the interval $[a,b]$.
We call such a property of $f$, if it holds,
the quadratic $\lambda$-variation of $f$. 
In a similar manner we modify a Stieltjes type integral.
Namely, we consider two variants of a $\lambda$-integral 
defined to be a limit, if it exists, of special type Riemann-Stieltjes sums along 
a sequence $\lambda$.
The two integrals, when they exist, are called the Left Cauchy and Right Cauchy
$\lambda$-integrals, respectively.
Existence of the quadratic $\lambda$-variation of a function $f$ 
implies the existence of the two $\lambda$-integrals of a composition
$\phi{\circ}f$ with respect to $f$ provided $\phi$ is a smooth
function.
In Section \ref{uniqueness}, the extended Dol\'eans exponential is shown 
to be the unique solution to the augmented linear $\lambda$-integral 
equation.
The main result is in Section \ref{e-representation}, where the evolution 
representation problem is solved for a class of functions having the quadratic 
$\lambda$-variation.

\section{The quadratic $\lambda$-variation and $\lambda$-covariation}

\paragraph*{The quadratic $\lambda$-variation: definition and examples.}
We start with a notation related to partitions of intervals.

\begin{notat}\label{lambda}
{\rm Let $J$ be a bounded nondegenerate interval of real numbers,  open or closed 
at either end.
A \emph{partition} of $J$ is a finite increasing sequence $\{t_0,t_1,\dots,t_n\}$ 
of points of $J$, where $t_0$ if the left endpoint of $J$ if it is left closed
and $t_n$ is the right endpoint of $J$ if it is right closed.
The class of all sequences $\lambda=\{\lambda_m\colon\,m\geq 1\}$ of 
partitions $\lambda_m=\{t_i^m\colon\,i=0,\dots,n(m)\}$ of $J$
such that $\lambda_m\subset\lambda_{m+1}$ for $m=1,2,\dots$ and
$\cup_m\lambda_m$ is dense in the closure of $J$ is denoted by 
$\Lambda (J)$.}
\end{notat}

Next we define a partition of an extended subinterval $\lei u,v\rei$,
$(u,v)\in S\lei a,b\rei$, induced by a partition of an ordinary 
interval $[a,b]$.

\begin{notat}\label{lmuv}
{\rm Let $\kappa=\{t_i\colon\,i=0,\dots,n\}$, $n\geq 2$, be a partition of $[a,b]$, 
and let $(u,v)\in S\lei a,b\rei$. 
First suppose that $u\in\{x,x+\}$ and $v\in\{y-,y\}$ for some
$a\leq x < y\leq b$.
If between $u$ and $v$ there are no points of $\kappa$, then let
$\kuv:=\{u,v\}$.
Otherwise denote by $i(u)$ and $i(v)$ the two indices in $\{1,\dots,n\}$ 
such that $t_{i(u)-1}\leq u<t_{i(u)}$ and $t_{i(v)-1}<v\leq t_{i(v)}$,
respectively.
Denoting $x_{i(u)-1}:=u$, $x_{i(v)}:=v$ and $x_i:=t_i$ for $i(u)\leq 
i<i(v)$, let 
$$
\kuv:=\{u,t_{i(u)},\dots,t_{i(v)-1},v\}=\{x_{i(u)-1},\dots,x_{i(v)}\}.
$$
Now suppose that $u=x-$ and $v\in\{y-,y\}$ for some $a<x\leq y\leq b$.
If $x\in\kappa$ then $\kuv$ is defined as before.
If $x\not\in\kappa$ then let $\kuv:=(\kappa\cup\{x\})\Cap\loc u,v\roc$.
Finally suppose that $v=y+$ for some $a\leq y<b$.
If $y\in\kappa$ then $\kuv$ is defined as before.
If $y\not\in\kappa$ then let $\kuv:=(\kappa\cup\{y\})\Cap\loc u,v\roc$.
Then $\kuv$ is called the \emph{trace partition of $\lei u,v \rei$ induced 
by $\kappa$}.}
\end{notat}

Notice that the trace partition of $\loc u,v\roc$
consists of ordinary points except possibly for the endpoints
$u$ and $v$.
If $(u,v)\in S[a,b]$, that is, if $[u,v]$ is an ordinary interval then 
we write $\kappa\Cap [u,v]$ instead of $\kuv$.
Recalling notation (\ref{s_p}) with $p=2$,  for a partition 
$\kappa=\{x_i\colon\,i=0,\dots,n\}$ and any $a\leq t\leq b$, we have
$$
s_2(f;\kappa\Cap [a,t])=\sum_{i=1}^n\big [f(x_i\wedge t)-f(x_{i-1}
\wedge t)\big ]^2.
$$
Also by the definition of the trace partition, we have
$$
\left\{ \begin{array}{l}
s_2(f;\kappa\Cap \lei s-,s\rei)=s_2(f;\{s-,s\})=(\Delta^{-}f(s))^2,\\
s_2(f;\kappa\Cap \lei t,t+\rei)=s_2(f;\{t,t+\})=(\Delta^{+}f(t))^2
\end{array}\right.
$$
for any partition $\kappa$ of $[a,b]$ and
for any $a\leq t<s\leq b$.

\begin{defn}\label{qv}
{\rm Let $f$ be a regulated function on $[a,b]$, and let
$\lambda=\{\lambda_m\colon\,m\geq 1\}\in\Lambda [a,b]$.
We say that $f$ has the \emph{quadratic $\lambda$-variation on $[a,b]$}
if there exists an additive upper continuous function
$\alpha_{\lambda}(f)$ on $S\lei a,b\rei$ such that for each
$(u,v)\in S\lei a,b\rei$,
\beq\label{1alpha}
\alpha_{\lambda}(f;u,v)=\lim_{m\to\infty}s_2(f;\lmuv).
\eeq
If $f$ has the quadratic $\lambda$-variation on $[a,b]$ then
the right distribution function of $\alpha_{\lambda}(f)$ restricted
to $[a,b]$ is called the \emph{bracket function} of $f$ and is denoted by
$[f]_{\lambda}$, that is, $[f]_{\lambda}(t)=\alpha_{\lambda} (f;a,t)$ for 
$t\in [a,b]$ (cf.\ (\ref{RandLdf})).}
\end{defn}

Alternatively the quadratic $\lambda$-variation on an interval $[a,b]$
can be defined solely in terms of a regulated function on $[a,b]$ 
as follows:

\begin{prop}\label{alpha}
Let $f$ be a regulated function on $[a,b]$, and let
$\lambda=\{\lambda_m\colon\,m\geq 1\}\in\Lambda [a,b]$.
The function $f$ has the quadratic $\lambda$-variation on $[a,b]$ if 
and only if there exists a regulated function $H=H_{\lambda}$ on $[a,b]$ 
such that $H(a)=0$ and for each $(s, t)\in S[a,b]$, 
\beq\label{3alpha}
H(t)-H(s)=\lim_{m\to\infty}s_2(f;\lambda_m\Cap [s,t]), 
\eeq
\beq\label{4alpha}
\Delta^{-}H(t)=(\Delta^{-}f(t))^2\qquad\mbox{and}\qquad
\Delta^{+}H(s)=(\Delta^{+}f(s))^2. 
\eeq
If the two statements hold then $[f]_{\lambda}=H$.
Moreover, condition {\rm (\ref{3alpha})} for each $(s,t)\in S[a,b]$,
can be replaced by the condition
\beq\label{5alpha}
H(t)=\lim_{m\to\infty}s_2(f;\lambda_m\Cap [a,t]), 
\eeq
for each $t\in [a,b]$, provided two-sided discontinuity points of $f$ are
accessible by $\lambda$, that is {\rm (\ref{jumps-access})} holds.
\end{prop}

\begin{proof}
First suppose that there exists a regulated function $H$ on $[a,b]$,
equal to zero at $a$, and such that (\ref{3alpha}) and (\ref{4alpha}) hold.
Define the function $\alpha_{\lambda}(f)$ on $S\lei a,b\rei $
by $\alpha_{\lambda}(f;u,v):=H(v)-H(u)$ for $(u,v)\in S\lei a,b\rei $.
Then $\alpha_{\lambda}(f)$ is additive and upper continuous
function on $S\lei a,b\rei $ by Theorem \ref{variation}.\ref{interv}.
We prove (\ref{1alpha}) only for the cases $( u,v)=( a,t-)$ 
and $( u,v)=( a,t+)$  because the proofs for the other cases 
are similar.
To this aim let 
$\lambda_m=\{x_i^m\colon\,i=0,\dots,n(m)\}$, $m=1,2,\dots$,
and let $t\in (a,b]$.
For each $m\geq 1$, there is an index $i(t)\in\{1,\dots,n(m)-1\}$
such that $x_{i(t)}^m<t\leq x_{i(t)+1}^m$.
Then by (\ref{3alpha}) and (\ref{4alpha}), we have
$$
\lim_{m\to\infty}s_2(f;\lambda_m\Cap\lei a,t-\rei )
=\lim_{m\to\infty}\Big\{s_2(f;\lambda_m\Cap [a,t])
-[f(t)-f(x_{i(t)}^m)]^2+[f(t-)-f(x_{i(t)}^m)]^2\Big\}
$$
$$
=H(t)-(\Delta^{-}f(t))^2=H(t-)=\alpha_{\lambda}(f;a,t-),
$$
proving $(\ref{1alpha})$ for $(u,v)= (a,t-)\in S\lei a,b\rei $.
Now let $t\in [a,b)$.
By the definition of the trace partition, we have for each $m\geq 1$,
$$
s_2(f;\lambda_m\Cap\lei a,t+\rei )
=s_2(f;\lambda_m\Cap [a,t]) +\big (\Delta^{+}f(t)\big )^2,
$$
proving $(\ref{1alpha})$ for $(u,v)= (a,t+)\in S\lei a,b\rei $.

To prove the converse implication, let $H$ be the restriction to $[a,b]$
of the right distribution function of $\alpha_{\lambda} (f)$, that is,
$H(t):=\alpha_{\lambda}(f;a,t)$ for $t\in [a,b]$.
By statement $(iv)$ of Theorem \ref{variation}.\ref{interv} and 
(\ref{1alpha}), $H$ is regulated on $[a,b]$, is zero at $a$, 
and satisfies properties (\ref{3alpha}), (\ref{4alpha}). 

To prove the last part of the conclusion, let $(s,t)\in S[a,b]$
be such that $s<t$, and let (\ref{5alpha}) holds for each $t\in [a,b]$.
For each $m\geq 1$, let $i(s)=i_m(s)\in\{1,\dots,n(m)\}$ be the index such that
$x_{i(s)-1}^m<s\leq x_{i(s)}^m$.
Then we have
\begin{eqnarray*}
s_2(f;\lambda_m\Cap [s,t])&=&s_2(f;\lambda_m\Cap [a,t])-s_2(f;\lambda_m 
\Cap [a,s])\\
& &+[f(x_{i(s)}^m)-f(s)]^2+[f(s)-f(x_{i(s)-1}^m)]^2
-[f(x_{i(s)}^m)-f(x_{i(s)-1}^m)]^2.
\end{eqnarray*}
Since a two-sided discontinuity point eventually is a partition point, 
by (\ref{5alpha}), the right side converges to
$H(t)-H(s)$ as $m\to\infty$, proving (\ref{3alpha}).
The proof of Proposition \ref{alpha} is now complete.
\qed\end{proof}

Suppose that for some $\lambda\in\Lambda [a,b]$, a regulated function $f$ 
on $[a,b]$ has the quadratic $\lambda$-variation.
By (\ref{3alpha}), the bracket function $H=[f]_{\lambda}$ is nondecreasing
on $[a,b]$, and so it has a decomposition into a sum of a pure jump function
and a continuous function.
By (\ref{4alpha}), the pure jump function of $[f]_{\lambda}$ is given by 
$\sum_{(a,t]}\{\Delta^{-}f\}^2+\sum_{[a,t)}\{\Delta^{+}f\}^2$, $t\in [a,b]$.
The difference between $[f]_{\lambda}$ and its pure jump function is
nondecreasing continuous function on $[a,b]$, and is denoted by
$[f]_{\lambda}^c$, so that for each $a\leq t\leq b$,
\beq\label{cnt-part}
[f]_{\lambda}(t)=[f]_{\lambda}^c(t)+\sum_{(a,t]}\{\Delta^{-}f\}^2
+\sum_{[a,t)}\{\Delta^{+}f\}^2.
\eeq

The notion of the quadratic $\lambda$-variation is similar to the one
defined by F\"ollmer \cite{HFa} for a regulated and right-continuous function 
on $[0,\infty)$.
We recall his definition in the case when a regulated and right-continuous 
function $f$ is defined on a closed interval $[a,b]$.
For a sequence $\lambda=\{\lambda_m\colon\,m\geq 1\}$
of partitions $\lambda_m=\{x_i^m\colon\,i=0,\dots,n(m)\}$ of $[a,b]$ 
such that the mesh tends to zero as $m\to\infty$, let
$$
\xi_m:=\sum_{i=1}^{n(m)}\big [f(x_i^m)-f(x_{i-1}^m)]^2
\epsilon_{x_{i-1}^m},
$$
where $\epsilon_x$ is a point mass at $x$.
The function $f$ has the quadratic variation on $[a,b]$ in the sense of 
F\"ollmer provided a Radon measure $\xi$ exists on Borel sets of $[a,b]$
such that point measures  $\xi_m$ converge weakly to $\xi$ as 
$m\to\infty$, and $\xi (\{t\})=(\Delta^{-}f(t))^2$ for each $a<t\leq b$.
By Theorem \ref{intqv} below, if for such a function, the quadratic 
$\lambda$-variation exists then F\"ollmer's
quadratic variation also exists, and the measure $\xi$ is the 
Lebesgue-Stieltjes measure on $[a,b]$ generated
by the bracket function $[f]_{\lambda}$.

\begin{exmp}\label{qv-for-W2}
{\rm If $f$ has the local $2$-variation on $[a,b]$ and the two-sided discontinuity
points of $f$ are accessible by $\lambda\in \Lambda [a,b]$,
then $f$ has the quadratic $\lambda$-variation, and its bracket function 
$[f]_{\lambda}$ is the pure jump function given by
\beq\label{2exmp}
[f]_{\lambda}(t)=\sum_{(a,t]}\big\{\Delta^{-}f\big\}^2
+\sum_{[a,t)}\big\{\Delta^{+}f\big\}^2
\eeq
for $a\leq t\leq b$ (see Corollary \ref{p-var-qv-var} below).
If $f$ is either right-continuous or left-continuous then assumption 
(\ref{jumps-access}) always holds because in such cases the set on the left side 
is empty.}
\end{exmp}

Relation (\ref{2exmp}) may be compared with the property (\ref{0exmp}) 
defining the Wiener class ${\cal W}_2^{\ast}$.
Given $\lambda\in\Lambda [0,1]$, almost every sample function 
$B(\cdot,\omega)$ of a standard Brownian motion $B$ has the 
quadratic $\lambda$-variation on $[0,1]$ with the bracket function
$[B(\cdot,\omega)]_{\lambda}(t)=t$, $0\leq t\leq 1$,
and since it is continuous $B(\cdot,\omega)$ cannot have
the local $2$-variation (see Theorem \ref{variation}.\ref{liminf}).
Next we give an example of a non-random function with the same
properties.

Let $f$ be a real-valued function defined on $[0,1]$, let $0<H<1$,
and let $r\geq 4$ be an integer.
Following K\^ono \cite{NK86} and \cite{NK88},
the function $f$ is called \emph{self-affine with the scale parameter
$H$ to base $r$} if for any $i=0,\dots,r^m-1$, $m=1,2,\dots$, there are
$\epsilon_{i,m}\in\{1,-1\}$ and the relation
\beq\label{self-affine}
f(ir^{-m}+h)-f(ir^{-m})=\epsilon_{i,m}r^{-mH}f(r^mh)
\eeq
holds for each $0\leq h< r^{-m}$.
A bounded self-affine function is continuous if and only if 
(\ref{self-affine}) holds for all $0\leq h\leq r^{-m}$.
For any rational number $0<H<1$, 
one can construct a continuous self-affine function of the scale 
parameter $H$  satisfying $0\leq f\leq 1$, $f(0)=0$ and $f(1)=1$.
Typical examples of self-affine functions, such as the coordinate
functions of the Peano curve, can be found in \cite{NK86}.
By the definition, a continuous self-affine function $f$ of the
scale parameter $H$ is H\"older continuous of order $H$, and
for each $t, s\in [0,1]$,
\beq\label{self-affine1}
\big |f(t)-f(s)\big |\leq 5M|t-s|^H,
\eeq
where $M:=\sup_{0\leq t\leq 1}|f(t)|$.
Also, a non-trivial self-affine function with the scale parameter
$H$ does not satisfy the H\"older continuity of order
$H'>H$ at any point (see Theorem 1 in \cite{NK86}).
A more general class of self-affine functions were considered
by Kamae \cite{TK86} and Bertoin \cite{JB88}.

\begin{exmp}\label{NK}
{\rm Let $f$ be a continuous self-affine function on $[0,1]$ with the scale 
parameter $1/2$ to base $r$ for some integer $r\geq 4$, such that $f(1)=1$, 
and let $\lambda_m:=\{ir^{-m}\colon\,i=0,\dots,r^m\}$.
Then $\lambda:=\{\lambda_m\colon\,m\geq 1\}\in\Lambda [0,1]$,
and $f$ has the quadratic $\lambda$-variation with the bracket function
$[f]_{\lambda}(t)=t$ for $0\leq t\leq 1$.
By (\ref{self-affine1}) and Theorem \ref{variation}.\ref{liminf},
it then follows that $f\in {\cal W}_2\setminus {\cal W}_2^{\ast}$.}
\end{exmp}

Let $X=\{X(t)\colon\,t\geq 0\}$ be a square integrable continuous martingale
on a filtered complete probability space  
$(\Omega,{\cal F},\Pr,\{{\cal F}_t\colon\,t\geq 0\})$ satisfying
the usual hypotheses.
Meyer \cite{PAM62},\cite{PAM63} proved that there exists a unique increasing 
stochastic process $[X]=\{[X](t)\colon\,t\geq 0\}$ such that
almost surely,
$$
E\big\{[X(t)-X(s)]^2|{\cal F}_s\big\}=E\big\{[X](t)|{\cal F}_s\big\}
-[X](s)
$$
for $0<s<t<\infty$.
It can be shown that for $0<T<\infty$, almost every sample function of $X$ 
has the quadratic $\lambda$-variation on $[0,T]$ for some 
$\lambda\in\Lambda [0,T]$ as follows. 
We say that a sequence $\{\tau_i\colon\,i=1,2,\dots\}$ of stopping
times is a {\it random partition of $[0,\infty)$
composed by stopping times} if almost surely,
$\tau_0=0$, $\tau_i\leq \tau_{i+1}$ and $\tau_i\to +\infty$ with
$i\to\infty$.
Kunita and Watanabe \cite[Theorem 1.3]{KW67} proved that there exists
a sequence $\lambda=\{\lambda_m\colon\,m\geq 1\}$ of random partitions of 
$[0,\infty)$ composed by stopping times $\lambda_m=\{\tau_i^m\colon\,
i=1,2,\dots\}$ such that almost surely,
$\lambda_m\subset\lambda_{m+1}$, $\max_i|\tau_{i+1}^m-\tau_i^m|
\to 0$ as $m\to\infty$, and
$$
\lim_{m\to\infty}\sum_{i=0}^{\infty}\big [X(\tau_{i+1}^m\wedge t)
-X(\tau_i^m\wedge t)\big ]^2=[X](t)-[X](0)\qquad
\mbox{for}\quad t\geq 0.
$$
For $\omega\in\Omega$ and $0<T<\infty$, we have
$\lambda_m(\omega)\Cap [0,T]=\{\tau_i^m(\omega)\wedge T\colon\,
i=0,1,\dots\}$, $m\geq 1$, and $\lambda (\omega,T):=\{\lambda_m(\omega)
\Cap [0,T] \colon\,m\geq 1\}\in\Lambda [0,T]$.
Therefore by the second part of Proposition \ref{alpha}, we have:

\begin{exmp}
{\rm Let $X$ be a square integrable continuous martingale on
a complete probability space $(\Omega,{\cal F},\Pr)$.
There exist a sequence $\lambda$ of random partitions of
$[0,\infty)$ composed by stopping times and a subset $\Omega_0
\subset\Omega$ such that $\Pr (\Omega_0)=1$, and for each $0<T<\infty$
and $\omega\in\Omega_0$, $X(\cdot,\omega)$ has the
quadratic $\lambda (\omega,T)$-variation on $[0,T]$ with the bracket
function
$$
[X(\cdot,\omega)]_{\lambda (\omega,T)}(t)=[X](t,\omega)-[X](0,\omega)
\qquad\mbox{for}\quad 0\leq t\leq T.
$$}
\end{exmp}

Let $k$ be a positive integer, and let $B$ be a standard Brownian 
motion.
A stochastic process $X=\{X(t)\colon\,0\leq t\leq 1\}$ is called
a {\em weak Brownian motion of order $k$} if for every $k$-tuple
$(t_1,t_2,\dots,t_k)$, the equality
$$
(X(t_1),X(t_2),\dots,X(t_k))=(B(t_1),B(t_2),\dots,B(t_k))
$$
holds in distribution.
If $k\geq 4$ the process $X$ has the same quadratic variation as a 
standard Brownian motion $B$, even so $X$ may not be a semimartingle.
For this and many other properties of a weak Brownian motion see
F\"ollmer, Wu and Yor \cite{FWandY}.

\begin{exmp}
{\rm Let $X=\{X(t)\colon\,0\leq t\leq 1\}$ be a weak Brownian motion
of order $k\geq 4$, and let $\lambda\in\Lambda [0,1]$.
By Proposition 2.1 of \cite{FWandY} and its proof, $X$ admits
a continuous version $\bar X$ on a complete probability space 
$(\Omega,{\cal F},\Pr)$ such that for almost all 
$\omega\in\Omega$, $\bar X(\cdot,\omega)$
has the quadratic $\lambda$-variation with the bracket function
$[\bar X(\cdot,\omega)]_{\lambda}(t)=t$ for $0\leq t\leq 1$.}
\end{exmp}

\paragraph*{Integrals with respect to the bracket function.}
Since the bracket function $[f]_{\lambda}$, if exists,
is nondecreasing, any continuous and bounded function is Riemann-Stieltjes
integrable with respect to $[f]_{\lambda}$.
We show next that any regulated function is Left Young and Right Young
integrable with respect to a function $V$ of bounded variation, and 
prove an integral representation corresponding to the decomposition of $V$ 
into the continuous and pure jump parts (see Definition 
\ref{variation}.\ref{LYandRY} for the definitions of the Left and Right Young 
integrals).

\begin{lem}\label{bv}
Let $\psi$ be a regulated function on $[a,b]$, and
let $V$ be a function on $[a,b]$ of bounded variation such that for 
$a\leq x\leq b$,
$$
V(x)=V^c(x)+\Delta^{+}V(a)+\sum_{(a,x)}\{\Delta^{-}V+\Delta^{+}V\}
+\Delta^{-}V(x), \quad a\leq x\leq b,
$$          
where $V^c$ is a continuous function.
Then $\psi$ is Left Young integrable with respect to $V$, it is 
Riemann-Stieltjes integrable with respect to $V^c$, and 
\beq\label{1bv}
(LY)\int_a^b\psi\,dV=(RS)\int_a^b\psi\,dV^c+\sum_{(a,b]}\psi_{-}
\Delta^{-}V+\sum_{[a,b)}\psi\Delta^{+}V.
\eeq
Also, $\psi$ is Right Young integrable with respect to $V$, it is 
Riemann-Stieltjes integrable with respect to $V^c$, and
\beq\label{3bv}
(RY)\int_a^b\psi\,dV=(RS)\int_a^b\psi\,dV^c+\sum_{(a,b]}\psi
\Delta^{-}V+\sum_{[a,b)}\psi_{+}\Delta^{+}V.
\eeq
\end{lem}

Remark. 
We notice that using the Central Young integral 
in the preceding statement instead of the Left or Right Young integrals,
gives the usual representation:
$$
(CY)\int_a^b\psi\,dV=(RS)\int_a^b\psi\,dV^c+\sum_{(a,b]}\psi\Delta^{-}V
+\sum_{[a,b)}\psi\Delta^{+}V.
$$

\begin{proof}
Due to boundedness of the variation of $V$, the sum in the
definition of $(LY)\smallint_a^b\psi\,dV$ converges absolutely,
so that it is enough to prove the existence of the $RRS$ integral of 
$\psi_{-}^{(a)}$ with respect to $V_{+}^{(b)}$.
Recall that $\Delta_a^{-}V$ and $\Delta_b^{+}V$ are defined on $[a,b]$ so
that $\Delta_a^{-}V(a)=0$ with $\Delta_a^{-}V=\Delta^{-}V$ on $(a,b]$,
and $\Delta_b^{+}V(b)=0$ with $\Delta_b^{+}V=\Delta^{+}V$
on $[a,b)$, respectively.
Then $V_{+}^{(b)}=V^c+U$ on $[a,b]$, where
$U(x):=\sum_{[a,x]}\big\{\Delta_a^{-}V+\Delta_b^{+}V\big\}$ for
$a\leq x\leq b$.
First, we prove that $\psi_{-}^{(a)}$ is $RRS$ integrable with respect 
to $U$ and
\beq\label{2bv}
(RRS)\int_a^b\psi_{-}^{(a)}\,dU=\sum_{(a,b]}\psi_{-}\big\{\Delta^{-}V
+\Delta_b^{+}V\big\}=:A.
\eeq
To this aim let $\epsilon >0$.
There exists a partition $\mu=\{z_j\colon\,j=0,\dots,m\}$
of $[a,b]$ such that
$$
\|\psi\|_{\infty}\sum_{j=1}^m\sum_{(z_{j-1},z_j)}\big\{|\Delta^{-}V|
+|\Delta^{+}V|\big\}<\epsilon\quad\mbox{and}\quad
\Big |A-\sum_{\mu\setminus\{a\}}\psi_{-}\big\{\Delta^{-}V
+\Delta_b^{+}V\big\}\Big |<\epsilon.
$$
Then choose $\{u_j\colon\,j=1,\dots,m\}\subset (a,b)$ such that
$z_{j-1}<u_j<z_j$ for $j=1,\dots,m$ and
$$
\max_{1\leq j\leq m}\sup_{u_j\leq x<z_j}\big |\psi_{-}(x)-\psi_{-}(z_j)
\big |<\epsilon/v_1(V).
$$
Let $\tau=\{([x_{i-1},x_i],y_i)\colon\,i=1,\dots,n\}$ be a tagged
partition of $[a,b]$ such that $\{u_j,z_j\colon\,j=1,\dots,m\}
\subset\{x_i\colon\,i=1,\dots,n\}$.
For $j=1,\dots,m$,
let $i(j)$ be the index in $\{1,\dots,n\}$ such that $x_{i(j)}=z_j$,
and let $I(j):=\{i(j-1)+1,\dots,i(j)\}$.
Then we have
$$
\Big |S_{RS}(\psi_{-}^{(a)},U;\tau)-A\Big |\leq\epsilon
+\Big |\sum_{j=1}^m\sum_{i\in I(j)}\psi_{-}(y_i)\sum_{(x_{i-1},x_i]}
\big\{\Delta^{-}V+\Delta_b^{+}V\big\}-\sum_{\mu\setminus\{a\}}
\psi_{-}\big\{\Delta^{-}V+\Delta_b^{+}V\big\}\Big |
$$
$$
\leq\epsilon +\|\psi\|_{\infty}\sum_{j=1}^m\sum_{(z_{j-1},z_j)}
\big\{|\Delta^{-}V|+|\Delta^{+}V|\big\}+\sum_{j=1}^m\big |\psi_{-}
(y_{i(j)})-\psi_{-}(z_j)\big |
\big |\Delta^{-}V(z_j)+\Delta^{+}V(z_j)\big |<3\epsilon.
$$
Therefore $\psi_{-}^{(a)}$ is $RRS$ integrable with respect to $U$,
and (\ref{2bv}) holds.
Since $\psi$ is regulated, we have that $(RS)\smallint_a^b\psi\,dV^c$ 
exists and  $(RS)\smallint_a^b\Delta_{a}^{-}\psi\,dV^c=0$.
This in conjunction with (\ref{2bv}) yields that $\psi_{-}^{(a)}$
is $RRS$ integrable with respect to $V_{+}^{(b)}$ and
$$
(RRS)\int_a^b\psi_{-}^{(a)}\,dV_{+}^{(b)}=(RS)\int_a^b\psi\,dV^c
+\sum_{(a,b]}\psi_{-}\big\{\Delta^{-}V+\Delta_b^{+}V\big\}.
$$
Adding the preceding sum with the sum from the definition of
the Left Young integral we get
$$
\sum_{(a,b]}\psi_{-}\big\{\Delta^{-}V+\Delta_b^{+}V\big\}
+\big\{\psi\Delta^{+}V\big\}(a)+\sum_{(a,b)}\Delta^{-}\psi\Delta^{+}V
=\sum_{(a,b]}\psi_{-}\Delta^{-}V+\sum_{[a,b)}\psi\Delta^{+}V.
$$
Therefore $\psi$ is Left Young integrable with respect to $V$,
and (\ref{1bv}) holds.

The proof of the second part of the lemma is similar.
Indeed, using the same notation as before we have $V_{-}^{(a)}=V^c+U$, 
where $U(a):=0$ and $U(x):=\sum_{[a,x)}\{\Delta_a^{-}V+\Delta_b^{+}V\}$
for $a<x\leq b$.
Then one can show similarly that $\psi_{+}^{(b)}$ is $RRS$ integrable with 
respect to  $U$, and
$$
(RRS)\int_a^b\psi_{+}^{(b)}\,dU=\sum_{[a,b)}\psi_{+}\big\{
\Delta_a^{-}V+\Delta^{+}V\big\}.
$$
Again, adding the preceding sum with the sum from the definition of
the Right Young integral we get
$$
\sum_{[a,b)}\psi_{+}\big\{\Delta_a^{-}V+\Delta^{+}V\big\}
+\big\{\psi\Delta^{-}V\big\}(b)-\sum_{(a,b)}\Delta^{+}\psi\Delta^{-}V
=\sum_{(a,b]}\psi\Delta^{-}V+\sum_{[a,b)}\psi_{+}\Delta^{+}V.
$$
This yields the existence of the Right Young integral and relation
(\ref{3bv}).
The proof of Lemma \ref{bv} is now complete.
\qed\end{proof}

To prove an approximation of an integral with respect to a bracket function,
we use the following convergence theorem; its proof is given in
\cite[Theorem II.15.4]{THH63}.

\begin{lem}\label{conv-thm}
Let $g$ and $\{g_m\colon\,m\geq 1\}$ be nondecreasing functions 
on $[a,b]$ such that for each $a\leq x\leq b$, $g_m(x)\to g(x)$ as $m\to\infty$.
If $(RRS)\smallint_a^bh\,dg$ and $(RRS)\smallint_a^bh\,dg_m$, $m\geq 1$, 
all exist, then 
$$
\lim_{m\to\infty}(RRS)\int_a^bh\,dg_m=(RRS)\int_a^bh\,dg.
$$
\end{lem}

The next approximation theorem extends  property (\ref{3alpha})
in the definition of the quadratic $\lambda$-variation.

\begin{thm}\label{intqv}
Let $\lambda\in\Lambda [a,b]$,  let $f$ have the quadratic 
$\lambda$-variation on $[a,b]$, and let $h$ be a regulated function on $[a,b]$.
If the right discontinuity points of $f$ are accessible by $\lambda$, then
\beq\label{1intqv}
\lim_{m\to\infty}\sum_{i=1}^{n(m)}h(x_{i-1}^m)\big [f(x_i^m)
-f(x_{i-1}^m)\big ]^2=(LY)\int_a^bh\,d[f]_{\lambda},
\eeq
where $\lambda_m=\{x_i^m\colon\,i=0,\dots,n(m)\}$, $m\geq 1$,
constitute $\lambda$.
Also with the same notation,
if the left discontinuity points of $f$ are accessible by $\lambda$, then
\beq\label{6intqv}
\lim_{m\to\infty}\sum_{i=1}^{n(m)}h(x_{i}^m)\big [f(x_i^m)
-f(x_{i-1}^m)\big ]^2=(RY)\int_a^bh\,d[f]_{\lambda}.
\eeq
\end{thm}

\begin{proof}
To prove (\ref{1intqv}) let $\lambda\in\Lambda [a,b]$ be such that $f$ has 
the quadratic $\lambda$-variation and the right discontinuity points of $f$ 
are accessible by $\lambda$, that is, (\ref{right-jumps-access}) holds.
Since the bracket function $[f]:=[f]_{\lambda}$ has bounded
variation, the Left Young integral in (\ref{1intqv})  exists by Lemma \ref{bv}.
We show that each sum on the left side of (\ref{1intqv}) is the Left Young integral 
with respect to a suitable nondecreasing step function, and use Lemma
\ref{conv-thm} to prove (\ref{1intqv}).
For each $m\geq 1$, let $F_m$ be the nondecreasing step function on $[a,b]$
defined by
$$
F_m(y):=\left\{ \begin{array}{ll}
    0, &\mbox{if $y=a$}\\
    \sum_{i=1}^k\big [f(x_i^m)-f(x_{i-1}^m)\big ]^2, &\mbox{if
           $y\in (x_{k-1}^m,x_k^m]$ for some $k=1,\dots,n(m)$.}
\end{array}
\right.
$$
For each $m\geq 1$ and $1\leq k\leq n(m)$,  the function 
$(F_m)_{+}^{(x_k^m)}$ on $[x_{k-1}^m,x_k^m]$ (see (\ref{rmodification})
for the notation) is constant,
and $\Delta^{+}F_m=0$ on $(x_{k-1}^m,x_k^m)$.
Thus $h$ is Left Young integrable with respect to $F_m$
on $[x_{k-1}^m,x_k^m]$, and the equality
\begin{eqnarray*}
(LY)\int_{x_{k-1}^m}^{x_k^m}h\,dF_m&=&(RRS)\int_{x_{k-1}^m}^{x_k^m}
h_{-}^{(x_{k-1}^m)}\,d(F_m)_{+}^{(x_k^m)}\\
& &\, +\{h\Delta^{+}F_m\}(x_{k-1}^m)+\sum_{(x_{k-1}^m,x_k^m)}
\Delta^{-}h\Delta^{+}F_m\\
&=& h(x_{k-1}^m)\big [f(x_k^m)-f(x_{k-1}^m)]^2
\end{eqnarray*}
holds for each $m\geq 1$ and $1\leq k\leq n(m)$.
Since the Left Young integral is additive over consecutive intervals
(Theorem \ref{variation}.\ref{additivity}), $h$ is Left Young integrable 
with respect to $F_m$ over $[a,b]$, and the equality
$$
\sum_{i=1}^{n(m)}h(x_{i-1}^m)\big [f(x_i^m)-f(x_{i-1}^m)\big ]^2
=(LY)\int_a^bh\,dF_m
$$
holds for each $m\geq 1$.
Therefore to prove (\ref{1intqv}) it is enough to show that
\beq\label{3intqv}
\lim_{m\to\infty}(LY)\int_a^bh\,dF_m=(LY)\int_a^bh\,d[f].
\eeq
First, we show that for each $a\leq y\leq b$,
\beq\label{2intqv}
\lim_{m\to\infty}F_m(y)=[f](y).
\eeq
Indeed this follows just by the definitions of the functions in
the case $y\in\cup_m\lambda_m$.
Otherwise, if $y\in [a,b]\setminus\cup_m\lambda_m$, then $x_{k-1}^m<y<x_k^m$
for some $k=k(m)$, and
\begin{eqnarray*}
\limsup_{m\to\infty}\big |F_m(y)-[f](y)\big |
&\leq &\limsup_{m\to\infty}\big |[f(x_k^m)-f(x_{k-1}^m)]^2-
[f(y)-f(x_{k-1}^m)]^2\big |\\
&=&\left |[f(y+)-f(y-)]^2-[f(y)-f(y-)]^2\right |=0
\end{eqnarray*}
because $f$ is right-continuous at $y$ in this case by (\ref{right-jumps-access}).
Therefore (\ref{2intqv}) holds for each $a\leq y\leq b$.
Second, in addition to (\ref{2intqv}), for $a\leq y<b$,
\beq\label{4intqv}
\lim_{m\to\infty}F_m(y+)=[f](y+).
\eeq
Indeed, if $y\in\cup_m\lambda_m$ then $y=x_k^m$ for all large enough $m$
and some $1\leq k<n(m)$, so that $F_m(x_k^m+)=F_m(x_k^m)+
[f(x_{k+1}^m)-f(x_k^m)]^2\to [f](y)+(\Delta^{+}f(y))^2
=[f](y+)$ as $m\to\infty$.
If $y\in [a,b)\setminus\cup_m\lambda_m$ then $x_{k-1}^m<y<x_k^m$ and
$F_m(y+)=F_m(y)\to [f](y)=[f](y+)$ as $m\to\infty$ because
$f$ is right-continuous at $y$ in this case by (\ref{right-jumps-access}).
Therefore (\ref{4intqv}) holds for each $y\in [a,b)$.
This in conjunction with Lemma \ref{conv-thm} yields that
$$
\lim_{m\to\infty}(RRS)\int_a^bh_{-}^{(a)}\,d(F_m)_{+}^{(b)}
=(RRS)\int_a^bh_{-}^{(a)}\,d[f]_{+}^{(b)}.
$$
Since $\lim_{m\to\infty}\Delta^{+}F_m(a)=\Delta^{+}[f](a)$,
we are left to show that
\beq\label{5intqv}
\lim_{m\to\infty}\sum_{(a,b)}\Delta^{-}h\Delta^{+}F_m=\sum_{(a,b)}
\Delta^{-}h\Delta^{+}[f].
\eeq
Given $\epsilon >0$ choose a finite set $\nu\subset (a,b)$ such that
$$
\Big |\sum_{(a,b)}\Delta^{-}h\Delta^{+}[f]-\sum_{\nu}
\Delta^{-}h\Delta^{+}[f]\Big |<\epsilon
\quad\mbox{and}\quad
|\Delta^{-}h(y)|<\epsilon  
$$
for each $y\not\in \nu\cap(\cup_m\lambda_m)$.
For each $m\geq 1$, let $\mu_m:=\{x_k^m\colon\,k=1,\dots,n(m)-1\}
\setminus\nu$ and $\nu_m:=\{x_k^m\colon\,k=1,\dots,n(m)-1\}
\cap\nu$.
Then we have
$$
\limsup_{m\to\infty}\Big |\sum_{(a,b)}\Delta^{-}h\Delta^{+}F_m
-\sum_{(a,b)}\Delta^{-}h\Delta^{+}[f]\Big |
$$
$$
\leq\epsilon +\limsup_{m\to\infty}\Big |\sum_{x_k^m\in\mu_m}
\Delta^{-}h(x_k^m)\big [f(x_{k+1}^m)-f(x_k^m)\big ]^2\Big |
$$
$$
+\limsup_{m\to\infty}\Big |\sum_{x_k^m\in\nu_m}\Delta^{-}h(x_k^m)
\Big [\big [f(x_{k+1}^m)-f(x_k^m)\big ]^2-\big\{\Delta^{+}f(x_k^m)
\big\}^2\Big ]\Big |\leq \epsilon +\epsilon [f](b).
$$
Since $\epsilon$ is arbitrary (\ref{5intqv}), and hence (\ref{3intqv}) holds.
The proof of (\ref{1intqv}) is complete.

Next we prove (\ref{6intqv}).
For each $m\geq 1$, let $G_m$ be  a nondecreasing step function 
on $[a,b]$ defined by
$$
G_m(y):=\left\{ \begin{array}{lll}
    0, &\mbox{if $y\in [a,y_1)$}\\
    \sum_{i=1}^{k-1}\big [f(x_i^m)-f(x_{i-1}^m)\big ]^2, &\mbox{if
           $y\in [x_{k-1}^m,x_k^m)$ for some $k=2,\dots,n(m)$}\\
    \sum_{i=1}^{n(m)}\big [f(x_i^m)-f(x_{i-1}^m)\big ]^2, &\mbox{if
           $y=b$.}\\
\end{array}
\right.
$$
For each $m\geq 1$ and $1\leq k\leq n(m)$, the function 
$(G_m)_{-}^{(x_{k-1}^m)}$ on $[x_{k-1}^m,x_k^m]$ (see (\ref{lmodification})
for the notation) is constant,
and $\Delta^{-}G_m=0$ on $(x_{k-1}^m,x_k^m)$.
Thus $h$ is Right Young integrable with respect to $G_m$
on $[x_{k-1}^m,x_k^m]$, and the equality
\begin{eqnarray*}
(RY)\int_{x_{k-1}^m}^{x_k^m}h\,dG_m&=&(RRS)\int_{x_{k-1}^m}^{x_k^m}
h_{+}^{(x_{k}^m)}\,d(G_m)_{-}^{(x_{k-1}^m)}\\
& &\, -\sum_{(x_{k-1}^m,x_k^m)}\Delta^{+}h\Delta^{-}G_m
+\big (h\Delta^{-}G_m\big )(x_{k}^m)\\
&=& h(x_{k}^m)\big [f(x_k^m)-f(x_{k-1}^m)]^2
\end{eqnarray*}
holds for each $m\geq 1$ and $1\leq k\leq n(m)$.
Since the Right Young integral is additive over consecutive intervals
(Theorem \ref{variation}.\ref{additivity}), $h$ is Right Young integrable 
with respect to $G_m$ over $[a,b]$, and the equality
$$
\sum_{i=1}^{n(m)}h(x_{i}^m)\big [f(x_i^m)-f(x_{i-1}^m)\big ]^2
=(RY)\int_a^bh\,dG_m
$$
holds for each $m\geq 1$.
Therefore to prove (\ref{6intqv}) it is enough to show that
\beq\label{7intqv}
\lim_{m\to\infty}(RY)\int_a^bh\,dG_m=(RY)\int_a^bh\,d[f].
\eeq
First, we show that for each $a\leq y\leq b$,
\beq\label{8intqv}
\lim_{m\to\infty}G_m(y)=[f](y).
\eeq
Indeed this follows just by the definitions of the functions in
the case $y\in\cup_m\lambda_m$.
Otherwise, if $y\in [a,b]\setminus\cup_m\lambda_m$, then $x_{k-1}^m<y<x_k^m$
for some $k=k(m)$, and
$$
\limsup_{m\to\infty}\big |G_m(y)-[f](y)\big |
\leq \limsup_{m\to\infty}\big |[f(y)-f(x_{k-1}^m)]^2
=(\Delta^{-}f(y))^2=0
$$
because $f$ is left-continuity at $y$ in this case by (\ref{left-jumps-access}).
Therefore (\ref{8intqv}) holds for each $a\leq y\leq b$.
Second, in addition to (\ref{8intqv}), for $a< y\leq b$,
\beq\label{9intqv}
\lim_{m\to\infty}G_m(y-)=[f](y-).
\eeq
Indeed, if $y\in\cup_m\lambda_m$ then $y=x_k^m$ for all large enough $m$
and some $1\leq k<n(m)$, so that $G_m(x_k^m-)=G_m(x_k^m)-
[f(x_{k}^m)-f(x_{k-1}^m)]^2\to [f](y)-(\Delta^{-}f(y))^2
=[f](y-)$ as $m\to\infty$.
If $y\in (a,b]\setminus\cup_m\lambda_m$ then $x_{k-1}^m<y<x_k^m$ and
$G_m(y-)=G_m(y)\to [f](y)=[f](y-)$ as $m\to\infty$ because
$f$ is left-continuity at $y$ in this case by (\ref{left-jumps-access}).
Therefore (\ref{9intqv}) holds for each $a\leq y<b$.
This in conjunction with Lemma \ref{conv-thm} yields that
$$
\lim_{m\to\infty}(RRS)\int_a^bh_{+}^{(b)}\,d(G_m)_{-}^{(a)}
=(RRS)\int_a^bh_{+}^{(b)}\,d[f]_{-}^{(a)}.
$$
Since $\lim_{m\to\infty}\Delta^{-}G_m(b)=\Delta^{-}[f](b)$,
we are left to show that
\beq\label{10intqv}
\lim_{m\to\infty}\sum_{(a,b)}\Delta^{+}h\Delta^{-}G_m=\sum_{(a,b)}
\Delta^{+}h\Delta^{-}[f].
\eeq
Given $\epsilon >0$ choose a finite set $\nu\subset (a,b)$ such that
$$
\Big |\sum_{(a,b)}\Delta^{+}h\Delta^{-}[f]-\sum_{\nu}
\Delta^{+}h\Delta^{-}[f]\Big |<\epsilon
\quad\mbox{and}\quad
|\Delta^{+}h(y)|<\epsilon  
$$
for each $y\not\in \nu\cap(\cup_m\lambda_m)$.
For each $m\geq 1$, let $\mu_m:=\{x_k^m\colon\,k=1,\dots,n(m)-1\}
\setminus\nu$ and $\nu_m:=\{x_k^m\colon\,k=1,\dots,n(m)-1\}
\cap\nu$.
Then we have
$$
\limsup_{m\to\infty}\Big |\sum_{(a,b)}\Delta^{+}h\Delta^{-}G_m
-\sum_{(a,b)}\Delta^{+}h\Delta^{-}[f]\Big |
$$
$$
\leq\epsilon +\limsup_{m\to\infty}\Big |\sum_{x_k^m\in\mu_m}
\Delta^{+}h(x_k^m)\big [f(x_{k}^m)-f(x_{k-1}^m)\big ]^2\Big |
$$
$$
+\limsup_{m\to\infty}\Big |\sum_{x_k^m\in\nu_m}\Delta^{+}h(x_k^m)
\Big [\big [f(x_{k}^m)-f(x_{k-1}^m)\big ]^2-\Delta^{-}f^2(x_k^m)
\Big ]\Big |\leq \epsilon +\epsilon [f](b).
$$
Since $\epsilon$ is arbitrary (\ref{10intqv}), 
and hence (\ref{7intqv}) holds.
The proof of Theorem \ref{intqv} is complete.
\qed\end{proof}

\paragraph*{The quadratic $\lambda$-covariation.}
Let $f_1$ and $f_2$ be two real-valued functions defined
on an interval $[a,b]$.
For a partition $\kappa=\{x_i\colon\,i=0,\dots,n\}$ of $[a,b]$, let
$$
C(f_1,f_2;\kappa):=\sum_{i=1}^{n}
\big [f_1(x_i)-f_1(x_{i-1})\big ]
\big [f_2(x_i)-f_2(x_{i-1})\big ].
$$
Next we extend the notion of the quadratic $\lambda$-variation of
a single function to a pair of functions by replacing  
sums $s_2$ of squares of increments in its definition by sums 
$C$ just defined.

\begin{defn}\label{8bracket}
{\rm Let $f_1$, $f_2$ be regulated functions on $[a,b]$, 
and let $\lambda=\{\lambda_m\colon\,m\geq 1\}\in\Lambda [a,b]$.
We say that the $2$-vector function $(f_1,f_2)$ has the \emph{quadratic
$\lambda$-covariation on $[a,b]$} if there exists an additive
upper continuous function $\alpha_{\lambda}(f_1,f_2)$ on 
$S\loc a,b\roc$ such that for each $(u,v)\in S\loc a,b\roc$
$$
\alpha_{\lambda}(f_1,f_2;u,v)=\lim_{m\to\infty}C(f_1,f_2;
\lambda_m\Cap \loc u,v\roc ).
$$
If the $2$-vector function $(f_1,f_2)$ has the quadratic $\lambda$-covariation 
on $[a,b]$ then the right distribution function of $\alpha_{\lambda}(f_1,f_2)$
restricted to $[a,b]$ is called the \emph{bracket function} of 
$(f_1,f_2)$, and is denoted by $[f_1,f_2]_{\lambda}$.}
\end{defn}

As in the case of the quadratic $\lambda$-variation, the quadratic
$\lambda$-covariation on $[a,b]$ can be defined solely in terms of a 
regulated function on $[a,b]$.
The proof of the following statement is similar to the proof
of the Proposition \ref{alpha}, and therefore is omitted.

\begin{prop}\label{ext-alpha}
Let $f_1$, $f_2$ be regulated functions on $[a,b]$,
and let $\lambda=\{\lambda_m\colon\,m\geq 1\}\in\Lambda [a,b]$.
The $2$-vector function $(f_1,f_2)$ has the quadratic $\lambda$-covariation
on $[a,b]$ if and only if there exists a regulated function $B=B_{\lambda}$
on $[a,b]$ such that $B(a)=0$ and for each $(s,t)\in S [a,b]$,
\beq\label{11bracket}
B(t)-B(s)=\lim_{m\to\infty}C(f_1,f_2;\lambda_m \Cap [s,t]),
\eeq
\beq\label{10bracket}
\Delta^{-}B(t)=\{\Delta^{-}f_1\Delta^{-}f_2\}(t)
\qquad\mbox{and}\qquad
\Delta^{+}B(s)=\{\Delta^{+}f_1\Delta^{+}f_2\}(s).
\eeq
If the two statements hold then $[f_1,f_2]_{\lambda}=B$.
Moreover, condition {\rm (\ref{11bracket})} for each $(s,t)\in S[a,b]$,
can be replaced by the condition
$$
B(t)=\lim_{m\to\infty}C(f_1,f_2;\lambda_m \Cap [a,t])
$$
for each $t\in [a,b]$, if in addition, $\lambda$ is such that
\beq\label{2qcforpq}
\{x\in (a,b)\colon\,\big (\Delta^{-}f_1\Delta^{+}f_2 + \Delta^{+}f_1
\Delta^{-}f_2\big )(x)\not =0\}\subset\cup_m\lambda_m.
\eeq
\end{prop}

If in the preceding definition $f_1=f_2=f$, then the $2$-vector
function $(f,f)$ has the quadratic $\lambda$-covariation if and only if
$f$ has the quadratic $\lambda$-variation,
and then $[f,f]_{\lambda}=[f]_{\lambda}$.
In Section \ref{vector} we extend the notion
of the quadratic $\lambda$-variation of a real-valued function to 
a vector valued function.
The notion of the quadratic $\lambda$-covariation of a $2$-vector
function is used to define the quadratic $\lambda$-variation of a 
$d$-vector valued function with $d\geq 2$.

For notation simplicity let ${\cal W}_{\infty}^{\ast}[a,b]:=
{\cal R}[a,b]$.

\begin{prop}\label{qcforpq}
For $p^{-1}+q^{-1}=1$ and $p\geq 1$,
let $f_1\in {\cal W}_p[a,b]$ and $f_2\in {\cal W}_q^{\ast}[a,b]$, and
let $\lambda=\{\lambda_m\colon\,m\geq 1\}\in\Lambda [a,b]$ be 
such that {\rm (\ref{2qcforpq})} holds.
Then the $2$-vector function $(f_1,f_2)$ has the quadratic $\lambda$-covariation
on $[a,b]$, and its bracket function $[f_1,f_2]_{\lambda}$ is a pure jump function
given by
\beq\label{1qcforpq}
[f_1,f_2]_{\lambda}(y)=\sum_{(a,y]}\Delta^{-}f_1\Delta^{-}f_2
+\sum_{[a,y)}\Delta^{+}f_1\Delta^{+}f_2,\quad
\mbox{for $a\leq y\leq b$.}
\eeq
\end{prop}

\begin{proof}
By H\"older's inequality, the two sums $\sum_{(a,b]}\Delta^{-}f_1
\Delta^{-}f_2$ and $\sum_{[a,b)}\Delta^{+}f_1\Delta^{+}f_2$
converge absolutely.
Therefore the function on the right side of (\ref{1qcforpq}) is of 
bounded variation, and its jumps agree with those of the bracket
function of $(f_1,f_2)$ if it exists (cf.\ (\ref{10bracket})).
Let $\lambda_m=\{x_i^m\colon\,i=0,\dots,n(m)\}$ for each $m\geq 1$,
and let $a\leq u<v\leq b$.
We have to prove that
\beq\label{3qcforpq}
\lim_{m\to\infty}C(f_1,f_2; \lambda_m\Cap [u,v])
=\sum_{(u,v]}\Delta^{-}f_1\Delta^{-}f_2
+\sum_{[u,v)}\Delta^{+}f_1\Delta^{+}f_2.
\eeq
Let $\epsilon >0$.
There exists a partition $\kappa=\{z_k\colon\,k=0,\dots,l\}$ of $[u,v]$ 
such that
$$
\Big |\sum_{\kappa\cap (u,v]}\Delta^{-}f_1\Delta^{-}f_2
+\sum_{\kappa\cap [u,v)}\Delta^{+}f_1\Delta^{+}f_2
-\sum_{(u,v]}\Delta^{-}f_1\Delta^{-}f_2
-\sum_{[u,v)}\Delta^{+}f_1\Delta^{+}f_2\Big |<\epsilon
$$
and either $\max_{1\leq k\leq l}\osc (f_2;(z_{k-1},z_k))
<\epsilon /v_1(f_1)$ if $p=1$, or
$$
\sum_{k=1}^lv_q(f_2;(z_{k-1},z_k))< \epsilon^q/v_p(f_1)^{q/p}
\qquad\mbox{if $p>1$}.
$$
The last relation holds by Lemma \ref{variation}.\ref{lemma1}.
Let $I:=\{k=1,\dots,l-1\colon\,z_k\in\cup_m\lambda_m\}$. 
If $k\in I$ then $z_k=x_{i(k)}^m$ for some 
$i(k)\in\{0,\dots,n(m)\}$ and for all sufficiently large $m$.
Thus for $k\in I$, we have
$$
\lim_{m\to\infty}\big [\Delta_{i(k)}^mf_1\Delta_{i(k)}^mf_2
+\Delta_{i(k)+1}^mf_1\Delta_{i(k)+1}^mf_2\big ]=\big (\Delta^{-}f_1
\Delta^{-}f_2\big )(z_k)+\big (\Delta^{+}f_1\Delta^{+}f_2\big )(z_k),
$$
where $\Delta_j^m\rho:=\rho (x_j^m)-\rho (x_{j-1}^m)$ if $x_j^m$
is defined, and $:=0$ otherwise.
If $k\in\{1,\dots,l-1\}\setminus I$ then $x_{i(k)-1}^m<z_k<x_{i(k)}^m$
for some $i(k)\in\{1,\dots,n(m)\}$ and for all sufficiently large $m$.
By assumption (\ref{2qcforpq}), it follows that
$$
\big [f_1(z_k+)-f_1(z_k-)\big ]\big [f_2(z_k+)-f_2(z_k-)\big ]
=\big (\Delta^{-}f_1\Delta^{-}f_2\big )(z_k)
+\big (\Delta^{+}f_1\Delta^{+}f_2\big )(z_k).
$$
Thus for $k\in\{1,\dots,l-1\}\setminus I$, we have
$$
\lim_{m\to\infty}\Delta_{i(k)}^mf_1\Delta_{i(k)}^mf_2
=\big (\Delta^{-}f_1\Delta^{-}f_2\big )(z_k)
+\big (\Delta^{+}f_1\Delta^{+}f_2\big )(z_k).
$$
For $k=0$ and $k=l$, that is at the endpoints of $[u,v]$,
recalling notation $\lmuv=\{u,x_{i(u)}^m,\dots,x_{i(v)}^m,v\}$, we have
$$
\lim_{m\to\infty}\big [f_1(x_{i(u)}^m)-f_1(u)\big ]
\big [f_2(x_{i(u)}^m)-f_2(u)\big ]
=\big (\Delta^{+}f_1\Delta^{+}f_2\big )(u)
$$
and
$$
\lim_{m\to\infty}\big [f_1(v)-f_1(x_{i(v)}^m)\big ]
\big [f_2(v)-f_2(x_{i(v)}^m)\big ]
=\big (\Delta^{-}f_1\Delta^{-}f_2\big )(v).
$$
Let $m_0$ be the first integer such that $\{z_k\colon\,k\in I\}
\subset\lambda_{m_0}$
and each pair $(z_{k-1},z_k)$ of points from $\kappa$ is separated
by at least two points from $\lambda_{m_0}$.
Let $x_{i(u)-1}:=u$ and $x_{i(v)+1}:=v$.
For $m\geq m_0$, let $I_1(m)$ be the set of all indices $i(k), i(k)+1$
with $k\in I$, all indices $i(k)$ with $k\in\{1,\dots,l-1\}\setminus I$,
and the indices $i(u)$, $i(v)+1$.
Then we have
$$
\lim_{m\to\infty}\sum_{i\in I_1(m)}\Delta_i^mf_1\Delta_i^mf_2
=\sum_{\kappa\cap (u,v]}\Delta^{-}f_1\Delta^{-}f_2
+\sum_{\kappa\cap [u,v)}\Delta^{+}f_1\Delta^{+}f_2.
$$
For $m\geq m_0$, let $I_2(m):=\{i(u),\dots,i(v)+1\}\setminus I_1(m)$.
Then by H\"older's inequality, we have
$$
\Big |\sum_{i\in I_2(m)}\Delta_i^mf_1\Delta_i^mf_2
\Big |\leq\Big (\sum_{k=1}^lv_q(f_2;(z_{k-1},z_k))\Big )^{1/q}
v_p(f_1)^{1/p}<\epsilon
$$
for all $m\geq m_0$ if $p>1$.
If $p=1$ then we simply have the bound
$$
\Big |\sum_{i\in I_2(m)}\Delta_i^mf_1\Delta_i^mf_2
\Big |\leq\max_{1\leq k\leq l}\osc (f_2;(z_{k-1},z_k))
v_1(f_1)<\epsilon
$$
for all $m\geq m_0$.
Therefore
$$
\limsup_{m\to\infty}\Big |C(f_1,f_2;\lambda_m\Cap [u,v])
-\sum_{(a,y]}\Delta^{-}f_1\Delta^{-}f_2
-\sum_{[a,y)}\Delta^{+}f_1\Delta^{+}f_2\Big |\leq 2\epsilon.
$$
Since $\epsilon>0$ is arbitrary (\ref{3qcforpq}) holds,
proving Proposition \ref{qcforpq}.
\qed\end{proof}

In Proposition \ref{qcforpq} taking $f_1=f_2$ we get:

\begin{cor}\label{p-var-qv-var}
Let $\lambda\in\Lambda [a,b]$.
If $f$ has the local $2$-variation on $[a,b]$ and two-sided discontinuity
points of $f$ are accessible by $\lambda$,
then $f$ has the quadratic $\lambda$-variation, and its bracket function
$[f]_{\lambda}$ is a pure jump function given by
$$
[f]_{\lambda}(y)=\sum_{(a,y]}\big (\Delta^{-}f\big )^2
+\sum_{[a,y)}\big (\Delta^{+}f\big )^2,\qquad
\mbox{for $a\leq y\leq b$.}
$$
\end{cor}

For the next corollary recall Definition \ref{introduction}.\ref{decomposable}
of the class ${\cal D}_{\lambda,p}$ of $(\lambda,p)$-decomposable functions.

\begin{cor}\label{qvclassQ}
Let $f\in {\cal D}_{\lambda,p}[a,b]$ for some $\lambda\in\Lambda [a,b]$
and $1\leq p<2$.
If $(g,h)\in D_{\lambda,p}^{+}(f)\cup D_{\lambda,p}^{-}(f)\not =\emptyset$
then $f$ has the quadratic $\lambda$-variation,
and its bracket function is given by
\beq\label{1qvclassQ}
[f]_{\lambda}(x)=[g]_{\lambda}^c(x)+\sum_{(a,x]}\big\{\Delta^{-}f\big\}^2
+\sum_{[a,x)}\big\{\Delta^{+}f\big\}^2,\quad
\mbox{for $a\leq x\leq b$.}
\eeq
\end{cor}

\begin{proof}
Let $\lambda=\{\lambda_m\colon\,m\geq 1\}$ and $1\leq p<2$ be such that
$f\equiv C+g+h$ for some $(g,h)\in D_{\lambda,p}^{+}(f)$. 
Since
\beq\label{2qvclassQ}
\{x\in (a,b)\colon\,(\Delta^{-}h\Delta^{+}h)(x)\not =0\}
\subset\{x\in (a,b)\colon\,(\Delta^{+}h)(x)\not =0\}
\subset\cup_m\lambda_m,
\eeq
by the preceding corollary, $h$ has the quadratic $\lambda$-variation with the 
bracket function
$$
[h]_{\lambda}(x)=\sum_{(a,x]}\{\Delta^{-}h\}^2+\sum_{[a,x)}
\{\Delta^{+}h\}^2,\quad
\mbox{for $a\leq x\leq b$.}
$$
Since $g\in {\cal W}_q^{\ast}[a,b]$ with $p^{-1}+q^{-1}=1$, and
\begin{eqnarray}
\lefteqn{\{x\in (a,b)\colon\,\big (\Delta^{-}g\Delta^{+}h+\Delta^{+}g\Delta^{-}h
\big )(x)\not =0\}}\label{3qvclassQ}\\[2mm]
&\subset&\{x\in (a,b)\colon\,\big (\Delta^{+}g\big )(x)\not =0\}
\cup\{x\in (a,b)\colon\,\big (\Delta^{+}h\big )(x)\not =0\}
\subset\cup_m\lambda_m,\nonumber
\end{eqnarray}
by Proposition \ref{qcforpq} applied to $f_1=g$ and $f_2=h$, 
the $2$-vector function $(g,h)$ has the quadratic $\lambda$-covariation 
with the bracket function
$$
[g,h]_{\lambda}(x)=\sum_{(a,x]}\Delta^{-}g\Delta^{-}h
+\sum_{[a,x)}\Delta^{+}g\Delta^{+}h,\quad
\mbox{for $a\leq x\leq b$.}
$$
Let $H$ be the function on $[a,b]$ defined by the right side
of (\ref{1qvclassQ}).
Since $g$ has the quadratic $\lambda$-variation, it follows
that for $(s,t)\in S[a,b]$,
\begin{eqnarray*}
\lefteqn{\lim_{m\to\infty}s_2(f;\lambda_m\cap [s,t])}\\[2mm]
&=&[g]_{\lambda}(t)-[g]_{\lambda}(s)+2\Big\{\sum_{(s,t]}\Delta^{-}g\Delta^{-}h
+\sum_{[s,t)}\Delta^{+}g\Delta^{+}h\Big\}
+\sum_{(s,t]}\{\Delta^{-}h\}^2+\sum_{[s,t)}\{\Delta^{+}h\}^2\\[2mm]
&=&[g]_{\lambda}^c(t)-[g]_{\lambda}^c(s)
+\sum_{(s,t]}\{\Delta^{-}(g+h)\}^2+\sum_{[s,t)}\{\Delta^{+}(g+h)\}^2
=H(t)-H(s).
\end{eqnarray*}
Also since $\Delta^{-}H(t)=\{\Delta^{-}f(t)\}^2$
for $a<t\leq b$, and $\Delta^{+}H(s)=\{\Delta^{+}f(s)\}^2$
for $a\leq s<b$, by Proposition \ref{alpha},
$f$ has the quadratic $\lambda$-variation with the bracket
function (\ref{1qvclassQ}) in the case $(g,h)\in D_{\lambda}^{+}(f)\not =
\emptyset$.
If $(g,h)\in D_{\lambda}^{-}(f)\not =\emptyset$ then the proof is the same
because (\ref{2qvclassQ}) and (\ref{3qvclassQ}) hold in this case
with $\Delta^{+}$ in the middle terms replaced by $\Delta^{-}$.
The proof of Corollary \ref{qvclassQ} is complete.
\qed\end{proof}

In (\ref{1qvclassQ}) we have that $[g]_{\lambda}^c$ does not depend
on a $(\lambda,p)$-decomposition $f\equiv C+g+h$ such that 
$(g,h)\in D_{\lambda,p}^{+}(f)\cup D_{\lambda,p}^{-}(f)$.
Therefore for $f\in {\cal D}_{\lambda,p}[a,b]$ with $\lambda\in\Lambda [a,b]$
and $1\leq p<2$,
such that  $D_{\lambda,p}^{+}(f)\cup D_{\lambda,p}^{-}(f)\not =\emptyset$,
the continuous part of the bracket function $[f]_{\lambda}$ is uniquely
defined.

\section{The Left and Right $\lambda$-integrals}\label{l-integrals}

In this section we define two integrals with respect to a 
$(\lambda,p)$-decomposable function.
Let $f\equiv C+g+h$, where $(g,h)$ is a $(\lambda,p)$-dual pair for some
$\lambda\in\Lambda [a,b]$ and $1\leq p<2$.
The two integrals with respect to $f$ are defined to be a sum
of an integral with respect to $g$ and an integral with respect
to $h$.
Since $h$ has bounded $p$-variation with $1\leq p<2$, we use the Left 
Young and Right Young integrals to integrate with respect to $h$.
To integrate with respect to a function having the quadratic $\lambda$-variation,
we introduce the Left Cauchy and Right Cauchy $\lambda$-integrals
which are extensions of the Left Young and Right Young integrals,
respectively.

\paragraph*{The Left and Right Cauchy $\lambda$-integrals.}
Let $\phi$ and $g$ be real-valued functions defined on a
closed interval $[a,b]$.
For a partition $\kappa=\{x_i\colon\,i=0,\dots,n\}$ of $[a,b]$,
let
$$
S_{LC}(\phi,g;\kappa):=
\sum_{i=1}^{n}\phi (x_{i-1})\big [g(x_i)-g(x_{i-1})\big ]
\quad\mbox{and}\quad
S_{RC}(\phi,g;\kappa):=
\sum_{i=1}^{n}\phi (x_{i})\big [g(x_i)-g(x_{i-1})\big ].
$$

\begin{defn}\label{LCint}
{\rm Let $g$ be a regulated function on $[a,b]$, let $\phi$ be a function
on $[a,b]$ such that a left limit $\phi (t-)$ exists at $t\in (a,b]$
whenever $\Delta^{-}g(t)\not =0$, and 
let $\lambda=\{\lambda_m\colon\,m\geq 1\}\in\Lambda [a,b]$.
We say that the \emph{Left Cauchy $\lambda$-integral $(LC)\smallint
\phi\,d_{\lambda}g$ is defined on $[a,b]$} if there exists a regulated 
function $\Phi$ on $[a,b]$ such that $\Phi (a)=0$ and for each 
$a\leq u<v\leq b$,
\beq\label{1LCint}
\Phi (v)-\Phi (u)=\lim_{m\to\infty}S_{LC}(\phi,g;\lambda_m\Cap [u,v]),
\eeq
\beq\label{2LCint}
\Delta^{-}\Phi (v)=\left\{\begin{array}{ll}
\big (\phi_{-}\Delta^{-}g\big )(v) &\mbox{if $(\Delta^{-}g)(v)\not =0$,}\\
0 &\mbox{if $(\Delta^{-}g)(v)=0$}\end{array}\right.
\quad\mbox{and}\quad
\Delta^{+}\Phi (u)=\big (\phi\Delta^{+}g\big )(u).
\eeq
If a regulated function $\Phi$ on $[a,b)$ exists and satisfies the same conditions
except that (\ref{1LCint}) and (\ref{2LCint}) do not hold for $v=b$,
then we say that the \emph{Left Cauchy $\lambda$-integral $(LC)\smallint
\phi\,d_{\lambda}g$ is defined on $[a,b)$}. }
\end{defn}

In Theorem \ref{modelling}.\ref{BandS} below the Left
Cauchy $\lambda$-integral $(LC)\smallint \phi\,d_{\lambda}g$ 
on $[a,b]$ is used when $g$ is a suitable continuous function on $[a,b]$
and $\phi$ is a function on $[a,b]$ having no the left limit at the
right-endpoint $b$.
In other places we use this integral when both the integrator and integrand
are regulated in which case a simple characteriztion of the integral holds.
But first  we define the Right Cauchy $\lambda$-integral.

 \begin{defn}\label{RCint}
{\rm Let $g$ be a regulated function on $[a,b]$, let $\phi$ be a function
on $[a,b]$ such that a right limit $\phi (t+)$ exists at $t\in [a,b)$
whenever $\Delta^{+}g(t)\not =0$, and 
let $\lambda=\{\lambda_m\colon\,m\geq 1\}\in\Lambda [a,b]$.
We say that the \emph{Right Cauchy $\lambda$-integral $(RC)\smallint
\phi\,d_{\lambda}g$ on $[a,b]$} exists if there exists a regulated function 
$\Psi$ on $[a,b]$ such that $\Psi (a)=0$ and for each $a\leq u<v\leq b$,
$$
\Psi (v)-\Psi (u)=\lim_{m\to\infty}S_{RC}(\phi,g;\lambda_m\Cap [u,v]),
$$
$$
\Delta^{-}\Psi (v)=\big (\phi\Delta^{-}g\big )(v)
\quad\mbox{and}\quad
\Delta^{+}\Psi (u)=\left\{\begin{array}{ll}
\big (\phi_{+}\Delta^{+}g\big )(u) &\mbox{if $(\Delta^{+}g)(u)\not =0$,}\\
0 &\mbox{if $(\Delta^{+}g)(u)=0$.}\end{array}\right.
$$}
\end{defn}

The Left Cauchy and Right Cauchy $\lambda$-integrals, if exist, are
additive functions on $S[a,b]$.
For notation simplicity, the following is proved under stronger assumptions
than necessary. 

\begin{prop}\label{property1}
Let $g$, $\phi$ be regulated functions on $[a,b]$, and
let $\lambda=\{\lambda_m\colon\,m\geq 1\}\in\Lambda [a,b]$.
The Left Cauchy $\lambda$-integral $(LC)\smallint \phi\,d_{\lambda}g$ 
is defined on $[a,b]$ if and only if there exists an additive upper continuous function 
$\mu_L$ on $S\loc a,b\roc$  such that for each $(u,v)\in S\loc a,b\roc$,
\beq\label{1property1}
\mu_L(u,v)=\lim_{m\to\infty}S_{LC}(\phi,g;\lmuv).
\eeq
Also, the Right Cauchy $\lambda$-integral $(RC)\smallint \phi\,d_{\lambda}g$ 
is defined on $[a,b]$  if and only if there exists an additive upper continuous function $\mu_R$ on 
$S\loc a,b\roc$ such that for each $(u,v)\in S\loc a,b\roc$,
$$
\mu_R(u,v)=\lim_{m\to\infty}S_{RC}(\phi,g;\lmuv).
$$
\end{prop}

\begin{proof}
We prove only the first part of the proposition because the proof of the second 
part is symmetric.
First suppose that there exists a regulated function $\Phi$ on $[a,b]$,
equal zero at $a$ and satisfying (\ref{1LCint}) and (\ref{2LCint}).
Define a function $\mu_{L}$ on $S\loc a,b\roc$ by 
$\mu_L(u,v):=\Phi (v)-\Phi(u)$ for  $(u,v)\in S\loc a,b\roc$.
Then $\mu_L$ is additive and upper continuous function on $S\lei a,b\rei $ 
by Theorem \ref{variation}.\ref{interv}.
We prove (\ref{1property1}) only for the cases $( u,v ) =(a,t-)$ and $( u,v) =(a,t+)$ 
because the proofs for the other cases are similar.
To this aim let 
$\lambda_m=\{x_i^m\colon\,i=0,\dots,n(m)\}$, $m=1,2,\dots$,
and first let $t\in (a,b]$.
For each $m\geq 1$, there is an index $i(t)=i_m(t)\in\{1,\dots,n(m)-1\}$
such that $x_{i(t)}^m<t\leq x_{i(t)+1}^m$.
Then by (\ref{1LCint}) and (\ref{2LCint}), we have
\begin{eqnarray*}
\lefteqn{\lim_{m\to\infty}S_{LC}(\phi,g;\lambda_m\Cap\lei a,t-\rei )}\\[2mm]
&=&\lim_{m\to\infty}\Big\{S_{LC}(\phi,g;\lambda_m\Cap [a,t])
-\phi (x_{i(t)}^m)[g(t)-g(x_{i(t)}^m)]+\phi (x_{i(t)}^m)[g(t-)-g(x_{i(t)}^m)]\Big\}
\\[2mm]
&=&\Phi (t)-\phi (t-)\Delta^{-}g(t)=\Phi (t-)=\mu_{L}(a,t-),
\end{eqnarray*}
proving (\ref{1property1}) for $(u,v)= (a,t-)\in S\lei a,b\rei $.
Now let $t\in [a,b)$.
By the definition of the trace partition (Notation \ref{lmuv}), 
we have for each $m\geq 1$, 
$$
S_{LC}(\phi,g;\lambda_m\Cap\lei a,t+\rei )
=S_{LC}(\phi,g;\lambda_m\Cap [ a,t])+(\phi\Delta^{+}g)(t).
$$
Letting $m\to\infty$, it follows that (\ref{1property1}) holds for
$(u,v)=(a,t+)\in S\loc a,b\roc$. 

To prove the converse implication, let $\Phi$ be the restriction to $[a,b]$
of the right distribution function of $\mu_{L} $, that is,
$\Phi (t):=\mu_{L}(a,t)$ for $t\in [a,b]$.
By statement $(iv)$ of Theorem \ref{variation}.\ref{interv} and 
(\ref{1property1}), $\Phi$ is regulated on $[a,b]$, is zero at $a$, 
and satisfies properties (\ref{1LCint}), (\ref{2LCint}).
The proof of Proposition \ref{property1} is complete. 
\qed\end{proof}

The next property identifies a difference between values of the Left Cauchy
and Right Cauchy $\lambda$-integrals, provided the integrand and integrator
possess the quadratic $\lambda$-covariation.

\begin{prop}\label{property4}
Let a $2$-vector function $(f_1,f_2)$ on $[a,b]$ have the quadratic 
$\lambda$-covariation for some $\lambda\in\Lambda [a,b]$.
The $\lambda$-integral $(LC)\smallint f_1\,d_{\lambda}f_2$ on $[a,b]$ 
exists if and only if so does the $\lambda$-integral $(RC)\smallint f_1\,d_{\lambda}f_2$ 
on $[a,b]$, and if the two $\lambda$-integrals exist then for 
$a\leq x\leq b$,
$$
[f_1,f_2]_{\lambda}(x)=(RC)\int_a^xf_1\,d_{\lambda}f_2
-(LC)\int_a^xf_1\,d_{\lambda}f_2.
$$
\end{prop}

\begin{proof}
The conclusion follows from the relation
$C(f_1,f_2;\kappa)=S_{RC}(f_1,f_2;\kappa)-S_{LC}(f_1,f_2;\kappa)$,
valid for any partition $\kappa$.
\qed\end{proof}

The following extends an integration by parts formula for the
Left Cauchy and Right Cauchy $\lambda$-integrals, provided
the integrand and integrator possess the quadratic $\lambda$-covariation.

\begin{prop}\label{property5}
For $\lambda\in\Lambda [a,b]$, let a $2$-vector function $(f_1,f_2)$ on 
$[a,b]$ have the quadratic $\lambda$-covariation.
For $\sharp$ equal to $LC$ or $RC$, the $\lambda$-integral 
$(\sharp)\smallint f_1\,d_{\lambda}f_2$ on $[a,b]$ exists if and only if so does
the integral $(\sharp)\smallint f_2\,d_{\lambda}f_1$ on $[a,b]$,
and if the four $\lambda$-integrals exist then for $a\leq x\leq b$
\begin{eqnarray*}
[f_1,f_2]_{\lambda}(x)&=&(f_1f_2)(x)-(f_1f_2)(a)-(LC)\int_a^xf_1\,
d_{\lambda}f_2-(LC)\int_a^xf_2\,d_{\lambda}f_1\\[2mm]
&=&-\Big [(f_1f_2)(x)-(f_1f_2)(a)-(RC)\int_a^xf_1\,
d_{\lambda}f_2-(RC)\int_a^xf_2\,d_{\lambda}f_1\Big ].
\end{eqnarray*}
\end{prop}

\begin{proof}
The conclusion follows from the relations
\begin{eqnarray*}
C(f_1,f_2;\kappa)&=&(f_1f_2)(x_n)-(f_1f_2)(x_0)-S_{LC}(f_1,f_2;\kappa)
-S_{LC}(f_1,f_2;\kappa)\\[2mm]
&=&-\Big [(f_1f_2)(x_n)-(f_1f_2)(x_0)-S_{RC}(f_1,f_2;\kappa)
-S_{RC}(f_1,f_2;\kappa)\Big ],
\end{eqnarray*}
valid for any partition $\kappa=\{x_i\colon\,i=0,\dots,n\}$.
\qed\end{proof}

\paragraph*{The Left and Right $\lambda$-integrals.}
For $\lambda\in\Lambda [a,b]$ and $1\leq p<2$, let $f$ be a
$(\lambda,p)$-decomposable function on $[a,b]$, that is, $f\equiv C+g+h$
for a $(\lambda,p)$-dual pair $(g,h)$.
Next two integrals with respect to $f$ are defined to be a sum of two
integrals with respect to $g$ and $h$, respectively.
Since the component $h$ is less erratic than the component $g$, 
one may expect that an integral with respect to $h$ may exist 
in a stronger sense for a large enough class of integrands.
Therefore we use the Left (Right) Young integral to integrate with respect to $h$.
Under the conditions of Theorem \ref{variation}.\ref{approximation}, 
the Left (Right) Young integral exists and equals to the
Left (Right) Cauchy $\lambda$-integral for any $\lambda\in\Lambda [a,b]$
such that the rigt (left) discontinuity points of $h$ are accessible by $\lambda$.
This type of independence on a sequence of partitions $\lambda$ 
makes the Left (Right) Young integral preferable whenever it exists.
Moreover, any result concerning the Left (Right) Cauchy $\lambda$-integral with 
respect to such a function is an \emph{extension} of the corresponding result 
concerning the Left (Right) Young integral.
Since a $(\lambda,p)$-decomposition of an integrator may be non-unique we 
have to insure that the sum of the two integrals with respect to different
decompositions give the same value.

\begin{defn}\label{LandR}
{\rm Let $F$ be a regulated function on $[a,b]$,
let $f\in {\cal D}_{\lambda,p}[a,b]$ for some $\lambda\in\Lambda [a,b]$
and $1\leq p<2$, and let $A\subset D_{\lambda,p}(f)$.

I. We say that the \emph{Left $\lambda$-integral $(L)\smallint F\,d_{\lambda}f$
on $[a,b]$ is defined with respect to $A$} if for each $(g,h)\in A$ the sum
$(LC)\smallint F\,d_{\lambda}g+(LY)\smallint F\,dh$ is defined on
$[a,b]$ and its values for each $(u,v)\in S[a,b]$ do not depend on $(g,h)\in A$. 
If the Left $\lambda$-integral $(L)\smallint F\,d_{\lambda}f$ is defined 
on $[a,b]$ with respect to $A$ then let
$$
(L)\int_u^vF\,d_{\lambda}f:=
(LC)\int_u^vF\,d_{\lambda}g+(LY)\int_u^vF\,dh
$$
for $(u,v)\in S[a,b]$ and for any $(g,h)\in A$.}

{\rm II. We say that the \emph{Right $\lambda$-integral 
$(R)\smallint F\,d_{\lambda}f$ on $[a,b]$  
 is defined with respect to $A$} if for each $(g,h)\in A$ the sum
$(RC)\smallint F\,d_{\lambda}g+(RY)\smallint F\,dh$ is defined on
$[a,b]$ and its values for each $(u,v)\in S[a,b]$ do not depend on $(g,h)\in A$. 
If the Right $\lambda$-integral $(R)\smallint F\,d_{\lambda}f$ is defined 
on $[a,b]$ with respect to $A$ then let
$$
(R)\int_u^vF\,d_{\lambda}f:=
(RC)\int_u^vF\,d_{\lambda}g+(RY)\int_u^vF\,dh
$$
for $(u,v)\in S[a,b]$ and for any $(g,h)\in A$.}
\end{defn}

Under the assumptions of the following proposition, 
the Left and Right $\lambda$-integrals with respect to a set $A$
are defined whenever they are defined for a single element of $A$.

\begin{prop}\label{property3}
For $\lambda\in\Lambda [a,b]$ and $1\leq p<2$, let $f$ be a 
$(\lambda,p)$-decomposable function on $[a,b]$, 
and let $F\in \dual({\cal W}_p)[a,b]$.
\begin{enumerate}
\item[$(a)$] If $(LC)\smallint F\,d_{\lambda}g$ is defined on $[a,b]$
for some $(g,h)\in D_{\lambda,p}^{+}(f)$ then $(L)\smallint F\,d_{\lambda}f$
is defined on $[a,b]$ with respect to $D_{\lambda,p}^{+}(f)$.
\item[$(b)$] If $(RC)\smallint F\,d_{\lambda}g$ is defined on $[a,b]$
for some $(g,h)\in D_{\lambda,p}^{-}(f)$ then $(R)\smallint F\,d_{\lambda}f$
is defined on $[a,b]$ with respect to $D_{\lambda,p}^{-}(f)$.
\end{enumerate}
\end{prop}

\begin{proof}
We prove only statement $(a)$ because a proof of statement $(b)$
is symmetric.
If the set $D_{\lambda,p}^{+}(f)$ has a single element then 
there is nothing to prove.
Suppose that there are two different elements $(g_1,h_1)$ and $(g_2,h_2)$
in $D_{\lambda,p}^{+}(f)$.
By the $(\lambda,p)$-decomposability of $f$, 
$g_1-g_2=-(h_1-h_2)\in {\cal W}_p$.
In the case $1<p<2$, $F\in {\cal W}_q$ with $p^{-1}+q^{-1}>1$, and the Left 
Young integral of $F$ with respect to $g_1-g_2$ exists by 
Theorem \ref{variation}.\ref{LYineq}.
While in the case $p=1$, this integral exists by Lemma \ref{bv}.
Since the $\lambda$-integrals $(LC)\smallint F\,d_{\lambda}g_1$
and $(LC)\smallint F\,d_{\lambda}g_2$ are defined, it follows that
$(LC)\smallint F\,d_{\lambda}(g_1-g_2)$ exists and equals to the difference
of the two $\lambda$-integrals.
If $\Delta^{+}(g_1-g_2)(x)\not =0$ for some $x\in [a,b)$ then either
$\Delta^{+}g_1(x)\not =0$ or $\Delta^{+}g_2(x)\not =0$, and so
$x\in\cup_m\lambda_m$.
By Theorem \ref{variation}.\ref{approximation}, we then have
$$
(LC)\int F\,d_{\lambda}g_1-(LC)\int F\,d_{\lambda}g_2
=(LC)\int F\,d_{\lambda}(g_1-g_2)
$$
$$
=(LY)\int F\,d(g_1-g_2)=(LY)\int F\,dh_2-(LY)\int F\,dh_1.
$$
Thus the sum $(LC)\smallint F\,d_{\lambda}g+(LY)\smallint F\,dh$
does not depend on whether $(g,h)=(g_1,h_1)$ or $(g,h)=(g_2,h_2)$.
The proof of Proposition \ref{property3} is complete.
\qed\end{proof}

Thus existence of the Left and Right $\lambda$-integrals essentially
reduces to the problem of existence of the Left and Right Cauchy
$\lambda$-integrals.

\paragraph*{Chain rules.}
Next we turn to establishing sufficient conditions for existence
of the $\lambda$-integrals $(L)\smallint (\psi{\circ}f)\,d_{\lambda}f$
and $(R)\smallint (\psi{\circ}f)\,d_{\lambda}f$ when $f$ has the quadratic
$\lambda$-variation and $\psi$ is a smooth function.
Also, we show that each such $\lambda$-integral satisfies a relation to be
called a \emph{chain rule}.

\begin{thm}\label{chainrule}
Let $f\in {\cal D}_{\lambda,p}[a,b]$ for some $\lambda\in \Lambda[a,b]$
and $1\leq p<2$, and let $\phi$ be a $C^2$ class function 
on an open set containing the range of $f$.
Then the following two statements hold{\rm :}
\begin{enumerate}
\item[$(a)$] If $D_{\lambda,p}^{+}(f)\not =\emptyset$ then the Left 
$\lambda$-integral $(L)\smallint (\phi'{\circ}f)\,d_{\lambda}f$ is defined 
on $[a,b]$ with respect to $D_{\lambda,p}^{+}(f)$, and for $a\leq z<y\leq b$,
\begin{eqnarray}\label{7HF}
(\phi{\circ}f)(y)&=&(\phi{\circ}f)(z)
+(L)\int_z^y(\phi'{\circ}f)\,d_{\lambda}f
+\frac{1}{2}(RS)\int_z^y(\phi''{\circ}f)\,d[f]_{\lambda}^c\\[2mm]
& &+\sum_{(z,y]}\Big\{\Delta^{-}(\phi{\circ}f)-(\phi'{\circ}f)_{-}
\Delta^{-}f\Big\}+\sum_{[z,y)}\Big\{\Delta^{+}(\phi{\circ}f)-
(\phi'{\circ}f)\Delta^{+}f\Big\},\nonumber
\end{eqnarray}
where the two sums converge unconditionally.
\item[$(b)$] If $D_{\lambda,p}^{-}(f)\not =\emptyset$, then the Right 
$\lambda$-integral  $(R)\smallint (\phi'{\circ}f)\,d_{\lambda}f$  is defined 
on $[a,b]$ with respect to $D_{\lambda,p}^{-}(f)$, and for $a\leq z<y\leq b$,
\begin{eqnarray}\label{17HF}
(\phi{\circ}f)(y)&=&(\phi{\circ}f)(z)
+(R)\int_z^y(\phi'{\circ}f)\,d_{\lambda}f
-\frac{1}{2}(RS)\int_z^y(\phi''{\circ}f)\,d[f]_{\lambda}^c\\[2mm]
& &+\sum_{(z,y]}\Big\{\Delta^{-}(\phi{\circ}f)-(\phi'{\circ}f)
\Delta^{-}f\Big\}+\sum_{[z,y)}\Big\{\Delta^{+}(\phi{\circ}f)-
(\phi'{\circ}f)_{+}\Delta^{+}f\Big\},\nonumber
\end{eqnarray}
where the two sums converge unconditionally.
\end{enumerate}
\end{thm}

\begin{proof}
We prove statement $(a)$ and indicate necessary changes needed
to prove statement $(b)$.
Let $f\equiv C+g+h$  for some $(g,h)\in D_{\lambda,p}^{+}(f)$.
Then $g$ has the quadratic $\lambda$-variation and
has bounded $q$-variation for some $2<q<p/(p-1)$ if $1<p<2$.
In the case $1<p<2$, the function $f$ has the same $q$-variation
property as the function $g$ because $p<2<q$.
Since $\phi'$ is Lipschitz function on the range of $f$, 
the composition $\phi'{\circ}f\in {\cal W}_q[a,b]$ if $1<p<2$,
and it is a regulated function on $[a,b]$ if $p=1$.
Thus the Left Young integral $(LY)\smallint_a^b(\phi'{\circ}f)\,dh$ 
exists by Theorem \ref{variation}.\ref{approximation} in the case $1<p<2$ 
and by Lemma \ref{bv} in the case $p=1$.
By Proposition \ref{property3} with $A=D_{\lambda,p}^{+}(f)$,
it is enough to prove that the Left Cauchy $\lambda$-integral
$(LC)\smallint (\phi'{\circ}f)\,d_{\lambda}g$ on $[a,b]$
is defined and relation (\ref{7HF}) holds for each $a\leq z<y\leq b$ .
In fact, it is enough to prove that  for each $a\leq z<y\leq b$, the limit
\beq\label{6HF}
(LC)\int_z^y(\phi'{\circ}f)\,d_{\lambda}g:=\lim_{m\to\infty}
S_{LC}(\phi'{\circ}f,g;\lambda_m\Cap [z,y])
\eeq
exists and relation (\ref{7HF}) holds.
Indeed, let $\Phi (y):=(LC)\smallint_a^y(\phi'{\circ}f)\,d_{\lambda}g$
for $a\leq y\leq b$.
The Left Young integral and the Riemann-Stieltjes integral are additive 
over adjacent intervals (Theorem \ref{variation}.\ref{additivity}). 
While the unconditional sums are additive over disjoint intervals
(Theorem \ref{convergence}.\ref{sum-add}).
Therefore using (\ref{7HF}), it follows that 
$\Phi (y)-\Phi (z)=(LC)\smallint_z^y(\phi'{\circ}f)\,d_{\lambda}g$
for $a\leq z<y\leq b$, that is, condition (\ref{1LCint}) of  
Definition \ref{LCint} holds.
To show that condition (\ref{2LCint}) holds, for $a\leq y<x\leq b$, we have
\begin{eqnarray}
\Phi (x)-\Phi (y)=&&(LC)\int_y^x\big (\phi'{\circ}f\big )\,
d_{\lambda}g=\big [(\phi{\circ}f)(x)-(\phi{\circ}f)(y)\big ]
\label{2qvofint}\\
&&-(LY)\int_x^y(\phi'{\circ}f)\,dh
-\frac{1}{2}(RS)\int_x^y(\phi''{\circ}f)\,d[f]_{\lambda}^c
\nonumber\\[2mm]
&&-\sum_{(x,y]}\Big\{\Delta^{-}(\phi{\circ}f)-(\phi'{\circ}f)_{-}
\Delta^{-}f\Big\}-\sum_{[x,y)}\Big\{\Delta^{+}(\phi{\circ}f)-
(\phi'{\circ}f)\Delta^{+}f\Big\}\nonumber\\[2mm]
\stackrel{y\uparrow x}{\longrightarrow}&&
\Delta^{-}(\phi{\circ}f)(x)-\big\{(\phi'{\circ}f)_{-}\Delta^{-}h\big\}
(x)-\big\{\Delta^{-}(\phi{\circ}f)-(\phi'{\circ}f)_{-}\Delta^{-}f
\big \}(x)\nonumber\\
=&&(\phi'{\circ}f)_{-}(x)\Delta^{-}g(x).\nonumber
\end{eqnarray}
For the last step we used Proposition \ref{variation}.\ref{LYjumps}.
Likewise for $a\leq x<y\leq b$, we have
\begin{eqnarray*}
\Phi (y)-\Phi (x)&\stackrel{y\downarrow x}{\longrightarrow}&
\Delta^{+}(\phi{\circ}f)(x)-\big\{(\phi'{\circ}f)\Delta^{+}h\big\}
(x)-\big\{\Delta^{+}(\phi{\circ}f)-(\phi'{\circ}f)\Delta^{+}f
\big \}(x)\\
&=&(\phi'{\circ}f)(x)\Delta^{+}g(x).
\end{eqnarray*}
Therefore $\Phi$ is a regulated function on $[a,b]$, is equal zero at $a$,
and has jumps satisfying the relations
\beq\label{3qvofint}
\Delta^{-}\Phi =(\phi'{\circ}f)_{-}\Delta^{-}g
\qquad\mbox{and}\qquad
\Delta^{+}\Phi =(\phi'{\circ}f)\Delta^{+}g.
\eeq
That is, (\ref{2LCint}) of Definition \ref{LCint} holds.
Hence statement $(a)$ of the theorem will be proved once
we show that for each $a\leq z<y\leq b$,  the limit (\ref{6HF})
exists and relation (\ref{7HF}) holds.

For notation simplicity, the existence of the limit (\ref{6HF})
and relation (\ref{7HF}) will be proved  when $z=a$ and $y=b$.
Let $A_1:=(LY)\smallint_a^b(\phi'{\circ}f)\,dh$.
By Lemma \ref{bv}, there exists the Riemann--Stieltjes integral 
$A_2:=(RS)\smallint_a^b (\phi''{\circ}f)\,d[f]_{\lambda}^c$.
By the mean value theorem, for any $y\in (a,b]$,
$$
\big |\Delta^{-}(\phi{\circ}f)(y)-(\phi'{\circ}f)(y-)\Delta^{-}f(y)\big |
\leq K\big [\Delta^{-}f(y)\big ]^2,
$$
where $K$ is the Lipschitz constant of $\phi'$ on the range of $f$.
Thus the first sum in (\ref{7HF}) converges unconditionally to a value 
$A_3$, because $\sum_{(a,b]}[\Delta^{-}f]^2<\infty$.
Likewise, it follows  that the second sum in (\ref{7HF}) converges 
unconditionally to a value $A_4$.
For each $m\geq 1$, let $\lambda_m=\{x_i^m\colon\,i=0,\dots,n(m)\}$, 
and let $\Delta_i^m\rho:=\rho (x_i^m)-\rho (x_{i-1}^m)$, $i=1,\dots,n(m)$, 
for a function $\rho$ on $[a,b]$.
Let $\epsilon >0$.
Since $(g,h)\in D_{\lambda,p}^{+}(f)$, and so
$\{x\in (a,b)\colon\,\Delta^{+}h(x)\not =0\}\subset\cup_m\lambda_m$,
by Theorem \ref{variation}.\ref{approximation}, there exists an integer $m_1$ 
such that for all $m\geq m_1$,
\beq\label{9HF}
\Big |\sum_{i=1}^{n(m)}\phi'(f(x_{i-1}^m))\Delta_i^mh
-A_1\Big |<\epsilon.
\eeq
Since $\{x\in (a,b)\colon\,\Delta^{+}g(x)\not =0\}\subset\cup_m\lambda_m$,
by  Theorem \ref{intqv} and Lemma \ref{bv},
there exists an integer $m_2$ such that for all $m\geq m_2$,
\beq\label{8HF}
\Big |\sum_{i=1}^{n(m)}\phi''(f(x_{i-1}^m))\big\{\Delta_i^mg\big\}^2
-A_2-A_5-A_6\Big |<\epsilon,
\eeq
where 
$$
A_5:=\sum_{(a,b]}\big (\phi''{\circ}f\big )_{-}\big\{\Delta^{-}g\big\}^2
\quad\mbox{and}\quad
A_6:=\sum_{[a,b)}\big (\phi''{\circ}f\big )\big\{\Delta^{+}g\big\}^2.
$$
Since $\phi''$ is uniformly continuous on the range of $f$,
there is a $\delta >0$ such that $|\phi''(u)-\phi''(v)|<\epsilon$
whenever $|u-v|<\delta$ and $|u|\vee |v|\leq\|f\|_{\infty}$.
Therefore and recalling  Lemma \ref{variation}.\ref{lemma1}, we conclude that
there exists a partition $\mu=\{y_j\colon\,j=0,\dots,k\}$ of $[a,b]$ such that
for any $\nu\supset\mu$,
\beq\label{10HF}
\Big |\sum_{\nu\setminus\{a\}}\Big\{\Delta^{-}(\phi{\circ}f)
-(\phi'{\circ}f)_{-}\Delta^{-}f\Big\}-A_3\Big |\vee
\Big |\sum_{\nu\setminus\{a\}}(\phi''{\circ}f)_{-}\big\{\Delta^{-}g\big\}^2
-A_5\Big |<\epsilon,
\eeq
\beq\label{12HF}
\Big |\sum_{\nu\setminus\{b\}}\Big\{\Delta^{+}(\phi{\circ}f)
-(\phi'{\circ}f)\Delta^{+}f\Big\}-A_4\Big |\vee
\Big |\sum_{\nu\setminus\{b\}}(\phi''{\circ}f)\big\{\Delta^{+}g\big\}^2
-A_6\Big |<\epsilon,
\eeq
\beq\label{11HF}
\max_{1\leq j\leq k}Osc (f;(y_{j-1},y_j))<\delta
\quad\mbox{and}\quad
\sum_{j=1}^k v_2(h;(y_{j-1},y_j))<\epsilon.
\eeq
Let $m_0$ be the minimal integer such that $\mu\subset\lambda_{m_0}$
and each pair $\{y_{j-1},y_j\}$ is separated by at least two other
points from $\lambda_{m_0}$.
For each $j=1,\dots,k-1$ and $m\geq m_0$, there exists an integer 
$i(j)=i_m(j)\in\{1,\dots,n(m)-1\}$ such that $x_{i(j)}^m=y_j$.
For each $m\geq m_0$, let $I_1(m):=\{x_1^m,x_{n(m)}^m,x_{i(j)}^m,
x_{i(j)+1}^m\colon\,j=1,\dots,k-1\}$ and 
$I_2(m):=\{1,\dots,n(m)\}\setminus I_1(m)$.
By (\ref{10HF}) and (\ref{12HF}), there exists an integer $m_3\geq m_0$ 
such that for all $m\geq m_3$
\beq\label{13HF}
\Big |\sum_{i\in I_1(m)}\Big\{\Delta_i^m(\phi{\circ}f)
-(\phi'{\circ}f)(x_{i-1}^m)\Delta_i^mf\Big\}-A_3-A_4\Big |<4\epsilon
\eeq
and
\beq\label{14HF}
\Big |\sum_{i\in I_1(m)}\phi''(f(x_{i-1}^m))\big\{\Delta_i^mg\big\}^2
-A_5-A_6\Big |<4\epsilon.
\eeq
By Taylor's theorem with Lagrange's form of the remainder:
for real numbers $u, v$, there is a $\theta=\theta (u,v)\in (0,1)$
such that 
\beq\label{3ps}
\phi(v)=\phi (u)+\phi'(u)(v-u)+\frac{1}{2}\phi''(u+\theta(v-u))[v-u]^2.
\eeq
Then by a telescoping sum, for $m\geq m_0$, we have
$$
(\phi{\circ}f)(b)-(\phi{\circ}f)(a)
=\sum_{i\in I_2(m)}\Big\{\phi'(f(x_{i-1}^m))\Delta_i^mf+\frac{1}{2}
\phi''(y_i^m)\{\Delta_i^mf\}^2\Big\}
+\sum_{i\in I_1(m)}\Delta_i^m(\phi{\circ}f)
$$
\beq\label{16HF}
=\sum_{i=1}^{n(m)}\phi'(f(x_{i-1}^m))\big\{\Delta_i^mg+\Delta_i^mh\big\}
+\frac{1}{2}\sum_{i=1}^{n(m)}\phi''(f(x_{i-1}^m))\big\{\Delta_i^mg\big\}^2
\eeq
$$
-\frac{1}{2}\sum_{i\in I_1(m)}\phi''(f(x_{i-1}^m))\big\{\Delta_i^mg\big\}^2
+\sum_{i\in I_1(m)}\Big\{\Delta_i^m(\phi{\circ}f)
-(\phi'{\circ}f)(x_{i-1}^m)\Delta_i^mf\Big\}+R_m,
$$
where   $y_i^m:=f(x_{i-1}^m)+\theta_i^m\Delta_i^mf$ with $\theta_i^m\in
(0,1)$ for $i=1,\dots,n(m)$, and
$$
R_m:=\frac{1}{2}\sum_{i\in I_2(m)}\Big\{\phi''(y_i^m)
\big [2\Delta_i^mg\Delta_i^mh +\{\Delta_i^mh\}^2\big ]
+\big [\phi''(y_i^m)-\phi''(f(x_{i-1}^m))]\{\Delta_i^mg\}^2\Big\}.
$$
Let 
$
A:=(\phi{\circ}f)(b)-(\phi{\circ}f)(a)-A_1-\frac{1}{2}A_2-A_3-A_4.
$
Then by (\ref{9HF}), (\ref{8HF}), (\ref{13HF}) and (\ref{14HF}), for all 
$m\geq m_1\vee m_2\vee m_3$,
$$
\big |S_{LC}(\phi'{\circ}f,g;\lambda_m)-A\big |
\leq\Big |\sum_{i=1}^{n(m)}\phi'(f(x_{i-1}^m))\Delta_i^mh-A_1\Big |
$$
$$
+\frac{1}{2}\Big |\sum_{i=1}^{n(m)}\phi''(f(x_{i-1}^m))
\big\{\Delta_i^mg\big\}^2-(A_2+A_5+A_6)\Big |
+\frac{1}{2}\Big |\sum_{i\in I_1(m)}\phi''(f(x_{i-1}^m))
\big\{\Delta_i^mg\big\}^2-(A_5+A_6)\Big |
$$
$$
+\Big |\sum_{i\in I_1(m)}\Big\{\Delta_i^m(\phi{\circ}f)
-(\phi'{\circ}f)(x_{i-1}^m)\Delta_i^mf\Big\}-(A_3+A_4)\Big |
+|R_m|
\leq \frac{15}{2}\epsilon +|R_m|.
$$
To bound $|R_m|$ notice that for each $i\in I_2(m)$ with $m\geq m_0$, 
$[x_{i-1}^m,x_i^m]\subset (y_{j-1},y_j)$ for some $j=1,\dots,k$.
By (\ref{11HF}) and H\"older's inequality, it then follows that
\begin{eqnarray*}
\sup_{m\geq m_0}|R_m|&\leq &\frac{\|\phi''{\circ}f\|_{\infty}}{2}
\Big\{2\Big ([g]_{\lambda}(b)\sum_{j=1}^kv_2(h;(y_{j-1},y_j))\Big )^{1/2}
+\sum_{j=1}^kv_2(h;(y_{j-1},y_j))\Big\}+\epsilon [g]_{\lambda}(b)\\
&\leq &\sqrt{\epsilon}\big\{(\|\phi''{\circ}f\|_{\infty}/2)
\big (2\sqrt{[g]_{\lambda}(b)}+1\big )+[g]_{\lambda}(b)\big\}.
\end{eqnarray*}
Since $\epsilon\in (0,1)$ is arbitrary, the limit (\ref{6HF}) exists
and relation (\ref{7HF}) holds for $z=a$ and $y=b$.
The proof of statement $(a)$ is complete.

The proof of statement $(b)$ is similar.
Instead of (\ref{16HF}), in this case we use the relations
$$
(\phi{\circ}f)(b)-(\phi{\circ}f)(a)
=\sum_{i\in I_2(m)}\Big\{\phi'(f(x_{i}^m))\Delta_i^mf-\frac{1}{2}
\phi''(y_i^m)\{\Delta_i^mf\}^2\Big\}
+\sum_{i\in I_1(m)}\Delta_i^m(\phi{\circ}f)
$$
$$
=\sum_{i=1}^{n(m)}\phi'(f(x_{i}^m))\big\{\Delta_i^mg+\Delta_i^mh\big\}
-\frac{1}{2}\sum_{i=1}^{n(m)}\phi''(f(x_{i}^m))\big\{\Delta_i^mg\big\}^2
$$
$$
+\frac{1}{2}\sum_{i\in I_1(m)}\phi''(f(x_{i}^m))\big\{\Delta_i^mg\big\}^2
+\sum_{i\in I_1(m)}\Big\{\Delta_i^m(\phi{\circ}f)
-(\phi'{\circ}f)(x_{i}^m)\Delta_i^mf\Big\}+R_m,
$$
valid for $m\geq m_0$,
where   $y_i^m:=f(x_{i}^m)-\theta_i^m\Delta_i^mf$ with $\theta_i^m\in
(0,1)$ for $i=1,\dots,n(m)$, and
$$
R_m:=-\frac{1}{2}\sum_{i\in I_2(m)}\Big\{\phi''(y_i^m)
\big [2\Delta_i^mg\Delta_i^mh +\{\Delta_i^mh\}^2\big ]
+\big [\phi''(y_i^m)-\phi''(f(x_{i}^m))]\{\Delta_i^mg\}^2\Big\}.
$$
Here we used the rearranged Taylor's theorem (\ref{3ps}).
The other steps of the proof are the same as in the proof of
statement $(a)$ except that now we approximate the integrals by Right 
Cauchy sums.
The proof of Theorem \ref{chainrule} is complete.
\qed\end{proof}

The preceding chain rules for $(\lambda,p)$-decomposable functions
are used to characterize functions having the quadratic 
$\lambda$-variation in terms of the existence of the $\lambda$-integrals 
of a special form.
A statement of this result for the Left Cauchy $\lambda$-integral
is given in the introductory section.
The case of the Right Cauchy $\lambda$-integral is symmetric and is treated next.

\begin{prop}\label{1chqv}
Let $f$ be a regulated function on $[a,b]$, and let 
$\lambda\in\Lambda [a,b]$ be such that 
the left discontinuity points of $f$ are accessible by $\lambda$.
The following three statements about $f$ are equivalent{\rm :}
\begin{enumerate}
\item[$(a)$] $f$ has the quadratic $\lambda$-variation on $[a,b]${\rm ;}
\item[$(b)$] for a $C^1$ class function $\psi$, the composition $\psi{\circ}f$
is Right Cauchy $\lambda$-integrable on $[a,b]$ with respect to $f$, and 
for any $a\leq u<v\leq b$,
\begin{eqnarray*}
(\Psi{\circ}f)(v)&=&(\Psi{\circ}f)(u)
+(RC)\int_u^v(\psi{\circ}f)\,d_{\lambda}f
-\frac{1}{2}(RS)\int_u^v(\psi'{\circ}f)\,d[f]_{\lambda}^c\\[2mm]
& &+\sum_{(u,v]}\Big\{\Delta^{-}(\Psi{\circ}f)-(\psi{\circ}f)
\Delta^{-}f\Big\}+\sum_{[u,v)}\Big\{\Delta^{+}(\Psi{\circ}f)-
(\psi{\circ}f)_{+}\Delta^{+}f\Big\},
\end{eqnarray*}
where the two sums converge unconditionally
and $\Psi (u)=\Psi (0)+\smallint_0^u\psi (x)\,dx$, $u\in\RR${\rm ;}
\item[$(c)$] $f$ is Right Cauchy $\lambda$-integrable on $[a,b]$ with respect to
itself.
\end{enumerate}
If any of the three statements holds, then for $a\leq x\leq b$,
\beq\label{2chqv}
(RC)\int_a^xf\,d_{\lambda}f=\frac{1}{2}\left\{f^2(x)-f^2(a)
+[f]_{\lambda}(x)\right\}.
\eeq
\end{prop}

\begin{proof}
Since $f$ is $(\lambda,1)$-decomposable and 
$(f-f(a),0)\in D_{\lambda,1}^{+}(f)$,
implication $(a)\Rightarrow (b)$ follows from statement $(b)$
of Theorem \ref{chainrule}.
Taking $\psi (u)=u$ in statement $(b)$ we get statement $(c)$.
To prove $(c)\Rightarrow (a)$, let 
$H(x):=-f^2(x)+f^2(a)+2(RC)\smallint_a^xf\,d_{\lambda}f$
for $x\in [a,b]$, and notice that
for each $a\leq u<v\leq b$ and any partition $\kappa$ of $[a,b]$,
we have
$$
s_2(f;\kappa\Cap [u,v])=-f^2(v)+f^2(u)+2S_{RC}(f,f;\kappa\Cap
[u,v]).
$$
This identity also yields relation (\ref{2chqv}).
The proof of Proposition \ref{1chqv} is complete.
\qed\end{proof}

Let $f$ be a $(\lambda,p)$-decomposable function on $[a,b]$
for some $\lambda\in\Lambda [a,b]$ and $1\leq p<2$ with a 
decomposition $f\equiv C+g+h$.
By Theorem \ref{chainrule}, if $\phi$ is a function of the class
$C^2$, and if $(g,h)\in D_{\lambda,p}^{+}(f)$, then the composition 
$\phi{\circ}f$ has the decomposition
\beq\label{10qvofint}
\phi{\circ}f(x)=\phi{\circ}f(a)+\Phi (x)+U(x),
\qquad a\leq x\leq b,
\eeq
where
\beq\label{8qvofint}
\Phi (x):=(LC)\int_a^x(\phi'{\circ}f)\,d_{\lambda}g
\eeq
and
\begin{eqnarray}\label{6qvofint}
U(x)&:=&(LY)\int_a^x(\phi'{\circ}f)\,dh
+\frac{1}{2}(RS)\int_a^x\big (\phi''{\circ}f\big )\,d[g]_{\lambda}^c\\
& &+\sum_{(a,x]}\Big\{\Delta^{-}(\phi{\circ}f)-(\phi'{\circ}f)_{-}\Delta^{-}
f\Big\}+\sum_{[a,x)}\Big\{\Delta^{+}(\phi{\circ}f)-(\phi'{\circ}f)
\Delta^{+}f\Big\}.\nonumber
\end{eqnarray}
Moreover, for the same function $\phi$, if $(g,h)\in D_{\lambda,p}^{-}(f)$
then the composition $\phi{\circ}f$ has the decomposition
\beq\label{11qvofint}
\phi{\circ}f(x)=\phi{\circ}f(a)+\Psi (x)+V(x),
\qquad a\leq x\leq b,
\eeq
where
\beq\label{5qvofint}
\Psi (x):=(RC)\int_a^x(\phi'{\circ}f)\,d_{\lambda}g
\eeq
and
\begin{eqnarray}\label{7qvofint}
V(x)&:=&(RY)\int_a^x(\phi'{\circ}f)\,dh
-\frac{1}{2}(RS)\int_a^x\big (\phi''{\circ}f\big )\,d[g]_{\lambda}^c\\
& &+\sum_{(a,x]}\Big\{\Delta^{-}(\phi{\circ}f)-(\phi'{\circ}f)\Delta^{-}
f\Big\}+\sum_{[a,x)}\Big\{\Delta^{+}(\phi{\circ}f)-(\phi'{\circ}f)_{+}
\Delta^{+}f\Big\}.\nonumber
\end{eqnarray}
Using Theorem \ref{variation}.\ref{LYineq}, one can show that the functions 
$U$ and $V$ have bounded $p$-variation.
Also, the functions $\Phi$ and $\Psi$ have the same $p$-variation
property as the function $g$.
The latter follows from the corresponding chain rule formulas.
Next we show that $\Phi$ and $\Psi$ have the quadratic $\lambda$-variation.
Thus the composition $\phi{\circ}f$ is a function of
the class ${\cal D}_{\lambda,p}[a,b]$ by Corollary \ref{invar} below.
This result is reminiscent of the important fact of Stochastic Analysis
that the class of semimartingales is closed under
taking a composition with a $C^2$ class function,
which is a consequence of It\^o's formula and the fact that
the stochastic integral with respect to a martingale
is also a martingale.

\begin{thm}\label{qvofint}
For $\lambda\in\Lambda [a,b]$ and $1\leq p<2$,
let $f\equiv C+g+h$ with $(g,h)\in D_{\lambda,p}(f)$, and 
let $\phi$ be a $C^2$ class function.
Then the following two statements hold{\rm :}
\begin{enumerate}
\item[$(a)$] If $(g,h)\in D_{\lambda,p}^{+}(f)$ then the indefinite Left 
Cauchy $\lambda$-integral $\Phi$ defined by {\rm (\ref{8qvofint})} has 
the quadratic $\lambda$-variation with the bracket function
\beq\label{1qvofint}
[\Phi]_{\lambda}(x)=(LY)\int_a^x\big (\phi'{\circ}f\big )^2\,
d[g]_{\lambda},\qquad x\in [a,b],
\eeq
and the $2$-vector function $(\Phi,g)$ has the quadratic 
$\lambda$-covariation with the bracket function
\beq\label{1qcofint}
[\Phi,g]_{\lambda}(x)=(LY)\int_a^x\big (\phi'{\circ}f\big )\,
d[g]_{\lambda},\qquad x\in [a,b].
\eeq
\item[$(b)$] If $(g,h)\in D_{\lambda,p}^{-}(f)$ then the indefinite Right 
Cauchy $\lambda$-integral $\Psi$ defined by {\rm (\ref{5qvofint})}
has the quadratic $\lambda$-variation with the bracket function
$$
[\Psi]_{\lambda}(x)=(RY)\int_a^x\big (\phi'{\circ}f\big )^2\,
d[g]_{\lambda},\qquad x\in [a,b],
$$
and the $2$-vector function $(\Psi,g)$ has the quadratic 
$\lambda$-covariation with the bracket function
$$
[\Psi,g]_{\lambda}(x)=(RY)\int_a^x\big (\phi'{\circ}f\big )\,
d[g]_{\lambda},\qquad x\in [a,b].
$$
\end{enumerate}
\end{thm}

\begin{proof}
We prove only statement $(a)$ because the proof of statement $(b)$ is symmetric.
Let $(g,h)\in D_{\lambda,p}^{+}(f)$ be such that $f\equiv C+g+h$.
By Definition \ref{LCint}, the indefinite Left Cauchy $\lambda$-integral $\Phi$ 
is a regulated function on $[a,b]$.
Statement $(a)$ will be proved once we show that for $a\leq u <v\leq b$,
\beq\label{9qvofint}
\lim_{m\to\infty}s_2(\Phi;\lambda_m\Cap [ u,v] )
=(LY)\int_u^v\big (\phi'{\circ}f\big )^2\,d[g]_{\lambda}
\eeq
and
\beq\label{2qcofint}
\lim_{m\to\infty}C(\Phi,g;\lambda_m\Cap [u,v] )
=(LY)\int_u^v\big (\phi'{\circ}f\big )\,d[g]_{\lambda}.
\eeq
Indeed, let $[\Phi]_{\lambda}(x):=(LY)\smallint_a^x(\phi'{\circ}g)^2
\,d[g]_{\lambda}$ for $a\leq x\leq b$.
Then by statement $(a)$ of Proposition \ref{variation}.\ref{LYjumps},
by (\ref{4alpha}) with $H=[g]_{\lambda}$, and by (\ref{2LCint}), we have
$$
\Delta^{-}[\Phi]_{\lambda}=\big (\phi'{\circ}f\big )_{-}^2
\Delta^{-}[g]_{\lambda}=\big\{\Delta^{-}\Phi\big\}^2
\quad\mbox{and}\quad
\Delta^{+}[\Phi]_{\lambda}=\big (\phi'{\circ}f\big )^2
\Delta^{+}[g]_{\lambda}=\big\{\Delta^{+}\Phi\big\}^2.
$$
Thus the indefinite Left Cauchy integral $\Phi$ has the 
quadratic $\lambda$-variation with the bracket function 
$[\Phi]_{\lambda}$ satisfying (\ref{1qvofint}) provided (\ref{9qvofint})
holds.
Likewise, let $[\Phi,g]_{\lambda}(x):=(LY)\smallint_a^x(\phi'{\circ}g)
\,d[g]_{\lambda}$ for $a\leq x\leq b$.
By the same arguments, we have
$$
\Delta^{-}[\Phi,g]_{\lambda}=\big (\phi'{\circ}f\big )_{-}
\Delta^{-}[g]_{\lambda}=\Delta^{-}\Phi\Delta^{-}g
\quad\mbox{and}\quad
\Delta^{+}[\Phi,g]_{\lambda}=\big (\phi'{\circ}f\big )
\Delta^{+}[g]_{\lambda}=\Delta^{+}\Phi\Delta^{+}g.
$$
Thus the $2$-vector function $(\Phi,g)$ has the quadratic 
$\lambda$-covariation with the bracket function 
$[\Phi,g]_{\lambda}$ satisfying (\ref{1qcofint}) provided (\ref{2qcofint})
holds.

Let $\lambda_m=\{x_i^m\colon\,i=0,\dots,n(m)\}$ for $m=1,2,\dots$.
For notation simplicity, we prove  (\ref{9qvofint}) and (\ref{2qcofint}) only 
for $u=a$ and $v=b$.
Using notation (\ref{6qvofint}) for $U$, 
by the chain rule formula (\ref{7HF}), $\Phi (x)=(\phi{\circ}f)(x)
-(\phi{\circ}f)(a)-U(x)$ for $a\leq x\leq b$.
For $i=1,\dots,n(m)$ and $m=1,2,\dots$,
by the mean value theorem, there is a $\theta_i^m\in (0,1)$ such that
$$
\Delta_i^m\Phi =\Delta_i^m(\phi{\circ}f)-\Delta_i^mU
=(\phi'{\circ}f)(x_{i-1}^m)\Delta_i^mg+B_i^m,
$$
where $\Delta_i^m\rho:=\rho(x_i^m)-\rho (x_{i-1}^m)$ for a function
$\rho$ on $[a,b]$, 
$$
B_i^m:=(\phi'{\circ}f)(x_{i-1}^m)\Delta_i^mh+\big [\phi'(y_i^m)
-\phi'(f(x_{i-1}^m))\big ]\Delta_i^mf-\Delta_i^mU
$$
and $y_i^m:=f(x_{i-1}^m)+\theta_i^m\Delta_i^mf$.
Let $\mu=\{z_j\colon\,j=1,\dots,s-1\}$ be a non-empty set of points in the
open interval $(a,b)$.
Then either $I(\mu\cap\lambda):=
\{j=1,\dots,s-1\colon\,z_j\in\cup_m\lambda_m\}\not =\emptyset$
or $I(\mu\setminus\lambda):=\{j=1,\dots,s-1\}\setminus 
I(\mu\cap\lambda)\not =\emptyset$.
If $j\in I(\mu\cap\lambda)$ then $z_j=x_{i(j)}^m$ for some index 
$i(j)=i_m(j)\in\{1,\dots,n(m)-1\}$ and for all sufficiently large $m$.
If $j\in I(\mu\setminus\lambda)$ then $x_{i(j)-1}^m<z_j<x_{i(j)}^m$ for some
index $i(j)=i_m(j)\in\{1,\dots,n(m)\}$ and for all $m$.
Let $m_0$ be the first integer such that in between each pair 
$\{z_{j-1},z_j\}$ there are at least two different points from $\lambda_{m_0}$
and $\{z_j\colon\,j\in I(\mu\cap\lambda)\}\subset\lambda_{m_0}$
provided $I(\mu\cap\lambda)\not =\emptyset$.
For $m\geq m_0$, let 
\beq\label{12qvofint}
I_1(m):=\{x_1^m,x_{n(m)}^m\}\cup\{x_{i(j)}^m,
x_{i(j)+1}^m\colon\,j\in I(\mu\cap\lambda)\}\cup\{x_{i(j)}^m\colon\,
j\in I(\mu\setminus\lambda)\}
\eeq
and let  $I_2(m):=\{1,\dots,n(m)\}\setminus I_1(m)$.
Then we have
\begin{eqnarray*}
\lefteqn{s_2(\Phi;\lambda_m)=\sum_{i\in I_1(m)}\big\{\Delta_i^m\Phi\big\}^2
+\sum_{i\in I_2(m)}\big\{(\phi'{\circ}f)(x_{i-1}^m)\Delta_i^mg
+B_i^m\big\}^2}\\[2mm]
&=&\sum_{i\in I_1(m)}\Big [\big\{\Delta_i^m\Phi\big\}^2-
\big\{(\phi'{\circ}f)(x_{i-1}^m)\Delta_i^mg\big\}^2\Big ]
+\sum_{i=1}^{n(m)}\big\{(\phi'{\circ}f)(x_{i-1}^m)\Delta_i^mg\big\}^2
+R_m^{(1)}
\end{eqnarray*}
for each $m\geq m_0$, where
$$
R_m^{(1)}:=2\sum_{i\in I_2(m)}\big\{(\phi'{\circ}f)(x_{i-1}^m)\Delta_i^mg\big\}
B_i^m+\sum_{i\in I_2(m)}\big\{B_i^m\big\}^2.
$$
Likewise, for $m\geq m_0$, we have
\begin{eqnarray*}
\lefteqn{C(\Phi,g;\lambda_m)=\sum_{i\in I_1(m)}\Delta_i^m\Phi\Delta_i^mg
+\sum_{i\in I_2(m)}\big\{(\phi'{\circ}f)(x_{i-1}^m)\Delta_i^mg
+B_i^m\big\}\Delta_i^mg}\\[2mm]
&=&\sum_{i\in I_1(m)}\Big [\Delta_i^m\Phi\Delta_i^mg-
(\phi'{\circ}f)(x_{i-1}^m)\big\{\Delta_i^mg\big\}^2\Big ]
+\sum_{i=1}^{n(m)}(\phi'{\circ}f)(x_{i-1}^m)\big\{\Delta_i^mg\big\}^2
+R_m^{(2)},
\end{eqnarray*}
where $R_m^{(2)}:=\sum_{i\in I_2(m)}B_i^m\Delta_i^mg$.
Suppose that given $\epsilon>0$ one can find a set $\mu$ for (\ref{12qvofint})
such that
$|R_m^{(1)}|\vee |R_m^{(2)}|<C\epsilon$ for all sufficiently large $m$.
Since $(g,h)\in D_{\lambda,p}^{+}(f)$, and so
$\Delta^{+}g(z_j)=0$ for $j\in I(\mu\setminus\lambda)$, we have
$$
\lim_{m\to\infty}\sum_{i\in I_1(m)}\Big [\big (\Delta_i^m\Phi\big )^2-
\big ((\phi'{\circ}f)(x_{i-1}^m)\Delta_i^mg\big )^2\Big ]
$$
$$
=\sum_{\mu\cup\{b\}}\Big\{\big (\Delta^{-}\Phi\big )^2-\big (
(\phi'{\circ}f)_{-}\Delta^{-}g\big )^2\Big\}
+\sum_{\{a\}\cup\mu}\Big\{\big (\Delta^{+}\Phi\big )^2-\big (
(\phi'{\circ}f)\Delta^{+}g\big )^2\Big\}=0
$$
by condition (\ref{2LCint}) of Definition \ref{LCint} for any given $\mu$.
Likewise for any given $\mu$, we have
$$
\lim_{m\to\infty}\sum_{i\in I_1(m)}\Big [\Delta_i^m\Phi\Delta_i^mg-
(\phi'{\circ}f)(x_{i-1}^m)\big (\Delta_i^mg\big )^2\Big ]
$$
$$
=\sum_{\mu\cup\{b\}}\Big\{\Delta^{-}\Phi\Delta^{-}g-
(\phi'{\circ}f)_{-}\big (\Delta^{-}g\big )^2\Big\}
+\sum_{\{a\}\cup\mu}\Big\{\Delta^{+}\Phi\Delta^{+}g-
(\phi'{\circ}f)\big (\Delta^{+}g\big )^2\Big\}=0.
$$
Since $\{x\in (a,b)\colon\,\Delta^{+}g(x)\not =0\}\subset \cup_m
\lambda_m$, we can and do apply Theorem \ref{intqv} which then yields
$$
\lim_{m\to\infty}s_2(\Phi;\lambda_m)=\lim_{m\to\infty}\sum_{i=1}^{n(m)}
(\phi'{\circ}f)^2(x_{i-1}^m)\big (\Delta_i^mg\big )^2
=(LY)\int_a^b(\phi'{\circ}f)^2\,d[g]_{\lambda}
$$
and
$$
\lim_{m\to\infty}C(\Phi,g;\lambda_m)=\lim_{m\to\infty}\sum_{i=1}^{n(m)}
(\phi'{\circ}f)(x_{i-1}^m)\big (\Delta_i^mg\big )^2
=(LY)\int_a^b(\phi'{\circ}f)\,d[g]_{\lambda}.
$$
Thus (\ref{9qvofint}) and (\ref{2qcofint}) hold for $u=a$ and $u=b$
provided given $\epsilon \in (0,1)$ one can find a set $\mu\subset (a,b)$ 
for (\ref{12qvofint}) such that $|R_m^{(1)}|\vee |R_m^{(2)}|<C\epsilon$ 
for all sufficiently large $m$.
We show this for the remainder $R_m^{(1)}$ because
the bound of the second remainder $R_m^{(2)}$ is similar.
Let $\epsilon \in (0,1)$.
Since $\phi'$ is uniformly continuous on the range of $f$, 
there exists a $\delta >0$ such that
$\big |\phi'(u)-\phi'(v)\big |<\epsilon $
whenever $|u|\vee |v|\leq\|f\|_{\infty}$ and $|u-v|<\delta$.
In the case $p\in (1,2)$ there is a $q<p/(p-1)$ such that
$f\in {\cal W}_q^{\ast}$.
Let $z_0:=a$, $z_s:=b$,
$\Delta^{-}S:=\Delta^{-}(\phi{\circ}f)-(\phi'{\circ}f)_{-}\Delta^{-}f$
and $\Delta^{+}S:=\Delta^{+}(\phi{\circ}f)-(\phi'{\circ}f)\Delta^{+}f$.
Then choose a set $\mu=\{z_j\colon\,j=1,\dots,s-1\}\subset (a,b)$
such that
$$
\max_{1\leq j\leq s}\osc (f;(z_{j-1},z_j))<\delta,\qquad
\max_{1\leq j\leq s}\osc (g;(z_{j-1},z_j))<\epsilon,
$$
$$
\max_{1\leq j\leq s}\sup_{z_{j-1}\leq u<v\leq z_j}(RS)\int_u^v
(\phi''{\circ}f)\,d[g]^c<\epsilon,\qquad
\sum_{(a,b)\setminus\mu}\big\{|\Delta^{-}S|+|\Delta^{+}S|\big\}
<\epsilon,
$$
and in addition for the case $p\in (1,2)$,
$$
\sum_{j=1}^sv_q(f;(z_{j-1},z_j))<\epsilon\quad\mbox{by Lemma 
\ref{variation}.\ref{lemma1}, and}\quad
\sum_{(a,b)\setminus\mu}|\Delta^{-}f|^q <\epsilon.
$$
To bound $|R_m^{(1)}|$, recall that
\begin{eqnarray*}
B_i^m&=&\big [\phi'(y_i^m)-\phi'(f(x_{i-1}^m))\big ]\Delta_i^mf
-\Big [(LY)\int_{x_{i-1}^m}^{x_i^m}(\phi'{\circ}f)\,dh-(\phi'{\circ}f)
(x_{i-1}^m)\Delta_i^mh\Big ]\\
&&-\frac{1}{2}(RS)\int_{x_{i-1}^m}^{x_i^m}(\phi''{\circ}f)\,d[g]_{\lambda}^c
-\sum_{(x_{i-1}^m,x_i^m]}\Delta^{-}S-\sum_{[x_{i-1}^m,x_i^m)}
\Delta^{+}S
\end{eqnarray*}
for $i=1,\dots,n(m)$.
By Definition \ref{LYandRY} of the $LY$ integral,
the absolute value of the second term on the right side
is equal to
$$
D_i^m:=
\Big |(RRS)\int_{x_{i-1}^m}^{x_i^m}\Big\{(\phi'{\circ}f)_{-}^{(x_{i-1}^m)}
-(\phi'{\circ}f)(x_{i-1}^m)\Big\}\,dh_{+}^{(x_i^m)}
+\sum_{(x_{i-1}^m,x_i^m)}\Delta^{-}(\phi'{\circ}f)\Delta^{+}h\Big |.
$$
In the case $p\in (1,2)$, by inequality (\ref{1LCY})
and H\"older's inequality, we have the bound
\begin{eqnarray*}
D_i^m &\leq& KV_q\big ((\phi'{\circ}f)_{-}^{(x_{i-1}^m)};[x_{i-1}^m,x_i^m]
\big )V_p\big (h_{+}^{(x_i^m)};[x_{i-1}^m,x_i^m]\big )\\[2mm]
& &+\Big (\sum_{(x_{i-1}^m,x_i^m)}\big |\Delta^{-}(\phi'{\circ}f)\big |^q
\Big )^{1/q}\Big (\sum_{(x_{i-1}^m,x_i^m)}\big |\Delta^{+}h\big |^p
\Big )^{1/p}\\[2mm]
&\leq &K_1V_q\big (f;[x_{i-1}^m,x_i^m]\big )V_p\big (h;
[x_{i-1}^m,x_i^m]\big )
+K_2\Big (\sum_{(x_{i-1}^m,x_i^m)}\big |\Delta^{-}f\big |^q
\Big )^{1/q}\Big (\sum_{(x_{i-1}^m,x_i^m)}\big |\Delta^{+}h\big |^p
\Big )^{1/p}
\end{eqnarray*}
for some finite constants $K_1$ and $K_2$.
In the case $p=1$, again by inequality (\ref{1LCY}), we have the bound
$$
D_i^m\leq K_3\osc\,(f;[x_{i-1}^m,x_i^m])\Big\{V_1(h;[x_{i-1}^m,x_i^m])
+\sum_{(x_{i-1}^m,x_i^m)}\big |\Delta^{+}h\big |\Big\}
$$
for some finite constant $K_3$.
Let $m_0$ be the integer defined in the first part of the proof.
Then in the case $p\in (1,2)$, for all $m\geq m_0$, we have
\begin{eqnarray*}
|R_m^{(1)}|&\leq&2\|\phi'{\circ}f\|_{\infty}\sum_{i\in I_2(m)}|\Delta_i^mg|
|B_i^m|+\sum_{i\in I_2(m)}\big\{B_i^m\big\}^2\\[2mm]
&\leq&\|\phi'{\circ}f\|_{\infty}\Big\{\max_{i\in I_2(m)}|\phi'(
y_i^m)-\phi'(f(x_{i-1}^m))|\big (3[g]_{\lambda}(b)+v_2(h)\big )\\[2mm]
&&+2K_1V_p(h)V_q(g)\Big (\sum_{j=1}^sv_q(f;(z_{j-1},z_j))
\Big )^{1/q}+2K_2\sigma_p(h)^{1/p}V_q(g)
\Big (\sum_{(a,b)\setminus\mu}|\Delta^{-}f|^q\Big )^{1/q}\\[2mm]
&&+\max_{i\in I_2(m)}|\Delta_i^mg|\Big ((RS)\int_a^b\phi''{\circ}f\,
d[g]_{\lambda}^c+2\sum_{(a,b)}\big\{|\Delta^{-}S|
+|\Delta^{+}S|\big\}\Big )\Big\}\\[2mm]
&&+4\max_{i\in I_2(m)}\big [\phi'(y_i^m)-\phi'(f(x_{i-1}))\big ]^2
\big ([g]_{\lambda}(b)+v_2(h)\big )\\[2mm]
&&+2K_1^2V_p(h)^2\Big (\sum_{j=1}^sv_q(f;(z_{j-1},z_j)
\Big )^{2/q}+2K_2\sigma_p(h)^{2/p}\Big (\sum_{(a,b)\setminus\mu}
|\Delta^{-}f|^q\Big )^{2/q}\\[2mm]
&&+\frac{1}{2}(RS)\int_a^b\phi''{\circ}f\,d[g]_{\lambda}^c\max_{i\in I_2(m)}
(RS)\int_{x_{i-1}^m}^{x_i^m}\!\phi''{\circ}f\,d[g]_{\lambda}^c
+2\Big (\sum_{(a,b)\setminus\mu}\big\{|\Delta^{-}S|+
|\Delta^{+}S|\big\}\Big )^{2}\\
&\leq& C\epsilon
\end{eqnarray*}
for a constant $C$ independent of $m$ by the choice of the set 
$\mu$.
The same bound with a different constant $C$ follows in the case
$p=1$.
The proof of Theorem \ref{qvofint} is complete.
\qed\end{proof}

Finally we are ready to prove Theorem \ref{introduction}.\ref{invar}.

\medskip
\noindent
{\bf Proof of Theorem \ref{introduction}.\ref{invar}.}
Suppose that $(g,h)\in D_{\lambda,p}^{+}(f)$ is such that $f\equiv C+g+h$.
By Theorem \ref{chainrule}, decomposition (\ref{10qvofint}) holds.
By the preceding Theorem \ref{qvofint}, $\Phi$ has the quadratic
$\lambda$-variation, and
by Theorem \ref{variation}.\ref{LYineq}, the function $U$ has bounded 
$p$-variation.
If $1<p<2$ then since $\phi$ is uniformly Lipschitz function
on the range of $f$, it follows that $\Phi$ has bounded
$q$-variation for some $q<p/(p-1)$, and so $g\in\dual ({\cal W}_p)[a,b]$.
Therefore (\ref{10qvofint}) is a $(\lambda,p)$-decomposition of $\phi{\circ}f$.
If $(g,h)\in D_{\lambda,p}^{-}(f)$ then it follows similarly that
(\ref{11qvofint}) is a $(\lambda,p)$-decomposition of $\phi{\circ}f$.
The proof is complete.
\qed

\section{Vector valued functions}\label{vector}

\begin{defn}\label{vectorqv}
{\rm Let $f_k$, $k=1,\dots,d$, $d\geq 2$, be regulated functions on $[a,b]$, 
and let $\lambda\in\Lambda [a,b]$.
We say that the $d$-vector function $\boldf=(f_1,\dots,f_d)$ has the 
\emph{quadratic $\lambda$-variation on $[a,b]$} if each $2$-vector
function $(f_k,f_l)$, $k,l=1,\dots,d$ has the quadratic 
$\lambda$-covariation on $[a,b]$.
The matrix valued function $[\boldf]_{\lambda}:= 
([f_k,f_l]_{\lambda})_{k,l=1,\dots,d}$ is called the bracket function 
of $\boldf$.}
\end{defn}

In the case $d=2$, a $2$-vector function $\boldf =(f_1,f_2)$
has the quadratic $\lambda$-variation if $f_1$, $f_2$ each has
the quadratic $\lambda$-variation and the $2$-vector function $(f_1,f_2)$
has the quadratic $\lambda$-covariation, and then the bracket
function of $\boldf$ is given by $2\times 2$-matrix function
$$
[\boldf]_{\lambda}=
\left ( \begin{array}{cc}
[f_1]_{\lambda} & [f_1,f_2]_{\lambda}\cr 
[f_2,f_1]_{\lambda} & [f_2]_{\lambda}
\end{array}\right ).
$$
To illustrate the preceding definition we reformulate
Theorem \ref{qvofint} as follows.
Let $f\equiv C+g+h$ with $(g,h)\in D_{\lambda,p}^{+}(f)$ for some 
$\lambda\in\Lambda [a,b]$ and $1\leq p<2$, 
let $\phi$ be a function of the class $C^2$
and let $\Phi$ be the indefinite Left Cauchy $\lambda$-integral
given by (\ref{8qvofint}).
Then the $2$-vector function $(\Phi,g)$ has
the quadratic $\lambda$-variation with the bracket function
$$
[(\Phi,g)]_{\lambda}(x)=
\left ( \begin{array}{cc}
[\Phi]_{\lambda} & [\Phi,g]_{\lambda}\cr 
[g,\Phi]_{\lambda} & [g]_{\lambda}
\end{array}\right ) (x)=(LY)\int_a^x
\left ( \begin{array}{cc}
(\phi'{\circ}f)^2 & \phi'{\circ}f\cr 
\phi'{\circ}f & 1
\end{array}\right )\,d[g]_{\lambda}
$$
for $a\leq x\leq b$.
Likewise for the indefinite Right Cauchy integral $\Psi$
given by (\ref{5qvofint}), if $(g,h)\in D_{\lambda,p}^{-}(f)$ then the 
$2$-vector function $(\Psi,g)$ 
has the quadratic $\lambda$-variation with the bracket function
$$
[(\Psi,g)]_{\lambda}(x)=
\left ( \begin{array}{cc}
[\Psi]_{\lambda} & [\Psi,g]_{\lambda}\cr 
[g,\Psi]_{\lambda} & [g]_{\lambda}
\end{array}\right ) (x)=(RY)\int_a^x
\left ( \begin{array}{cc}
(\phi'{\circ}f)^2 & \phi'{\circ}f\cr 
\phi'{\circ}f & 1
\end{array}\right )\,d[g]_{\lambda}
$$
for $a\leq x\leq b$.

\begin{prop}\label{equivqv}
Let $f_k\in {\cal R}[a,b]$, $k=1,\dots,d$, and let $\lambda\in
\Lambda [a,b]$.
The following two statements are equivalent{\rm :}
\begin{enumerate}
\item[$(a)$] the $d$-vector function $(f_1,\dots,f_d)$ has
the quadratic $\lambda$-variation on $[a,b]${\rm ;}
\item[$(b)$] for each pair $(k,l)$ of indices $k,l=1,\dots,d$,
the functions $f_k+f_l$ and $f_k-f_l$ have the quadratic $\lambda$-variation 
on $[a,b]$.
\end{enumerate}
If either of the two statements hold, then for each pair $(k,l)$, $k,l=1,\dots,d$, 
\beq\label{9bracket}
[f_k,f_l]_{\lambda}=\frac{1}{4}\Big \{[f_k+f_l]_{\lambda}
-[f_k-f_l]_{\lambda}\Big\}.
\eeq
\end{prop}

\begin{proof}
$(a)\Rightarrow (b)$:
Let $(k,l)$ be a pair of indices $k,l\in\{1,\dots,d\}$ such that
$k\not =l$.
By Definition \ref{vectorqv}, there exist the bracket functions
$[f_k]_{\lambda}$, $[f_l]_{\lambda}$ and $[f_k,f_l]_{\lambda}$.
Let $[f_k+f_l]_{\lambda}:=[f_k]_{\lambda}+2[f_k,f_l]_{\lambda}
+[f_l]_{\lambda}$ and $[f_k-f_l]_{\lambda}:=[f_k]_{\lambda}
-2[f_k,f_l]_{\lambda}+[f_l]_{\lambda}$.
It is easy to check that $[f_k+f_l]_{\lambda}$ and 
$[f_k-f_l]_{\lambda}$ are the bracket functions of $f_k+f_l$ and
$f_k-f_l$, respectively.

$(b)\Rightarrow (a)$:
Given a pair $(k,l)$ of indices, let 
$[f_k,f_l]_{\lambda}:=(1/4)\big\{[f_k+f_l]_{\lambda}
-[f_k-f_l]_{\lambda}\big\}$.
It is then easy to check that the $2$-vector function
$(f_k,f_l)$ has the quadratic $\lambda$-covariation with the bracket 
function  $[f_k,f_l]_{\lambda}$. 
The proof of Proposition \ref{equivqv} is complete.
\qed\end{proof}

Suppose that a $2$-vector function $(f_1,f_2)$ on $[a,b]$ has the 
quadratic $\lambda$-variation for some $\lambda\in\Lambda [a,b]$.
Let
$$
[f_1,f_2]_{\lambda}^c:=\frac{1}{4}\Big\{[f_1+f_2]_{\lambda}^c
-[f_1-f_2]_{\lambda}^c\Big\}.
$$
Then the following is true.

\begin{prop}\label{contqc}
Let $(f_1,f_2)$ be a $2$-vector function on $[a,b]$ 
having the quadratic $\lambda$-variation for some $\lambda\in
\Lambda [a,b]$.
Then for each $a\leq x\leq b$,
\beq\label{1contqc}
[f_1,f_2]_{\lambda}(x)=[f_1,f_2]_{\lambda}^c(x)
+\sum_{(a,x]}\Delta^{-}f_1\Delta^{-}f_2+\sum_{[a,x)}
\Delta^{+}f_1\Delta^{+}f_2.
\eeq
\end{prop}

\begin{proof}
Since $[f_1+f_2]_{\lambda}^c$ and $[f_1-f_2]_{\lambda}^c$ are
continuous parts of the bracket functions $[f_1+f_2]_{\lambda}$
and $[f_1-f_2]_{\lambda}$, respectively, it follows that for each
$a\leq x\leq b$,
$$
\frac{1}{4}\Big\{[f_1+f_2]_{\lambda}^c(x)-[f_1-f_2]_{\lambda}^c(x)\Big\}
=\frac{1}{4}\Big\{[f_1+f_2]_{\lambda}(x)-[f_1-f_2]_{\lambda}(x)\Big\}
$$
$$
-\frac{1}{4}\Big\{\sum_{(a,x]}\big\{\Delta^{-}(f_1+f_2)\big\}^2
-\sum_{(a,x]}\big\{\Delta^{-}(f_1-f_2)\big\}^2\Big\}
-\frac{1}{4}\Big\{\sum_{[a,x)}\big\{\Delta^{+}(f_1+f_2)\big\}^2
-\sum_{[a,x)}\big\{\Delta^{+}(f_1-f_2)\big\}^2\Big\}
$$
$$
=[f_1,f_2]_{\lambda}(x)-\sum_{(a,x]}\Delta^{-}f_1\Delta^{-}f_2
-\sum_{[a,x)}\Delta^{+}f_1\Delta^{+}f_2,
$$
where the last equality follows by
linearity of the unconditional sums and by the preceding
Proposition \ref{equivqv}.
Therefore (\ref{1contqc}) holds by the definition of $[f_1,f_2]_{\lambda}^c$.
\qed\end{proof}

\begin{prop}\label{covardec}
For $\lambda\in\Lambda [a,b]$, $1\leq p<2$ and $k=1,2$, let
$f_k$ be a function on $[a,b]$ which is $(\lambda,p)$-decomposable by 
$(g_k,h_k)\in D_{\lambda,p}^{+}(f_k)\cup D_{\lambda,p}^{-}(f_k)$.
The $2$-vector function $(f_1,f_2)$ has the quadratic $\lambda$-covariation 
if and only if
the $2$-vector function $(g_1,g_2)$ has the quadratic $\lambda$-covariation.
If both $2$-vector functions have the quadratic $\lambda$-covariation then
for $a\leq x\leq b$,
\beq\label{2covardec}
[f_1,f_2]_{\lambda}(x)=[g_1,g_2]_{\lambda}^c(x)
+\sum_{(a,x]}\Delta^{-}f_1\Delta^{-}f_2
+\sum_{[a,x)}\Delta^{+}f_1\Delta^{+}f_2.
\eeq
\end{prop}

\begin{proof}
For $k=1,2$,  let $(g_k,h_k)\in D_{\lambda,p}^{+}(f_k)$.
By Proposition \ref{qcforpq}, the three $2$-vector functions
$(g_1,h_2)$, $(g_2,h_1)$ and $(h_1,h_2)$ have the quadratic
$\lambda$-covariations, and for $a\leq x\leq b$,
$$
[g_1,h_2]_{\lambda}(x)=\sum_{(a,x]}\Delta^{-}g_1\Delta^{-}h_2
+\sum_{[a,x)}\Delta^{+}g_1\Delta^{+}h_2,
$$
$$
[g_2,h_1]_{\lambda}(x)=\sum_{(a,x]}\Delta^{-}g_2\Delta^{-}h_1
+\sum_{[a,x)}\Delta^{+}g_2\Delta^{+}h_1,
$$
$$
[h_1,h_2]_{\lambda}(x)=\sum_{(a,x]}\Delta^{-}h_1\Delta^{-}h_2
+\sum_{[a,x)}\Delta^{+}h_1\Delta^{+}h_2.
$$
If $(g_k,h_k)\in D_{\lambda,p}^{-}(f_k)$ for $k=1,2$, then the same conclusion
holds.
Let $\lambda=\{\lambda_m\colon\,m\geq 1\}$.
Then for each $m\geq 1$, we have
$$
C(f_1,f_2;\lambda_m)-C(g_1,g_2;\lambda_m)
=C(g_1,h_2;\lambda_m)+C(g_2,h_1;\lambda_m)+C(h_1,h_2;\lambda_m).
$$
Since the right side has a limit as $m\to\infty$, the first part of the
conclusion holds.
By linearity of the uncondition sums, it follows that for $a\leq x\leq b$,
\begin{eqnarray*}
[f_1,f_2]_{\lambda}(x)&=&[g_1,g_2]_{\lambda}(x)+[g_1,h_2]_{\lambda}(x)
+[g_2,h_1]_{\lambda}(x)+[h_1,h_2]_{\lambda}(x)\\[2mm]
&=&[g_1,g_2]_{\lambda}^c(x)+\sum_{(a,x]}\Big\{\Delta^{-}g_1\Delta^{-}g_2
+\Delta^{-}g_1\Delta^{-}h_2+\Delta^{-}g_2\Delta^{-}h_2
+\Delta^{-}h_1\Delta^{-}h_2\Big\}\\[2mm]
& &\qquad\qquad+\sum_{[a,x)}\Big\{\Delta^{+}g_1\Delta^{+}g_2
+\Delta^{+}g_1\Delta^{+}h_2+\Delta^{+}g_2\Delta^{+}h_2
+\Delta^{+}h_1\Delta^{+}h_2\Big\},
\end{eqnarray*}
proving (\ref{2covardec}).
The proof is complete.
\qed\end{proof}

The next statement extends Theorem \ref{intqv} to bracket functions
of the quadratic $\lambda$-covariation, and hence to bracket functions
of the quadratic $\lambda$-variation of vector functions.

\begin{lem}\label{intqc}
Let $(f_1,\dots,f_d)$ have the quadratic $\lambda$-variation on $[a,b]$  
for some sequence
$\lambda=\{\lambda_m\colon\,m\geq 1\}\in\Lambda [a,b]$ of partitions
$\lambda_m=\{x_i^m\colon\,i=0,\dots,n(m)\}$, 
and let $h\in {\cal R}[a,b]$.
Then the following two statements hold{\rm :}
\begin{enumerate}
\item[$(a)$] If $\cup_k\{x\in (a,b)\colon\,\Delta^{+}f_k(x)\not =0\}
\subset\cup_m\lambda_m$, then
\beq\label{1intqc}
\lim_{m\to\infty}\sum_{i=1}^{n(m)}h(x_{i-1}^m)\big [f_k(x_i^m)
-f_k(x_{i-1}^m)\big ]\big [f_l(x_i^m)-f_l(x_{i-1}^m)\big ]
=(LY)\int_a^bh\,d[f_k,f_l]_{\lambda}.
\eeq
\item[$(b)$] If $\cup_k\{x\in (a,b)\colon\,\Delta^{-}f_k(x)\not =0\}
\subset\cup_m\lambda_m$, then
\beq\label{2intqc}
\lim_{m\to\infty}\sum_{i=1}^{n(m)}h(x_{i}^m)\big [f_k(x_i^m)
-f_k(x_{i-1}^m)\big ]\big [f_l(x_i^m)-f_l(x_{i-1}^m)\big ]
=(RY)\int_a^bh\,d[f_k,f_l]_{\lambda}.
\eeq
\end{enumerate}
\end{lem}

\begin{proof}
Since $[f_k,f_l]_{\lambda}$, $k,l=1,\dots,d$, are functions of bounded
variation on $[a,b]$, the Left Young and Right Young integrals in 
(\ref{1intqc}) and (\ref{2intqc}), respectively, exist by Lemma \ref{bv}.
For $(a)$ suppose that $\cup_k\{x\in (a,b)\colon\,\Delta^{+}f_k(x)\not =0\}
\subset\cup_m\lambda_m$.
If $k=l$ then (\ref{1intqc}) holds by Theorem \ref{intqv}.
If $k\not =l$ then $\{x\in (a,b)\colon\,\Delta^{+}(f_k+f_l)(x)\not =0\}
\subset\cup_k\{x\in (a,b)\colon\,\Delta^{+}f_k(x)\not =0\}$
and $\{x\in (a,b)\colon\,\Delta^{+}(f_k-f_l)(x)\not =0\}
\subset\cup_k\{x\in (a,b)\colon\,\Delta^{+}f_k(x)\not =0\}$.
Letting $\Delta_i^k\rho :=\rho (x_i^k)-\rho(x_{i-1}^k)$
for a function $\rho$ on $[a,b]$, by (\ref{9bracket}), by Theorem \ref{intqv}
and by the bilinearity of the Left Young integral, we have
$$
\sum_{i=1}^{n(m)}h(x_{i-1}^m)\Delta_i^mf_k\Delta_i^mf_l
=\frac{1}{4}\sum_{i=1}^{n(m)}h(x_{i-1}^m)\Big\{(\Delta_i^m(f_k+f_l))^2
-(\Delta_i^m(f_k-f_l))^2\Big\}
$$
$$
\stackrel{m\to\infty}{\longrightarrow}
\frac{1}{4}\Big\{(LY)\int_a^bh\,d[f_k+f_l]_{\lambda}
-(LY)\int_a^bh\,d[f_k-f_l]_{\lambda}\Big\}
=(LY)\int_a^bh\,d[f_k,f_l]_{\lambda},
$$
proving (\ref{1intqc}).
The proof of statement $(b)$ is analogous and is omitted.
\qed\end{proof}

Next we extend Definitions \ref{LCint} and \ref{RCint} of the two 
$\lambda$-integrals for vector-valued functions.
Recall the notation of the sums $S_{LC}$ and $S_{RC}$ given
before Definition \ref{LCint}.

\begin{defn}\label{mdLC}
{\rm Let $\boldff =(F_1,\dots,F_d)$ and $\boldf =(f_1,\dots,f_d)$ be two
$d$-vector functions  such that for each $k=1,\dots,d$
$f_k\in {\cal R}[a,b]$, $F_k(t-)$ exists at $t\in (a,b]$ whenever 
$\Delta^{-}f_k(t)\not =0$, 
and let $\lambda=\{\lambda_m\colon\,m\geq 1\}\in\Lambda [a,b]$.
We say that the \emph{Left Cauchy $\lambda$-integral 
$(LC)\smallint \langle\boldff,d_{\lambda}\boldf\rangle$
is defined on $[a,b]$} if there exists a regulated function $\Phi$ on $[a,b]$ 
such that $\Phi (a)=0$ and for each $a\leq u< v\leq b$,
\beq\label{1mdLC}
\Phi (v)-\Phi (u)=\lim_{m\to\infty}\sum_{k=1}^d
S_{LC}(F_k,f_k;\lambda_m\Cap [u,v]),
\eeq
\beq\label{2mdLC}
\Delta^{-}\Phi (v)=\sum_{k=1}^dF_k(v-)\Delta^{-}f_k(v)
\quad\mbox{and}\quad
\Delta^{+}\Phi (u)=\sum_{k=1}^dF_k(u)\Delta^{+}f_k(u),
\eeq
where $F_k(v-):=0$ whenever $\Delta^{-}f_k(v)=0$.
For $a\leq u<v\leq b$, let
$(LC)\smallint_u^v\langle\boldff,d_{\lambda}\boldf\rangle
:=\Phi (v)-\Phi (u)$.
If a regulated function $\Phi$ on $[a,b)$ exists and satisfies the same conditions
except that (\ref{1mdLC}) and (\ref{2mdLC}) do not hold for $v=b$,
then we say that the \emph{Left Cauchy $\lambda$-integral $(LC)\smallint
\langle\boldff,d_{\lambda}\boldf\rangle$ is defined on $[a,b)$}. }
\end{defn}

\begin{defn}\label{mdLRint}
{\rm Let $\boldff =(F_1,\dots,F_d)$ and $\boldf =(f_1,\dots,f_d)$ be two
$d$-vector functions  such that for each $k=1,\dots,d$
$f_k\in {\cal R}[a,b]$, $F_k(t+)$ exists at $t\in [a,b)$ whenever 
$\Delta^{+}f_k(t)\not =0$, 
and let $\lambda=\{\lambda_m\colon\,m\geq 1\}\in\Lambda [a,b]$.
We say that the \emph{Right Cauchy $\lambda$-integral 
$(RC)\smallint \langle\boldff,d_{\lambda}\boldf\rangle$
is defined on $[a,b]$} if there exists a regulated function $\Psi$ on $[a,b]$ 
such that $\Psi (a)=0$ and for each $a\leq u< v\leq b$,
$$
\Psi (y)-\Psi (z)=\lim_{m\to\infty}\sum_{k=1}^d
S_{RC}(F_k,f_k;\lambda\Cap [u,v]),
$$
$$
\Delta^{-}\Psi (v)=\sum_{k=1}^dF_k(v)\Delta^{-}f_k(v)
\quad\mbox{and}\quad
\Delta^{+}\Psi (u)=\sum_{k=1}^dF_k(u+)\Delta^{+}f_k(u),
$$
where $F_k(u+):=0$ whenever $\Delta^{+}f_k(u)=0$.
For $a\leq u<v\leq b$, let
$(RC)\smallint_u^v\langle\boldff,d_{\lambda}\boldf\rangle
:=\Psi (v)-\Psi (u)$.}
\end{defn}

Next chain rule formulas (\ref{7HF}) and (\ref{17HF})
of Theorem \ref{chainrule} are extended in two directions.
First, the increment of the composition on the left side
of each formula is replaced by an
integral with respect to the composition.
Second, the real valued functions $f$, $g$, $h$ are replaced by
$d$-vector functions $\boldf=(f_1,\dots,f_d)$,
$\boldg=(g_1,\dots,g_d)$ and $\boldh=(h_1,\dots,h_d)$ so that
for each $k=1,\dots,d$, $f_k$, $g_k$ and $h_k$ satisfy
the conditions of the case $d=1$.
For $\lambda\in\Lambda [a,b]$ and $1\leq p<2$, we say that
\emph{$\boldf$ is $(\lambda,p)$-decomposable by $(\boldg,\boldh)$}
if $\boldf=\boldC+\boldg+\boldh$ for some constant $\boldC$ and
for each $k=1,\dots,d$, $(g_k,h_k)$ is a $(\lambda,p)$-dual pair. 
If in addition (cf.\ (\ref{right-jumps-access})),
\beq\label{mdaccess}
\cup_{k=1}^d\Big (N_{(a,b)}(\Delta^{+}g_k)\cup N_{(a,b)}(\Delta^{+}h_k)
\Big )\subset\cup\lambda
\eeq
then we write $(\boldg,\boldh)\in D_{\lambda,p}^{+}(\boldf)$.
Likewise,  if in addition, 
(\ref{mdaccess}) holds with $\Delta^{+}$ replaced by $\Delta^{-}$
then we write $(\boldg,\boldh)\in D_{\lambda,p}^{-}(\boldf)$.

\begin{thm}\label{mdchrule}
Let $\boldf=(f_1,\dots,f_d)$, $\boldg=(g_1,\dots,g_d)$,
$\boldh=(h_1,\dots,h_d)$, $\lambda\in\Lambda [a,b]$ and
$1\leq p<2$ be such that $\boldf$ is $(\lambda,p)$-decomposable
by $(\boldg,\boldh)$.
Let $\phi\colon\,\RR^d\mapsto\RR$ be a function of the class $C^2$,
and let $\psi\in {\rm \dual} ({\cal W}_p)[a,b]$.
Then the following two statements hold{\rm :}
\begin{enumerate}
\item[$(a)$]
If $(\boldg,\boldh)\in D_{\lambda,p}^{+}(\boldf)$
then the Left Cauchy $\lambda$-integral
$(LC)\smallint\big\langle\big (\psi,-\psi (\nabla\phi{\circ}\boldf)
\big ),d_{\lambda}\big (\phi{\circ}\boldf,\boldg )\big\rangle$
is defined on $[a,b]$, and  for each $a\leq u< v\leq b$,
\begin{eqnarray}
\lefteqn{(LC)\int_u^v\big\langle\big (\psi,-\psi (\nabla\phi{\circ}\boldf)
\big ),d_{\lambda}\big (\phi{\circ}\boldf,\boldg )\big\rangle}
\label{1mdchrule}\\[2mm]
&=&\sum_{k=1}^d(LY)\int_u^v\psi (\phi_k'{\circ}\boldf)\,dh_k
+\frac{1}{2}\sum_{l,k=1}^d(RS)\int_u^v\psi \big 
(\phi_{lk}''{\circ}\boldf\big )\,d[g_l,g_k]_{\lambda}^c\nonumber\\[2mm]
& &+\sum_{(u,v]}\psi_{-}\Big\{\Delta^{-}(\phi{\circ}\boldf)
-\sum_{k=1}^d(\phi_k'{\circ}\boldf)_{-}\Delta^{-}f_k\Big\}
+\sum_{[u,v)}\psi\Big\{\Delta^{+}(\phi{\circ}\boldf)
-\sum_{k=1}^d(\phi_k'{\circ}\boldf)\Delta^{+}f_k\Big\},\nonumber
\end{eqnarray}\noindent
where the two sums converge unconditionally.
\item[$(b)$]
If $(\boldg,\boldh)\in D_{\lambda,p}^{-}(\boldf)$
 then the Right Cauchy $\lambda$-integral
$(RC)\smallint\big\langle\big (\psi,-\psi (\nabla\phi{\circ}\boldf)
\big ),d_{\lambda}\big (\phi{\circ}\boldf,\boldg )\big\rangle$
is defined on $[a,b]$, and  for each $a\leq u< v\leq b$,
\begin{eqnarray}\label{11mdchrule}
\lefteqn{(RC)\int_u^v\big\langle\big (\psi,-\psi (\nabla\phi{\circ}\boldf)
\big ),d_{\lambda}\big (\phi{\circ}\boldf,\boldg )\big\rangle}\\[2mm]
&=&\sum_{k=1}^d(RY)\int_u^v\psi (\phi_k'{\circ}\boldf)\,dh_k
-\frac{1}{2}\sum_{l,k=1}^d(RS)\int_u^v\psi \big 
(\phi_{lk}''{\circ}\boldf\big )\,d[g_l,g_k]_{\lambda}^c\nonumber\\[2mm]
& &+\sum_{(u,v]}\psi\Big\{\Delta^{-}(\phi{\circ}\boldf)
-\sum_{k=1}^d(\phi_k'{\circ}\boldf)\Delta^{-}f_k\Big\}
+\sum_{[u,v)}\psi_{+}\Big\{\Delta^{+}(\phi{\circ}\boldf)
-\sum_{k=1}^d(\phi_k'{\circ}\boldf)_{+}\Delta^{+}f_k\Big\},\nonumber
\end{eqnarray}\noindent
where the two sums converge unconditionally.
\end{enumerate}
\end{thm}

\begin{proof}
We prove statement $(a)$ and indicate necessary changes needed to
prove statement $(b)$.
Suppose that (\ref{mdaccess}) holds.
In the case $1<p<2$ there exists $2<q<p/(p-1)$ such that
$g_k$ has bounded $q$-variation for $k=1,\dots,d$.
We can and do assume that in this case $\psi$ has bounded $q$-variation 
as well, and $\psi$ is a regulated function in the case $p=1$.
Let $\phi\colon\,\RR^d\mapsto\RR$ be a function of the class $C^2$.
Statement $(a)$ will be proved once we show
that the limit
\begin{eqnarray}
\lefteqn{(LC)\int_u^v\big\langle\big (\psi,-\psi (\nabla\phi{\circ}\boldf)
\big ),d_{\lambda}\big (\phi{\circ}\boldf,\boldg )\big\rangle}
\label{9mdchrule}\\[2mm]
&:=&\lim_{m\to\infty}\left [S_{LC}(\psi,\phi{\circ}\boldf;
\lambda\Cap [u,v])-\sum_{k=1}^dS_{LC}(\psi (\phi_k'{\circ}
\boldf),g_k;\lambda\Cap [u,v])\right ]\nonumber
\end{eqnarray}
exists for each $a\leq u<v\leq b$ and satisfies relation
(\ref{1mdchrule}).
Indeed, let 
$$
\Phi (x):=(LC)\int_a^x\big\langle\big (\psi,-\psi 
(\nabla\phi{\circ}\boldf)
\big ),d_{\lambda}\big (\phi{\circ}\boldf,\boldg )\big\rangle,
\qquad a\leq x\leq b.
$$
Clearly, condition (\ref{1mdLC}) holds by definition of $\Phi$, while 
condition (\ref{2mdLC})  follows from relation (\ref{1mdchrule})
as in the proof of Theorem \ref{chainrule}.

For notation simplicity, we prove that (\ref{9mdchrule}) exists and 
(\ref{1mdchrule}) holds for $u=a$ and $v=b$.
If $1<p<2$ then $f_k\in {\cal W}_q$ for $2<q<p/(p-1)$ and $k=1,\dots,d$.
Since $\phi'$ is Lipschitz function on the range of $\boldf$,
and since $\WW_q$ is a Banach algebra,
$\psi (\phi_k'{\circ}\boldf)\in {\cal W}_q$ for $k=1,\dots,d$. 
By Theorem \ref{variation}.\ref{LYineq} if $1<p<2$
and by Lemma \ref{bv} if $p=1$, it then follows that there exist the Left Young
integrals $A_k:=(LY)\smallint_a^b\psi (\phi_k'{\circ}\boldf)\,dh_k$
for $k=1,\dots,d$.
The Riemann-Stieltjes integrals $B_{l,k}:=(RS)\smallint_a^b
\psi (\phi_{lk}''{\circ}\boldf)\,d[g_l,g_k]_{\lambda}^c$ for $l, k=
1,\dots,d$ exist by Lemma \ref{bv}.
The mean value theorem assure that the two sums in 
(\ref{1mdchrule}) converge absolutely because $\max_k\sigma_2(f_k)
<\infty$, and thus converge unconditionally to values, say
$C^{-}$ and $C^{+}$, respectively.
Let $\lambda_m=\{x_i^m\colon\,i=0,\dots,n(m)\}$ for $m=1,2,\dots$.
For each $m\geq 1$ and $i=1,\dots,n(m)$, let $\Delta_i^m\rho
:=\rho (x_i^m)-\rho (x_{i-1}^m)$ for a function $\rho$ on $[a,b]$.
Let $\epsilon >0$.
Since $\cup_k\{x\in (a,b)\colon\,\Delta^{+}h_k(x)\not =0\}\subset
\cup_m\lambda_m$, by Theorem \ref{variation}.\ref{approximation},
there exists an integer $m_1$ such that for all $m\geq m_1$,
\beq\label{2mdchrule}
\left |A_k-S_{LC}(\psi (\phi_k'{\circ}\boldf),h_k;\lambda_m)\right |
<\epsilon
\eeq
for $k=1,\dots,d$.
Since $\cup_k\{x\in (a,b)\colon\,\Delta^{+}g_k(x)\not =0\}\subset
\cup_m\lambda_m$, by Lemmas \ref{intqc} and \ref{bv}, there exists 
an integer $m_2$ such that for all $m\geq m_2$,
\beq\label{3mdchrule}
\Big |B_{l,k}+B_{l,k}^{-}+B_{l,k}^{+}-\sum_{i=1}^{n(m)}\big (\psi
(\phi_{lk}''{\circ}\boldf)\big )(x_{i-1}^m)\Delta_i^mg_l\Delta_i^mg_k
\Big |<\epsilon
\eeq
for $l, k=1,\dots,d$, where
$$
B_{l,k}^{-}:=\sum_{(a,b]}\psi_{-}(\phi_{lk}''{\circ}\boldf)_{-}\Delta^{-}
g_l\Delta^{-}g_k\quad\mbox{and}\quad
B_{l,k}^{+}:=\sum_{[a,b)}\psi (\phi_{lk}''{\circ}\boldf)\Delta^{+}
g_l\Delta^{+}g_k.
$$
Since $\phi_{lk}''$ is uniformly continuous on the range
of $\boldf$, there is a $\delta >0$ such that
$|\phi_{lk}''(u)-\phi_{lk}''(v)|<\epsilon$ for $l,k=1,\dots,d$
whenever $\max_r|u_r-v_r|<\delta$ and $|u_r|\vee |v_r|\leq
\|f_r\|_{\infty}$ for $r=1,\dots,d$ and $u=(u_1,\dots,u_d)$,
$v=(v_1,\dots,v_d)$.
Therefore recalling Lemma \ref{variation}.\ref{lemma1}, we conclude that there 
exists a finite set $\mu=\{z_j\colon\,j=1,\dots,s-1\}\subset (a,b)$ 
such that  for each $l,k=1,\dots,d$,
\beq\label{4mdchrule}
\Big |C^{-}-\sum_{\mu\cup\{b\}}\psi_{-}\Big\{\Delta^{-}(\phi
{\circ}\boldf)-\sum_{k=1}^d(\phi_k'{\circ}\boldf)_{-}\Delta^{-}f_k
\Big\}\Big |\vee\Big |B_{lk}^{-}-\sum_{\mu\cup\{b\}}\psi_{-}(\phi_{lk}''
{\circ}\boldf)_{-}\Delta^{-}g_l\Delta^{-}g_k\Big |<\epsilon,
\eeq
\beq\label{5mdchrule}
\Big |C^{+}-\sum_{\{a\}\cup\mu}\psi\Big\{\Delta^{+}(\phi
{\circ}\boldf)-\sum_{k=1}^d(\phi_k'{\circ}\boldf)\Delta^{+}f_k
\Big\}\Big |\vee\Big |B_{lk}^{+}-\sum_{\{a\}\cup\mu}\psi (\phi_{lk}''{\circ}
\boldf)\Delta^{+}g_l\Delta^{+}g_k\Big |<\epsilon,
\eeq
\beq\label{6mdchrule}
\max_{1\leq k\leq d}\max_{1\leq j\leq s}\osc (f_k;(z_{j-1},z_j))
<\delta\quad\mbox{and}\quad
\max_{1\leq k\leq d}\sum_{j=1}^sv_2(h_k;(z_{j-1},z_j))<\epsilon,
\eeq
where $z_0:=a$ and $z_s:=b$.
Let $I(\mu\cap\lambda):=\{j=1,\dots,s-1\colon\,z_j\in\cup_m\lambda_m\}$.
If $j\in I(\mu\cap\lambda)$ then $z_j=x_{i(j)}^m$ for some index 
$i(j)\in\{1,\dots,n(m)-1\}$ and for all sufficiently large $m$.
If $j\in I(\mu\setminus\lambda):=\{j=1,\dots,s-1\}\setminus 
I(\mu\cap\lambda)$ then $x_{i(j)-1}^m<z_j<x_{i(j)}^m$ for some
index $i(j)\in\{1,\dots,n(m)\}$ and for all $m$.
Let $m_0$ be the first integer such that $\{z_j\colon\,j\in
I(\mu\cap\lambda)\}\subset\lambda_{m_0}$ and in between each pair $(z_{j-1},
z_j)$ there are at least two different points from $\lambda_{m_0}$.
For $m\geq m_0$, let $I_1(m):=\{x_1^m,x_{n(m)}^m\}\cup\{x_{i(j)}^m,
x_{i(j)+1}^m\colon\,j\in I(\mu\cap\lambda)\}\cup\{x_{i(j)}^m\colon\,
j\in I(\mu\setminus\lambda)\}$ and let 
$I_2(m):=\{1,\dots,n(m)\}\setminus I_1(m)$.
Since $\Delta^{+}g_k(z_j)=\Delta^{+}h_k(z_j)=0$ for $z_j\not\in\cup_m
\lambda_m$ and $k=1,\dots,d$, we have that
$$
\lim_{m\to\infty}
\sum_{i\in I_1(m)}\psi (x_{i-1}^m)\Big\{\Delta_i^m\phi{\circ}
\boldf-\sum_{k=1}^d(\phi_k'{\circ}\boldf)(x_{i-1}^m)\Delta_i^mf_k\Big\}
$$
$$
= \sum_{\mu\cup\{b\}}\psi_{-}\Big\{\Delta^{-}(\phi
{\circ}\boldf)-\sum_{k=1}^d(\phi_k'{\circ}\boldf)_{-}\Delta^{-}f_k
\Big\}+ \sum_{\{a\}\cup\mu}\psi\Big\{\Delta^{+}(\phi
{\circ}\boldf)-\sum_{k=1}^d(\phi_k'{\circ}\boldf)\Delta^{+}f_k
\Big\}
$$
and
$$
\lim_{m\to\infty}\sum_{i\in I_1(m)}\big (\psi (\phi_{lk}''{\circ}\boldf)\big )
(x_{i-1}^m)\Delta_i^mg_l\Delta_i^mg_k
$$
$$
=\sum_{\mu\cup\{b\}}\psi_{-}(\phi_{lk}''
{\circ}\boldf)_{-}\Delta^{-}g_l\Delta^{-}g_k
+\sum_{\{a\}\cup\mu}\psi (\phi_{lk}''{\circ}
\boldf)\Delta^{+}g_l\Delta^{+}g_k.
$$
Therefore by (\ref{4mdchrule}) and (\ref{5mdchrule}), there exists an integer 
$m_3\geq m_0$ such that for all $m\geq m_3$
\beq\label{7mdchrule}
\Big |\sum_{i\in I_1(m)}\psi (x_{i-1}^m)\Big\{\Delta_i^m\phi{\circ}
\boldf-\sum_{k=1}^d(\phi_k'{\circ}\boldf)(x_{i-1}^m)\Delta_i^mf_k\Big\}
-C^{-}-C^{+}\Big |<4\epsilon
\eeq
and
\beq\label{8mdchrule}
\Big |\sum_{i\in I_1(m)}\big (\psi (\phi_{lk}''{\circ}\boldf)\big )
(x_{i-1}^m)\Delta_i^mg_l\Delta_i^mg_k -B_{lk}^{-}-B_{lk}^{+}\Big |
<4\epsilon
\eeq
for each $l,k =1,\dots,d$.
We use Taylor's theorem with Lagrange's form of the remainder:
for vectors $u=(u_1,\dots,u_d)$, $v=(v_1,\dots,v_d)$ there is
a $\theta=\theta (u,v)\in (0,1)$ such that
$$
\phi (v)=\phi (u)+\sum_{k=1}^d\phi_k'(u)[v_k-u_k]+\frac{1}{2}
\sum_{l,k=1}^d\phi_{lk}''(u+\theta (v-u))[v_l-u_l][v_k-u_k],
$$
where $\phi_l'(y)=\frac{\partial \phi}{\partial y_l}(y)$
and $\phi_{lk}''(y)=\frac{\partial^2\phi}{\partial y_l\partial y_k}(y)$
for $y=(y_1,\dots,y_d)$.
Then by a telescoping sum, for $m\geq m_0$,
it follows that
\begin{eqnarray}\label{10mdchrule}
\lefteqn{S_{LC}(\psi,\phi{\circ}\boldf;\lambda_m)
=\sum_{i\in I_1(m)}\psi(x_{i-1}^m)\Delta_i^m(\phi{\circ}\boldf)}\\[2mm]
& &+\sum_{i\in I_2(m)}\psi (x_{i-1}^m)\Big\{\sum_{k=1}^d(\phi_k'{\circ}\boldf)
(x_{i-1}^m)\Delta_i^mf_k+\frac{1}{2}\sum_{l,k=1}^d
\phi_{lk}''(y_i^m)\Delta_i^mf_l\Delta_i^mf_k\Big\}\nonumber\\[2mm]
&=&\sum_{k=1}^dS_{LC}(\psi(\phi_k'{\circ}\boldf),g_k;\lambda_m)
+\sum_{k=1}^dS_{LC}(\psi(\phi_k'{\circ}\boldf),h_k;\lambda_m)
\nonumber\\[2mm]
& &+\frac{1}{2}\sum_{l,k=1}^d\left (\sum_{i=1}^{n(m)}-\sum_{i\in I_1(m)}
\right )\big (\psi (\phi_{lk}''{\circ}\boldf)\big )
(x_{i-1}^m)\Delta_i^mg_l\Delta_i^mg_k\nonumber\\[2mm]
& &+\sum_{i\in I_1(m)}\psi (x_{i-1}^m)\Big\{\Delta_i^m
(\phi{\circ}\boldf)-\sum_{k=1}^d(\phi_k'{\circ}\boldf)(x_{i-1}^m)
\Delta_i^mf_k\Big\}+R_m,\nonumber
\end{eqnarray}
where $y_i^m:=\boldf(x_{i-1}^m)+\theta_i^m\Delta_i^m\boldf$, 
$\theta_i^m\in (0,1)$ for $i=1,\dots,n(m)$, and
\begin{eqnarray*}
R_m&:=&\frac{1}{2}\sum_{l,k=1}^d\sum_{i\in I_2(m)}\Big\{
\psi (x_{i-1}^m)\phi_{lk}''
(y_i^m)\big [\Delta_i^mg_l\Delta_i^mh_k+\Delta_i^mh_l\Delta_i^mg_k
+\Delta_i^mh_l\Delta_i^mh_k\big ]\\[2mm]
& &\qquad\qquad\qquad+\psi (x_{i-1}^m)\big [\phi_{lk}''(y_i^m)-
\phi_{lk}''(\boldf(x_{i-1}^m))\big ]\Delta_i^mg_k\Delta_i^mg_l\Big\}.
\end{eqnarray*}
Then by (\ref{2mdchrule}), (\ref{3mdchrule}), (\ref{7mdchrule}) and
(\ref{8mdchrule}), for all $m\geq m_1\vee m_1\vee m_3$,
$$
\Big |S_{LC}(\psi,\phi{\circ}\boldf;\lambda_m)
-\sum_{k=1}^dS_{LC}(\psi(\phi_k'{\circ}\boldf),g_k;\lambda_m)
-\sum_{k=1}^dA_k-\frac{1}{2}\sum_{l,k=1}^dB_{l,k}-C^{-}-C^{+}\Big |
$$
$$
\leq\sum_{k=1}^d\Big |A_k-
S_{LC}(\psi(\phi_k'{\circ}\boldf),h_k;\lambda_m)\Big |
+\frac{1}{2}\sum_{l,k=1}^d\Big |B_{l,k}-\sum_{i\in I_2(m)}\big (\psi
(\phi_{lk}''{\circ}\boldf)\big )(x_{i-1}^m)\Delta_i^mg_l\Delta_i^mg_k\Big |
$$
$$
+\Big |C^{-}+C^{+}
-\!\sum_{i\in I_1(m)}\psi (x_{i-1}^m)\Big\{\Delta_i^m(\phi{\circ}\boldf)
-\sum_{k=1}^d(\phi_k'{\circ}\boldf)(x_{i-1}^m)\Delta_i^mf_k\Big\}\Big |
+|R_m|
\leq d\epsilon+(5/2)d^2\epsilon +4\epsilon +|R_m|.
$$
To bound $|R_m|$ notice that for each $i\in I_2(m)$, $[x_{i-1}^m,
x_i^m]\subset (z_{j-1},z_j)$ for some $j=1,\dots,s$.
By (\ref{6mdchrule}) and H\"older's inequality, it then follows
$$
\limsup_{m\to\infty}|R_m|\leq
\frac{\sqrt{\epsilon}}{2}\sum_{l,k=1}^d\|\phi_{lk}''{\circ}
\boldf\|_{\infty}\big (\sqrt{[g_l](b)}+\sqrt{[g_k](b)}+\sqrt{\epsilon}
\big )+\frac{\epsilon}{2}\sum_{l,k=1}^d\sqrt{[g_k](b)[g_l](b)}.
$$
Since $\epsilon >0$ is arbitrary (\ref{9mdchrule}) exists and 
(\ref{1mdchrule}) holds with $z=a$ and $y=b$.
The proof of statement $(a)$ is complete.

The proof of statement $(b)$ is analogous to the proof of statement $(a)$.
Instead of representation (\ref{10mdchrule}), in this case
we have a telescoping representation
\begin{eqnarray*}
\lefteqn{S_{RC}(\psi,\phi{\circ}\boldf;\lambda_m)
=\sum_{i\in I_1(m)}\psi(x_{i}^m)\Delta_i^m(\phi{\circ}\boldf)}\\[2mm]
& &+\sum_{i\in I_2(m)}\psi (x_{i-1}^m)\Big\{\sum_{k=1}^d(\phi_k'{\circ}
\boldf)(x_{i}^m)\Delta_i^mf_k-\frac{1}{2}\sum_{l,k=1}^d\phi_{lk}''(y_i^m)
\Delta_i^mf_l\Delta_i^mf_k\Big\}\\[2mm]
&=&\sum_{k=1}^dS_{RC}(\psi(\phi_k'{\circ}\boldf),g_k;\lambda_m)
+\sum_{k=1}^dS_{RC}(\psi(\phi_k'{\circ}\boldf),h_k;\lambda_m)\\[2mm]
& &-\frac{1}{2}\sum_{l,k=1}^d\left (\sum_{i=1}^{n(m)}-\sum_{i\in I_1(m)}
\right )\big (\psi (\phi_{lk}''{\circ}\boldf)\big )
(x_{i}^m)\Delta_i^mg_l\Delta_i^mg_k\\[2mm]
& &+\sum_{i\in I_1(m)}\psi (x_{i}^m)\Big\{\Delta_i^m(\phi{\circ}\boldf)
-\sum_{k=1}^d(\phi_k'{\circ}\boldf)(x_{i}^m)\Delta_i^mf_k\Big\}
+R_m
\end{eqnarray*}
for $m\geq m_0$,
where $y_i^m:=\boldf(x_{i}^m)-\theta_i^m\Delta_i^m\boldf$, 
$\theta_i^m\in (0,1)$ for $i=1,\dots,n(m)$, and
\begin{eqnarray*}
R_m&:=&-\frac{1}{2}\sum_{l,k=1}^d\sum_{i\in I_2(m)}\Big\{\psi (x_i^m)
\phi_{lk}''(y_i^m)\big [\Delta_i^mg_l\Delta_i^mh_k+\Delta_i^mh_l\Delta_i^mg_k
+\Delta_i^mh_l\Delta_i^mh_k\big ]\\[2mm]
& &\qquad\qquad\qquad+\psi (x_i^m(\big [\phi_{lk}''(y_i^m)-\phi_{lk}''
(\boldf(x_{i}^m))\big ]\Delta_i^mg_k\Delta_i^mg_l\Big\}.
\end{eqnarray*}
This representation allows to prove the existence of the 
Right Cauchy $\lambda$-integral in (\ref{11mdchrule}) and
relation (\ref{11mdchrule}).
The proof of Theorem \ref{mdchrule} is complete.
\qed\end{proof}

We finish this section with a chain rule which is used in Chapter \ref{modelling}
to prove the Black-Scholes option pricing formula.
Let $U\subset\RR$ be an open set and let $\phi\colon\,[a,b]\times U
\mapsto\RR$ be a continuous function.
We say that $\phi$ is a $C^2$ class function on $[a,b)\times U$
if it is a $C^2$ class function on $(a,b)\times U$ and has continuous
one-sided first and second derivatives at $a$.
For such function $\phi$ and a function $f\colon\,[a,b]\mapsto U$, let
\beq\label{Fdt}
\Phi f(t):=\phi (t,f(t)),\quad t\in [a,b],\qquad
\Phi_t'f(t):=\left\{ \begin{array}{ll}
\frac{\partial\phi}{\partial t}(t,f(t))\ &\mbox{if $t\in [a,b)$, }\\ 
0 &\mbox{if $t=b$,} 
\end{array} \right. 
\eeq
\beq\label{Fdx}
\Phi_x'f(t):=\left\{ \begin{array}{ll}
\frac{\partial\phi}{\partial x}(t,f(t))\ &\mbox{if $t\in [a,b)$, }\\ 
0 &\mbox{if $t=b$,} 
\end{array} \right. \qquad
\Phi_{xx}''f(t):=\left\{ \begin{array}{ll}
\frac{\partial^2\phi}{\partial x^2}(t,f(t))\ &\mbox{if $t\in [a,b)$, }\\ 
0 &\mbox{if $t=b$.} 
\end{array} \right. 
\eeq
The functions $\Phi_t'f$, $\Phi_x'f$ and $\Phi_{xx}''f$ so defined on $[a,b]$,
may be discontinuous at $b$ and unbounded.

\begin{prop}\label{chrule-HK}
For an open set $U\subset\RR$ and $\lambda\in\Lambda [a,b]$,
let $f\colon\,[a,b]\mapsto U$ be a continuous function having the quadratic 
$\lambda$-variation on $[a,b]$, and let 
$\phi \colon\,[a,b]\times U\mapsto\RR$ be a continuous function such
that its restriction to $[a,b)\times U$ is a $C^2$ class function.
Suppose that there exist the Henstock-Kurzweil integrals
$(HK)\smallint_a^b\Phi_t'f$ and $(HK)\smallint_a^b\Phi_{xx}''f\,d[f]_{\lambda}$.
Then the Left Cauchy $\lambda$-integral $(LC)\smallint \Phi_x'f\,d_{\lambda}f$
is defined on $[a,b)$ and there exists the limit
\beq\label{1chrule-HK}
\lim_{u\uparrow b}(LC)\int_a^u\Phi_x'f\,d_{\lambda}f
=\Phi f(b)-\Phi f(a)-(HK)\int_a^b\Phi_t'f-\frac{1}{2}(HK)\int_a^b
\Phi_{xx}''f\,d[f]_{\lambda}.
\eeq
\end{prop}

\begin{proof}
By the proof of Theorem \ref{mdchrule} with $d=2$, $\psi\equiv 1$, 
$g_1\equiv 0$, $h_1(x)\equiv x$, $g_2\equiv f$ and $h_2\equiv 0$, it follows that
the Left Cauchy $\lambda$-integral $(LC)\smallint \Phi_x'f\,d_{\lambda}f$
is defined on $[a,b)$, and for each $a\leq u<b$,
\beq\label{2chrule-HK}
\Phi f(u)-\Phi f(a)=(R)\int_a^u\Phi_t'f +(LC)\int_a^u\Phi_x'f\,d_{\lambda}f
+\frac{1}{2}(RS)\int_a^u\Phi_{xx}''f\,d[f]_{\lambda}.
\eeq
By the Hake theorem (see e.g.\  Section 7.3 in \cite{RMM}),
for a continuous function $h$ with bounded variation on $[a,b]$,
a function $g$ is Henstock-Kurzweil integrable over $[a,b]$ with
respect to $h$ if and only if for each $a\leq u<b$,
$g$ is Henstock-Kurzweil integrable over $[a,u]$ with respect to $h$
and there exists the limit $\lim_{u\uparrow b}(HK)\smallint_a^ug\,dh$.
In this case, 
$$
\lim_{u\uparrow b}(HK)\int_a^ug\,dh=(HK)\int_a^bg\,dh.
$$
Since the Henstock-Kurzweil integral extends the Riemann-Stieltjes integral,
the second part of the conclusion follows.
The proof is complete.
\qed\end{proof}

At this writting we are not able to prove under the same conditions
that the Left Cauchy $\lambda$-integral in Proposition \ref{chrule-HK}
is defined on the whole interval $[a,b]$.
The method explored in this paper requires a convergence theorem
for the Henstock-Kurzweil integral similar to the convergence theorem
of Lemma  \ref{conv-thm}, which cannot hold by Example 10 in 
\cite[p.\ 213]{RMM}.

\section{The product $\lambda$-integral}

Let $f$ be a real-valued function on $[a,b]$.
Recall that the product integral on $[a,b]$ with respect to $f$ 
is defined to be the limit
$$
\prodi_a^b (1+df):=\lim_{\kappa,\prtn }P(f;[a,b],\kappa)
$$
provided it exists, where $P(f;[a,b],\kappa)$ is defined by
(\ref{PandS}).
As it is mentioned earlier,
by Theorem 4.4 of \cite[Part II]{DNa}, if the product integral with
respect to $f$ exists then $f\in{\cal W}_2^{\ast}$.
In this section we extend the definition of the product
integral, so that the extended product integral
exists with respect to functions having the quadratic $\lambda$-variation 
(Definition \ref{lambda-prod-int} below).
But first we give a new proof of  Theorem 4.4 of \cite[Part II]{DNa}
for regulated functions which is based on a new characterization of 
the Wiener class ${\cal W}_2^{\ast}$  (Theorem 
\ref{variation}.\ref{liminf} above).

\paragraph*{The product integral.}
The product integral is more often used for non-commutative vector valued 
functions when no explicit formulas can be given.
The product integral with respect to a \emph{real-valued}
function, or more generaly with respect to a function with values
in a commutative algebra, if exists, can be computed explicitly. 
Let $f$ be a regulated function on $[a,b]$, and let $(u,v)\in S\loc a,b\roc$
be such that $u<v$.
For each $\delta >0$, let $D(u,v;\delta):=D_{-}(u,v;\delta)\cup 
D_{+}(u,v;\delta)$, where
$$
\left\{\begin{array}{ll}
D_{-}(u,v;\delta):=\{x\in (a,b]\colon\,u<x\leq v,\,\,
|\Delta^{-}f(x)|>\delta\}\\
D_{+}(u,v;\delta):=\{x\in [a,b)\colon\,u\leq x<v,\,\,
|\Delta^{+}f(x)|>\delta\}.
\end{array}
\right.
$$
Here we use the order on the extended interval defined in Section \ref{add&mult}.
Also, let
$$
\prod_{D(u,v;\delta)}(1+\Delta f)e^{-\Delta f}
:=\prod_{D_{-}(u,v;\delta)}(1+\Delta^{-} f)e^{-\Delta^{-} f}
\prod_{D_{+}(u,v;\delta)}(1+\Delta^{+} f)e^{-\Delta^{+} f}.
$$
Define the function $\gamma (f)$ on $S\loc a,b\roc$ by
\beq\label{1bvofV}
\gamma (f;u,v):=
\prod_{(u,v)}(1+\Delta f)e^{-\Delta f}:=\lim_{\delta\downarrow 0}
\prod_{D(u,v;\delta)}(1+\Delta f)e^{-\Delta f}
\eeq
provided $(u,v)\in S\loc a,b\roc$, $u<v$, and the infinite product
converges absolutely.
If $u=v$ then let $\gamma (f;u,u):=1$.
Now the necessary and sufficient conditions for the existence of
the product integral can be formulated as follows:

\begin{thm}\label{prodint}
The following two statements about a regulated function $f$ on
$[a,b]$ are equivalent{\rm :}
\begin{enumerate}
\item[$(a)$] the product integral $\prodi_a^b (1+df)$ exists and
is non-zero{\rm ;}
\item[$(b)$] $f\in {\cal W}_2^{\ast}[a,b]$ and
\beq\label{nondg}
\Delta^{-}f(x)\not =-1\not =\Delta^{+}f(y)\qquad
\mbox{ for each $a\leq y<x\leq b$.}
\eeq
\end{enumerate}
Moreover, if any one of the two statements holds then
the product integral extends to the function 
\beq\label{1prodint}
S\loc a,b\roc\ni (u,v)\mapsto
\prodi_u^v(1+df)=\exp\{f(v)-f(u)\}\prod_{(u,v)}(1+\Delta f)
e^{-\Delta f},
\eeq
which is nondegenerate, bounded, multiplicative and upper continuous 
function on $S\loc a,b\roc$.
\end{thm}

The function $\gamma (f)$ defined by (\ref{1bvofV}) is a discontinuous part 
of values (\ref{1prodint}) and it exists under the weaker condition 
$\sigma_2(f)<\infty$ (see (\ref{wiener0}) for its definition)
than above condition $(b)$. 

\begin{lem}\label{gamma}
Let $f$ be a regulated function on $[a,b]$ such that $\sigma_2(f)<\infty$.
The function $\gamma (f)$ on $S\loc a,b\roc$ is well defined, 
multiplicative, bounded, and for any $(u,v)\in S\loc a,b\roc$ such that 
$f$ has no jumps on $\loc u,v\roc$ exceeding $1/2$,
\beq\label{1gamma}
\big |\gamma (f;u,v)-1\big |\leq 2e^{2\sigma_2 (f)}\sum_{\loc u,v\roc}
\big [\Delta f\big ]^2.
\eeq
Moreover, $\gamma (f)$ is nondegenerate if and only if
{\rm (\ref{nondg})} holds.
If $\gamma (f)$ is nondegenerate then it is upper continuous.
\end{lem}

\begin{proof}
First we show that $\gamma (f)$ is well defined.
By a Taylor series expansion with remainder, we have that
for $|u|\leq 1/2$,
\beq\label{7bvofV}
\xi (u):=(1+u)e^{-u}=1-\frac{\theta (u)}{2}u^2\,,
\eeq
where $1/(2\sqrt{e})\leq\theta (u)\leq 3\sqrt{e}/2$.
Thus 
$$
\sum_{[a,b)}|1-\xi (\Delta^{-} f)|1_{\{|\Delta^{-}f|\leq 1/2\}}
+\sum_{(a,b]}|1-\xi (\Delta^{+} f)|1_{\{|\Delta^{+}f|\leq 1/2\}}
\leq \frac {3\sqrt{e}}{4}\sigma_2(f)^2<\infty\,.
$$
Since $f$ is regulated there exists a partition $\{z_k\colon\,
k=0,\dots,l\}$ of $[a,b]$ such that 
\beq\label{4bvofV}
\osc (f;(z_{k-1},z_k))\leq 1/2\qquad\mbox{ for $k=1,\dots,l$.}
\eeq
Let $(u,v)\in S\loc a,b\roc$.
By the preceding bounds, each product $\gamma (f;z_{k-1}+,z_k-)$, 
$k=1,\dots,l$, converges absolutely, and hence so does the 
product $\gamma (f;u,v)$.
Therefore the function $\gamma (f)$ is well defined.

For each $\delta >0$, the map
$$
S\loc a,b\roc\ni (u,v)\mapsto\prod_{D(u,v;\delta)}(1+\Delta f)
e^{-\Delta f}
$$
is multiplicative.
Thus the function $\gamma (f)$ is multiplicative as a limit
of multiplicative functions.
Next we show that $\gamma (f)$ is bounded.
A Taylor series expansion with remainder gives for $|u|\leq 1/2$,
\beq\label{5bvofV}
\log (1+u)=u-\frac{\theta (u)}{2}u^2\,,
\eeq
where $4/9\leq\theta (u)\leq 4$.
Let $(u,v)\in S\loc a,b\roc$ be such that jumps of
$f$ on $\loc u,v\roc$ do not exceed $1/2$, and let $0<\delta \leq 1/2$.
By (\ref{5bvofV}), we have
\begin{eqnarray}\label{6bvofV}
\prod_{D(u,v;\delta)}(1+\Delta f)e^{-\Delta f}
&=&\exp\Big\{-\frac{1}{2}\sum_{D(u,v;\delta)}\big [\theta (\Delta^{-}f)
(\Delta^{-}f)^2+\theta (\Delta^{+}f)(\Delta^{+}f)^2\big ]\Big\}
\nonumber \\[2mm]
&\leq & \exp\Big\{-\frac{2}{9}\;\sum_{D(u,v;\delta)}
\big [(\Delta^{-}f)^2+(\Delta^{+}f)^2\big ]\Big\}.
\end{eqnarray}
Letting $\delta\downarrow 0$ on both sides, we obtain the upper bound
$$
\gamma (f;u,v)\leq \exp\Big\{-\frac{2}{9}
\;\sum_{\loc u,v\roc}\big [(\Delta^{-}f)^2+(\Delta^{+}f)^2\big ]\Big\}\leq 1.
$$
For any $(u,v)\in S\loc a,b\roc$, applying these bounds to 
the intersections $\loc u,v\roc\cap \loc z_{k-1}+,z_k-\roc$, where 
$\{z_k\colon\,k=0,\dots,l\}$ is defined by (\ref{4bvofV}), and using 
notation (\ref{7bvofV}) we get
$$
\big\|\gamma (f)\big\|_{\infty}\leq\prod_{k=1}^l\Big [1\vee\big |
\xi (\Delta^{+}f(z_{k-1}))\big |\Big ]\Big [1\vee\big |\xi (
\Delta^{-}f(z_k))\big |\Big ]<+\infty,
$$
proving that $\gamma (f)$ is bounded.
Let $(u,v)\in S\loc a,b\roc$ be such that jumps of
$f$ on $\loc u,v\roc$ do not exceed $1/2$, and let $0<\delta \leq 1/2$.
An application of the inequality $|e^x-1|\leq |x|e^{|x|}$
for $x\in\RR$ and (\ref{5bvofV}) imply that
\begin{eqnarray*}
\Big |\prod_{D(u,v;\delta)}(1+\Delta f)e^{-\Delta f}-1\Big |
&=&\Big |\exp\Big\{\sum_{D(u,v;\delta)}\big [\log (1+\Delta f)-\Delta f
\big ]\Big\}-1\Big |\\[2mm]
&\leq &2e^{2\sigma_2(f)^2}\;\sum_{\loc u,v\roc}\big [(\Delta^{-}f)^2
+(\Delta^{+}f)^2\big ].
\end{eqnarray*}
Letting $\delta\downarrow 0$ on the left side, we obtain that
(\ref{1gamma}) holds.

Second we show that $\gamma (f)$ is nondegenerate on $[a,b]$ if and only
if (\ref{nondg}) holds.
Clearly if $\gamma (f)$ is nondegenerate on $[a,b]$ then (\ref{nondg}) 
holds.
To show the converse let $(x,y)\subset [a,b]$ be an open interval such 
that $f$ have no jumps on $(x,y)$ bigger than $1/2$, 
and let $0<\delta \leq 1/2$.
In the left side of (\ref{6bvofV}) using the upper bound for $\theta$ 
(see (\ref{5bvofV})), and then letting $\delta\downarrow 0$ on both sides 
of the resulting inequality, we get the lower bound
$$
\gamma (f;x+,y-)\geq \exp\Big\{-2
\;\sum_{(x,y)}\big [(\Delta^{-}f)^2+(\Delta^{+}f)^2\big ]\Big\}
$$
for any open interval $(x,y)$ which does not contain discontinuities
of $f$ with jumps bigger than $1/2$.
Applying this bound to each interval $(z_{k-1},z_k)$ satisfying
(\ref{4bvofV}), and using notation (\ref{7bvofV}) we get
$$
\big |\gamma (f;a,b)\big |\geq\Big |\prod_{k=1}^l\xi (\Delta^{+}
f(z_{k-1}))\xi (\Delta^{-}f(z_k))\Big |\exp\Big\{-2
\sum_{k=1}^l\sigma_2(f;(z_{k-1},z_k))^2\Big \}.
$$
Since $\sigma_2(f)<\infty$ and (\ref{nondg}) holds,
$\gamma (f;a,b)\not =0$.

Finally we show that $\gamma (f)$ is upper
continuous provided it is nondegenerate.
To this aim we show that statement $(vii)$ of Theorem 
\ref{variation}.\ref{interv} holds.
For $a\leq y<x\leq b$ with $|\Delta^{-}f(x)|\leq 1/2$ and
$|\Delta^{+}f(y)|\leq 1/2$, by (\ref{7bvofV}),
$$
\big |\gamma (f;x-,x)-1\big |\leq \frac{3\sqrt{e}}{8}
(\Delta^{-}f(x))^2
\quad\mbox{and}\quad
\big |\gamma (f;y,y+)-1\big |\leq \frac{3\sqrt{e}}{8}
(\Delta^{+}f(y))^2.
$$
Thus (\ref{singletons}) holds for $\mu =\gamma (f)$ by the analogous
property for the regulated function $f$.
Suppose that extended open intervals 
$\loc x_k+,y_k-\roc\downarrow\loc u,u\roc$ as
$k\to\infty$, for some $u\in \loc a,b\roc$.
Then either $y_k=y$ and $x_k\uparrow y$ for some $y\in (a,b]$ and all 
sufficiently large $k$, or $x_k=x$ and $y_k\downarrow x$ for some 
$x\in [a,b)$ and all sufficiently large $k$.
In either case for all sufficiently large $k$, $\loc x_k+,y_k-\roc$ does not
contain discontinuity points with jumps bigger than $1/2$.
Since $\sigma_2(f)<\infty$, by (\ref{1gamma}), it then follows that
$$
\limsup_{k\to\infty}\left |\gamma (f;x_k+,y_k-)-1\right |
\leq 2e^{2\sigma_2(f)^2}\limsup_{k\to\infty}
\sum_{(x_k,y_k)}\big [(\Delta^{-}f)^2+(\Delta^{+}f)^2\big ]=0.
$$
Thus statement $(vii)$ of Theorem \ref{variation}.\ref{interv} holds for 
$\mu=\gamma (f)$, and hence $\gamma (f)$ is upper continuous.
The proof of Lemma \ref{gamma} is now complete.
\qed\end{proof}

\medskip
\noindent
{\bf Proof of Theorem \ref{prodint}.}
We start with a useful representation of the product $P(f;\kappa)$
defined by (\ref{PandS}).
Let $\zeta=\{z_j\colon\,j=0,\dots,k\}$ be a partition of $[a,b]$
to be chosen later,
and let $\kappa=\{x_i\colon\,i=0,\dots,n\}$ be a refinement of $\zeta$
such that each intersection $(z_{j-1},z_j)\cap\kappa$, $j=1,\dots,k$,
contains at least two different points.
For $j=0,\dots,k$, let $i(j)\in\{0,\dots,n\}$ be such that 
$x_{i(j)}=z_j$.
Then let 
$$
I_1(\zeta,\kappa):=\{i(j-1)+1,i(j)\colon\,j=1,\dots,k\}\subset\{1,\dots,n\}
\quad\mbox{and}\quad
I_2(\zeta,\kappa):=\{1,\dots,n\}\setminus I_1(\zeta,\kappa).
$$
For a function $f$ on $[a,b]$, let $\Delta_if:=f(x_i)-f(x_{i-1})$,
$i=1,\dots,n$.
If $|\Delta_if|\leq 1/2$ for each $i\in I_2(\zeta,\kappa)$, then
by (\ref{5bvofV}),
\begin{eqnarray}
P(f;\kappa)&=&
\exp\Big\{\sum_{i\in I_2(\zeta,\kappa)}\Delta_if
-\frac{1}{2}\sum_{i\in I_2(\zeta,\kappa)}\theta_i(\Delta_if)^2\Big\}
\prod_{i\in I_1(\zeta,\kappa)}\big (1+\Delta_if\big )\nonumber\\[2mm]
&=&\exp\Big\{f(b)-f(a)-\frac{1}{2}\sum_{i\in I_2(\zeta,\kappa)}
\theta_i(\Delta_if)^2\Big\}
\prod_{i\in I_1(\zeta,\kappa)}\Big\{(1+\Delta_if)e^{-\Delta_if}\Big\},
\label{2prodint}
\end{eqnarray}
where $\theta_i:=\theta (\Delta_if)$ for $i\in I_2(\zeta,\kappa)$
and $\theta (u)=(1+v)^{-2}$ for some $v$ between $u$ and $0$.

Suppose that statement $(a)$ holds.
First we show that (\ref{nondg}) holds.
Since $\prodi_a^b(1+df)\not =0$, there exist $\epsilon >0$ and 
a partition $\zeta=\{z_j\colon\,j=0,\dots,k\}$ of $[a,b]$ such that
$D(a,b;1/2)\subset\zeta$ and $|P(f;\kappa)|\geq\epsilon$ for each
refinement $\kappa$ of $\zeta$.
Suppose that $\Delta^{-}f(x)=-1$ for some $x\in (a,b]$.
Then $x=z_j\in\zeta$ for some $j\in \{1,\dots,k\}$.
Let $\kappa$ be a refinement of $\zeta$ such that $z_{j-1}<
x_{i(j)-1}<z_j (=x_{i(j)})$.
Letting $x_{i(j)-1}\uparrow z_j$ and keeping other $x_i\in\kappa$ unchanged,
we have that $(1+\Delta_{i(j)-1}f)(1+\Delta_{i(j)}f)\to 0$.
Thus $P(f;\kappa)\to 0$, a contradiction.
Similarly, the assumption that $\Delta^{+}f(y)=-1$ for some $y\in [a,b)$
leads to a contradiction.
Therefore (\ref{nondg}) must hold.

Next we show that $v_2(f)<\infty$.
This is so if $\sum_{j=1}^kv_2(f;(z_{j-1},z_j))<\infty$ for
some partition $\{z_j\colon\,j=0,\dots,k\}$ of $[a,b]$.
Again, since $f$ is regulated and
$\prodi_a^b(1+df)\not =0$, there exist $\epsilon >0$
and a partition $\zeta=\{z_j\colon\,j=0,\dots,k\}$ of $[a,b]$
such that $D(a,b;1/2)\subset\zeta$, $\osc\, (f;(z_{j-1},z_j))\leq 
1/2$ for $j=1,\dots,k$, and $|P(f;\kappa)|\geq\epsilon$
for each refinement $\kappa$ of $\zeta$.
Since $|f|$ is bounded by a finite constant $M$, by (\ref{2prodint}), 
we have that 
$$
0<\epsilon\leq \exp\Big\{f(b)-f(a)-\frac{2}{9}\sum_{i\in I_2(\zeta,\kappa)}
(\Delta_if)^2\Big\} \Big [(1+2M)e^{2M}\Big ]^{2k}
$$
for all refinements $\kappa$ of $\zeta$.
If $\sum_{j=1}^kv_2(f;(z_{j-1},z_j))=+\infty$ then the right
side of the preceding bound can be made arbitrarily small.
The contradiction yields that $v_2(f)<\infty$.

To prove that $f\in {\cal W}_2^{\ast}[a,b]$, let $\epsilon\in (0,1/2]$.
Since $f$ is regulated and since the product integral with respect to $f$ is defined,
there exists a partition $\kappa_0=\{y_r\colon\,r=0,\dots,l\}$ of $[a,b]$ 
such that
$$
\max_{1\leq r\leq l}\osc\, (f;(y_{r-1},y_r))<\epsilon
\quad\mbox{and}\quad
\Big |\log\prodi_a^b(1+df)-\log P(f;\kappa)\Big |<\epsilon 
$$
for each refinement $\kappa$ of $\kappa_0$.
By Lemma \ref{variation}.\ref{wiener}, $\sigma_2(f)\leq v_2(f)<+\infty$.
Thus by Lemma \ref{gamma}, there exists a $\delta_0\in (0,1/2]$ such that
\beq\label{3prodint}
\Big |\log \gamma (f;a,b)-\log\Big\{ \prod_{D(a,b;\delta)}
(1+\Delta f)e^{-\Delta f}\Big\}\Big |<\epsilon 
\eeq
for each $0<\delta\leq\delta_0$.
Let $\delta$ be the minimum of the subset of non-zero values in the set
\beq\label{5prodint}
\{\delta_0,|\Delta^{+}f(a)|,|\Delta^{-}f(b)|,|\Delta^{-}f(x)|,
|\Delta^{+}f(x)|\colon\,x\in (a,b)\cap\kappa_0\},
\eeq
and let $\zeta=\{z_j\colon\,j=0,\dots,k\}:=\kappa_0\cup D(a,b;\delta)$.
Then choose $\{u_{j-1},v_j\colon\,j=1,\dots,k\}\subset (a,b)$ such that
for $j=1,\dots,k$, $z_{j-1}<u_{j-1}<v_j<z_j$,
$$
\osc\,\big ([f-f(z_{j-1})]^2;(z_{j-1},u_{j-1}]\big )<\epsilon /2k,
\qquad
\osc\,\big ([f(z_j)-f]^2;[v_j,z_j)\big )<\epsilon /2k,
$$
and for each refinement $\kappa$ of $\kappa_1:=\zeta\cup\{u_{j-1},v_j
\colon\,j=1,\dots,k\}$,
\beq\label{6prodint}
\Big |\log\Big\{\prod_{D(a,b;\delta)}\big (1+\Delta f\big )e^{-\Delta f}
\Big\}-\log\Big\{\prod_{i\in I_1(\zeta,\kappa)}\big (1+\Delta_if\big )
e^{-\Delta_if}\Big\}\Big |<\epsilon.
\eeq
Here we use the notation introduced at the beginning of the proof.
Letting
\beq\label{7prodint}
A:=f(b)-f(a)+\log\gamma (f;a,b)-\log\prodi_a^b(1+df)
\eeq
and using (\ref{2prodint}), we then have the bound
\begin{eqnarray*}
\lefteqn{\Big |A-\frac{1}{2}\sum_{i\in I_2(\zeta,\kappa)}
\theta_i(\Delta_if)^2\Big |}\nonumber\\[2mm]
&=&\Big |A-\big [f(b)-f(a)+\log\Big\{\prod_{i\in I_1(\zeta,\kappa)}
(1+\Delta_if)e^{-\Delta_if}\Big\}-\log P(f;\kappa)\big ]\Big |\\[2mm]
&\leq&\Big |\log\gamma (f;[a,b])-\log\Big\{\prod_{i\in I_1(\zeta,\kappa)}
\big (1+\Delta_if\big )e^{-\Delta_if}\Big\}\Big |+\Big |\log
\prodi_a^b(1+df)-\log P(f;\kappa)\Big |<3\epsilon
\end{eqnarray*}
for each refinement $\kappa$ of $\kappa_1$.
Since $|\Delta_if|<\epsilon$ for $i\in I_2(\zeta,\kappa)$,
it then follows that
$$
2(1-\epsilon )^2(A-3\epsilon) 
<\sum_{i\in I_2(\zeta,\kappa)}(\Delta_if)^2
<2(1+\epsilon )^2(A+3\epsilon)
$$
for each refinement $\kappa$ of $\kappa_1$.
Thus for any two refinements $\kappa'=\{x_i'\}$ and $\kappa''=
\{x_i''\}$ of $\kappa_1$,
$$
\big |s_2(f;\kappa')-s_2(f;\kappa'')\big |
\leq\sum_{j=1}^k\Big |[f(x_{i(j-1)+1}')-f(z_{j-1})]^2
-[f(x_{i(j-1)+1}'')-f(z_{j-1})]^2\Big |
$$
$$
+\sum_{j=1}^k\Big |[f(z_j)-f(x_{i(j)-1}')]^2
-[f(z_j)-f(x_{i(j)-1}'')]^2\Big |
+\Big |\sum_{i\in I_2(\zeta,\kappa')}(\Delta_i'f)^2
-\sum_{i\in I_2(\zeta,\kappa'')}(\Delta_i''f)^2\Big |
$$
$$
<\epsilon +8\epsilon A+12\epsilon (1+\epsilon^2)
\leq\epsilon +8\epsilon (3+A),
$$
where $\Delta_i'f:=f(x_i')-f(x_{i-1}')$ and 
$\Delta_i''f:=f(x_i'')-f(x_{i-1}'')$.
Hence
$$
0\leq v_2^{\ast}(f)^2-\sigma_2^{\ast}(f)^2
\leq\sup_{\kappa\supset\kappa_1}s_2(f;\kappa)
-\inf_{\kappa\supset\kappa_1}s_2(f;\kappa)
\leq\epsilon +8\epsilon (3+A).
$$
Since $\epsilon$ is arbitrary, $f\in {\cal W}_2^{\ast}[a,b]$.
Thus $(b)$ follows from $(a)$.

Now suppose that statement $(b)$ holds.
The following proof of existence of the product integral is a
variation of the preceding argument showing that $A$ defined by 
(\ref{7prodint}) actualy is zero.
Let $\epsilon\in (0,1/4]$.
By Lemma \ref{variation}.\ref{lemma1}, there exists a partition 
$\kappa_0=\{y_r\colon\,r=0,\dots,l\}$ of $[a,b]$ such that
$$
\sum_{r=1}^lv_2\big (f;(y_{r-1},y_r)\big )<\epsilon.
$$
In particular, we have that $\osc\,(f;(y_{r-1},y_r))<\sqrt{\epsilon}
\leq 1/2$ for each $r=1,\dots,l$.
By Lemma \ref{gamma}, $\gamma (f;a,b)\not =0$.
Using Lemma \ref{gamma} again we find a $\delta_0 >0$ such that
(\ref{3prodint}) holds, and then choose a $\delta\in (0,\delta_0]$
to be the minimum of the subset of non-zero values in the set
(\ref{5prodint}).
As before, let 
$\zeta :=\{z_j\colon\,j=0,\dots,k\}:=\kappa_0\cup D(a,b;\delta)$.
Then choose $\{u_{j-1},v_j\colon\,j=1,\dots,k\}\subset (a,b)$
such that $z_{j-1}<u_{j-1}<v_j<z_j$ for $j=1,\dots,k$, and
(\ref{6prodint}) holds for each refinement $\kappa$ of
$\kappa_1:=\zeta\cup\{u_{j-1},v_j\colon\,j=1,\dots,k\}$.
By (\ref{2prodint}), it then follows that
$$
\Big |\log P(f;\kappa)-\big [f(b)-f(a)+\log\gamma (f;a,b)\big ]\Big |
$$
$$
\leq\Big |\log\Big\{\prod_{i\in I_1(\zeta,\kappa)}(1+\Delta_if)
e^{-\Delta_if}\Big\}-\log\gamma (f;a,b)\Big |
+\frac{2}{9}\sum_{i\in I_2(\zeta,\kappa)}(\Delta_if)^2
<2\epsilon +\frac{2}{9}\epsilon
$$
for each refinement of $\kappa_1$.
Thus $\prodi_a^b(1+df)$ exists and equals to the right side of
(\ref{1prodint}) with $[u,v]=[a,b]$.
Since $\gamma (f;a,b)\not =0$, $(a)$ follows from $(b)$.

Since the interval function $\mu (f)$ on $S\loc a,b\roc$ defined by
(\ref{muf}) is additive and upper continuous, by Lemma \ref{gamma},
the map
$$
S\loc a,b\roc\ni (u,v)\mapsto \exp\big\{\mu (f;u,v)\big\}
\gamma (f;u,v)=:\pi (f;u,v)
$$
is nondegenerate, bounded, multiplicative and upper continuous 
function on $S\loc a,b\roc$.
The second part of the theorem will follow once we show that 
for each $(u,v)\in S\loc a,b\roc$, the product
integral $\prodi_u^v(1+df)$ exists and equals $\pi (f;u,v)$.
For the intervals $\loc u,v\roc =\loc x-,x\roc$ and $\loc u,v\roc=\loc y,y+\roc$, 
we have
$$
\prodi_{x-}^x(1+df)=1+\Delta^{-}f(x)=\pi (f;x-,x)
\quad\mbox{and}\quad
\prodi_y^{y+}(1+df)=1+\Delta^{+}f(y)=\pi (f;y,y+).
$$
Due to multiplicativity of the product integral it is enough
to prove the existence of $\prodi_u^v(1+df)$ when $\loc u,v\roc =\loc x+,y-\roc$ 
for some $a\leq x<y\leq b$.
The proof in this case is the same as for $[u,v]=[a,b]$
except that a partition $\zeta=\{z_j\colon\,j=0,\dots,k\}$ is such that
$z_0>x$ if $u=x+$ and $z_k<y$ if $v=y-$.
Also, $\kappa_1$ is the same as before except that in the cases
$x<z_0$ and $z_k<y$ we include additional points $u_{-1}\in (x,z_0)$
and $v_{k+1}\in (z_k,y)$, respectively.
By the proof it also follows that 
the product integral $\prodi_u^v(1+df)$ has the value $\pi (f;u,v)$.
Thus $\prodi (1+df)$ is the interval function on $S\loc a,b\roc$ which agree 
with the interval function $\pi (f)$. 
The proof of Theorem \ref{prodint} is now complete.
\qed

\paragraph*{The product $\lambda$-integral.}
By Theorem \ref{prodint}, if the product integral with respect to a regulated
function $f$ exists then $f\in {\cal W}_2^{\ast}$.
Thus Example \ref{NK} shows that the product integral may not exist with
respect to a function having the quadratic $\lambda$-variation.
Next we extend the product integral so as to accommodate a product
integration with such a situation.
Recall that $P(f;\kappa)$ is defined by (\ref{PandS})
and the trace partition $\kappa\Cap\lei u,v\rei$ is defined in
Notation  \ref{lmuv}.

\begin{defn}\label{lambda-prod-int}
{\rm Let $f$ be a regulated function on $[a,b]$,
and let $\lambda=\{\lambda_m\colon\,m\geq 1\}\in\Lambda [a,b]$.
We say that the {\em product $\lambda$-integral} 
$\prodi (1+d_{\lambda}f)$ is defined on $[a,b]$ if there exists 
a nondegenerate, bounded, multiplicative and upper continuous  
function $\pi_{\lambda}(f)$ on $S\lei a,b\rei$ such that for each 
$(u,v)\in S\lei a,b\rei$,
\beq\label{1lambda-prod-int}
\pi_{\lambda} (f;u,v)=\lim_{m\to\infty}P(f;\lambda_m\Cap \lei u,v\rei).
\eeq
The right distribution function of $\pi_{\lambda}(f)$ restricted to $[a,b]$,
that is, ${\cal P}_{\lambda}f(x):=\pi (f;a,x)$ for $a\leq x\leq b$, is
called the \emph{indefinite product $\lambda$-intergal}.}
\end{defn}

Alternatively the product $\lambda$-integral can be defined
solely in terms of a regulated function on $[a,b]$ as follows:

\begin{prop}\label{prod-int2}
Let $f$ be a regulated function on $[a,b]$, and let
$\lambda=\{\lambda_m\colon\,m\geq 1\}\in\Lambda [a,b]$.
The product $\lambda$-integral $\prodi (1+d_{\lambda}f)$ is defined 
on $[a,b]$ if and only if there exists a regulated function $G$ on $[a,b]$ such 
that $G(a)=1$, $|G|\gg 0$ and for $a\leq s<t\leq b$,
\beq\label{1prod-int2}
G(t)/G(s)=\lim_{m\to\infty}P(f;\lambda_m\Cap [s,t]),
\eeq
\beq\label{2prod-int2} 
G(t)/G(t-)=1+(\Delta^{-}f)(t)\qquad\mbox{and}\qquad
G(s+)/G(s)=1+(\Delta^{+}f)(s). 
\eeq
If the two equivalent assertions hold then $\pi_{\lambda}(f;u,v)
=G(v)/G(u)$ for each $(u,v)\in S\lei a,b\rei $.
\end{prop}

\begin{proof}
First suppose that there exists a regulated function $G$ on $[a,b]$
equal $1$ at $a$, $|G|\gg 0$, and 
which satisfy both (\ref{1prod-int2}) and (\ref{2prod-int2}).
Define a function $\pi_{\lambda}(f)$ on $S\lei a,b\rei$
by $\pi_{\lambda}(f;u,v):=G(v)/G(u)$ for $(u,v)\in S\lei a,b\rei$.
Then $\pi_{\lambda}(f)$ is nondegenerate, bounded, multiplicative and 
upper continuous function on $\lei a,b\rei$ by Theorem 
\ref{variation}.\ref{interv}.
We prove (\ref{1lambda-prod-int}) only for the cases $( u,v) =( a,t-)$ and 
$(u,v) = (a,t+)$ because the proofs for the other cases are similar.
To this aim let 
$\lambda_m=\{x_i^m\colon\,i=0,\dots,n(m)\}$, $m=1,2,\dots$,
and first let $t\in (a,b]$.
For each $m\geq 1$, there is an index $i(t)=i_m(t)\in\{1,\dots,n(m)-1\}$
such that $x_{i(t)}^m<t\leq x_{i(t)+1}^m$.
Then by  (\ref{1prod-int2}) and (\ref{2prod-int2}), we have
\begin{eqnarray*}
\lim_{m\to\infty}P(f;\lambda_m\Cap\lei a,t-\rei)
&=&\lim_{m\to\infty}\Big\{P(f;\lambda_m\Cap [a,t])
\frac{[1+f(t-)-f(x_{i(t)}^m)]}{[1+f(t)-f(x_{i(t)}^m)]}\Big\}\\[2mm]
&=&\frac{G(t)}{1+(\Delta^{-}f)(t)}=G(t-)=\pi_{\lambda}(f;a,t-),
\end{eqnarray*}
proving $(\ref{1lambda-prod-int})$ for $(u,v) =(a,t-)$.
Now let $t\in [a,b)$.
By the definition of the trace partition, we have for each $m\geq 1$,
$$
P(f;\lambda_m\Cap \lei a,t+\rei )=P(f;\lambda_m\Cap [ a,t] )
(1+\Delta^{+}f(t)),
$$
proving $(\ref{1lambda-prod-int})$ for $(u,v) =(a,t+)$.

To prove the converse implication let $G(t):=\pi_{\lambda}(f;a,t)$
for $t\in [a,b]$.
By the statement $(iv)$ of Theorem \ref{variation}.\ref{interv} and 
(\ref{1lambda-prod-int}), $G$ is regulated on $[a,b]$, it is $1$ at $a$, 
$|G(t)|\geq C$ with $C:=\inf\{|\pi_{\lambda}(f;a,t)|\colon\,t\in [a,b]\}
>0$, and satisfies (\ref{1prod-int2}) and (\ref{2prod-int2}). 
The proof of Proposition \ref{prod-int2} is complete.
\qed\end{proof}

Next we show that the product $\lambda$-integral extends
the product integral for $\lambda\in\Lambda$ such that 
(\ref{jumps-access}) holds.

\begin{prop}\label{prod->l-prod}
Let $f$ be a regulated function on $[a,b]$, and let $\lambda\in\Lambda [a,b]$
be such that the two-sided discontinuity points of $f$ are accessible by 
$\lambda$.
If the product integral $\prodi (1+df)$ exists and is non-zero on $[a,b]$
then the product $\lambda$-integral $\prodi (1+d_{\lambda}f)$
is defined on $[a,b]$, and both have the same values on $S\lei a,b\rei $.
\end{prop}

\begin{proof}
By the second part of Theorem \ref{prodint} and by Proposition
\ref{prod-int2}, it is enough to prove that
for each $(s,t)\in S[a,b]$,
$$
\prodi_s^t(1+df)=\lim_{m\to\infty}P(f;\lambda_m\Cap [s,t]).
$$
Given $(s,t)\in S[a,b]$, let $\epsilon >0$.
Then there is a partition $\kappa_0$ of $[s,t]$
such that
$$
\Big |\prodi_s^t(1+df)-P(f;\kappa)\Big |<\epsilon
$$
for each refinement $\kappa$ of $\kappa_0$.
If $\kappa_0\subset\cup_m\lambda_m$ then $P(f;\lambda_m\Cap
[s,t])$ is within $\epsilon$ from $\prodi_s^t(1+df)$ for all
sufficiently large $m$.
Suppose that $\zeta=\{z_j\colon\,j=1,\dots,k\}:=
\kappa_0\setminus (\cup_m\lambda_m\cup\{s,t\})$ is non-empty.
For $j=1,\dots,k$ and all $m$ such that each $(z_{j-1},z_j)\cap
\lambda_m$ is non-empty, let $i(j)=i(j,m)\in\{1,\dots,n(m)\}$ be such that
$x_{i(j)-1}^m<z_j<x_{i(j)}^m$.
Also let $\Delta_{i(j)}^mf:=f(x_{i(j)}^m)-f(x_{i(j)-1}^m)$,
$\Delta_{i(j)}^{m,1}f:=f(z_j)-f(x_{i(j)-1}^m)$ and
$\Delta_{i(j)}^{m,2}f:=f(x_{i(j)}^m)-f(z_j)$.
By a telescoping sum
$$
\prod_{j=1}^ka_j-\prod_{j=1}^kb_j
=\sum_{j=1}^k\Big (\prod_{i=j+1}^{k}a_i\Big )(a_j-b_j)
\Big (\prod_{i=1}^{j-1}b_i\Big ),
$$
we have the relation
$$
P(f;(\lambda_m\Cap [s,t])\cup\zeta)-P(f;\lambda_m\Cap [s,t])
$$
$$
=P(f;(\lambda_m\Cap [s,t])\cup\zeta)\sum_{j=1}^k\Big\{
\frac{\Delta_{i(j)}^{m,1}f\Delta_{i(j)}^{m,2}f}{[1+\Delta_{i(j)}^{m,1}f]
[1+\Delta_{i(j)}^{m,2}f]}\prod_{l=1}^{j-1}\frac{1+\Delta_{i(l)}^mf}
{[1+\Delta_{i(l)}^{m,1}f][1+\Delta_{i(l)}^{m,2}f]}\Big\}
$$
for all sufficiently large $m$, where the products over the empty set
of indices are equal to $1$.
By Theorem \ref{prodint}$(b)$, $\Delta^{-}f(z_j)\not =-1\not 
=\Delta^{+}f(z_j)$. 
By (\ref{jumps-access}), $\Delta_{i(j)}^{m,1}f\Delta_{i(j)}^{m,2}f
\to 0$ as $m\to\infty$.
Thus each term of the sum on the right side tends to $0$ as $m\to\infty$.
Therefore since for all sufficiently large $m$,
$P(f;(\lambda_m\Cap [s,t])\cup\zeta)$ is within $\epsilon$ from
the value of the product integral over $[s,t]$,
$$
\limsup_{m\to\infty}\Big |\prodi_s^t(1+df)-P(f;
\lambda_m\Cap [s,t])\Big |\leq\epsilon,
$$
proving the proposition.
\qed\end{proof}

Suppose that a function $f$ on $[a,b]$ has the quadratic $\lambda$-variation 
for some $\lambda\in\Lambda [a,b]$.
Thus by Definition \ref{qv}, there is an additive and upper
continuous function $\alpha_{\lambda}(f)$ on $S\lei a,b\rei $
such that $\alpha_{\lambda}(f;u,v)=[f]_{\lambda}(v)-
[f]_{\lambda}(u)$ for $(u,v)\in S\lei a,b\rei $.
Define the additive upper continuous function $\alpha_{\lambda}^c(f)$ by
\beq\label{2alpha}
\alpha_{\lambda}^c(f;u,v):=[f]_{\lambda}^c(v)-[f]_{\lambda}^c
(u),\qquad (u,v)\in S\lei a,b\rei .
\eeq
Then by (\ref{cnt-part}),
for $(u,v)\in S[a,b]$, we have the decomposition
\beq\label{alpha-decomp}
\alpha_{\lambda}(f;u,v)=\alpha_{\lambda}^c(f;u,v)
+\sum_{(u,v]}\{\Delta^{-}f\}^2+\sum_{[u,v)}\{\Delta^{+}f\}^2.
\eeq
Define the function $\beta_{\lambda}(f)$ on $S\lei a,b\rei $ by
\beq\label{beta1}
\beta_{\lambda}(f;u,v):=\exp\big\{\mu (f;u,v)-\frac{1}{2}\alpha_{\lambda}^c
(f;u,v)\big\}\gamma (f;u,v),\qquad (u,v)\in S\lei a,b\rei ,
\eeq
where $\mu (f)$, $\alpha_{\lambda}^c(f)$ and $\gamma (f)$ are defined 
by (\ref{muf}), (\ref{2alpha}) and (\ref{1bvofV}), respectively.
By the definition of $\gamma (f)$, it follows that for $a\leq s<t\leq b$,
\beq\label{beta2}
\beta_{\lambda}(f;t-,t)=1+\Delta^{-}f(t)
\qquad\mbox{and}\qquad
\beta_{\lambda}(f;s,s+)=1+\Delta^{+}f(s).
\eeq
The following is a consequence of Lemma \ref{gamma}.

\begin{cor}\label{beta}
If $f$ has the quadratic $\lambda$-variation for some $\lambda\in\Lambda [a,b]$
and {\rm (\ref{nondg})} holds, then $\beta_{\lambda}(f)$ is
nondegenerate, bounded,
multiplicative and upper continuous function on $S\lei a,b\rei$.
\end{cor}

\begin{thm}\label{ExistLprod}
Let $f$ be a regulated function on $[a,b]$ having the quadratic 
$\lambda$-variation for some $\lambda\in\Lambda [a,b]$ such that 
the two-sided discontinuity points of $f$ are accessible by $\lambda$. 
If also {\rm (\ref{nondg})} holds then the product $\lambda$-integral
$\prodi (1+d_{\lambda}f)$ is defined on $[a,b]$,
and equals $\beta_{\lambda}(f)$.
\end{thm}

\begin{proof}
For $x\in [a,b]$, let $G(x):=\beta_{\lambda}(f;a,x)$.
By Theorem   \ref{variation}.\ref{interv} and Corollary \ref{beta}, 
$G$ is a regulated function on $[a,b]$ such that $G(a)=1$ and 
(\ref{2prod-int2}) holds for each $a\leq x<y\leq b$.
Thus by Proposition \ref{prod-int2}, it is enough to prove that
for each $a\leq x<y\leq b$,
\beq\label{1ExistLprod}
\lim_{m\to\infty}\log P(f;\lambda_m\Cap [x,y])
=\mu (f;x,y)-\frac{1}{2}\alpha_{\lambda}^c(f;x,y)
+\log\gamma (f;x,y).
\eeq
To begin the proof, let $a\leq x<y\leq b$ and let $\epsilon >0$.
Since $f$ is regulated there exists a partition $\kappa=\{t_r\colon\,
r=0,\dots,l\}$ of $[x,y]$ such that
\beq\label{2ExistLprod}
\max_{1\leq r\leq l}\osc (f;(t_{r-1},t_r))<\epsilon.
\eeq
Since $f$ has the quadratic $\lambda$-variation, by Lemma \ref{gamma},
there exists $\delta_0>0$ such that for each $0<\delta\leq\delta_0$
and $\sigma := D(x,y;\delta)$,
\beq\label{3ExistLprod}
\Big |\log \gamma(f;x,y)-\log\Big\{\prod_{\sigma}(1+\Delta f)
e^{-\Delta f}\Big\}\Big |<\epsilon 
\quad\mbox{and}\quad
\Big |\sum_{[x,y]}\big (\Delta f\big )^2-\sum_{\sigma}
\big (\Delta f\big )^2\Big |<\epsilon.
\eeq
Let $\delta$ be the minimal number among the positive values in the
set $\{\delta_0,|\Delta^{-}f(t_r)|,|\Delta^{+}f(t_{r-1})|\colon
r=1,\dots,l\}$, and let $\zeta =\{z_j\colon\,j=0,\dots,k\}:=
\kappa\cup D(x,y;\delta)$.
Thus if $z_j\in\zeta$ and $f$ has a jump at $z_j$ then $z_j\in 
D(x,y;\delta)$.
Let $\zeta_0:=\zeta\setminus\cup_m\lambda_m$, and let
$m_0$ be the minimal integer such that $\zeta\setminus\zeta_0
\subset\lambda_{m_0}$ and each $(z_{j-1},z_j)\cap\lambda_{m_0}$
contains at least two points.
By (\ref{jumps-access}), if $z_j\in\zeta_0$ then $f$ cannot have
jumps from the two sides at $z_j$.
Suppose that each $\lambda_m\Cap [x,y]=\{t_i^m\colon\,i=i_m(x)-1,\dots,
i_m(y)\}$, and let $\Delta_i^mf:=f(t_i^m)-f(t_{i-1}^m)$ for each 
$i\in I(m):=\{i_m(x),\dots,i_m(y)\}$ and $m\geq m_0$.
For each $m\geq m_0$, if $z_j\in\zeta_0$ then there is
$i(j)\in I(m)$ such that $t_{i(j)-1}^m<z_j<t_{i(j)}^m$, and
if $z_j\in\zeta\setminus\zeta_0$ then there is $i(j)\in I(m)$
such that $z_j=t_{i(j)}^m$.
For $m\geq m_0$, let $I_1(m)$ be the set of all indices $i(j)$
with $j$ such that $z_j\in\zeta_0$, and all pairs of indices
$i(j), i(j)+1$ with $j$ such that $z_j\in\zeta\setminus\zeta_0$,
and let $I_2(m):=I(m)\setminus I_1(m)$.
By (\ref{jumps-access}), if $z_j\in\zeta_0$ then $\lim_{m\to\infty}
\Delta_{i(j)}^mf$ exists and equals either to $\Delta^{-}f(z_j)$ or to
$\Delta^{+}f(z_j)$.
Otherwise if $z_j\in\zeta\setminus\zeta_0$ then $\lim_{m\to\infty}
\Delta_{i(j)}^mf=\Delta^{-}f(z_j)$ and $\lim_{m\to\infty}
\Delta_{i(j)+1}^mf=\Delta^{+}f(z_j)$.
Again since $f$ is regulated and has the quadratic $\lambda$-variation, 
there is an $m_1\geq m_0$ such that for each $m\geq m_1$,
\beq\label{4ExistLprod}
\big |s_2(f;\lambda_m\Cap [x,y])-\alpha_{\lambda}(f;x,y)\big |
<\epsilon,\quad
\Big |\sum_{i\in I_1(m)}\big (\Delta_i^mf\big )^2-\sum_{\sigma}
\big (\Delta f\big )^2\Big |<\epsilon
\eeq
and
\beq\label{5ExistLprod}
\Big |\log\Big\{\prod_{i\in I_1(m)}(1+\Delta_i^m f)e^{-\Delta_i^m f}
\Big\}-\log\Big\{\prod_{\sigma}(1+\Delta f)e^{-\Delta f}\Big\}\Big |
<\epsilon,
\eeq
where $\sigma= D(x,y;\delta)$.
By Taylor's theorem with remainder, for $|z|\leq 1/2$, 
\beq\label{Taylor3}
\log (1+z)=z-z^2/2+3\theta z^3,
\eeq
where $|\theta |=|\theta (z)|\leq 1$.
Thus for each $m\geq m_1$, we have
\begin{eqnarray*}
\lefteqn{\log P(f;\lambda_m\Cap [x,y])}\\[2mm]
&=&\log\prod_{i\in I_1(m)}(1+\Delta_i^mf)
+\sum_{i\in I_2(m)}\Delta_i^mf-\frac{1}{2}
\sum_{i\in I_2(m)}\big (\Delta_i^mf\big )^2
+3\sum_{i\in I_2(m)}\theta_i\big (\Delta_i^mf\big )^3\\[2mm]
&=&\log\Big\{\prod_{i\in I_1(m)}(1+\Delta_i^mf)e^{-\Delta_i^mf}
\Big\}+\mu (f;x,y)\\[2mm]
& &-\frac{1}{2}\Big\{s_2(f;\lambda_m\Cap [x,y])
-\sum_{i\in I_1(m)}\big (\Delta_i^mf\big )^2\Big\}
+3\sum_{i\in I_2(m)}\theta_i\big (\Delta_i^mf\big )^3,
\end{eqnarray*}
where $\theta_i:=\theta (\Delta_i^mf)$.
By  (\ref{alpha-decomp}), (\ref{4ExistLprod}) and by the second inequality in 
(\ref{3ExistLprod}), we have for each $m\geq m_1$,
$$
\Big | s_2(f;\lambda_m\Cap [x,y])-\sum_{i\in I_1(m)}
\big (\Delta_i^mf\big )^2-\alpha_{\lambda}^c (f;x,y)\Big |
$$
$$
\leq \big |s_2(f;\lambda_m\Cap [x,y])-\alpha_{\lambda}(f;x,y)\big |
+\Big |\sum_{[x,y]}\big (\Delta f\big )^2-\sum_{i\in I_1(m)}
\big (\Delta_i^mf\big )^2\Big |<3\epsilon.
$$
By (\ref{5ExistLprod}) and by the first inequality in 
(\ref{3ExistLprod}), we have for each $m\geq m_1$,
$$
\Big |\log\Big\{\prod_{i\in I_1(m)}(1+\Delta_i^m f)e^{-\Delta_i^m f}
\Big\}-\log \gamma (f;x,y)\Big |<2\epsilon.
$$
This in conjunction with (\ref{2ExistLprod}) yields that the bound
\begin{eqnarray*}
\lefteqn{\Big |\log P(f;\lambda_m\Cap [x,y])-\log\gamma (f;x,y)
-\mu (f;x,y)+\frac{1}{2}\alpha_{\lambda}^c(f;x,y)\Big |}\\[2mm]
&\leq&\Big |\log\Big\{\prod_{i\in I_1(m)}(1+\Delta_i^mf)
e^{-\Delta_i^mf}\Big\}-\log\gamma (f;x,y)\Big |\\[2mm]
& &+\frac{1}{2}\Big |s_2(f;\lambda_m\Cap [x,y])
-\sum_{i\in I_1(m)}\big (\Delta_i^mf\big )^2
-\alpha_{\lambda}^c(f;x,y)\Big |
+3\sum_{i\in I_2(m)}|\theta_i|\,\big |\Delta_i^mf\big |^3\\[2mm]
&\leq& 2\epsilon +(3/2)\epsilon +3\epsilon s_2(f;\lambda_m\Cap [x,y])
\end{eqnarray*}
holds for all $m\geq m_1$.
Since $s_2(f;\lambda_m\Cap [x,y])$ is bounded in $m$ and since 
$\epsilon$ is arbitrary, (\ref{1ExistLprod}) holds for any
$a\leq x<y\leq b$.
The proof of Theorem \ref{ExistLprod} is complete.
\qed\end{proof}


\section{Extended Dol\'eans exponentials}\label{doleansexp}

Let $f$ be a regulated function on $[a,b]$ having the quadratic 
$\lambda$-variation for some $\lambda\in\Lambda [a,b]$.
Define the forward \emph{Dol\'eans exponential} ${\cal E}_a(f)\equiv
{\cal E}_{\lambda,a}(f)$ on  $[a,b]$ by
\beq\label{Dol1}
{\cal E}_a(f;x):=\left\{ \begin{array}{ll}
\exp\big\{f(x)-f(a)-\frac{1}{2}[f]_{\lambda}^c(x)\big\}
\prod_{[a,x]}(1+\Delta f)e^{-\Delta f}&
                          \mbox{if $x\in (a,b]$,}\\
1 &\mbox{if $x=a$.}
\end{array}\right. 
\eeq
Here the product is defined by (\ref{1bvofV}) and it exists by Lemma \ref{gamma}.
Similarly, define the backward \emph{Dol\'eans exponential} 
${\cal E}_b(f)\equiv{\cal E}_{\lambda,b}(f)$ on $[a,b]$ by
$$
{\cal E}_b(f;x):=\left\{ \begin{array}{ll}
\exp\big\{f(b)-f(x)-\frac{1}{2}[f]_{\lambda}^c(b)+\frac{1}{2}
[f]_{\lambda}^c(x)\big\}
\prod_{[x,b]}(1+\Delta f)e^{-\Delta f}&
                          \mbox{if $x\in [a,b)$,}\\
1 &\mbox{if $x=b$.}
\end{array}\right. 
$$
The forward Dol\'eans exponential is a pathwise variant of the stochastic
Dol\'eans exponential which is the unique solution to the linear
stochastic differential equation (Dol\'eans-Dade \cite{DD}).
We show in the next section that the two Dol\'eans exponentials are unique
solutions to the augmented forward and backward linear $\lambda$-integral 
equations, respectively.
If a real valued function $f$ is in the class ${\cal W}_2^{\ast}[a,b]$,
so that $[f]_{\lambda}^c\equiv 0$, then the two Dol\'eans exponentials 
agree with the values of the product integral with respect to
$f$ over the intervals $[a,x]$ and $[x,b]$, respectively (see
Theorem \ref{prodint}).
In this section we show that the $2$-vector functions $({\cal E}_{\lambda,a}
(f),f)$ and $({\cal E}_{\lambda,b}(f),f)$ have the quadratic $\lambda$-variations.

It is easy to see that the Dol\'eans exponentials ${\cal E}_a(f;x)
=\beta_{\lambda}(f;a,x)$ and ${\cal E}_b(f;x)=\beta_{\lambda}(f;x,b)$ for 
$x\in [a,b]$, where $\beta_{\lambda}(f)$ is defined by (\ref{beta1}).
The two Dol\'eans exponentials are regulated functions on $[a,b]$.
Jumps of ${\cal E}_a(f)$ satisfy the relations:
\begin{eqnarray}
(\Delta^{-}{\cal E}_a(f))(x)
&=&\beta_{\lambda} (f;a,x)-\beta_{\lambda} (f;a,x-)
={\cal E}_a(f;x-)\big [\beta_{\lambda} (f;x-,x)-1\big ]\nonumber\\[2mm]
&=&{\cal E}_a(f;x-)\Delta^{-}f(x),\quad x\in (a,b],\label{-Ea}
\end{eqnarray}
and 
\begin{eqnarray}
(\Delta^{+}{\cal E}_a(f))(x)
&=&\beta_{\lambda} (f;a,x+)-\beta_{\lambda} (f;a,x)
={\cal E}_a(f;x)\big [\beta_{\lambda} (f;x,x+)-1\big ]\nonumber\\[2mm]
&=&{\cal E}_a(f;x)\Delta^{+}f(x),\quad x\in [a,b).\label{+Ea}
\end{eqnarray}
While jumps of  ${\cal E}_b(f)$ satisfy the relations:
\begin{eqnarray}
(\Delta^{-}{\cal E}_b(f))(x)
&=&\beta_{\lambda} (f;x,b)-\beta_{\lambda} (f;x-,b)
={\cal E}_b(f;x)\big [1-\beta_{\lambda} (f;x-,x)\big ]\nonumber\\[2mm]
&=&-{\cal E}_b(f;x)\Delta^{-}f(x),\quad x\in (a,b],\label{-Eb}
\end{eqnarray}
and 
\begin{eqnarray}
(\Delta^{+}{\cal E}_b(f))(x)
&=&\beta_{\lambda} (f;x+,b)-\beta_{\lambda} (f;x,b)
={\cal E}_b(f;x+)\big [1-\beta_{\lambda} (f;x,x+)\big ]\nonumber\\[2mm]
&=&-{\cal E}_b(f;x+)\Delta^{+}f(x),\quad x\in [a,b).\label{+Eb}
\end{eqnarray}
Next we show that the distribution functions 
$R_{\gamma (f)}$ and $L_{\gamma (f)}$ of the function
$\gamma (f)$ restricted to $[a,b]$ (see (\ref{RandLdf}))
 are pure jump functions of bounded variation.
Therefore let
$$
\mbox{$V_a(x):=\gamma (f;a,x)$ and $V_b(x):=\gamma (f;x,b)$
for $x\in [a,b]$}.
$$
By Theorem \ref{variation}.\ref{interv}, if $f$ is a regulated function such that
$\sigma_2(f)<\infty$ and (\ref{nondg}) holds, then
$V_a(x-)=\gamma (f;a,x-)$, $V_b(x-)=\gamma (f;x-,b)$ for
$x\in (a,b]$ and $V_a(x+)=\gamma (f;a,x+)$, $V_b(x+)=\gamma (f;
x+,b)$ for $x\in [a,b)$.
The multiplicativity of $\gamma (f)$ yields the following
relations: for $x\in (a,b]$,
\beq\label{2bvofV}
V_a(x)=\gamma (f;a,x)=\gamma (f;a,x-)\gamma (f;x-,x)
=V_a(x-)\big [1+\Delta^{-}f(x)\big ]e^{-\Delta^{-}f(x)},
\eeq
\beq\label{8bvofV}
V_b(x-)=\gamma (f;x-,b)=\gamma (f;x-,x)\gamma (f;x,b)
=V_b(x)\big [1+\Delta^{-}f(x)\big ]e^{-\Delta^{-}f(x)},
\eeq
and for $x\in [a,b)$,
\beq\label{3bvofV}
V_a(x+)=\gamma (f;a,x+)=\gamma (f;a,x)\gamma (f;x,x+)
=V_a(x)\big [1+\Delta^{+}f(x)\big ]e^{-\Delta^{+}f(x)},
\eeq
$$
V_b(x)=\gamma (f;x,b)=\gamma (f;x,x+)\gamma (f;x+,b)
=V_b(x+)\big [1+\Delta^{+}f(x)\big ]e^{-\Delta^{+}f(x)}.
$$

\begin{lem}\label{bvofV}
The functions $V_a$ and $V_b$ are pure jump functions on $[a,b]$
of bounded variation.
\end{lem}

\begin{proof}
We prove that $V:=V_a$ is a pure jump function of bounded variation.
The proof of the same property for $V_b$ is similar, and hence is
omitted.
If $\Delta^{-}f(x)=-1$ for some $x\in (a,b]$ then $V(y)=0$ for all
$y\geq x$, and if $\Delta^{+}f(x)=-1$ for some $x\in [a,b)$ then 
$V(y)=0$ for all $y > x$.
In either case the claim of the lemma holds for $V$ if it holds
for its restriction to the interval $[a,x)$ if $\Delta^{-}f(x)=-1$,
or for its restriction to $[a,x]$ if $\Delta^{+}f(x)=1$.
For notation simplicity we assume that (\ref{nondg}) holds.
By (\ref{2bvofV}) and (\ref{7bvofV}), for any $a<x\leq b$, we have
$$
\left |\Delta^{-}V(x)\right|
=\left |V(x-)\right |\left |(1+\Delta^{-}f(x))e^{-\Delta^{-}f(x)}
-1\right |\leq\frac{3\sqrt{e}}{2}\|\gamma (f)\|_{\infty}\left [\Delta^{-}f(x)
\right ]^2
$$
provided $|\Delta^{-}f(x)|\leq 1/2$.
The same bound holds for $|\Delta^{+}V(y)|$, $a\leq y<b$, with
$\Delta^{-}f$ replaced by $\Delta^{+}f$, provided
$|\Delta^{+}f(y)|\leq 1/2$.
Therefore $\sigma_1(V;(a,b))<\infty$,
and hence there exists the limit
$$
\sum_{(a,x)}\big\{\Delta^{-}V+\Delta^{+}V\big\}
=\lim_{\kappa,\prtn}\sum_{\kappa\cap (a,x)}
\big\{\Delta^{-}V+\Delta^{+}V\big\}
$$
for each $a<x\leq b$.
Let $a<x\leq b$.
For a partition $a=x_0<x_1<\dots <x_n=x$ of $[a,x]$, by a
telescoping sum, we have
$$
V(x)-1=\Delta^{+}V(a)+\sum_{i=1}^n\left\{V(x_i-)-V(x_{i-1}+)
\right\}+\sum_{i=1}^{n-1}\left\{\Delta^{-}V(x_i)+\Delta^{+}V(x_i)
\right\}+\Delta^{-}V(x).
$$
By multiplicativity of $\gamma (f)$ and by (\ref{1gamma}),
for each $i=1,\dots,n$, we get 
\begin{eqnarray*}
\big |V(x_i-)-V(x_{i-1}+)\big |&=&\big |V(x_{i-1}+)\big |
\Big |\prod_{(x_{i-1}+,x_i-)}(1+\Delta f)e^{-\Delta f}-1\Big |\\
&\leq& 2e^{2\sigma_2(f)}\|\gamma (f)\|_{\infty}\sigma_2
(f;(x_{i-1},x_i))
\end{eqnarray*}
provided $(x_{i-1},x_i)$ does not contain discontinuities of $f$
with jumps bigger than $1/2$.
Therefore given $\epsilon >0$, one can choose a partition of $[a,x]$ such that
for each its refinement $\{x_i\colon\,i=0,\dots,n\}$,
$$
\Big |\sum_{(a,x)}\big\{\Delta^{-}V+\Delta^{+}V\big\}
-\sum_{i=1}^n\big\{\Delta^{-}V(x_i)+\Delta^{+}V(x_i)\big\}\Big |<\epsilon
\quad\mbox{and}\quad
\Big |\sum_{i=1}^n\big\{V(x_i-)-V(x_{i-1}+)\big\}\Big |<\epsilon.
$$
Since $\epsilon >0$ is arbitrary, the relation
$$
V(x)=1+\Delta^{+}V(a)+\sum_{(a,x)}\big\{\Delta^{-}V+\Delta^{+}V
\big\}+\Delta^{-}V(x)
$$
holds for each $a<x\leq b$.
Also, for a partition $\kappa=\{x_i\colon\,i=0,\dots,n\}$ of $[a,b]$,
we have
$$
s_1(V;\kappa)\leq\sum_{i=1}^n\Big\{|\Delta^{+}V(x_{i-1})|
+\sum_{(x_{i-1},x_i)}\big\{|\Delta^{-}V|+|\Delta^{+}V|\big\} 
+|\Delta^{-}V(x_i)|\Big\}
\leq\sigma_1(V;[a,b])<\infty.
$$
Since $\kappa$ is arbitrary, 
the bound implies that $V=V_a$ is of bounded variation.
The proof of Lemma \ref{bvofV} is now complete.
\qed\end{proof}

\begin{thm}\label{qvofEa}
Let  $f$ be a regulated function on $[a,b]$ having the quadratic 
$\lambda$-variation for some $\lambda\in \Lambda[a,b]$ such 
that right discontinuity points of $f$ are accessible by $\lambda$.
Then the Dol\'eans exponential ${\cal E}_a(f)$ has the quadratic 
$\lambda$-variation on $[a,b]$, and its bracket function is given by
\beq\label{1qvofE}
\big [{\cal E}_a(f)\big ]_{\lambda}(x)=(LY)\int_a^x{\cal E}_a(f)^2\,
d[f]_{\lambda},\qquad a\leq x\leq b.
\eeq
Also, the $2$-vector function $({\cal E}_a(f),f)$ has the
quadratic $\lambda$-covariation on $[a,b]$, and its bracket function
is given by 
\beq\label{2qvofE}
\big [{\cal E}_a(f),f\big ]_{\lambda}(x)=(LY)\int_a^x{\cal E}_a(f)\,
d[f]_{\lambda},\qquad a\leq x\leq b.
\eeq
\end{thm}

\begin{proof}
Since the Dol\'eans exponential ${\cal E}_a(f)$ is a regulated function on 
$[a,b]$ and the bracket function $[f]_{\lambda}$ is nondecreasing,
the two Left Young integrals in (\ref{1qvofE})
and (\ref{2qvofE}) exist by Lemma \ref{bv}.
Also by Proposition \ref{variation}.\ref{LYjumps}, 
the two indefinite Left Young integrals are regulated functions with jumps
$$
\left\{{\cal E}_a(f)_{-}\Delta^{-}f\right\}^2,\,\,
\left\{{\cal E}_a(f)\Delta^{+}f\right\}^2
\quad\mbox{and}\quad
{\cal E}_a(f)_{-}\left\{\Delta^{-}f\right\}^2,\,\,
{\cal E}_a(f)\left\{\Delta^{+}f\right\}^2,
$$
respectively.
Hence the jumps of both sides of (\ref{1qvofE}) and (\ref{2qvofE})
agree by (\ref{-Ea}) and (\ref{+Ea}). 
By Propositions \ref{alpha} and \ref{ext-alpha},
it is enough to prove the relations 
\beq\label{3qvofE}
\lim_{m\to\infty}C({\cal E}_a(f),f;\lambda_m\Cap [a,x])
=(LY)\int_a^x{\cal E}_a(f)\,d[f]_{\lambda}
\eeq
and
\beq\label{4qvofE}
\lim_{m\to\infty}s_2({\cal E}_a(f);\lambda_m\Cap [a,x])
=(LY)\int_a^x{\cal E}_a(f)^2\,d[f]_{\lambda}
\eeq
for $a\leq x\leq b$.
For notation simplicity, we prove only that (\ref{3qvofE}) and
(\ref{4qvofE}) hold for $x=b$.
The proof for $x<b$ is the same.

Let $\lambda_m=\{x_i^m\colon\,i=0,\dots,n(m)\}$, $m\geq 1$, be partitions
of $[a,b]$ constituting the sequence $\lambda=\{\lambda_m\colon\,m\geq 1
\}\in\Lambda [a,b]$ such that $f$ has the quadratic $\lambda$-variation 
and the right discontinuity points of $f$ are accessible by $\lambda$.
For a real valued function $\rho$ on $[a,b]$, let $\Delta_i^m\rho:=
\rho (x_i^m)-\rho (x_{i-1}^m)$ for $i=1,\dots,n(m)$ and $m=1,2,\dots $.
Let $\phi (u,v):=e^uv$, $u, v\in\RR$, and let $\boldf :=(\bar f,V_a)$,
where $\bar f=f-f(a)-(1/2)[f]_{\lambda}^c$ and $V_a=\gamma (f;a,\cdot )$ 
is defined by (\ref{1bvofV}).
Then ${\cal E}_a(f)=\phi{\circ}\boldf$.
For the $2$-vector function $\boldf$, let $\Delta_i^m\boldf 
:=(\Delta_i^m\bar f,\Delta_i^mV_a)$.
By the mean value theorem in the Lagrange form, for each $m\geq 1$
and $i=1,\dots,n(m)$, there is a $\theta_i^m\in (0,1)$ such that
\beq\label{11qvofE}
\Delta_i^m\left ({\cal E}_a(f)\right )
=\Delta_i^m\left (\phi{\circ}\boldf\right )
=\left\langle\nabla\phi (\boldy_i^m),\Delta_i^m\boldf\right\rangle
\eeq
where $\nabla\phi (u,v)=(\frac{\partial\phi}{\partial u},
\frac{\partial\phi}{\partial v})(u,v)=(e^uv,e^u)$, $u,v\in\RR$,
and $\boldy_i^m:=\boldf (x_{i-1}^m)+\theta_i^m\Delta_i^m\boldf$.
Given a partition $\mu=\{z_j\colon\,j=0,\dots,m\}$ of $[a,b]$
let $J:=\{j=1,\dots,m-1\colon\,z_j\in\cup_m\lambda_m\}$.
There exists the minimal integer $m_0\geq 1$ such that 
$\{z_j\colon\,j\in J\}\subset\lambda_{m_0}$ and $(z_{j-1},z_j)
\cap\lambda_{m_0}\not =\emptyset$ for $j=1,\dots,m$.
For each $m\geq m_0$, let $i(j)$ be the index in
$\{1,\dots,n(m)\}$ such that $x_{i(j)-1}^m<z_j\leq x_{i(j)}^m$, 
$I_1(m):=\{1,i(j)\colon\,j=1,\dots,m\}\cup\{i(j)+1\colon\,
j\in J\}$ and let $I_2(m):=\{1,\dots,n(m)\}\setminus I_1(m)$.

First consider relation (\ref{3qvofE}).
For each  $m\geq m_0$, we then have
$$
C({\cal E}_a(f),f;\lambda_m)
=\sum_{i\in I_1(m)}\Delta_i^m{\cal E}_a(f)\Delta_i^mf
+\sum_{i\in I_2(m)}\Delta_i^m\big (\phi{\circ}\boldf\big )\Delta_i^mf
$$
$$
=\sum_{i=1}^{n(m)}{\cal E}_a(f)(x_{i-1}^m)\big\{\Delta_i^mf\big\}^2
+\sum_{i\in I_1(m)}\Big [\Delta_i^m{\cal E}_a(f)\Delta_i^mf
-{\cal E}_a(f)(x_{i-1}^m)\big\{\Delta_i^m f\big\}^2\Big ]+R_m^{(1)},
$$
where 
$$
R_m^{(1)}\!:=\!\sum_{i\in I_2(m)}\!\left\{\Big [\frac{\partial\phi}{\partial u}
(\boldy_i^m)-\frac{\partial\phi}{\partial u}(\boldf (x_{i-1}^m))\Big ]
\{\Delta_i^mf\}^2
+\Big [\frac{\partial\phi}{\partial v}(\boldy_i^m)\Delta_i^mV_a
-\frac{1}{2}\frac{\partial\phi}{\partial u}
(\boldy_i^m)\Delta_i^m[f]_{\lambda}^c\Big ]\Delta_i^mf\right\}.
$$
Suppose that given $\epsilon >0$ one can find a partition $\mu$ such that
$|R_m^{(1)}|<C\epsilon$ for all sufficiently large $m$.
Since $\Delta^{+}f(x)=0$ if $x\not\in\cup_m\lambda_m$, the limit
$$
\lim_{m\to\infty}
\sum_{i\in I_1(m)}\Big [\Delta_i^m{\cal E}_a(f)\Delta_i^mf
-{\cal E}_a(f)(x_{i-1}^m)\big\{\Delta_i^m f\big\}^2\Big ]
$$
$$
=\sum_{\mu\setminus a}\Big [\Delta^{-}{\cal E}_a(f)\Delta^{-}f
-{\cal E}_a(f)_{-}\{\Delta^{-}f\}^2\Big ]
+\sum_{\mu\setminus b}\Big [\Delta^{+}{\cal E}_a(f)\Delta^{+}f
-{\cal E}_a(f)\{\Delta^{+}f\}^2\Big ]=0
$$
for any given $\mu$ by (\ref{-Ea}) and (\ref{+Ea}).
An application of Lemma \ref{intqv} then yields
\beq\label{5qvofE}
\lim_{m\to\infty}C({\cal E}_a(f),f;\lambda_m)
=\lim_{m\to\infty}\sum_{i=1}^{n(m)}{\cal E}_a(f)(x_{i-1}^m)
\{\Delta_i^mf\}^2=(LY)\int_a^b{\cal E}_a(f)\,d[f]_{\lambda}.
\eeq
To choose a partition $\mu$, let $\epsilon \in (0,1)$ be given.
Since $\frac{\partial\phi}{\partial u}=\phi$ is uniformly continuous on 
the range of $\boldf$, there exists $\delta >0$ such that $|\phi (u_1,v_1)-
\phi (u_2,v_2)|<\epsilon$ if $|u_1|\vee |u_2|\leq\|\bar f\|_{\infty}$,
$|v_1|\vee |v_2|\leq\|V_a\|_{\infty}$ and $|u_1-u_2|\vee
|v_1-v_2|\leq\delta$.
Since $f$ and $V_a$ are regulated, there exists a partition
$\mu=\{z_j\colon\,j=0,\dots,m\}$ of $[a,b]$ such that
for each $j=1,\dots,m$,
$$
\osc (\bar f;(z_{j-1},z_j))<\delta,\qquad \osc(V_a;(z_{j-1},z_j))
<\delta,
$$
$$
\osc (f;(z_{j-1},z_j))<\epsilon\qquad\mbox{and}\qquad
\osc ([f]_{\lambda}^c;(z_{j-1},z_j))<\epsilon.
$$
Since $[\Delta_i^m]:=[x_{i-1}^m,x_i^m]\subset\cup_j(z_{j-1},z_j)$ 
for $i\in I_2(m)$, we then have the bound
$$
|R_m^{(1)}|\leq\max_{i\in I_2(m)}\osc (\phi{\circ}\boldf;[\Delta_i^m])
s_2(f;\lambda_m)
+\|\frac{\partial\phi}{\partial v}{\circ}\boldf\|_{\infty}
\max_{i\in I_2(m)}\osc (f;[\Delta_i^m])s_1(V_a;\lambda_m)
$$
$$
+\frac{1}{2}\|\phi{\circ}\boldf\|_{\infty}\Big (\max_{i\in I_2(m)}
\osc ([f]_{\lambda}^c;[\Delta_i^m])s_1([f]_{\lambda}^c;\lambda_m)
s_2(f;\lambda_m)\Big )^{1/2}\leq C\epsilon
$$ 
for all sufficiently large $m$.
Therefore (\ref{5qvofE}) holds.

Now consider relation (\ref{4qvofE}).
Again using the mean value theorem (\ref{11qvofE}) and the notation
following it, for each $m\geq m_0$, we have
$$
s_2({\cal E}_a(f);\lambda_m)
=\sum_{i\in I_1(m)}\left\{\Delta_i^m{\cal E}_a(f)\right\}^2
+\sum_{i\in I_2(m)}\left\{\Delta_i^m\big (\phi{\circ}\boldf \big )\right\}^2
$$
$$
=\sum_{i=1}^{n(m)}\left\{{\cal E}_a(f)(x_{i-1}^m)\Delta_i^mf\right\}^2
+\sum_{i\in I_1(m)}\Big [\left\{\Delta_i^m{\cal E}_a(f)\right\}^2
-\left\{{\cal E}_a(f)(x_{i-1}^m)\Delta_i^m f\right\}^2\Big ]+R_m^{(2)},
$$
where as before letting 
$\boldy_i^m:= \boldf (x_{i-1}^m)+\theta_i^m\Delta_i^m\boldf $,
$$
R_m^{(2)}:=\sum_{i\in I_2(m)}\Big\{\Big [\Big (\frac{\partial
\phi}{\partial u}( \boldy_i^m)\Big )^2-\Big (\frac{\partial
\phi}{\partial u}(\boldf (x_{i-1}^m))\Big )^2\Big ]
\big\{\Delta_i^mf\big\}^2
$$
$$
+ 2\frac{\partial\phi}{\partial u}(\boldy_i^m)
\Delta_i^mf\Big [\frac{\partial\phi}{\partial v}(\boldy_i^m)
\Delta_i^mV_a-\frac{1}{2}\frac{\partial\phi}{\partial u}
(\boldy_i^m )\Delta_i^m[f]^c\Big ]
+\Big [\frac{\partial\phi}{\partial v}(\boldy_i^m)\Delta_i^mV_a
-\frac{1}{2}\frac{\partial\phi}{\partial u}(\boldy_i^m )\Delta_i^m[f]^c
\Big ]^2\Big\}.
$$
As in the preceding case, given $\epsilon >0$ there exists a set $\mu$ 
such that $|R_m^{(2)}|<C\epsilon$ for all sufficiently large $m$.
Again since $\Delta^{+}f(x)=0$ if $x\not\in\cup_m\lambda_m$,
the limit
$$
\lim_{m\to\infty}
\sum_{i\in I_1(m)}\Big [\left\{\Delta_i^m{\cal E}_a(f)\right\}^2
-\left\{{\cal E}_a(f)(x_{i-1}^m)\Delta_i^m f\right\}^2\Big ]
$$
$$
\sum_{\mu\setminus a}\Big [\left\{\Delta^{-}{\cal E}_a(f)\right\}^2
-\left\{{\cal E}_a(f)_{-}\Delta^{-}f\right\}^2\Big ]
+\sum_{\mu\setminus b}\Big [\left\{\Delta^{+}{\cal E}_a(f)\right\}^2
-\left\{{\cal E}_a(f)\Delta^{+}f\right\}^2\Big ]=0
$$
by (\ref{-Ea}) and (\ref{+Ea}).
An application of Lemma \ref{intqv} then yields
$$
\lim_{m\to\infty}s_2({\cal E}_a(f);\lambda_m)=\lim_{m\to\infty}
\sum_{i=1}^{n(m)}\left\{{\cal E}_a(f)(x_{i-1}^m)\Delta_i^mf\right\}^2
=(LY)\int_a^b{\cal E}_a(f)^2\,d[f]_{\lambda}.
$$
The proof of Theorem \ref{qvofEa} is complete.
\qed\end{proof}

Next we prove analogous theorem with the forward Dol\'eans exponential 
${\cal E}_a(f)$ replaced by the backward Dol\'eans exponential ${\cal E}_b(f)$.

\begin{thm}\label{qvofEb}
Let  $f$ be a regulated function on $[a,b]$ having the quadratic 
$\lambda$-variation for some $\lambda\in \Lambda[a,b]$ such 
that left discontinuity points of $f$ are accessible by $\lambda$.
Then the Dol\'eans exponential ${\cal E}_b(f)$ has the quadratic
$\lambda$-variation on $[a,b]$, and its bracket function is given by
\beq\label{6qvofE}
\big [{\cal E}_b(f)\big ]_{\lambda}(x)=(RY)\int_a^x{\cal E}_b(f)^2\,
d[f]_{\lambda},\qquad a\leq x\leq b.
\eeq
Also, the $2$-vector function $({\cal E}_b(f),f)$ has the
quadratic $\lambda$-covariation on $[a,b]$, and its bracket function
is given by 
\beq\label{7qvofE}
\big [{\cal E}_b(f),f\big ]_{\lambda}(x)=-(RY)\int_a^x{\cal E}_b(f)\,
d[f]_{\lambda},\qquad a\leq x\leq b.
\eeq
\end{thm}

\begin{proof}
Since the Dol\'eans exponential ${\cal E}_b(f)$ is a regulated function 
on $[a,b]$ and the bracket function $[f]_{\lambda}$ is nondecreasing,
the two Right Young integrals in (\ref{6qvofE})
and (\ref{7qvofE}) are defined by Lemma \ref{bv}.
Also by Proposition \ref{variation}.\ref{LYjumps}, 
the two indefinite Right Young integrals are regulated functions with jumps
$$
\left\{{\cal E}_b(f)\Delta^{-}f\right\}^2,\,\,
\left\{{\cal E}_b(f)_{+}\Delta^{+}f\right\}^2
\quad\mbox{and}\quad
{\cal E}_b(f)\left\{\Delta^{-}f\right\}^2,\,\,
{\cal E}_b(f)_{+}\left\{\Delta^{+}f\right\}^2,
$$
respectively.
Hence the jumps of both sides of (\ref{6qvofE}) and (\ref{7qvofE})
agree by (\ref{-Eb}) and (\ref{+Eb}).
By Propositions \ref{alpha} and \ref{ext-alpha},
it is enough to prove the relations 
\beq\label{8qvofE}
\lim_{m\to\infty}C({\cal E}_b(f),f;\lambda_m\Cap [a,x])
=-(RY)\int_a^x{\cal E}_b(f)\,d[f]_{\lambda}
\eeq
and
\beq\label{9qvofE}
\lim_{m\to\infty}s_2({\cal E}_b(f);\lambda_m\Cap [a,x])
=(RY)\int_a^x{\cal E}_b(f)^2\,d[f]_{\lambda}
\eeq
for $a\leq x\leq b$.
For notation simplicity, we prove only that (\ref{8qvofE}) and
(\ref{9qvofE}) hold for $x=b$.

Let $\lambda_m=\{x_i^m\colon\,i=0,\dots,n(m)\}$, $m\geq 1$, be partitions
of $[a,b]$ constituting the sequence $\lambda=\{\lambda_m\colon\,m\geq 1
\}\in\Lambda [a,b]$ such that $f$ has the quadratic $\lambda$-variation 
and the left discontinuity points of $f$ are accessible by $\lambda$.
For a real valued function $\rho$ on $[a,b]$, let $\Delta_i^m\rho:=
\rho (x_i^m)-\rho (x_{i-1}^m)$ for $i=1,\dots,n(m)$ and $m=1,2,\dots $.
Let $\phi (u,v):=e^{-u}v$, $u, v\in\RR$, and let $\boldf :=(\bar f,V_b)$,
where $\bar f=f-f(b)+(1/2)[f]_{\lambda}^c(b)-(1/2)[f]_{\lambda}^c$ 
and $V_b=\gamma (f;\cdot,b)$ is defined by (\ref{1bvofV}).
Then ${\cal E}_b(f)=\phi{\circ}\boldf$.
For the $2$-vector function $\boldf$, let $\Delta_i^m\boldf 
:=(\Delta_i^m\bar f,\Delta_i^mV_b)$.
By the mean value theorem in the Lagrange form, for each $m\geq 1$
and $i=1,\dots,n(m)$, there is a $\theta_i^m\in (0,1)$ such that
\beq\label{12qvofE}
\Delta_i^m\left ({\cal E}_b(f)\right )
=\Delta_i^m\left (\phi{\circ}\boldf\right )
=\left\langle\nabla\phi (\boldy_i^m),\Delta_i^m\boldf\right\rangle
\eeq
where $\nabla\phi (u,v)=(\frac{\partial\phi}{\partial u},
\frac{\partial\phi}{\partial v})(u,v)=(-e^{-u}v,e^{-u})$, $u,v\in\RR$,
and $\boldy_i^m:=\boldf (x_{i-1}^m)+\theta_i^m\Delta_i^m\boldf$.
Given a partition $\mu=\{z_j\colon\,j=0,\dots,m\}$ of $[a,b]$
let $J:=\{j=1,\dots,m-1\colon\,z_j\in\cup_m\lambda_m\}$.
There exists the minimal integer $m_0\geq 1$ such that 
$\{z_j\colon\,j\in J\}\subset\lambda_{m_0}$ and $(z_{j-1},z_j)
\cap\lambda_{m_0}\not =\emptyset$ for $j=1,\dots,m$.
For each $m\geq m_0$, let $i(j)$ be the index in
$\{1,\dots,n(m)\}$ such that $x_{i(j)-1}^m<z_j\leq x_{i(j)}^m$, 
$I_1(m):=\{1,i(j)\colon\,j=1,\dots,m\}\cup\{i(j)+1\colon\,
j\in J\}$ and let $I_2(m):=\{1,\dots,n(m)\}\setminus I_1(m)$.

First consider relation (\ref{8qvofE}).
For each $m\geq m_0$, we then have
$$
C({\cal E}_b(f),f;\lambda_m)
=\sum_{i\in I_1(m)}\Delta_i^m{\cal E}_b(f)\Delta_i^mf
+\sum_{i\in I_2(m)}\Delta_i^m\big (\phi{\circ}\boldf\big )\Delta_i^mf
$$
$$
=-\sum_{i=1}^{n(m)}{\cal E}_b(f)(x_{i}^m)\big\{\Delta_i^mf\big\}^2
+\sum_{i\in I_1(m)}\Big [\Delta_i^m{\cal E}_b(f)\Delta_i^mf
+{\cal E}_b(f)(x_{i}^m)\big\{\Delta_i^m f\big\}^2\Big ]+R_m^{(1)},
$$
where
$$
R_m^{(1)}\!:=\!\sum_{i\in I_2(m)}\!\left\{\Big [\frac{\partial\phi}
{\partial u}(\boldy_i^m)
-\frac{\partial\phi}{\partial u}(\boldf (x_{i-1}^m))\Big ]\{\Delta_i^mf\}^2
+\Big [\frac{\partial\phi}{\partial v}(\boldy_i^m)\Delta_i^mV_b
-\frac{1}{2}\frac{\partial\phi}{\partial u}
(\boldy_i^m)\Delta_i^m[f]_{\lambda}^c\Big ]\Delta_i^mf\right\}.
$$
As in the proof of Theorem \ref{qvofEa}, given $\epsilon >0$ one 
can find a partition $\mu$ such that
$|R_m^{(1)}|<C\epsilon$ for all sufficiently large $m$.
Since $\Delta^{-}f(x)=0$ if $x\not\in\cup_m\lambda_m$, the limit
$$
\lim_{m\to\infty}
\sum_{i\in I_1(m)}\Big [\Delta_i^m{\cal E}_b(f)\Delta_i^mf
+{\cal E}_b(f)(x_{i}^m)\big\{\Delta_i^m f\big\}^2\Big ]
$$
$$
=\sum_{\mu\setminus a}\Big [\Delta^{-}{\cal E}_b(f)\Delta^{-}f
+{\cal E}_b(f)\{\Delta^{-}f\}^2\Big ]
+\sum_{\mu\setminus b}\Big [\Delta^{+}{\cal E}_b(f)\Delta^{+}f
+{\cal E}_b(f)_{+}\{\Delta^{+}f\}^2\Big ]=0
$$
for any given $\mu$ by (\ref{-Eb}) and (\ref{+Eb}).
An application of Lemma \ref{intqv} then yields
$$
\lim_{m\to\infty}C({\cal E}_b(f),f;\lambda_m)
=-\lim_{m\to\infty}\sum_{i=1}^{n(m)}{\cal E}_b(f)(x_{i}^m)
\{\Delta_i^mf\}^2=-(RY)\int_a^b{\cal E}_b(f)\,d[f]_{\lambda}.
$$
That is (\ref{8qvofE}) holds with $x=b$.

Now consider relation (\ref{9qvofE}).
Again using the mean value theorem (\ref{12qvofE}) and the notation
following it, for each $m\geq m_0$, we have
$$
s_2({\cal E}_b(f);\lambda_m)
=\sum_{i\in I_1(m)}\left\{\Delta_i^m{\cal E}_b(f)\right\}^2
+\sum_{i\in I_2(m)}\left\{\Delta_i^m\big (\phi{\circ}\boldf \big )\right\}^2
$$
$$
=\sum_{i=1}^{n(m)}\left\{{\cal E}_b(f)(x_{i}^m)\Delta_i^mf\right\}^2
+\sum_{i\in I_1(m)}\Big [\left\{\Delta_i^m{\cal E}_b(f)\right\}^2
-\left\{{\cal E}_b(f)(x_{i}^m)\Delta_i^m f\right\}^2\Big ]+R_m^{(2)},
$$
where  as before letting 
$\boldy_i^m:= \boldf (x_{i-1}^m)+\theta_i^m\Delta_i^m\boldf $,
$$
R_m^{(2)}:=\sum_{i\in I_2(m)}\Big\{\Big [\Big (\frac{\partial
\phi}{\partial u}( \boldy_i^m)\Big )^2-\Big (\frac{\partial
\phi}{\partial u}(\boldf (x_{i-1}^m))\Big )^2\Big ]
\big\{\Delta_i^mf\big\}^2
$$
$$
+ 2\frac{\partial\phi}{\partial u}(\boldy_i^m)
\Delta_i^mf\Big [\frac{\partial\phi}{\partial v}(\boldy_i^m)
\Delta_i^mV_b-\frac{1}{2}\frac{\partial\phi}{\partial u}
(\boldy_i^m )\Delta_i^m[f]^c\Big ]
+\Big [\frac{\partial\phi}{\partial v}(\boldy_i^m)\Delta_i^mV_b
-\frac{1}{2}\frac{\partial\phi}{\partial u}(\boldy_i^m )\Delta_i^m[f]^c
\Big ]^2\Big\}.
$$
As in the preceding theorem, given $\epsilon >0$ there exists a partition $\mu$ 
such that $|R_m^{(2)}|<C\epsilon$ for all sufficiently large $m$.
Again since $\Delta^{-}f(x)=0$ if $x\not\in\cup_m\lambda_m$,
the limit
$$
\lim_{m\to\infty}
\sum_{i\in I_1(m)}\Big [\left\{\Delta_i^m{\cal E}_b(f)\right\}^2
-\left\{{\cal E}_b(f)(x_{i}^m)\Delta_i^m f\right\}^2\Big ]
$$
$$
\sum_{\mu\setminus a}\Big [\left\{\Delta^{-}{\cal E}_b(f)\right\}^2
-\left\{{\cal E}_b(f)\Delta^{-}f\right\}^2\Big ]
+\sum_{\mu\setminus b}\Big [\left\{\Delta^{+}{\cal E}_b(f)\right\}^2
-\left\{{\cal E}_b(f)_{+}\Delta^{+}f\right\}^2\Big ]=0
$$
by (\ref{-Eb}) and (\ref{+Eb}).
An application of Lemma \ref{intqv} then yields
$$
\lim_{m\to\infty}s_2({\cal E}_b(f);\lambda_m)=\lim_{m\to\infty}
\sum_{i=1}^{n(m)}\left\{{\cal E}_b(f)(x_{i}^m)\Delta_i^mf\right\}^2
=(RY)\int_a^b{\cal E}_b(f)^2\,d[f]_{\lambda}.
$$
The proof of Theorem \ref{qvofEb} is complete.
\qed\end{proof}

\section{Augmented linear $\lambda$-integral equations}\label{uniqueness}

Let $f$ be a regulated function on $[a,b]$ which is $(\lambda,p)$-decomposable 
for some $\lambda\in\Lambda [a,b]$ and $ 1\leq p<2$, and let
$D_{\lambda,p}^{+}(f)$ be nonempty.
We consider the forward linear Left $\lambda$-integral equation
\beq\label{1LIE}
F(y)=1+(L)\int_a^yF\,d_{\lambda}f, \qquad a\leq y\leq b,
\eeq
provided the Left $\lambda$-integral
is defined on $[a,b]$ with respect to $D_{\lambda,p}^{+}(f)$.
It is not hard to check using a chain rule that the forward Dol\'eans 
exponential ${\cal E}_a(f)$ satisfy (\ref{1LIE}).
More challenging is the task to find a large enough class of functions $F$ 
which contains the Dol\'eans exponential ${\cal E}_a(f)$ and contains no 
other solutions of (\ref{1LIE}).
For the meantime we show that an augmented forward linear Left
$\lambda$-integral equation has the Dol\'eans exponential as the
unique solution in the class of all functions in $\dual ({\cal W}_p)$ which
have the quadratic $\lambda$-variation.

\begin{defn}\label{solution}
{\rm For $\lambda\in\Lambda [a,b]$ and $ 1\leq p<2$, let $f$ be 
a $(\lambda,p)$-decomposable function on $[a,b]$, and let
$D_{\lambda,p}^{+}(f)$ be nonempty.
We say that a function $F\in \dual({\cal W}_p)[a,b]\cap Q_{\lambda}[a,b]$ 
is a \emph{solution
of the augmented forward linear Left $\lambda$-integral equation}  
\beq\label{4LIE}
d_{\lambda}^{\to}F=F\,d_{\lambda}^{\to}f\qquad\mbox{on $[a,b]$}
\eeq
if $(a)$ and $(b)$ hold, where
\begin{enumerate}
\item[$(a)$]
for $(g,h)\in D_{\lambda,p}^{+}(f)$ and $\psi\in \dual ({\cal W}_p)[a,b]$,
$(LC)\smallint\big\langle (\psi,-\psi F),d_{\lambda}(F,g)\big\rangle$
is defined on $[a,b]$, and for $a \leq y\leq b$,
\beq\label{2LIE}
(LC)\int_a^y\big\langle (\psi,-\psi F),d_{\lambda}(F,g)\big\rangle
=(LY)\int_a^y\psi F\,dh;
\eeq
\item[$(b)$]
the $2$-vector function $(F,f)$ has the
quadratic $\lambda$-variation and for $a\leq y\leq b$,
\beq\label{3LIE}
\left ( \begin{array}{cc}
[F]_{\lambda} & [f,F]_{\lambda} \cr
[F,f]_{\lambda} & [f]_{\lambda}
\end{array} \right ) (y)=(LY)\int_a^y
\left ( \begin{array}{cc}
F^2 & F \cr
F & 1
\end{array} \right ) \,d[f]_{\lambda}.
\eeq
\end{enumerate}}
\end{defn}

Let $f$ be a $(\lambda,p)$-decomposable function on $[a,b]$
for some $\lambda\in\Lambda [a,b]$ and $1\leq p<2$, and let
$D_{\lambda,p}^{+}(f)$ be nonempty.
Suppose that $F\in\dual ({\cal W}_p)$ has the quadratic $\lambda$-variation
and is a solution of the augmented forward linear Left $\lambda$-integral 
equation (\ref{4LIE}).
Then in $(a)$ of Definition \ref{solution}
taking $\psi\equiv 1$, by Proposition \ref{property3}, it follows that the 
Left $\lambda$-integral $(L)\smallint F\,d_{\lambda}f$ is defined on $[a,b]$ 
with respect to $D_{\lambda,p}^{+}(f)$, and (\ref{1LIE}) holds.
The more general relation (\ref{2LIE}) contains a form of the substitution
rule for the Left $\lambda$-integral, analogous to the associativity property 
of the stochastic integral.
Indeed, assuming the following integrals exist 
and using (\ref{2LIE}), (\ref{1LIE}), we have
\begin{eqnarray*}
\lefteqn{(LC)\int_a^b\psi\,d_{\lambda}F-(LC)\int_a^b\psi F\,
d_{\lambda}g}\\[2mm]
&=&(LC)\int_a^b\big\langle (\psi,-\psi F),d_{\lambda}(F,g)\big\rangle
=(LY)\int_a^b\psi F\,dh\quad
\mbox{by Theorem \ref{variation}.\ref{srule}}\\[2mm]
&=&(LY)\int_a^b\psi\,d\Big ((LY)\int_a^{\cdot}F\,dh\Big )
=(LC)\int_a^b\psi\,d_{\lambda}F-(LC)\int_a^b\psi\,d_{\lambda}
\Big ((LC)\int_a^{\cdot}F\,d_{\lambda}g\Big ).
\end{eqnarray*}
Comparing the left and right sides of these relations we get 
the substitution rule for the Left $\lambda$-integral:
$$
(LC)\int_a^b\psi F\,d_{\lambda}g=(LC)\int_a^b\psi\,d_{\lambda}
\Big ((LC)\int_a^{\cdot}F\,d_{\lambda}g\Big ).
$$

The linear Left Young integral with respect to a function having
bounded $p$-variation for some $1\leq p<2$ was solved by
Dudley and Norvai\v sa \cite[Theorem II.5.21]{DNa}:

\begin{thm}\label{DN}
Let $f\in {\cal W}_p[a,b]$ for some $1\leq p<2$.
In the class {\rm $\dual ({\cal W}_p)[a,b]$} the Dol\'eans exponential 
${\cal E}_a(f)$ is the unique
solution of the forward linear Left Young integral equation
\beq\label{1DN}
F(y)=1+(LY)\int_a^yF\,df,\qquad a\leq y\leq b.
\eeq
\end{thm}

Dudley and Norvai\v sa considered equation (\ref{1DN}) with respect to
a function $f$ having values in a Banach algebra, and the unique solution
of their equation is the indefinite product integral with respect
to the function $f$.
In the case $f$ has real values, the indefinite product integral
with respect to $f$ has values which agree with the Dol\'eans 
exponential ${\cal E}_a(f)$ by Theorem II.4.4 of Dudley and 
Norvai\v sa \cite{DNa}.
The only difference with the cited result of Dudley
and Norvai\v sa is that  the uniqueness class in Theorem
\ref{DN} is larger in the case $p=1$, which is the class of
all regulated functions.

Now we are ready to solve the augmented forward linear
Left  $\lambda$-integral equation (\ref{4LIE}).

\begin{thm}\label{LIE}
Let $\lambda\in\Lambda [a,b]$, $1\leq p<2$ and let
$f$ be a $(\lambda,p)$-decomposable function on $[a,b]$.
If $D_{\lambda,p}^{+}(f)\not =\emptyset$ then in the class
{\rm $\dual ({\cal W}_p)[a,b]\cap Q_{\lambda}[a,b]$} the Dol\'eans 
exponential ${\cal E}_a(f)$ is the unique solution of the augmented 
forward linear Left $\lambda$-integral equation {\rm (\ref{4LIE})}.
\end{thm}

\begin{proof}
Let $(g,h)\in D_{\lambda,p}^{+}(f)$.
By Corollary \ref{qvclassQ}, $f$ has the quadratic $\lambda$-variation,
and therefore the Dol\'eans exponential ${\cal E}_a(f)$ is defined. 
Next we show that ${\cal E}_a(f)$ is a solution of the augmented forward
linear Left $\lambda$-integral equation (\ref{4LIE}).
By Theorem \ref{qvofEa}, the $2$-vector function $({\cal E}_a(f),f)$ 
has the quadratic $\lambda$-variation and (\ref{3LIE}) holds with 
$F={\cal E}_a(f)$.
Let $\bar f:=f-f(a)-(1/2)[g]_{\lambda}^c$, and $V_a:=\gamma (f;a,\,\cdot\,)$,
where $\gamma (f)$ is defined by (\ref{1bvofV}).
Then ${\cal E}_a(f)=e^{\bar f}V_a$ and $V_a$ is pure jump function of 
bounded variation by Lemma \ref{bvofV}.
Since the composition of the function $\bar f$ with a smooth function
does not change its $p$-variation, and since ${\cal W}_p$ is a Banach 
algebra, ${\cal E}_a(f)\in\dual ({\cal W}_p)[a,b]$, and so ${\cal E}_a(f)$
satisfies condition $(b)$ of Definition \ref{solution}.
To show that ${\cal E}_a(f)$ satisfies condition $(a)$, we apply 
the chain rule formula of Theorem \ref{mdchrule}.$(a)$ to the composition 
$\phi{\circ}\boldf={\cal E}_a(f)$, where the $2$-vector function
$\boldf$ is defined by
\beq\label{21LIE}
\boldf :=(\bar f,V_a)=(C-f(a),0)+(g,0)
+(h-\frac{1}{2}[g]_{\lambda}^c,V_a)\quad\mbox{and}\quad
\phi (u,v)=e^uv,\quad u,v\in\RR.
\eeq
Let $\boldg:= (g,0)$ and $\boldh :=(h-(1/2)[g]_{\lambda}^c,V_a)$.
Then $\boldf$ is $(\lambda,p)$-decomposable by $(\boldg,\boldh)$.
Since all discontinuity points of $V_a$ are the same as 
discontinuity points of $f$, the right
discontinuity points of $\boldf$ are accessible by $\lambda$.
Let $\psi\in\dual ({\cal W}_p)[a,b]$.
By statement $(a)$ of Theorem \ref{mdchrule}, the Left Cauchy
$\lambda$-integral
$$
(LC)\int\big\langle\big (\psi,-\psi(\nabla\phi{\circ}\boldf)\big ),
d_{\lambda}\big (\phi{\circ}\boldf,\boldg\big )\big\rangle 
=(LC)\int\big\langle\big (\psi,-\psi {\cal E}_a(f)\big ),
d_{\lambda}\big ({\cal E}_a(f),g\big )\big\rangle
$$
is defined on $[a,b]$,
where $\nabla\phi (u,v)=(\frac{\partial\phi}{\partial u},
\frac{\partial\phi}{\partial v})(u,v)=(e^uv,e^u)$ for $u,v\in\RR$,
and the equalities
\begin{eqnarray}
\lefteqn{(LC)\int_a^y\big\langle\big (\psi,-\psi {\cal E}_a(f)\big ),
d_{\lambda}\big ({\cal E}_a(f),g\big )\big\rangle}\nonumber\\[2mm]
&=&(LY)\int_a^y\psi {\cal E}_a(f)\,d(h-(1/2)[g]_{\lambda}^c)
+(LY)\int_a^y\psi e^{\bar f}\,dV_a
+\frac{1}{2}(RS)\int_a^y\psi {\cal E}_a(f)\,d[g]_{\lambda}^c\nonumber\\[2mm]
&+&\sum_{(a,y]}\psi_{-}\Big\{\Delta^{-}{\cal E}_a(f)-{\cal E}_a(f)_{-}
\Delta^{-}f-e^{\bar f_{-}}\Delta^{-}V_a\Big\}
+\sum_{[a,y)}\psi\Big\{\Delta^{+}{\cal E}_a(f)-{\cal E}_a(f)\Delta^{+}f
-e^{\bar f}\Delta^{+}V_a\Big\}\nonumber\\[2mm]
&=&(LY)\int_a^y\psi {\cal E}_a(f)\,dh
+(LY)\int_a^y\psi e^{\bar f}\,dV_a
-\Big\{\sum_{(a,y]}\psi_{-}e^{\bar f_{-}}\Delta^{-}V_a+\sum_{[a,y)}
\psi e^{\bar f}\Delta^{+}V_a\Big\}\label{9LIE}
\end{eqnarray}
hold for $a\leq y\leq b$.
The last equality follows by Lemma \ref{bv} and relations
(\ref{-Ea}), (\ref{+Ea}).
Let $\lambda=\{\lambda_m\colon\,m\geq 1\}$ with
$\lambda_m=\{x_i^m\colon\,i=0,\dots,n(m)\}$.
Since $V_a$ is pure jump function of bounded variation, and since 
$\{x\in (a,b)\colon\,\Delta^{+}V_a(x)\not =0\}\subset\cup_m\lambda_m$
by (\ref{3bvofV}), by Theorem \ref{variation}.\ref{approximation}, 
we have the equalities
\begin{eqnarray*}
(LY)\int_a^y\psi e^{\bar f}\,dV_a&=&\lim_{m\to\infty}\sum_{i=1}^{n(m)}
\psi (x_{i-1}^m)e^{\bar f(x_{i-1}^m)}\big [V_a(x_i^m\wedge y)
-V_a(x_{i-1}^m\wedge y)\big ]\\[2mm]
&=&\sum_{(a,y]}\psi_{-}e^{\bar f_{-}}\Delta^{-}V_a
+\sum_{[a,y)}\psi e^{\bar f}\Delta^{+}V_a
\end{eqnarray*}
for $a\leq y\leq b$.
Thus the last two terms in (\ref{9LIE}) give zero to the equality,
and hence the Dol\'eans exponential ${\cal E}_a(f)$ solves the augmented
forward linear left $\lambda$-integral equation (\ref{4LIE}).

Finally, suppose that $F\in\dual ({\cal W}_p)$, has the quadratic 
$\lambda$-variation and solves the augmented
forward linear Left $\lambda$-integral equation (\ref{4LIE}). 
We then show that  $F={\cal E}_a(f)$.
To this aim let $V:=e^{-\bar f}F$, where $\bar f:=f-f(a)-(1/2)
[g]_{\lambda}^c$.
It is enough to prove that $V$ is equal to the pure jump
function $V_a=\gamma (f;a,\,\cdot)$.
Let $\boldC:=(C-f(a),0)$, $\boldg:=(g,F)$, $\boldh:=(h-(1/2)
[g]_{\lambda}^c,0)$,
\beq\label{22LIE}
\boldf:=(\bar f,F)=\boldC+\boldg+\boldh
\qquad\mbox{and}\qquad\mbox{$\phi (u,v):=e^{-u}v\,\,\,$ 
for $u, v\in\RR$.}
\eeq
Then $V=\phi{\circ}\boldf$.
We wish to apply a chain rule to show that $V$ is the unique solution
to the linear Left Young integral equation with respect to a  
pure jump function of bounded variation.
To this aim we have to show that the $2$-vector function $\boldg =(g,F)$ 
has the quadratic $\lambda$-variation.
Recall that by condition $(b)$ of Definition \ref{solution}, the
$2$-vector function $(f,F)$ has the quadratic $\lambda$-variation.

\begin{lem}\label{QVandqc}
For $\lambda\in\Lambda [a,b]$
and $p\in [1,2)$, let $f$ be  a regulated function on $[a,b]$ which is
$(\lambda,p)$-decomposable by $(g,h)\in D_{\lambda,p}^{+}(f)$.
Let {\rm $F\in\dual ({\cal W}_p)[a,b]$} be such that 
\beq\label{3QVandqc}
\{x\in (a,b)\colon\, \Delta^{+}F(x)\not =0\}\subset
\{x\in (a,b)\colon\, \Delta^{+}f(x)\not =0\}.
\eeq
Then the $2$-vector function $(f,F)$ has the quadratic 
$\lambda$-covariation if and only if the $2$-vector function $(g,F)$
has the quadratic $\lambda$-covariation,
and if the two functions have the quadratic $\lambda$-covariation then 
$[f,F]_{\lambda}^c=[g,F]_{\lambda}^c$.
\end{lem}

\begin{proof}
To begin with we show that the $2$-vector function $(h,F)$
has the quadratic $\lambda$-covariation with the bracket function
\beq\label{1QVandqc}
[h,F]_{\lambda}(y)=\sum_{(a,y]}\Delta^{-}h\Delta^{-}F
+\sum_{[a,y)}\Delta^{+}h\Delta^{+}F,\qquad a\leq y\leq b.
\eeq
To this aim we apply Proposition \ref{qcforpq} with $f_1=h\in {\cal W}_p$ 
and $f_2=F\in\dual ({\cal W}_p)$.
If $(\Delta^{-}F\Delta^{+}h+\Delta^{-}h\Delta^{+}F)(x)\not =0$ for some 
$x\in (a,b)$, then either $\{x\in (a,b)\colon\, \Delta^{+}h(x)\not =0\}$
or $\{x\in (a,b)\colon\, \Delta^{+}F(x)\not =0\}$.
Since $(g,h)\in D_{\lambda,p}^{+}(f)$, by (\ref{3QVandqc}),
 (\ref{2qcforpq}) holds.
Therefore by Proposition \ref{qcforpq}, the function $(h,F)$ has the quadratic 
$\lambda$-covariation with the bracket function (\ref{1QVandqc}).
Letting $\lambda=\{\lambda_m\colon\,m\geq 1\}$, the conclusion then
follows from the relation
$$
C(f,F;\lambda_m\Cap [u,v])=C(g,F;\lambda_m\Cap [u,v])
+C(h,F;\lambda_m\Cap [u,v]),
$$
which holds for each $m$ and $a\leq u<v\leq v$,
and since $\Delta f\Delta F=\Delta g\Delta F+\Delta h\Delta F$.
The proof of Lemma \ref{QVandqc} is complete.
\qed\end{proof}

To apply Lemma \ref{QVandqc} to $f$, $g$, $h$ and $F$ as before
we have to check its assumption (\ref{3QVandqc}).
By condition $(a)$ of Definition \ref{solution} with $\psi\equiv 1$ and by 
Proposition \ref{property3}, the Left $\lambda$-integral 
$(L)\smallint F\,d_{\lambda}f$ is defined on $[a,b]$ 
with respect to $D_{\lambda,p}^{+}(f)$, and (\ref{1LIE}) holds.
By Definition \ref{LCint} of the Left Cauchy $\lambda$-integral in 
conjunction with Proposition \ref{variation}.\ref{LYjumps} concerning 
jumps of the indefinite Left Young integral, we have that
\beq\label{5LIE}
\Delta^{-}F(y)=\big (F_{-}\Delta^{-}f\big )(y)
\qquad\mbox{and}\qquad
\Delta^{+}F(x)=\big (F\Delta^{+}f\big )(x)
\eeq
for $a\leq y<x\leq b$.
In particular, it then follows that (\ref{3QVandqc}) holds.
Therefore by Lemma \ref{QVandqc}, the $2$-vector function
$\boldg=(g,F)$ has the quadratic $\lambda$-covariation, and thereby
has the quadratic $\lambda$-variation.
Since the pair $(\boldg,\boldh)$ is $(\lambda,p)$-dual and
(\ref{22LIE}) holds, $\boldf$ is $(\lambda,p)$-decomposable by  
$(\boldg,\boldh)$.
By (\ref{5LIE}), it follows that $(\boldg,\boldh)\in D_{\lambda,p}^{+}(\boldf)$.
Therefore all assumptions of Theorem \ref{mdchrule}.$(a)$ are satisfied
by $\boldf$ and $\phi$.
By the first part of statement $(a)$ of Theorem \ref{mdchrule} with 
$\psi\equiv 1$, the Left Cauchy $\lambda$-integrals
$$
(LC)\int\big\langle(1,-\nabla\phi{\circ}\boldf),d_{\lambda}
(\phi{\circ}\boldf,\boldg)\big\rangle
=\mu (\phi{\circ}\boldf)-
(LC)\int\langle\nabla\phi{\circ}\boldf,d_{\lambda}\boldg\rangle
$$
are defined on $[a,b]$,
where $\nabla\phi (u,v)=(\frac{\partial\phi}{\partial u},
\frac{\partial\phi}{\partial v})(u,v)=(-e^{-u}v,e^{-u})$ for
$u,v\in\RR$, and $\mu$ is defined by (\ref{muf}).
Further, by condition $(a)$ of Definition \ref{solution} 
with $\psi=\exp\{-\bar f\}$, we have the equality
\beq\label{6LIE}
(LC)\int \langle\nabla\phi{\circ}\boldf,d_{\lambda}\boldg\rangle
=(LC)\int \Big\langle \Big (-e^{-\bar f}F,e^{-\bar f}\Big ),
d_{\lambda}(g,F)\Big\rangle=(LY)\int e^{-\bar f}F\,dh.
\eeq
By condition $(b)$ of Definition \ref{solution}, by the second part of Lemma 
\ref{QVandqc}, and by (\ref{5LIE}), we have
$$
(LY)\int_a^yF\,d[f]_{\lambda}=
[f,F]_{\lambda}(y)=[g,F]_{\lambda}^c(y)+\sum_{(a,y]}F_{-}\big\{
\Delta^{-}f\big\}^2+\sum_{[a,y)}F\big\{\Delta^{+}f\big\}^2
$$
for $a\leq y\leq b$.
On the other hand, by the first part of Lemma \ref{qv} 
and since $[f]_{\lambda}^c=[g]_{\lambda}^c$, we have
$$
(LY)\int_a^yF\,d[f]_{\lambda}=(RS)\int_a^yF\,d[g]_{\lambda}^c
+\sum_{(a,y]}F_{-}\big\{
\Delta^{-}f\big\}^2+\sum_{[a,y)}F\big\{\Delta^{+}f\big\}^2
$$
for $a\leq y\leq b$.
Therefore $[g,F]_{\lambda}^c(y)=(RS)\smallint_a^yF\,d[g]_{\lambda}^c$
for $a\leq y\leq b$.
This in conjunction with the substitution theorem for the Riemann-Stieltjes
integral, yields the equality
\beq\label{7LIE}
(RS)\int_a^ye^{-\bar f}\,d[g,F]_{\lambda}^c
=(RS)\int_a^ye^{-\bar f}F\,d[g]_{\lambda}^c
\eeq
for $a\leq y\leq b$.
By the second part of Theorem \ref{mdchrule}.$(a)$ with $\psi\equiv 1$
applied to $V=\phi{\circ}\boldf$, we then have
\begin{eqnarray}
V(y)-V(a)&=&(LC)\int_a^y\langle\nabla\phi{\circ}\boldf,
d_{\lambda}\boldg\rangle-(LY)\int_a^ye^{-\bar f}F\,
d\big (h-\frac{1}{2}[g]_{\lambda}^c\big )\nonumber\\[2mm]
& &+\frac{1}{2}(RS)\int_a^ye^{-\bar f}F\,d[g]_{\lambda}^c-(RS)\int_a^y
e^{-\bar f}d[g,F]_{\lambda}^c\nonumber\\[2mm]
& &+\sum_{(a,y]}\Big\{\Delta^{-}V-\Big (e^{-\bar f}\Big )_{-}\Delta^{-}F
+\Big (e^{-\bar f}F\Big )_{-}\Delta^{-}f\Big\}\nonumber\\[2mm]
& &+\sum_{[a,y)}\Big\{\Delta^{+}V-e^{-\bar f}\Delta^{+}F
+e^{-\bar f}F\Delta^{+}f\Big\}\nonumber\\[2mm]
\mbox{by Lemma \ref{bv}}\quad 
&=&(LC)\int_a^y\langle\nabla\phi{\circ}\boldf,d_{\lambda}\boldg\rangle
-(LY)\int_a^ye^{-\bar f}F\,dh\nonumber\\[2mm]
\mbox{by (\ref{7LIE})}\quad
& &+\frac{1}{2}(RS)\int_a^ye^{-\bar f}F\,d[g]_{\lambda}^c
+\frac{1}{2}(RS)\int_a^ye^{-\bar f}F\,d[g]_{\lambda}^c
-(RS)\int_a^ye^{-\bar f}Fd[g]_{\lambda}^c\nonumber\\[2mm]
\mbox{by (\ref{5LIE})}\quad
& &+\sum_{(a,y]}\Big\{\Delta^{-}V-\Big (e^{-\bar f}\Big )_{-}F_{-}
\Delta^{-}f+\Big (e^{-\bar f}\Big )_{-}F_{-}\Delta^{-}f\Big\}
\nonumber\\[2mm]
& &+\sum_{[a,y)}\Big\{\Delta^{+}V-e^{-\bar f}F\Delta^{+}f
+e^{-\bar f}F\Delta^{+}f\Big\}\nonumber\\[2mm]
\mbox{by (\ref{6LIE})}\quad
&=&\sum_{(a,y]}\Delta^{-}V+\sum_{[a,y)}\Delta^{+}V\label{8LIE}
\end{eqnarray}
for all $a< y\leq b$.
Therefore $V$ is pure jump function of bounded variation.
Further by (\ref{5LIE}), we have
$$
\Delta^{-}V=\big (F_{-}+\Delta^{-}F\big )\left (e^{-\bar f}\right)_{-}
e^{-\Delta^{-}f}-\left (Fe^{-\bar f}\right )_{-}
=V_{-}\Big\{(1+\Delta^{-}f)e^{-\Delta^{-}f}-1\Big\}
$$
and
$$
\Delta^{+}V=\big (F+\Delta^{+}F\big )e^{-\bar f}e^{-\Delta^{+}f}
-Fe^{-\bar f}=V\Big\{(1+\Delta^{+}f)e^{-\Delta^{+}f}-1\Big\}.
$$
By Taylor's formula (\ref{7bvofV}), the function $A$ defined by
\beq\label{20LIE}
A(y):=\sum_{(a,y]}\Big\{(1+\Delta^{-}f)e^{-\Delta^{-}f}-1\Big\}
+\sum_{[a,y)}\Big\{(1+\Delta^{+}f)e^{-\Delta^{+}f}-1\Big\}
\eeq
for $a\leq y\leq b$, is pure jump function of bounded variation.
By Theorem \ref{variation}.\ref{approximation}, we have
$$
(LY)\int_a^yV\,dA=\lim_{m\to\infty}S_{LC}(V,A;\lambda_m\Cap [a,y])
=\sum_{(a,y]}V_{-}\Delta^{-}A+\sum_{[a,y)}V\Delta^{+}A
=\sum_{(a,y]}\Delta^{-}V+\sum_{[a,y)}\Delta^{+}V
$$
for $a\leq y\leq b$.
Thus by (\ref{8LIE}), the function $V$ is a solution to the linear
Left Young integral equation
$$
V(y)=1+(LY)\int_a^yV\,dA,\qquad a\leq y\leq b.
$$
By Theorem \ref{DN}, the solution is unique in the class of 
all regulated functions, which
certainly contains the class of functions of bounded variation.
Therefore
$$
V(y)=e^{A(y)-A(a)}\prod_{[a,y]}\big (1+\Delta A\big )e^{-\Delta A}
=\prod_{[a,y]}\big (1+\Delta f\big )e^{-\Delta f},
$$
for $a\leq y\leq b$. 
That is $V=V_a=\gamma (f;a,\,\cdot\,)$, and hence $F={\cal E}_a(f)$
as desired.
The proof of Theorem \ref{LIE} is complete.
\qed\end{proof}

Again for given $\lambda\in\Lambda [a,b]$ and $1\leq p<2$, 
let $f$ be a $(\lambda,p)$-decomposable function on $[a,b]$.
This time we consider the backward linear Right $\lambda$-integral equation
\beq\label{10LIE}
G(y)=1+(R)\int_y^bG\,d_{\lambda}f, \qquad a\leq y\leq b,
\eeq
provided the Right $\lambda$-integral is
defined on $[a,b]$ with respect to $D_{\lambda,p}^{-}(f)$.
Using a chain rule it is easy to see that the backward Dol\'eans exponential 
${\cal E}_b(f)$ satisfy (\ref{10LIE}).
Similarly to the forward case, we show that this function 
is the unique solution of an augmented equation.

\begin{defn}\label{Bsolution}
{\rm For $\lambda\in\Lambda [a,b]$ and $ 1\leq p<2$, let $f$ be 
a $(\lambda,p)$-decomposable function on $[a,b]$, and let
$D_{\lambda,p}^{-}(f)$ be nonempty.
We say that a function $G\in \dual({\cal W}_p)[a,b]\cap Q_{\lambda}[a,b]$ 
is a \emph{solution
of the augmented backward linear Right $\lambda$-integral equation}  
\beq\label{23LIE}
d_{\lambda}^{\gets}G=G\,d_{\lambda}^{\gets}f\qquad\mbox{on $[a,b]$}
\eeq
if $(a)$ and $(b)$ hold, where
\begin{enumerate}
\item[$(a)$]
for $(g,h)\in D_{\lambda,p}^{-}(f)$ and $\psi\in \dual ({\cal W}_p)[a,b]$,
$(RC)\smallint\big\langle (\psi,\psi G),d_{\lambda}(G,g)\big\rangle$
is defined on $[a,b]$, and for $a\leq y\leq b$,
\beq\label{13LIE}
(RC)\int_y^b\big\langle (\psi,\psi G),d_{\lambda}(G,g)\big\rangle
=-(RY)\int_y^b\psi G\,dh;
\eeq
\item[$(b)$]
the $2$-vector function $(G,f)$ has the
quadratic $\lambda$-variation and for $a\leq y\leq b$,
\beq\label{11LIE}
\left ( \begin{array}{cc}
[G]_{\lambda} & [f,G]_{\lambda} \cr
[G,f]_{\lambda} & [f]_{\lambda}
\end{array} \right ) (y)=(RY)\int_a^y
\left ( \begin{array}{cc}
G^2 & -G \cr
-G & 1
\end{array} \right ) \,d[f]_{\lambda}.
\eeq
\end{enumerate}}
\end{defn}

Let $f$ be a $(\lambda,p)$-decomposable function on $[a,b]$ for some
$\lambda\in\Lambda [a,b]$ and $1\leq p<2$, and let $D_{\lambda,p}^{-}(f)$
be nonempty.
Suppose that $G\in\dual ({\cal W}_p)$, has the quadratic $\lambda$-variation
and is a solution of the augmented backward linear Right $\lambda$-integral 
equation (\ref{23LIE}).
Then in $(a)$ of Definition \ref{Bsolution}
taking $\psi\equiv 1$, by Proposition \ref{property3}, it follows that the 
Right $\lambda$-integral $(R)\smallint G\,d_{\lambda}f$ is defined on $[a,b]$ 
with respect to $D_{\lambda,p}^{-}(f)$, and (\ref{10LIE}) holds.

\begin{thm}\label{BLIE}
Let $\lambda\in\Lambda [a,b]$, $1\leq p<2$ and let
$f$ be a $(\lambda,p)$-decomposable function on $[a,b]$.
If $D_{\lambda,p}^{-}(f)\not =\emptyset$
then in the class {\rm $\dual ({\cal W}_p)\cap Q_{\lambda}[a,b]$}
the Dol\'eans exponential ${\cal E}_b(f)$ is the unique solution 
of the augmented backward linear 
Right $\lambda$-integral equation {\rm (\ref{23LIE})}.
\end{thm}

\begin{proof}
Let $(g,h)\in D_{\lambda,p}^{-}(f)$.
By Corollary \ref{qvclassQ}, $f$ has the quadratic $\lambda$-variation,
and therefore the Dol\'eans exponential ${\cal E}_b(f)$ is defined. 
Next we show that ${\cal E}_b(f)$ is a solution of the augmented backward
linear Right $\lambda$-integral equation (\ref{23LIE}).
By Theorem \ref{qvofEb}, the $2$-vector function $({\cal E}_b(f),f)$ 
has the quadratic $\lambda$-variation and (\ref{11LIE}) holds with 
$G={\cal E}_b(f)$.
Let $\bar f:=f-f(b)+(1/2)[g]_{\lambda}^c(b)-(1/2)[g]_{\lambda}^c$,
and $V_b:=\gamma (f;\,\cdot\,,b)$,
where $\gamma (f)$ is defined by (\ref{1bvofV}).
Then ${\cal E}_b(f)=e^{-\bar f}V_b$ and $V_b$ is pure jump function of 
bounded variation by Lemma \ref{bvofV}.
Since the composition of the function $\bar f$ with a smooth function
does not change its $p$-variation, and since ${\cal W}_p$ is a Banach 
algebra, ${\cal E}_b(f)\in\dual ({\cal W}_p)[a,b]$, and so ${\cal E}_b(f)$
satisfies condition $(b)$ of Definition \ref{Bsolution}.
To show that ${\cal E}_b(f)$ satisfies condition $(a)$, we apply 
the chain rule formula of Theorem \ref{mdchrule}.$(b)$ to the composition 
$\phi{\circ}\boldf={\cal E}_b(f)$, where the $2$-vector function
$\boldf$ is defined by
$$
\boldf :=(\bar f,V^b)=(-g(b)-h(b)+\frac{1}{2}[g]_{\lambda}^c(b),0)+(g,0)
+(h-\frac{1}{2}[g]_{\lambda}^c,V^b)\quad\mbox{and}\quad
\phi (u,v)=e^{-u}v
$$
for $u, v\in\RR$.
Let $\boldg:= (g,0)$ and $\boldh :=(h-(1/2)[g]_{\lambda}^c,V_b)$.
Then $\boldf$ is $(\lambda,p)$-decomposable by $(\boldg,\boldh)$.
Since all discontinuity points of $V_b$ are the same as 
discontinuity points of $f$, the left
discontinuity points of $\boldf$ are accessible by $\lambda$.
Let $\psi\in\dual ({\cal W}_p)[a,b]$.
By statement $(b)$ of Theorem \ref{mdchrule}, the Right Cauchy
$\lambda$-integral
$$
(RC)\int\big\langle\big (\psi,-\psi(\nabla\phi{\circ}\boldf)\big ),
d_{\lambda}\big (\phi{\circ}\boldf,\boldg\big )\big\rangle
=(RC)\int\big\langle\big (\psi,\psi {\cal E}_b(f)\big ),
d_{\lambda}\big ({\cal E}_b(f),g\big )\big\rangle 
$$
is defined on $[a,b]$,
where $\nabla\phi (u,v)=(\frac{\partial\phi}{\partial u},
\frac{\partial\phi}{\partial v})(u,v)=(-e^{-u}v,e^{-u})$ for
$u,v\in\RR$, and the equalities
\begin{eqnarray}
\lefteqn{(RC)\int_y^b\big\langle\big (\psi,\psi {\cal E}_b(f)\big ),
d_{\lambda}\big ({\cal E}_b(f),g\big )\big\rangle}\nonumber\\[2mm]
&=&-(RY)\int_y^b\psi {\cal E}_b(f)\,d(h-(1/2)[g]^c)
+(RY)\int_y^b\psi e^{-\bar f}\,dV^b
-\frac{1}{2}(RS)\int_y^b\psi {\cal E}_b(f)\,d[g]^c\nonumber\\[2mm]
&+&\sum_{(y,b]}\psi\Big\{\Delta^{-}{\cal E}_b(f)+{\cal E}_b(f)
\Delta^{-}f-e^{-\bar f}\Delta^{-}V^b\Big\}\nonumber\\[2mm]
&+&\sum_{[y,b)}\psi_{+}\Big\{\Delta^{+}{\cal E}_b(f)+{\cal E}_b(f)_{+}
\Delta^{+}f-e^{-\bar f_{+}}\Delta^{+}V^b\Big\}
=(RY)\int_y^b\psi {\cal E}_b(f)\,dh\nonumber\\[2mm]
&+&(RY)\int_y^b\psi e^{-\bar f}\,dV^b
-\Big\{\sum_{(y,b]}\psi e^{-\bar f}\Delta^{-}V^b+\sum_{[y,b)}
\psi_{+} e^{-\bar f_{+}}\Delta^{+}V^b\Big\}\label{15LIE}
\end{eqnarray}
hold for $a\leq y\leq b$.
The last equality follows by Lemma \ref{bv} and relations
(\ref{-Eb}), (\ref{+Eb}).
Let $\lambda=\{\lambda_m\colon\,m\geq 1\}$ with
$\lambda_m=\{x_i^m\colon\,i=0,\dots,n(m)\}$.
Since $V_b$ is pure jump function of bounded variation, and since 
$\{x\in (a,b)\colon\,\Delta^{-}V_b(x)\not =0\}\subset\cup_m\lambda_m$
by (\ref{8bvofV}), by Theorem \ref{variation}.\ref{approximation}, 
we have the equalities
\begin{eqnarray*}
(RY)\int_y^b\psi e^{-\bar f}\,dV^b&=&\lim_{m\to\infty}\sum_{i=1}^{n(m)}
\psi (x_{i}^m)e^{-\bar f(x_{i}^m)}\big [V^b(x_i^m\vee y)
-V^b(x_{i-1}^m\vee y)\big ]\\[2mm]
&=&\sum_{(y,b]}\psi e^{-\bar f}\Delta^{-}V^b
+\sum_{[y,b)}\psi_{+} e^{-\bar f_{+}}\Delta^{+}V^b
\end{eqnarray*}
for  $a\leq y\leq b$.
Thus the last two terms in (\ref{15LIE}) give zero to the equality,
and hence the Dol\'eans exponential ${\cal E}_b(f)$ solves the augmented
backward linear Right $\lambda$-integral equation (\ref{23LIE}).

Finally, suppose that $G\in\dual ({\cal W}_p)\cap Q_{\lambda}[a,b]$
solves the augmented
backward linear Right $\lambda$-integral equation (\ref{23LIE}). 
We then show that  $G={\cal E}_b(f)$.
To this aim let $V:=e^{\bar f}G$, where 
$\bar f:=f-f(b)+(1/2)[g]_{\lambda}^c(b)-(1/2)[g]_{\lambda}^c$.
It is enough to prove that $V$ is equal to the pure jump
function $V_b=\gamma (f;\,\cdot,b)$.
Let 
\beq\label{24LIE}
\boldf:=(\bar f,G)=(C-f(b)+\frac{1}{2}[g]_{\lambda}^c(b),0)+(g,G)+
(h-\frac{1}{2}[g]_{\lambda}^c,0)
\quad\mbox{and}\quad \phi (u,v):=e^{u}v
\eeq
for $u, v\in\RR$.
Then $V=\phi{\circ}\boldf$.
We wish to apply a chain rule to show that $V$ is the unique solution
to the linear Right Young integral equation with respect to a  
pure jump function of bounded variation.
First, we claim that the $2$-vector function $\boldg :=(g,G)$ has 
the quadratic $\lambda$-variation.
By condition $(b)$ of Definition \ref{Bsolution}, $(f,G)$ has the quadratic 
$\lambda$-variation.
To apply Lemma \ref{QVandqc} we have to check its assumption
(\ref{3QVandqc}).
By condition $(a)$ of Definition \ref{Bsolution} with $\psi\equiv 1$,
and by Proposition \ref{property3}, the Right $\lambda$-integral 
$(R)\smallint G\,d_{\lambda}f$ is defined on $[a,b]$ 
with respect to $D_{\lambda,p}^{-}(f)$, and (\ref{10LIE}) holds.
By Definition \ref{RCint} of the Right Cauchy $\lambda$-integral in 
conjunction with Proposition \ref{variation}.\ref{LYjumps} concerning 
jumps of the indefinite Right Young integral, we have that
\beq\label{16LIE}
\Delta^{-}G(y)=-\big (G\Delta^{-}f\big )(y)
\qquad\mbox{and}\qquad
\Delta^{+}G(x)=-\big (G\Delta^{+}f\big )(x)
\eeq
for $a\leq y<x\leq b$.
In particular, it then follows that (\ref{3QVandqc}) holds.
Therefore by Lemma \ref{QVandqc}, the $2$-vector function
$\boldg=(g,G)$ has the quadratic $\lambda$-covariation, and thereby
has the quadratic $\lambda$-variation.
Since the pair $(\boldg,\boldh)$ is $(\lambda,p)$-dual and
(\ref{24LIE}) holds, $\boldf$ is $(\lambda,p)$-decomposable by  
$(\boldg,\boldh)$.
By (\ref{16LIE}), it follows that $(\boldg,\boldh)\in D_{\lambda,p}^{-}(\boldf)$.
Therefore all assumptions of Theorem \ref{mdchrule}.$(b)$ are satisfied
by $\boldf$ and $\phi$.
By the first part of statement $(b)$ of Theorem \ref{mdchrule} with 
$\psi\equiv 1$, the Right Cauchy $\lambda$-integrals
$$
(RC)\int\big\langle(1,-\nabla\phi{\circ}\boldf),d_{\lambda}
(\phi{\circ}\boldf,\boldg)\big\rangle
=\mu(\phi{\circ}\boldf)-
(RC)\int\langle\nabla\phi{\circ}\boldf,d_{\lambda}\boldg
\rangle
$$
are defined on $[a,b]$,
where $\nabla\phi (u,v)=(\frac{\partial\phi}{\partial u},
\frac{\partial\phi}{\partial v})(u,v)=(e^{u}v,e^{u})$ for
$u,v\in\RR$, and $\mu$ is defined by (\ref{muf}).
Further, by condition $(a)$ of Definition \ref{Bsolution} 
with $\psi=\exp\{\bar f\}$, we have the equality
\beq\label{17LIE}
(RC)\int\langle\nabla\phi{\circ}\boldf,d_{\lambda}\boldg\rangle
\equiv (RC)\int\Big\langle \Big (e^{\bar f}G,e^{\bar f}\Big ),
d_{\lambda}(g,G)\Big\rangle=-(RY)\int e^{\bar f}G\,dh.
\eeq
By condition $(b)$ of Definition \ref{Bsolution}, 
by the second part of Lemma \ref{QVandqc} and by (\ref{16LIE}),
we have
$$
(RY)\int_a^yG\,d[f]_{\lambda}=
-[f,G]_{\lambda}(y)=-[g,G]_{\lambda}^c(y)+\sum_{(a,y]}G\big\{
\Delta^{-}f\big\}^2+\sum_{[a,y)}G_{+}\big\{\Delta^{+}f\big\}^2
$$
for $a\leq y\leq b$.
On the other hand, by the second part of Lemma \ref{qv},
and since $[f]_{\lambda}^c=[g]_{\lambda}^c$, we have
$$
(RY)\int_a^yG\,d[f]_{\lambda}=(RS)\int_a^yG\,d[g]_{\lambda}^c
+\sum_{(a,y]}G\big\{\Delta^{-}f\big\}^2
+\sum_{[a,y)}G_{+}\big\{\Delta^{+}f\big\}^2
$$
for $a\leq y\leq b$.
Therefore $[g,G]_{\lambda}^c(y)=(RS)\smallint_a^yG\,d[g]_{\lambda}^c$
for $a\leq y\leq b$.
This in conjunction with the substitution theorem for the Riemann-Stieltjes
integral, yields the equality
\beq\label{18LIE}
(RS)\int_y^be^{\bar f}\,d[g,G]_{\lambda}^c
=-(RS)\int_y^be^{\bar f}G\,d[g]_{\lambda}^c
\eeq
for $a\leq y\leq b$.
By the second part of statement $(b)$ of Theorem \ref{mdchrule} 
with $\psi\equiv 1$, we then have
\begin{eqnarray}
V(b)-V(y)&=&(RC)\int_y^b\langle\nabla\phi{\circ}\boldf,
d_{\lambda}\boldg\rangle+(RY)\int_y^be^{\bar f}G\,
d\big (h-(1/2)[g]_{\lambda}^c\big )\nonumber\\[2mm]
& &-\frac{1}{2}(RS)\int_y^be^{\bar f}G\,d[g]_{\lambda}^c
-(RS)\int_y^be^{\bar f}d[g,G]_{\lambda}^c\nonumber\\[2mm]
& &+\sum_{(y,b]}\Big\{\Delta^{-}V-e^{\bar f}\Delta^{-}G
-e^{\bar f}G\Delta^{-}f\Big\}\nonumber\\[2mm]
& &+\sum_{[y,b)}\Big\{\Delta^{+}V-\Big (e^{\bar f}\Big )_{+}\Delta^{+}G
-\Big (e^{\bar f}G\Big )_{+}\Delta^{+}f\Big\}\nonumber\\[2mm]
\mbox{by Lemma \ref{bv}}\quad 
&=&(RC)\int_y^b\langle\nabla\phi{\circ}\boldf,d_{\lambda}\boldg\rangle
+(RY)\int_y^be^{\bar f}G\,dh\nonumber\\[2mm]
\mbox{by (\ref{18LIE})}\quad
& &-\frac{1}{2}(RS)\int_y^be^{\bar f}G\,d[g]_{\lambda}^c
-\frac{1}{2}(RS)\int_y^be^{\bar f}G\,d[g]_{\lambda}^c
+(RS)\int_y^be^{\bar f}Gd[g]_{\lambda}^c\nonumber\\[2mm]
\mbox{by (\ref{16LIE})}\quad
& &+\sum_{(y,b]}\Big\{\Delta^{-}V+e^{\bar f}G\Delta^{-}f
-e^{\bar f}G\Delta^{-}f\Big\}\nonumber\\[2mm]
& &+\sum_{[y,b)}\Big\{\Delta^{+}V+\Big (e^{\bar f}\Big )_{+}G_{+}\Delta^{+}f
-\Big (e^{\bar f}\Big )_{+}G_{+}\Delta^{+}f\Big\}\nonumber\\[2mm]
\mbox{by (\ref{17LIE})}\quad
&=&\sum_{(y,b]}\Delta^{-}V+\sum_{[y,b)}\Delta^{+}V\nonumber
\end{eqnarray}
for $a\leq y< b$.
Therefore $V$ is pure jump function of bounded variation
satisfying the relation
\beq\label{19LIE}
V(y)=1-\sum_{(y,b]}\Delta^{-}V-\sum_{[y,b)}\Delta^{+}V
\eeq
for $a\leq y<b$.
Further by (\ref{16LIE}), we have
$$
\Delta^{-}V=Ge^{\bar f}-\big (G-\Delta^{-}G\big )e^{\bar f}
e^{-\Delta^{-}f}
=V\Big\{1-(1+\Delta^{-}f)e^{-\Delta^{-}f}\Big\}=-V\Delta^{-}A
$$
and
$$
\Delta^{+}V=\Big (Ge^{\bar f}\Big )_{+}-\big (G_{+}-\Delta^{+}G
\big )\Big (e^{\bar f}\Big )_{+}e^{-\Delta^{+}f}
=V_{+}\Big\{1-(1+\Delta^{+}f)e^{-\Delta^{+}f}\Big\}
=-V_{+}\Delta^{+}A,
$$
where $A$ is the pure jump function of bounded variation defined
by (\ref{20LIE}).
By Theorem \ref{variation}.\ref{approximation}, we have
\begin{eqnarray*}
(RY)\int_y^bV\,dA&=&\lim_{m\to\infty}S_{RC}(V,A;\lambda_m\Cap [y,b])
=\sum_{(y,b]}V\Delta^{-}A+\sum_{[y,b)}V_{+}\Delta^{+}A\\[2mm]
&=&-\sum_{(y,b]}\Delta^{-}V-\sum_{[y,b)}\Delta^{+}V
\end{eqnarray*}
for $a\leq y\leq b$.
Thus by (\ref{19LIE}), the function $V$ is a solution to the linear
Right Young integral equation
$$
V(y)=1+(RY)\int_y^bV\,dA,\qquad a\leq y\leq b.
$$
By Theorem 5.22 of Dudley and Norvai\v sa \cite[Part II]{DNa}, the solution is unique 
in the class of functions of bounded variation.
Therefore
$$
V(y)=e^{A(b)-A(y)}\prod_{[y,b]}\big (1+\Delta A\big )e^{-\Delta A}
=\prod_{[y,b]}\big (1+\Delta f\big )e^{-\Delta f},
$$
for $a\leq y\leq b$. 
That is $V=V_b=\gamma (f;\,\cdot\,,b])$, and hence $G={\cal E}_b(f)$.
The proof of Theorem \ref{BLIE} is complete.
\qed\end{proof}

\section{The evolution representation}\label{e-representation}

Now we are prepared to prove the main result of the paper as formulated
by Theorem \ref{introduction}.\ref{evolution} in the introduction.
Recall that a family $U=\{U(t,s)\colon\,a\leq s\leq t\leq b\}$ of 
real numbers $U(t,s)$ is called an {\em evolution} on $[a,b]$ if 
$$
\left\{\begin{array}{ll}
U(r,t)\,U(t,s)=U(r,s) &\mbox{for $a\leq s\leq t\leq r\leq b$}\\
U(s,s)=1&\mbox{for $a\leq s\leq b$.}
\end{array}\right.
$$

\begin{defn}\label{regul-evol1}
{\rm An evolution $U$ on $[a,b]$ is called {\em regulated}
if there exists a multiplicative upper continuous bounded
and nondegenerate function $\pi=\pi_U$ on the extended simplex 
$S\lei a,b\rei $ such that $U(t,s)=\pi (s,t)$ for all
$(s,t)\in S[a,b]$. }
\end{defn}

Recall that for a function $f$ on $[a,b]$, $f\gg 0$ means that
there is a constant $C>0$ such that $f(t)\geq C$ for each
$t\in [a,b]$.

\begin{prop}\label{regul-evol2}
An evolution $U$ on $[a,b]$ is regulated if and only if
$U_a\equiv U(\cdot,a)$ is a regulated function on $[a,b]$ and $|U_a|\gg 0$.
\end{prop}

\begin{proof}
First, let an evolution $U$ on $[a,b]$ be regulated.
Since $U_a$ is the restriction of the right distribution function
of $\pi$ to $[a,b]$, $U_a$ is regulated on $[a,b]$ by statement
$(iv)$ of Theorem \ref{variation}.\ref{interv}.
If $\inf\{|U_a(u)|\colon\, u\in [a,b]\}=0$, then by compactness 
of $[a,b]$ it follows that $\pi (a,b)=0$, a contradiction proving that
$|U_a|\gg 0$.

Now, let $U_a$ be regulated on $[a,b]$ and let $|U_a|\gg 0 $.
For each $(u,v)\in S\lei a,b\rei$, let $\pi (u,v):=
U_a(u)/U_a(v)$.
It is easy to see that $\pi$ is a multiplicative, bounded and
nondegenerate function on $S\lei a,b\rei $.
By  statement $(iv)$ of Theorem \ref{variation}.\ref{interv}, it then
follows that $\pi$ is upper-continuous.
Since $U=\pi$ on $S[a,b]$, the conclusion follows.
The proof of the proposition is complete.
\qed\end{proof}

Similarly one can show that an evolution $U$ on $[a,b]$ is
regulated if and only if the function $U_b:=U(b,\cdot)$ is regulated
on $[a,b]$ and $|U_b|\gg 0$.

As it is mentioned in the introduction, we are interested in evolutions
satisfying property (\ref{Frechet}).
We formulate next its restricted variant.
For a function $\pi$ on $S\lei a,b\rei $ and
for a partition $\kappa=\{x_i\colon\,i=0,\dots,n\}$ of an
extended interval $\lei u,v\rei$, $a\leq u\leq v\leq b$, let
$$
S(\pi;\kappa):=\sum_{i=1}^n\big [\pi(x_{i-1},x_i)-1\big ].
$$
Recalling Notation \ref{function}.\ref{lmuv} of the trace
partition $\kuv$ we have the following:

\begin{defn}\label{regul-evol3}
{\rm Let $U$ be a regulated evolution on $[a,b]$, and let
$\lambda=\{\lambda_m\colon\,m\geq 1\}\in\Lambda [a,b]$.
We say that $U$ has the \emph{$\lambda$-generator} on $[a,b]$ if there exists
an additive and upper continuous function $\rho$ on $S\lei a,b\rei $
such that for each $(u,v)\in S\lei a,b\rei $,
\beq\label{3generator}
\rho (u,v)=\lim_{m\to\infty}S(\pi_U;\lambda_m\Cap\lei u,v\rei).
\eeq
A restriction to $[a,b]$ of the right distribution function of
$\rho$, that is the function $R_{\rho}$ defined by {\rm (\ref{RandLdf})}, 
is called the \emph{$\lambda$-generator of} $U$.}
\end{defn}

Also as it is mentioned in the introduction,
the $\lambda$-generator of an evolution generalizes the notion
of the generator of a one-parameter semigroup.
Alternatively the $\lambda$-generator on $[a,b]$ can be defined solely
in terms of a regulated function on $[a,b]$ as follows:

\begin{prop}\label{generator}
Let $U$ be a regulated evolution on $[a,b]$, and let
$\lambda=\{\lambda_m\colon\,m\geq 1\}\in\Lambda [a,b]$.
Then $U$ has the $\lambda$-generator on $[a,b]$ if and only if there exists
a regulated function $G$ on $[a,b]$ such that $G(a)=0$ and the
following hold for each $a\leq s<t\leq b${\rm :}
\beq\label{1generator}
G(t)-G(s)=\lim_{m\to\infty}S(U;\lambda_m\Cap [s,t]),
\eeq
\beq\label{2generator}
\Delta^{-}G(t)=\pi_U(t,t-)-1\qquad\mbox{and}\qquad
\Delta^{+}G(s)=\pi_U(s+,s)-1.
\eeq
If the two statements hold then $G$ is the $\lambda$-generator of $U$.
\end{prop}

\begin{proof}
First suppose that $U$ has the $\lambda$-generator, and let $G$ be
the $\lambda$-generator of $U$.
Then $G(a)=0$, $G$ is regulated by statement $(iv)$ of Theorem
\ref{variation}.\ref{interv}, 
$G$ satisfies (\ref{1generator}) and (\ref{2generator}),
and so the ``only if'' part of the conclusion holds.

For the converse, let $G$ be a regulated function on $[a,b]$ such
that $G(a)=0$, and (\ref{1generator}), (\ref{2generator}) hold.
For each $(u,v)\in S\lei a,b\rei $, let $\rho (u,v):=G(v)-G(u)$.
Then $\rho$ is additive and upper continuous function on 
$S\lei a,b\rei $ by Theorem \ref{variation}.\ref{interv}.
We prove (\ref{3generator}) only for the two cases $( u,v) =(a,t-)$ and 
$(u,v)=( a,t+)$ because the proofs for the other cases are similar.
To this aim let $\lambda_m=\{x_i^m\colon\,i=0,\dots,n(m)\}$,
$m=1,2,\dots$, and let $t\in (a,b]$.
For each $m\geq 1$, there is an index $i(t)\in\{1,\dots,n(m)-1\}$
such that $x_{i(t)}^m<t\leq x_{i(t)+1}^m$.
Then by (\ref{1generator}) and (\ref{2generator}), we have
\begin{eqnarray*}
\lim_{m\to\infty}S(\pi_U;\lambda_m\Cap \lei a,t-\rei )
&=&\lim_{m\to\infty}\Big\{S(U;\lambda_m\Cap [a,t])-\pi_U(t,x_{i(t)}^m)
+\pi_U(x_{i(t)}^m,t-)\Big\}\\
&=&G(t)-\pi_U(t,t-)+1=G(t-)=\rho (a,t-),
\end{eqnarray*}
proving (\ref{3generator}) for $( u,v) =( a,t-)\in S\lei a,b\rei$.
Now let $t\in [a,b)$.
By the definition of the trace partition, we have for each $m\geq 1$,
$$
S(\pi_U;\lambda\Cap\lei a,t+\rei )=S(\pi_U;\lambda\Cap [a,t])
+\Delta^{+}G(t),
$$
proving (\ref{3generator}) for $( u,v) =( a,t+)\in S\lei a,b\rei$.
The proof of Proposition \ref{generator} is complete.
\qed\end{proof}

Given $\lambda\in\Lambda [a,b]$ and a regulated evolution $U$, we prove the
existence of the $\lambda$-generator if $U$ is strictly positive and if
the right discontinuity points of $U_a\equiv U(\cdot,a)$ are accessible
by $\lambda$ in the sense of (\ref{right-jumps-access}), that is,
if $U_a$ is an element of the set $L_{\lambda}[a,b]$ defined by
(\ref{class-L}).
Also we show that the $\lambda$-generator $G$ of $U$ in this case
has discontinuities strictly larger than $-1$ and the right discontinuity points 
of $G$ are accessible by $\lambda$, that is, $G$ is an element of the set
$E_{\lambda}[a,b]$ defined by (\ref{class-E}).
Moreover, the $\lambda$-generator $G$ of $U$ has a Left Cauchy
$\lambda$-integral representation.

\begin{prop}\label{Elogarithm}
Let $f\in L_{\lambda}[a,b]$ for some $\lambda\in\Lambda [a,b]$.
Then the Left Cauchy $\lambda$-integral $(LC)\smallint f^{-1}
d_{\lambda}f$ is defined on $[a,b]$, and for $a\leq t\leq b$,
\begin{eqnarray}
\lefteqn{{\cal L}_{\lambda}f(t):=(LC)\int_a^tf^{-1}\,d_{\lambda}f}
\label{1Elogarithm}\\[2mm]
&=&\log\frac{f(t)}{f(a)}+\frac{1}{2}
(RS)\int_a^t\frac{d[f]_{\lambda}^c}{f^2}-\sum_{(a,t]}\Big (
\log\frac{f}{f_{-}}-\frac{\Delta^{-}f}{f_{-}}\Big )-\sum_{[a,t)}
\Big (\log\frac{f_{+}}{f}-\frac{\Delta^{+}f}{f}\Big ),\nonumber
\end{eqnarray}
where the two sums converge absolutely.
Moreover, ${\cal L}_{\lambda}f\in E_{\lambda}[a,b]$
with the bracket function
\beq\label{2Elogarithm}
\big [{\cal L}_{\lambda}f\big ]_{\lambda}(t)=(LY)\int_a^t
\frac{d[f]_{\lambda}}{f^2}\qquad a\leq t\leq b.
\eeq
\end{prop}

\begin{proof}
The Left Cauchy $\lambda$-integral $(LC)\smallint f^{-1}\,d_{\lambda}f$
is defined and (\ref{1Elogarithm}) holds by the chain rule of statement $(a)$
in Theorem \ref{chainrule} applied to 
$\phi{\circ}f$ with $\phi (u)=\log u$ for $u>0$.
By Theorem \ref{qvofint} with the same $\phi$, ${\cal L}_{\lambda}f$ 
has the quadratic $\lambda$-variation with the bracket function (\ref{2Elogarithm}).
By Definition \ref{LCint} of the Left Cauchy $\lambda$-integral
and by Proposition \ref{property1},
$1+\Delta^{-}{\cal L}_{\lambda}f\wedge\Delta^{+}{\cal L}_{\lambda}f\gg 0$
since $f\gg 0$,  $N_{(a,b)}(\Delta^{+}{\cal L}_{\lambda}f)\subset\cup\lambda$,
and so ${\cal L}_{\lambda}f\in E_{\lambda}[a,b]$.
The proof is complete.
\qed\end{proof}

Now we are ready to prove existence of the $\lambda$-generator for
a class of regulated evolutions. 

\begin{thm}\label{generator1}
Let $U$ be a regulated evolution on $[a,b]$, and let $\lambda\in\Lambda [a,b]$.
If $U_a\in L_{\lambda}[a,b]$ then on $[a,b]$, $\CAL_{\lambda}U_a$ is the 
$\lambda$-generator of $U$.
\end{thm}

\begin{proof}
Let $\lambda=\{\lambda_m\colon\,m\geq 1\}\in\Lambda [a,b]$ be such
$U_a\in L_{\lambda}[a,b]$.
The function $\CAL_{\lambda}U_a$ is defined by Proposition
\ref{Elogarithm}.
By Proposition \ref{property1}, for $a\leq s<t\leq b$, we have
\beq\label{6generator1}
1+\Delta^{-}(\CAL_{\lambda}U_a)(t)
=1+\frac{\Delta^{-}U_a(t)}{U_a(t-)}=\frac{U_a(t)}{U_a(t-)}
=\pi_U(t,t-)
\eeq
and
\beq\label{7generator1}
1+\Delta^{+}(\CAL_{\lambda}U_a)(s)
=1+\frac{\Delta^{+}U_a(s)}{U_a(s)}=\frac{U_a(s+)}{U_a(s)}
=\pi_U(s+,s).
\eeq
Thus by Proposition \ref{generator}, it is enough to prove that for 
$a\leq s<t\leq b$,
\beq\label{1generator1}
\lim_{m\to\infty}S(U;\lambda_m\Cap [s,t])=(LC)\int_s^t
U_a^{-1}\,d_{\lambda}U_a.
\eeq
Fix $a\leq s<t\leq b$, and suppose that $\lambda_m\Cap [s,t]=
\{x_i^m\colon\,i=0,\dots,n(m)\}$ for each $m\geq 1$.
For $i=1,\dots,n(m)$ and $m\geq 1$, letting
\beq\label{8generator1}
d(x_i^m,x_{i-1}^m):=U(x_i^m,x_{i-1}^m)-1-\big [\CAL_{\lambda}U_a(x_i^m)
-\CAL_{\lambda}U_a(x_{i-1}^m)\big ],
\eeq
by additivity of the Left Cauchy $\lambda$-integral, it amounts
to show that
\beq\label{10generator1}
\lim_{m\to\infty}\sum_{i=1}^{n(m)}d(x_i^m,x_{i-1}^m)=0.
\eeq
Let $\epsilon\in (0,1/2)$.
Since $U_a$ is regulated and $U_a\in L_{\lambda}[a,b]$,
there exists a partition $\zeta=\{z_j\colon\,j=0,\dots,k\}$ of 
$[s,t]$ such that
\beq\label{4generator1}
\sum_{j=1}^k\sum_{(z_{j-1},z_j)}\Big ((\Delta^{-}U_a)^2
+(\Delta^{+}U_a)^2\Big )<\epsilon
\quad\mbox{and}\quad
\max_{1\leq j\leq k}\osc\big (U_a;(z_{j-1},z_j)\big )<\epsilon C,
\eeq
where $C$ is the constant from the definition of $L_{\lambda}[a,b]$
such that $|U_a(t)|\geq C$ for $t\in [a,b]$.
Let $I$ be the set of indices $j\in\{1,\dots,k-1\}$ such that
$z_j\in\cup_m\lambda_m$, and let $m_0$ be the minimal integer such that
each intersection $\lambda_m\cap (z_{j-1},z_j)$, $j=1,\dots,k$,
contains at least two different points, and $\{z_j\colon\,j\in I\}
\subset\lambda_m$.
For each $m\geq m_0$ and $j=0,\dots,k$, let $i(j)=i_m(j)$ be the index
in $\{0,\dots,n(m)\}$ such that $x_{i(j)}^m=z_j$ if 
$z_j\in\cup_m\lambda_m$, or $x_{i(j)-1}^m<z_j<x_{i(j)}^m$
if $j\in\{1,\dots,k-1\}\setminus I$.
Then let
$$
I_1(\zeta,m):=\big\{i(j)\colon\,j=1,\dots,k\big\}\cup
\big\{1,i(j)+1\colon\,j\in I\big\}\quad\mbox{and}\quad
I_2(\zeta,m):=\{1,\dots,n(m)\}\setminus I_1(\zeta,m).
$$
By the definition of $L_{\lambda}[a,b]$,
$\Delta^{+}U_a(z_j)=0$ if $j\not\in I$.
Thus by (\ref{6generator1}), for each $j\in\{1,\dots,k\}$,  we have
$$
\lim_{m\to\infty}d(x_{i(j)}^m,x_{i(j)-1}^m)
=\pi_U(z_j,z_j-)-1-\Delta^{-}\big (\CAL_{\lambda}U_a\big )(z_j)=0.
$$
For each $j\in\{0\}\cup I$, by (\ref{7generator1}), we have
$$
\lim_{m\to\infty}d(x_{i(j)+1}^m,x_{i(j)}^m)
=\pi_U(z_j+,z_j)-1-\Delta^{+}\big (\CAL_{\lambda}U_a\big )(z_j)=0.
$$
Therefore
$$
\lim_{m\to\infty}\sum_{i\in I_1(\zeta,m)}d(x_i^m,x_{i-1}^m)=0.
$$
Since $U_a$ has the quadratic $\lambda$-variation and 
$\{t\in (a,b)\colon\,(\Delta^{+}U_a)(t)\not =0\}\subset\cup_m
\lambda_m$, by Lemmas \ref{intqv} and \ref{bv}, we have
\begin{eqnarray}
\lefteqn{\lim_{m\to\infty}\sum_{i=1}^{n(m)}
\Big (\frac{U_a(x_i^m)}{U_a(x_{i-1}^m)}-1\Big )^2
=(LY)\int_s^t\frac{d[U_a]_{\lambda}}{U_a^2}}\\[2mm]
&=&(RS)\int_s^t\frac{d[U_a]_{\lambda}^c}{U_a^2}
+\sum_{(s,t]}\Big (\frac{\Delta^{-}U_a}{(U_a)_{-}}\Big )^2
+\sum_{[s,t)}\Big (\frac{\Delta^{+}U_a}{U_a}\Big )^2.\label{3generator1}
\end{eqnarray}
For each $i\in I_2(\zeta,m)$, by the second relation in (\ref{4generator1}),
\beq\label{9generator1}
\big |1-U_a(x_i^m)/U_a(x_{i-1}^m)\big |
\leq C^{-1}\max_{1\leq j\leq k}\osc (U_a;(z_{j-1},z_j))<\epsilon.
\eeq
Thus by Taylor's theorem with remainder (\ref{Taylor3}), 
for $i\in I_2(\zeta,m)$
\beq\label{5generator1}
\log\Big (\frac{U_a(x_i^m)}{U_a(x_{i-1}^m)}\Big )
=\frac{U_a(x_i^m)}{U_a(x_{i-1}^m)}-1-\frac{1}{2}\Big (
\frac{U_a(x_i^m)}{U_a(x_{i-1}^m)}-1\Big )^2
+3\theta_i^m\Big (\frac{U_a(x_i^m)}{U_a(x_{i-1}^m)}-1\Big )^3
\eeq
with $|\theta_i^m|\leq 1$.
By (\ref{8generator1}) and (\ref{1Elogarithm}), for $i\in I_2(\zeta,m)$,
\begin{eqnarray}
d(x_i^m,x_{i-1}^m)&=&\frac{U_a(x_i^m)}{U_a(x_{i-1}^m)}
-1-\log\Big (\frac{U_a(x_i^m)}{U_a(x_{i-1}^m)}\Big )-\frac{1}{2}
(RS)\int_{x_{i-1}^m}^{x_i^m}\frac{d[U_a]_{\lambda}^c}{U_a^2}\nonumber\\[2mm]
& &+\sum_{(x_{i-1}^m,x_i^m]}\Big (\log\frac{U_a}{(U_a)_{-}}-
\frac{\Delta^{-}U_a}{(U_a)_{-}}\Big )+\sum_{[x_{i-1}^m,x_i^m)}
\Big (\log\frac{(U_a)_{+}}{U_a}-\frac{\Delta^{+}U_a}{U_a}\Big )\nonumber\\[2mm]
\mbox{by (\ref{5generator1})}\quad
&=&\frac{1}{2}\Big (\frac{U_a(x_i^m)}{U_a(x_{i-1}^m)}-1\Big )^2
-3\theta_i^m\Big (\frac{U_a(x_i^m)}{U_a(x_{i-1}^m)}-1\Big )^3
-\frac{1}{2}(LY)\int_{x_{i-1}^m}^{x_i^m}
\frac{d[U_a]_{\lambda}}{U_a^2}\nonumber\\[2mm]
\mbox{and (\ref{3generator1})}\quad
& &+\sum_{(x_{i-1}^m,x_i^m]}\Big (\log\frac{U_a}{(U_a)_{-}}-
\frac{\Delta^{-}U_a}{(U_a)_{-}}\Big )+\frac{1}{2}\sum_{(x_{i-1}^m,x_i^m]}
\Big (\frac{\Delta^{-}U_a}{(U_a)_{-}}\Big )^2\nonumber\\[2mm]
& &+\sum_{[x_{i-1}^m,x_i^m)}
\Big (\log\frac{(U_a)_{+}}{U_a}-\frac{\Delta^{+}U_a}{U_a}\Big )
+\frac{1}{2}\sum_{[x_{i-1}^m,x_i^m)}
\Big (\frac{\Delta^{+}U_a}{U_a}\Big )^2.\label{2generator1}
\end{eqnarray}
As for (\ref{8generator1}), we conclude that
$$
\lim_{m\to\infty}\sum_{i\in I_1(\zeta,m)}\Big\{
\Big (\frac{U_a(x_i^m)}{U_a(x_{i-1}^m)}-1\Big )^2
-(LY)\int_{x_{i-1}^m}^{x_i^m}\frac{d[U_a]_{\lambda}}{U_a^2}\Big\}=0.
$$
By (\ref{9generator1}), we have
$$
\limsup_{m\to\infty}\Big |\sum_{i\in I_2(\zeta,m)}
3\theta_i^m\Big (\frac{U_a(x_i^m)}{U_a(x_{i-1}^m)}-1\Big )^3\Big |
\leq 3\epsilon C^{-2}[U_a]_{\lambda}.
$$
By the Taylor series expansion (\ref{5bvofV}) and by the first relation 
in (\ref{4generator1}), we have
$$
\limsup_{m\to\infty}\Big |\sum_{i\in I_2(\zeta,m)}
\sum_{(x_{i-1}^m,x_i^m]}\Big (\log\frac{U_a}{(U_a)_{-}}-
\frac{\Delta^{-}U_a}{(U_a)_{-}}\Big )\Big |
\leq\frac{2}{C^2}\limsup_{m\to\infty}
\sum_{i\in I_2(\zeta,m)}\sum_{(x_{i-1}^m,x_i^m]}\big (\Delta^{-}
U_a\big )^2\leq\frac{2\epsilon}{C^2}.
$$
Similar bounds hold for the last three terms in (\ref{2generator1}).
All things considered, it then follows that the sum
\begin{eqnarray*}
\lefteqn{\sum_{i=1}^{n(m)}d(x_{i-1}^m,x_i^m)
=\sum_{i\in I_1(\zeta,m)}d(x_{i-1}^m,x_i^m)+\frac{1}{2}
\Big\{\sum_{i=1}^{n(m)}\Big (\frac{U_a(x_i^m)}{U_a(x_{i-1}^m)}-1\Big )^2
-(LY)\int_a^b\frac{d[U_a]_{\lambda}}{U_a^2}\Big\}}\\[2mm]
& &-\frac{1}{2}\sum_{i\in I_1(\zeta,m)}\Big\{
\Big (\frac{U_a(x_i^m)}{U_a(x_{i-1}^m)}-1\Big )^2
-(LY)\int_{x_{i-1}^m}^{x_i^m}\frac{d[U_a]_{\lambda}}{U_a^2}\Big\}
-\sum_{i\in I_2(\zeta,m)}
3\theta_i^m\Big (\frac{U_a(x_i^m)}{U_a(x_{i-1}^m)}-1\Big )^3\\[2mm]
& &+\sum_{i\in I_2(\zeta,m)}\Big\{
\sum_{(x_{i-1}^m,x_i^m]}\Big (\log\frac{U_a}{(U_a)_{-}}-
\frac{\Delta^{-}U_a}{(U_a)_{-}}\Big )+\frac{1}{2}\sum_{(x_{i-1}^m,x_i^m]}
\Big (\frac{\Delta^{-}U_a}{(U_a)_{-}}\Big )^2\Big\}\\[2mm]
& &+\sum_{i\in I_2(\zeta,m)}\Big\{\sum_{[x_{i-1}^m,x_i^m)}
\Big (\log\frac{(U_a)_{+}}{U_a}-\frac{\Delta^{+}U_a}{U_a}\Big )
+\frac{1}{2}\sum_{[x_{i-1}^m,x_i^m)}
\Big (\frac{\Delta^{+}U_a}{U_a}\Big )^2\Big\}
\end{eqnarray*}
is small as $m$ is large, proving (\ref{10generator1}).
The proof of Theorem \ref{generator1} is complete.
\qed\end{proof}

By Proposition \ref{Elogarithm},  given $\lambda\in\Lambda [a,b]$,
${\cal L}_{\lambda}$ is a mapping from $L_{\lambda}[a,b]$ to
$E_{\lambda}[a,b]$.
The following statement shows that this mapping is one-to-one
and the inverse mapping is the indefinite product $\lambda$-integral
${\cal P}_{\lambda}$.

\begin{thm}\label{ratio}
Let $\lambda\in\Lambda [a,b]$ and let $f\in L_{\lambda}[a,b]$. 
Then both ${\cal L}_{\lambda}f$ and $\prodi (1+d_{\lambda}
({\cal L}_{\lambda}f))$ are defined on $[a,b]$,
and for each $a\leq y\leq z\leq b$,
\beq\label{1ratio}
\prodi_y^z\big (1+d_{\lambda}\big ({\cal L}_{\lambda}f\big )\big )
=f(z)/f(y). 
\eeq
\end{thm}

\begin{proof}
Let $\lambda=\{\lambda_m\colon\,m\geq 1\}\in\Lambda [a,b]$ be such
that $f\in L_{\lambda}[a,b]$.
The Left Cauchy $\lambda$-integral $(LC)\smallint f^{-1}d_{\lambda}f$
is defined on $[a,b]$ by Proposition \ref{Elogarithm}.
To prove the conclusion concerning the product $\lambda$-integral
we use Proposition \ref{function}.\ref{prod-int2} with $G:=f/f(a)$.
It is clear that $G$ is a regulated function on $[a,b]$, $G(a)=1$
and $G(x)\not =0$ for $x\in (a,b]$.
Also by Proposition \ref{property1}, for $a\leq y<z\leq b$, we have 
\beq\label{7ratio}
\left\{ \begin{array}{ll}
1+\Delta^{-}\big ({\cal L}_{\lambda}f\big )(z)
=1+\Delta^{-}f(z)/f(z-)=G(z)/G(z-)\\
1+\Delta^{+}\big ({\cal L}_{\lambda}f\big )(y)
=1+\Delta^{+}f(y)/f(y)=G(y+)/G(y).
\end{array}\right.
\eeq
Therefore by Proposition \ref{function}.\ref{prod-int2}, 
it is enough to prove that for each $a\leq y<z\leq b$,
\beq\label{2ratio}
\lim_{m\to\infty}\big |\log\big\{P({\cal L}_{\lambda}f;\lambda_m\Cap 
[y,z])\big\}-\log \big\{f(z)/f(y)\big\}\big |=0.
\eeq
Fix $a\leq y<z\leq b$, and suppose that $\lambda_m\Cap [y,z]=
\{x_i^m\colon\,i=0,\dots,n(m)\}$.
By the second part of Proposition \ref{Elogarithm},  ${\cal L}_{\lambda}f$ 
has the quadratic $\lambda$-variation and
\begin{eqnarray}
\lefteqn{\lim_{m\to\infty}s_2({\cal L}_{\lambda}f;\lambda_m\Cap [y,z])
=(LY)\int_y^z\frac{d[f]_{\lambda}}{f^2}}\label{3ratio}\\[2mm]
&&=(RS)\int_y^z\frac{d[f]_{\lambda}^c}{f^2}+\sum_{(y,z]}\Big (
\frac{\Delta^{-}f}{f_{-}}\Big )^2+\sum_{[y,z)}\Big (
\frac{\Delta^{+}f}{f}\Big )^2,\label{4ratio}
\end{eqnarray}
where the second equality holds by Lemma \ref{bv}.
To begin the proof of (\ref{2ratio})  let $\epsilon >0$.
Since the sums in (\ref{1Elogarithm}) and (\ref{4ratio}) converge
unconditionally, and ${\cal L}_{\lambda}f$ is regulated,
there exists a partition $\zeta=\{z_j\colon\,j=0,\dots,k\}$ of
$[y,z]$ such that for each finite set $\mu\subset [y,z]$ disjoint from
$\zeta$,
\beq\label{5ratio}
\Big |\sum_{\mu}\Big (\log\frac{f}{f_{-}}
-\frac{\Delta^{-}f}{f_{-}}\Big )+\sum_{\mu}
\Big (\log\frac{f_{+}}{f}-\frac{\Delta^{+}f}{f}\Big )\Big |<\epsilon,
\eeq
\beq\label{6ratio}
\sum_{\mu}\Big (\frac{\Delta^{-}f}{f_{-}}\Big )^2
+\sum_{\mu}\Big (\frac{\Delta^{+}f}{f}\Big )^2 <\epsilon
\quad\mbox{and}\quad
\osc\,\big ({\cal L}_{\lambda}f;(z_{j-1},z_j)\big )<\epsilon
\wedge \frac 12.
\eeq
Let $I$ be the set of indices $j\in\{1,\dots,k-1\}$ such that 
$z_j\in\cup_m\lambda_m$, and let $m_0$ be the minimal integer such that 
$\{z_j\colon\,j\in I\}\subset\lambda_m$
and each intersection $\lambda_m\cap (z_{j-1},z_j)$, $j=1,\dots,k$,
contains at least two different points.
For each $m\geq m_0$ and $j=0,\dots,k$, let $i(j)=i_m(j)$ be the
index in $\{0,\dots,n(m)\}$ such that $x_{i(j)}^m=z_j$ if 
$z_j\in\cup_m\lambda_m$, or $x_{i(j)-1}^m<z_j<x_{i(j)}^m$ otherwise.
Then let 
$$
I_1(\zeta,m):=\{i(j)\colon\,j=0,\dots,k\}\cup\{1,i(j)+1\colon\,
j\in I\}\quad\mbox{and}\quad
I_2(\zeta,m):=\{0,\dots,n(m)\}\setminus I_1(\zeta,m).
$$
For each $m\geq m_0$, let $\Delta_i^m({\cal L}_{\lambda}f) :=
{\cal L}_{\lambda}f(x_i^m)-{\cal L}_{\lambda}f(x_{i-1}^m)
=(LC)\smallint_{x_{i-1}^m}^{x_i^m}f^{-1}d_{\lambda}f$, 
$i=1,\dots,n(m)$, and
\begin{eqnarray*}
\lefteqn{\log\big\{P({\cal L}_{\lambda}f;\lambda_m\Cap 
[y,z])\big\}-\log \big\{f(z)/f(y)\big\}}\\[2mm]
&=&\Big [\log\Big\{\prod_{i\in I_1(\zeta,m)}\big (1+\Delta_i^m\big (
{\cal L}_{\lambda}f\big )\big )\Big\}-\log\Big\{\prod_{i\in I_1(\zeta,m)}
\frac{f(x_i^m)}{f(x_{i-1}^m)}\Big\}\Big ]\\[2mm]
& &+\Big [\log\Big\{\prod_{i\in I_2(\zeta,m)}\big (1+\Delta_i^m\big (
{\cal L}_{\lambda}f\big )\big )\Big\}-\log\Big\{\prod_{i\in I_2(\zeta,m)}
\frac{f(x_i^m)}{f(x_{i-1}^m)}\Big\}\Big ]
=:D_1(\zeta,m)+D_2(\zeta,m).
\end{eqnarray*}
First consider the difference $D_1(\zeta,m)$.
Since $\Delta^{+}f(z_j)=0$ for $z_j\not\in\cup_m\lambda_m$, we have
\begin{eqnarray*}
\lim_{m\to\infty}D_1(\zeta,m)&=&
\log\Big\{\prod_{j=1}^k\big (1+\Delta^{-}({\cal L}_{\lambda}f)(z_j)
\big )\big (1+\Delta^{+}({\cal L}_{\lambda}f)(z_{j-1})\big )\Big\}\\[2mm]
\mbox{by {\rm (\ref{7ratio})}}
& &- \log\Big\{\prod_{j=1}^k\big (1+\Big (\frac{\Delta^{-}f}{f_{-}}
\Big )(z_j)\big )\big (1+\Big (\frac{\Delta^{+}f}{f}\Big )
(z_{j-1})\big )\Big\}=0.
\end{eqnarray*}
Now we examine the difference $D_2(\zeta,m)$. 
By Taylor's theorem with remainder (\ref{Taylor3}),
\begin{eqnarray*}
D_2(\zeta,m)
&=&\sum_{i\in I_2(\zeta,m)}\Big\{\Delta_i^m({\cal L}_{\lambda}f)
-\frac{1}{2}\big (\Delta_i^m({\cal L}_{\lambda}f)\big )^2
+3\theta_i\big (\Delta_i^m({\cal L}_{\lambda}f)\big )^3
-\log\Big (\frac{f(x_i^m)}{f(x_{i-1}^m)}\Big )\Big\}\\[2mm]
\mbox{by {\rm (\ref{1Elogarithm})}}
&=&\sum_{i\in I_2(\zeta,m)}\Big\{\frac{1}{2}(RS)\int_{x_{i-1}^m}^{x_i^m}
\frac{d[f]_{\lambda}^c}{f^2}-\sum_{(x_{i-1}^m,x_i^m]}\Big (\log
\frac{f}{f_{-}}-\frac{\Delta^{-}f}{f_{-}}\Big )\\[2mm]
& &-\sum_{[x_{i-1}^m,x_i^m)}
\Big (\log\frac{f_{+}}{f}-\frac{\Delta^{+}f}{f}\Big )
-\frac{1}{2}\big (\Delta_i^m({\cal L}_{\lambda}f)\big )^2
+3\theta_i\big (\Delta_i^m({\cal L}_{\lambda}f)\big )^3
\Big\}.
\end{eqnarray*}
By (\ref{4ratio}), we have
\begin{eqnarray*}
\lefteqn{\sum_{i\in I_2(\zeta,m)}\Big\{\big (\Delta_i^m({\cal L}_{\lambda}f
)\big )^2-(RS)\int_{x_{i-1}^m}^{x_i^m}\frac{d[f]_{\lambda}^c}{f^2}\Big\}}\\[2mm]
&=&\Big\{s_2({\cal L}_{\lambda}f;\lambda_m\Cap [y,z])-(LY)\int_y^z
\frac{d[f]_{\lambda}^c}{f^2}\Big\}+\sum_{i\in I_1(\zeta,m)}
(RS)\int_{x_{i-1}^m}^{x_i^m}\frac{d[f]_{\lambda}^c}{f^2}\\[2mm]
& &+\sum_{i\in I_1(\zeta,m)}\Big\{\sum_{(x_{i-1}^m,
x_i^m]}\Big (\frac{\Delta^{-}f}{f_{-}}\Big )^2+\sum_{[x_{i-1}^m,
x_i^m)}\Big (\frac{\Delta^{+}f}{f}\Big )^2
-\left (\Delta_i^m({\cal L}_{\lambda}f)\right )^2\Big\}\\[2mm]
& &+\sum_{i\in I_2(\zeta,m)}\Big\{\sum_{(x_{i-1}^m,
x_i^m]}\Big (\frac{\Delta^{-}f}{f_{-}}\Big )^2+\sum_{[x_{i-1}^m,
x_i^m)}\Big (\frac{\Delta^{+}f}{f}\Big )^2\Big\}
=:\sum_{j=1}^4A_j(m).
\end{eqnarray*}
By (\ref{3ratio}),  $A_1(m)\to 0$ as $m\to\infty$.
By continuity of $[f]_{\lambda}^c$,  $A_2(m)\to 0$ as $m\to\infty$.
Since $\Delta^{+}f(z_j)=0$ for $z_j\not\in\cup_m\lambda_m$, and
by (\ref{2LCint}), it follows that $A_3(m)\to 0$ as $m\to\infty$.
This in conjunction with (\ref{5ratio}) and (\ref{6ratio}) yields that
\begin{eqnarray*}
\lefteqn{\limsup_{m\to\infty}\big |D_2(\zeta,m)\big |
\leq 3\big [{\cal L}_{\lambda}f\big ]_{\lambda}(b)\limsup_{m\to\infty}
\max_{i\in I_2(\zeta,m)}\Big |\Delta_i^m({\cal L}_{\lambda}f)\Big |
+\frac{1}{2}\limsup_{m\to\infty}\big |A_4(m)\big |}\\[2mm]
&&+\limsup_{m\to\infty}\Big |
\sum_{i\in I_2(\zeta,m)}\Big\{\sum_{(x_{i-1}^m,x_i^m]}\Big (\log
\frac{f}{f_{-}}-\frac{\Delta^{-}f}{f_{-}}\Big )+\sum_{[x_{i-1}^m,x_i^m)}
\Big (\log\frac{f_{+}}{f}-\frac{\Delta^{+}f}{f}\Big )\Big\}\Big |\\[2mm]
&\leq &3\epsilon \big [{\cal L}_{\lambda}f\big ]_{\lambda}(b)
+\epsilon /2+\epsilon.
\end{eqnarray*}
All things considered, (\ref{2ratio}) holds.
The proof of Theorem \ref{ratio} is complete.                                     
\qed\end{proof}

The following shows that the mapping ${\cal L}_{\lambda}$ from
$L_{\lambda}[a,b]$ to $E_{\lambda}[a,b]$ is \emph{onto}.
Recall that the indefinite product $\lambda$-integral 
${\cal P}_{\lambda}$ is defined by ${\cal P}_{\lambda}f(t)=
\prodi_a^t(1+d_{\lambda}f)$ for $t\in [a,b]$.

\begin{thm}\label{Ulogarithm}
Let $\lambda\in\Lambda [a,b]$ and $f\in E_{\lambda}[a,b]$.
Then $\prodi (1+d_{\lambda}f)$ and $(LC)\smallint ({\cal P}_{\lambda}f
)^{-1}\,d_{\lambda}({\cal P}_{\lambda}f)$ are defined on $[a,b]$, and 
for all $a\leq y\leq z\leq b$,
\beq\label{1Ulog}
(LC)\int_y^z\big ({\cal P}_{\lambda}f\big )^{-1}d_{\lambda}
\big ({\cal P}_{\lambda}f\big )=f(z)-f(y).
\eeq
\end{thm}

\begin{proof}
Let $\lambda\in\Lambda [a,b]$ be such that $f\in E_{\lambda}[a,b]$.
The product $\lambda$-integral $\prodi (1+d_{\lambda}f)$ exists and
equals to $\beta_{\lambda}(f)$ by Theorem \ref{ExistLprod}.
The function $\beta_{\lambda}(f)$ on $S\loc a,b\roc$ defined by (\ref{beta1})
is nondegenerate, bounded, multiplicative and upper continuous by
Corollary \ref{beta}.
Therefore ${\cal P}_{\lambda}(f)$ exists and equals to the forward
Dol\'eans exponential ${\cal E}_{\lambda}(f):={\cal E}_{\lambda,a}(f)$.
By Theorem \ref{qvofEa}, ${\cal E}_{\lambda}(f)$ has the quadratic
$\lambda$-variation, and so ${\cal E}_{\lambda}(f)\in L_{\lambda}[a,b]$
since $f\in E_{\lambda}[a,b]$.
Therefore the Left Cauchy $\lambda$-integral in (\ref{1Ulog}) exists by 
Proposition \ref{Elogarithm}.
To prove (\ref{1Ulog})  notice that for $a\leq s<t\leq b$,
$$
\Delta^{-}{\cal L}_{\lambda}\big ({\cal P}_{\lambda}(f)\big )(t)
=\big ({\cal P}_{\lambda}(f)(t-)\big )^{-1}\Delta^{-}{\cal P}_{\lambda}
(f)(t)=\Delta^{-}f(t)
$$
and
$$
\Delta^{+}{\cal L}_{\lambda}\big ({\cal P}_{\lambda}(f)\big )(s)
=\big ({\cal P}_{\lambda}(f)(s)\big )^{-1}\Delta^{+}{\cal P}_{\lambda}
(f)(s)=\Delta^{+}f(s).
$$
Hence by Proposition \ref{property1}, it is enough to prove 
(\ref{1Ulog}) for ordinary points $a\leq y<z\leq b$.
Let $a\leq y<z\leq b$, $\lambda=\{\lambda_m\colon\,
m\geq 1\}$, and for each $m\geq 1$, let $\lambda_m\Cap [y,z]=\{x_i^m
\colon\,i=0,\dots,n(m)\}$.
Thus letting $\Delta_i^m\chi:=\chi (x_i^m)-\chi (x_{i-1}^m)$, it is enough to prove that
\beq\label{3Ulog}
\lim_{m\to\infty}\sum_{i=1}^{n(m)}\Big\{\beta_{\lambda}(f;x_{i-1}^m,
x_i^m)-1-\Delta_i^mf\Big\}=0.
\eeq
Let $\epsilon\in (0,1/2)$.
There exists a partition $\zeta=\{z_j\colon\,j=0,\dots,k\}$ of $[y,z]$
such that for each $j=1,\dots,k$,
\beq\label{2Ulog}
\osc (f;(z_{j-1},z_j))<\epsilon,\quad\osc ([f]_{\lambda}^c;[z_{j-1},z_j])
<\epsilon\quad\mbox{and}\quad
\sum_{j=1}^k\sum_{(z_{j-1},z_j)}(\Delta f)^2<\epsilon.
\eeq
Let $I$ be the set of indices $j\in\{1,\dots,k-1\}$ such that 
$z_j\in\cup_m\lambda_m$, and let $m_0$ be the minimal integer such that 
$\{z_j\colon\,j\in I\}\subset\lambda_m$
and each intersection $\lambda_m\cap (z_{j-1},z_j)$, $j=1,\dots,k$,
contains at least two different points.
For each $m\geq m_0$ and $j=0,\dots,k$, let $i(j)=i_m(j)$ be the
index in $\{0,\dots,n(m)\}$ such that $x_{i(j)}^m=z_j$ if 
$z_j\in\cup_m\lambda_m$, or $x_{i(j)-1}^m<z_j<x_{i(j)}^m$ otherwise.
Then let 
$$
I_1(\zeta,m):=\{i(j)\colon\,j=0,\dots,k\}\cup\{1,i(j)+1\colon\,
j\in I\}\quad\mbox{and}\quad
I_2(\zeta,m):=\{0,\dots,n(m)\}\setminus I_1(\zeta,m).
$$
Since $\Delta^{+}f(z_j)=0$ for $z_j\not\in\cup_m\lambda_m$,
by (\ref{beta2}), we have
$$
\lim_{m\to\infty}\sum_{i\in I_1(\zeta,m)}\Big\{\beta_{\lambda}(f;x_{i-1}^m,
x_i^m)-1-\Delta_i^mf\Big\}=0.
$$
For each $i\in\{1,\dots,n(m)\}$ and $m\geq m_0$, we have
\begin{eqnarray*}
\beta_{\lambda}(f;x_{i-1}^m,x_i^m)-1-\Delta_i^mf
&=&\exp\Big\{\Delta_i^mf-\frac{1}{2}\Delta_i^m[f]_{\lambda}^c\Big\}
\Big [\gamma (f;x_{i-1}^m,x_i^m)-1\Big ]\\[2mm]
&&+\exp\Big\{\Delta_i^mf-\frac{1}{2}\Delta_i^m[f]_{\lambda}^c\Big\}
-1-\Delta_i^mf.
\end{eqnarray*}
By (\ref{1gamma}) and (\ref{2Ulog}), there exists a finite constant $C=C(f)$ 
such that for each $m\geq m_0$,
$$
\Big |\sum_{i\in I_2(\zeta,m)}
\exp\Big\{\Delta_i^mf-\frac{1}{2}\Delta_i^m[f]_{\lambda}^c\Big\}
\Big [\gamma (f;x_{i-1}^m,x_i^m)-1\Big ]\Big |
\leq C\sum_{i\in I_2(\zeta,m)}\sum_{[x_{i-1}^m,x_i^m]}(\Delta f)^2
\leq \epsilon C.
$$
By Taylor's theorem for the exponential function, for each
$i\in\{1,\dots,n(m)\}$ and $m\geq m_0$, there is a  $\theta_i^m$,
uniformly bounded in $i$ and $m$, and such that
$$
d_i^m:=
\exp\Big\{\Delta_i^mf-\frac{1}{2}\Delta_i^m[f]_{\lambda}^c\Big\}
-1-\Delta_i^mf=\frac{1}{2}\Big\{\big (\Delta_i^mf\big )^2
-\Delta_i^m[f]_{\lambda}^c\Big\}+r_i^m,
$$
where
$$
r_i^m:=-\frac{1}{2}\Delta_i^mf\Delta_i^m[f]_{\lambda}^c+\frac{1}{8}
\big (\Delta_i^m[f]_{\lambda}^c\big )^2+\theta_i^m\max\big\{|
\Delta_i^mf|,\Delta_i^m[f]_{\lambda}^c\big\}\Big\{(\Delta_i^mf)^2
+\Delta_i^m[f]_{\lambda}^c\Big \}.
$$
For each $m\geq m_0$, we have
\begin{eqnarray*}
\lefteqn{\sum_{i\in I_2(\zeta,m)}\Big\{\big (\Delta_i^mf\big )^2
-\Delta_i^m[f]_{\lambda}^c\Big\}}\\[2mm]
&=&\Big\{s_2(f;\lambda_m\Cap [y,z])-\alpha_{\lambda}(f;y,z)\Big\}
+\sum_{i\in I_1(\zeta,m)}\Delta_i^m[f]_{\lambda}^c\\[2mm]
& &+\sum_{i\in I_1(\zeta,m)}\Big\{\sum_{(x_{i-1}^m,
x_i^m]}\big (\Delta^{-}f\big )^2+\sum_{[x_{i-1}^m,
x_i^m)}\big (\Delta^{+}f\big )^2-\big (\Delta_i^mf\big )^2\Big\}\\[2mm]
& &+\sum_{i\in I_2(\zeta,m)}\Big\{\sum_{(x_{i-1}^m,
x_i^m]}\big (\Delta^{-}f\big )^2+\sum_{[x_{i-1}^m,
x_i^m)}\big (\Delta^{+}f\big )^2\Big\}
=:\sum_{j=1}^4A_j(m).
\end{eqnarray*}
Since $f$ has the quadratic $\lambda$-variation,  $A_1(m)\to 0$ as $m\to\infty$.
By continuity of $[f]_{\lambda}^c$,  $A_2(m)\to 0$ as $m\to\infty$.
Since $\Delta^{+}f(z_j)=0$ for $z_j\not\in\cup_m\lambda_m$, 
it follows that $A_3(m)\to 0$ as $m\to\infty$.
Therefore by (\ref{2Ulog}), there exists a finite constant $C=C(f)$ such that 
$$
\limsup_{m\to\infty}\Big |\sum_{i\in I_2(\zeta,m)}d_i^m
\Big |\leq\limsup_{m\to\infty}\Big\{A_4(m)
+\sum_{i\in I_2(\zeta,m)}|r_i^m|\Big\}\leq\epsilon+\epsilon C.
$$
Since $\epsilon\in (0,1/2)$ is arbitrary, (\ref{3Ulog}) holds. 
The proof of Theorem \ref{Ulogarithm} is complete.
\qed\end{proof}

Summing up we have:

\begin{cor}
For $\lambda\in\Lambda [a,b]$, the mapping ${\cal L}_{\lambda}$ from 
$L_{\lambda}[a,b]$ to $E_{\lambda}[a,b]$ is one-to-one and onto,
with the inverse ${\cal P}_{\lambda}$.
\end{cor}

\section{Non-existence of $\lambda$-integrals}\label{nonex}

In this section  for a sequence $\lambda$ of dyadic partitions of 
an interval $[0,1]$, we show that the Left Cauchy and Right Cauchy
$\lambda$-integrals do not exist for almost all pairs of sample
functions of two jointly Gaussian but dependent Brownian motions.

\paragraph*{Integration with respect to itself of a Brownian motion.}
Let $B=\{B(t)\colon\,t\geq 0\}$ be a standard Brownian motion.
It is well-known that the integral
\beq\label{1BdB}
\int_0^TB(t)\,dB(t),\qquad 0<T<\infty,
\eeq
does not exist in the Riemann-Stieltjes sense.
To see this notice that the difference of the two Riemann-Stieltjes sums
evaluated at the same partition $\kappa=\{0=t_0<t_1<\cdots <t_n=T\}$
may be written as follows:
$$
\sum_{i=1}^nB(t_i)\big [B(t_i)-B(t_{i-1})\big ]
-\sum_{i=1}^nB(t_{i-1})\big [B(t_i)-B(t_{i-1})\big ]=s_2(B;\kappa).
$$
It follows from Th\'eor\`eme 9 of L\'evy \cite[p.\ 516]{PL40} that
with probability $1$,
$$
\limsup_{\epsilon\downarrow 0}\big\{s_2(B;\kappa)\colon\,
\kappa\in\Xi [0,T],\,\,|\kappa|\leq\epsilon\}=+\infty,
$$
where $|\kappa|$ denotes the mesh of $\kappa\in\Xi [0,T]$.
Thus the integral (\ref{1BdB}) does not exist in the Riemann-Stieltjes sense.

Let $\lambda=\{\lambda_m\colon\,m\geq 1\}\in\Lambda [0,T]$ and
let $\lambda_m=\{t_i^m\colon\,i=0,\dots,n(m)\}$ for $m\geq 1$.
Since $B(0)=0$, for each $m\geq 1$, we have
$$
S_{LC}(B,B;\lambda_m)
=\sum_{i=1}^{n(m)}B(t_{i-1}^m)\big [B(t_i^m)-B(t_{i-1}^m)\big ]
=\frac{1}{2}\big [B(T)^2-s_2(B;\lambda_m)\big ].
$$
Moreover, by Th\'eor\`eme 5 of L\'evy \cite[p.\ 510]{PL40}, the
limit $\lim_{m\to\infty}s_2(B;\lambda_m)=T$ exists with probability $1$.
Thus with probability $1$, there exists the limit
$$
(LC)\int_0^TB\,d_{\lambda}B=\lim_{m\to\infty}\sum_{i=1}^{n(m)}
B(t_{i-1}^m)
\big [B(t_i^m)-B(t_{i-1}^m)\big ]=\frac{1}{2}\big [B(T)^2-T\big ].
$$
We show next that such a limit may be infinite if a pair of
Brownian motions $(B,B)$ is replaced by a pair of two different
Brownian motions.

\paragraph*{A special construction.}
Let $L$ be a linear isometry between Hilbert spaces
$L^2[0,1]$ and $L^2(\Omega,\Pr)$.
Then $B(t):=L(1_{[0,t]})$, $0\leq t\leq 1$, is the Brownian motion
stochastic process. 
The trigonometric system 
\beq\label{trig}
1, \sqrt{2}\cos (2\pi t), \sqrt{2}\sin (2\pi t),
\dots, c_k(t):=\sqrt{2}\cos (2\pi kt), s_k(t):=\sqrt {2}
\sin (2\pi kt), \dots
\eeq
is the basis of $L^2[0,1]$. 
Thus $L$ maps the trigonometric system onto a sequence
of independent identically distributed
(iid) standard normal random variables $\xi_0$, $\xi_1$, $\eta_1$,
$\dots$, $\xi_k$, $\eta_k$, $\dots$. 
Let
$$
1_{[0,t]}=a_0(t)+\sum_{k=1}^{\infty}\Big\{a_k(t)c_k(t)+b_k(t)s_k(t)
\Big\}
$$
be the Fourier series representation of the indicator function
$1_{[0,t]}$ in $L^2[0,1]$,
so that $a_0(t)=\langle 1_{[0,t]},1\rangle =t$, 
$$
a_k(t)=\langle 1_{[0,t]},c_k\rangle=\frac{1}{2\pi k}s_k(t) 
\quad\mbox{and}\quad 
b_k(t)=\langle 1_{[0,t]},s_k\rangle =\frac{1}{2\pi k}(\sqrt{2}-c_k(t)).
$$
Then the Fourier-Wiener series representation of the
Brownian motion $B$ is given by
\begin{eqnarray}\label{bm}
B(t)&=&\xi_0t+\frac{1}{2\pi}\sum_{k=1}^{\infty}\frac{1}{k}
\Big\{\xi_ks_k(t)+\eta_k(\sqrt{2}-c_k(t))\Big\}\\
&=&\theta_0 t+\frac{\sqrt{2}}{\pi}\sum_{k=1}^{\infty}
\theta_k\frac{\sin (\pi kt)}{k},\qquad 0\leq t\leq 1\nonumber
\end{eqnarray}
where $\theta_0\equiv \xi_0$ and $\theta_k:=[\xi_k s_k(t)+\eta_k
(\sqrt{2}-c_k(t))]/2s_k(t/2)$, $k=1,2,\dots$ are iid standard
normal random variables.
Indeed, we then have
$$
EB(t)B(s)=ts+\frac{2}{\pi^2}\sum_{k=1}^{\infty}\frac{\sin (\pi kt)
\sin (\pi ks)}{k^2}=t\wedge s.
$$
The right side of (\ref{bm}) without the linear term $\xi_0t$ has the
distribution of a Brownian bridge.
The Fourier-Wiener series representation of a Brownian bridge was used 
in Paley, Wiener and Zygmund \cite{PWZ} and in Paley and Wiener 
\cite[Chapter IX]{PW} to prove several of its sample function properties.
For example, by Theorem XLIII of Paley and Wiener \cite{PW}, a subsequence of 
partial sums of (\ref{bm}) converges uniformly to a limit
almost surely which yields its almost sure continuity.
Almost sure uniform convergence of the sequence of partial sums
to its limit was proved by Hunt \cite{GAH} for a more general class
of random series.
Later on the almost sure uniform convergence of (\ref{bm}) was proved 
by several authors using different methods including general methods of 
highly developed theory of probability distributions on Banach spaces.
A new approach to proving this fact have been recently suggested
by Kwapie\'n and Woyczy\'nski \cite[Section 2.5]{KW}.
The series
$$
C(t):=\frac{1}{2\pi}\sum_{k=1}^{\infty}\frac{1}{k}
\Big\{\eta_ks_k(t)-\xi_k(\sqrt{2}-c_k(t))\Big\},\qquad 0\leq t\leq 1.
$$
is called the series conjugate to (\ref{bm}) because $C$ is the imaginary
part and $B$ is the real part of the complex series
$$
Z(t)=\zeta_0 t+\frac{1}{\pi}\sum_{k=1}^{\infty}\frac{\zeta_k}{ik}
\Big (1-e^{-2\pi ikt}\Big ),\qquad 0\leq t\leq 1,
$$
where $\zeta_k=(\xi_k+i\eta_k)/\sqrt{2}$, $k=0,1,\dots$.
As for $B$, one can check that $EC(t)C(s)=t\wedge s -ts$, that is
$C$ has the distribution of a Brownian bridge.
For $0\leq s,t\leq t$, we have
$$
EB(t)C(s)=\frac{1}{2\pi^2}\sum_{k=1}^{\infty}\frac{1}{k^2}
\Big\{\sin (2\pi k(t-s))-\big [\sin (2\pi kt)-\sin (2\pi ks)\big ]
\Big\}.
$$
Thus $B(t)$ and $C(s)$ are independent if either $t=s\in (0,1)$
or $t,s\in\{0,1\}$.

A behaviour of the Riemann-Stieltjes sums of $S_{RS}(B,C;\kappa)$
is defined by a behaviour of the Riemann-Stieltjes sums
$S_{RS}(Y,X;\kappa)$, where
$$
X(t):=B(t)-\xi_0t -\frac{1}{\sqrt{2}\pi}\sum_{k=1}^{\infty}
\frac{\eta_k}{k}=\frac{1}{2\pi}\sum_{k=1}^{\infty}\frac{1}{k}
\Big\{\xi_ks_k(t)-\eta_kc_k(t)\Big\},\qquad 0\leq t\leq 1.
$$
and
$$
Y(t):=C(t)+\frac{1}{\sqrt{2}\pi}\sum_{k=1}^{\infty}\frac{\xi_k}{k}
=\frac{1}{2\pi}\sum_{k=1}^{\infty}\frac{1}{k}\Big\{\eta_ks_k(t)
+\xi_kc_k(t)\Big\},\qquad 0\leq t\leq 1.
$$
It is easy to check that for any tagged partition $\kappa$ of $[0,1]$
\beq\label{1nonex}
S_{RS}(C,B;\kappa)=\frac{-\xi_0}{\sqrt{2}\pi}\sum_{k=1}^{\infty}
\frac{\xi_k}{k}+S_{RS}(Y,X;\kappa)
\quad\mbox{and}\quad
S_{RS}(Y,X;\kappa)=S^1(\kappa)+S^2(\kappa),
\eeq
where
\beq\label{3nonex}
S^1(\kappa):=\frac{1}{4\pi^2}\sum_{k,l=1}^{\infty}\Big\{
\frac{\xi_k\xi_l}{kl}S_{RS}(c_k,s_l;\kappa)
-\frac{\eta_k\eta_l}{kl}S_{RS}(s_k,c_l;\kappa)\Big\},
\eeq
\beq\label{4nonex}
S^2(\kappa):=\frac{1}{4\pi^2}\sum_{k,l=1}^{\infty}\Big\{
\frac{\eta_k\xi_l}{kl}S_{RS}(s_k,s_l;\kappa)
-\frac{\xi_k\eta_l}{kl}S_{RS}(c_k,c_l;\kappa)\Big\}.
\eeq
For $0\leq s,t\leq 1$, we have
$$
EX(t)Y(s)=\frac{1}{2\pi^2}\sum_{k=1}^{\infty}
\frac{\sin (2\pi k(t-s))}{k^2}.
$$
Thus the dependence structure between $X$ and $Y$ is 
the same as between $B$ and $C$.
Also, for $0\leq s,t\leq 1$, we have
$$
EX(t)X(s)=EY(t)Y(s)=\frac{1}{2\pi^2}\sum_{k=1}^{\infty}
\frac{\cos (2\pi k(t-s))}{k^2}.
$$
The random series $X$ and $Y$ are examples of the class of
Gaussian Fourier series
$$
\sum_{k=0}^{\infty}a_k\big\{\xi_k\cos(2\pi kt)+\eta_k\sin (2\pi kt)
\big\},\qquad 0\leq t\leq 1,
$$
for a sequence $\{a_k\colon\,k\geq 1\}$ of non-negative real
numbers (see Section 14 in Kahane \cite{J-PK85}).
By known results for Gaussian Fourier series, the series
$X$ and $Y$ are Fourier series of continuous functions almost surely.
Clearly, their sample functions have similar properties to sample
functions of a Brownian motion.
In particular, for a $\lambda\in \Lambda [0,1]$, almost all
sample functions of $X$ and $Y$ have the quadratic $\lambda$-variation and 
their bracket functions $[X]_{\lambda}(t)=[Y]_{\lambda}(t)=t$
for $0\leq t\leq 1$.

\begin{thm}\label{nonexistence}
Let $\lambda$ be the sequence of partitions
$\lambda_m=\{i2^{-m}\colon\,i=0,\dots,2^m\}$, $m\geq 1$, of $[0,1]$.
Then $S_{LC}(Y,X;\lambda_m)$ and $S_{RC}(Y,X;\lambda_m)$ 
do not converge in probability to finite limits as $m\to\infty$.
\end{thm}

\begin{proof}
More generally, we show that the statement of the theorem
holds for any subsequence of the sequence of partitions $\kappa_n
=\{i/n\colon\,i=0,\dots,n\}$, $n\geq 1$.
For $n\geq 1$, let 
$$
G_n:=\frac{n}{4\pi^2}\sum_{(k,l)\in C_{-}(n)\cup C_{+}(n)}
\Big\{\frac{\xi_k\xi_l}{kl}\sin \Big (\frac{2\pi l}{n}\Big )
+\frac{\eta_k\eta_l}{kl}\sin\Big (\frac{2\pi k}{n}\Big )\Big\}
$$
and
$$
F_n:=\frac{n}{4\pi^2}\sum_{(k,l)\in C(n)}\Big\{\frac{\xi_k\eta_l
+\eta_k\xi_l}{kl}
(1-\cos \Big (\frac{2\pi l}{n}\Big ))\Big\}
=\frac{n}{4\pi^2}\sum_{(k,l)\in C(n)}\Big\{\frac{\xi_k\eta_l
+\eta_k\xi_l}{kl}
(1-\cos \Big (\frac{2\pi k}{n}\Big ))\Big\},
$$
where
$C(n):=\{(k,l)\in\NN\times\NN\colon\, k+l\in n\NN\}$ and
\beq\label{spIto3}
\left\{ \begin{array}{l}
C_{-}(n):=\{(k,l)\in\NN\times\NN\colon\,\mbox{$k-l\in n\ZZ$
and $k+l\not\in n\NN$}\}\\
C_{+}(n):=\{(k,l)\in\NN\times\NN\colon\,\mbox{$k-l\not\in n\ZZ$
and $k+l\in n\NN$}\}.
\end{array}\right.
\eeq

\begin{lem}\label{representation}
For an integer $n\geq 1$,
$$
S_{LC}(Y,X;\kappa_n)=G_n+F_n\qquad \mbox{and}\qquad
S_{RC}(Y,X;\kappa_n)=G_n-F_n.
$$
\end{lem}

\begin{proof}
Let $n$ be a positive integer.
Let $\kappa_n^l=\{([(i-1)/n,i/n],(i-1)/n)\colon\,i=1,\dots,n\}$
and $\kappa_n^r=\{([(i-1)/n,i/n],i/n)\colon\,i=1,\dots,n\}$
be two tagged partitions of $[0,1]$.
Then by the second relation in (\ref{1nonex}), we have
\beq\label{2nonex}
S_{LC}(Y,X;\kappa_n)=S^1(\kappa_n^l)+S^2(\kappa_n^l)
\quad\mbox{and}\quad
S_{RC}(Y,X;\kappa_n)=S^1(\kappa_n^r)+S^2(\kappa_n^r),
\eeq
where functions $S^1$ and $S^2$ are defined by (\ref{3nonex})
and (\ref{4nonex}), respectively.
For functions $\phi$, $\psi$ defined on $[0,1]$,
let $\Sigma^n\phi :=\sum_{i=1}^n\phi (i/n)$, and $(\Delta\psi)(i/n)
:=\psi (i/n)-\psi ((i-1)/n)$ for $i=1,\dots,n$.
Notice that if $\phi$ is periodic: $\phi (0)=\phi (1)$, then
$\Sigma^n\phi=\sum_{i=1}^n\phi ((i-1)/n)$.
For each positive integer $l$, let $s_{l,n}:=s_l(1/n)=\sqrt{2}\sin (2\pi 
l/n)$ and $c_{l,n}:=\sqrt{2}-c_l(1/n)=\sqrt{2}(1-\cos (2\pi l/n))$.
We have by trigonometric identities for each $l$ and $i=1,\dots,n$
\beq\label{spIto4}
\Delta s_l\Big (\frac{i}{n}\Big )=
\left\{ \begin{array}{l}
\big [c_l\Big (\frac{i}{n}\Big )s_{l,n}
+s_l\Big (\frac{i}{n}\Big )c_{l,n}\big ]/\sqrt{2}\\
\big [c_l\Big (\frac{i-1}{n}\Big )s_{l,n}
-s_l\Big (\frac{i-1}{n}\Big )c_{l,n}\big ]/\sqrt{2}
\end{array}\right.
\eeq
and
\beq\label{spIto5}
\Delta c_l\Big (\frac{i}{n}\Big )=
\left\{ \begin{array}{l}
\big [c_l\Big (\frac{i}{n}\Big )c_{l,n}
-s_l\Big (\frac{i}{n}\Big )s_{l,n}\big ]/\sqrt{2}\\
-\big [c_l\Big (\frac{i-1}{n}\Big )c_{l,n}
+s_l\Big (\frac{i-1}{n}\Big )s_{l,n}\big ]/\sqrt{2}.
\end{array}\right.
\eeq
Recall the formulas
$$
D_n(x)=\frac{1}{2}+\sum_{i=1}^n\cos ix=\frac{\sin (n+1/2)x}{2\sin x/2}
\quad\mbox{and}\quad
\tilde D_n(x)=\sum_{i=1}^n\sin ix=\frac{\cos x/2 -\cos (n+1/2)x}{2\sin x/2}
$$
for Dirichlet's kernel $D_n$ and Dirichlet's conjugate kernel
$\tilde D_n$ respectively, valid for all $x\in (0,2\pi)$.
Then we have
\beq\label{spIto2}
\sum_{i=1}^n\cos (2\pi ji/n )= \left\{ \begin{array}{ll}
              0&\mbox{ if $j\not\in n\ZZ$}\\
              n&\mbox{ if $j\in n\ZZ$}
            \end{array}
\right. \quad\mbox{and}\quad 
\sum_{i=1}^n\sin (2\pi ji/n )=0\quad\mbox{ for any integer $j$.}
\eeq
By the preceding relation it follows that for each $k, l\in \NN$,
\beq\label{spIto7}
\Sigma^nc_ks_l=\sum_{i=1}^n\big [\sin (2\pi (l-k)i/n )+
\sin (2\pi (l+k)i/n )\big ]=0
\eeq
and
\beq\label{spIto8}
\Sigma^ns_kc_l=\sum_{i=1}^n\big [\sin (2\pi (k-l)i/n )+
\sin (2\pi (k+l)i/n )\big ]=0.
\eeq
Recalling definition (\ref{spIto3}) of $C_{-}=C_{-}(n)$ and $C_{+}
=C_{+}(n)$,
by (\ref{spIto4}) and (\ref{spIto7}), we have for each $k$ and $l$, 
$$
S_{LC}(c_k,s_l;\kappa_n)=S_{RC}(c_k,s_l;\kappa_n)
=\frac{s_{l,n}}{2}\Sigma^n\big (c_{k-l}+c_{k+l}\big )
=\left\{ \begin{array}{ll}
   n\frac{s_{l,n}}{\sqrt{2}}&\mbox{ if $(k,l)\in C_{-}\cup C_{+}$}\\
                    0&\mbox{ otherwise.}\end{array}
\right.
$$
The last equation holds by (\ref{spIto2}) and
because $s_{l,n}=0$ if $k-l\in n\ZZ$ and $k+l\in n\ZZ$.
Likewise by (\ref{spIto5}), (\ref{spIto8}) and (\ref{spIto2}), 
we have
$$
S_{LC}(s_k,c_l;\kappa_n)=S_{RC}(s_k,c_l;\kappa_n)
=-\frac{s_{l,n}}{2}\Sigma^n\big (c_{k-l}-c_{k+l}\big )
=\left\{ \begin{array}{ll}
   -n\frac{s_{k,n}}{\sqrt{2}}&\mbox{ if $(k,l)\in C_{-}\cup C_{+}$}\\
                    0&\mbox{ otherwise}\end{array}
\right.
$$
because $s_{l,n}=s_{k,n}$ if $(k,l)\in C_{-}$, and
$s_{l,n}=-s_{k,n}$ if $(k,l)\in C_{+}$.
Thus we arrive at formulas for the first terms on the right
sides of $(\ref{2nonex})$
\beq\label{5nonex}
S^1(\kappa_n^l)=S^1(\kappa_n^r)=G_n.
\eeq
Next by (\ref{spIto4}) and (\ref{spIto7}),
for each $k$ and $l$, we have
$$
S_{LC}(c_k,c_l;\kappa_n)
=-\frac{c_{l,n}}{2}\Sigma^n(c_{k-l}+c_{k+l})
=-S_{RC}(c_k,c_l;\kappa_n).
$$
Likewise, by (\ref{spIto5}) and (\ref{spIto8}),
we have
$$
S_{LC}(s_k,s_l;\kappa_n)
=-\frac{c_{l,n}}{2}\Sigma^n(c_{k-l}-c_{k+l})
=-S_{RS}(s_k,s_l;\kappa_n).
$$
Therefore by the first relation in (\ref{spIto2}), it follows that
\begin{eqnarray*}
S_2(\kappa_n^l)&=&
\frac{1}{4\pi^2}\sum_{k,l=1}^{\infty}\frac{\xi_k\eta_l
-\xi_l\eta_k}{kl}\frac{c_{l,n}}{2}\Sigma^nc_{k-l}+
\sum_{k,l=1}^{\infty}\frac{\xi_k\eta_l+\xi_l\eta_k}{kl}
\frac{c_{l,n}}{2}\Sigma^nc_{k+l}\\[2mm]
&=&\frac{n}{4\pi^2}\sum_{(k,l)\in C}\frac{\xi_k\eta_l+\xi_l\eta_k}{kl}
\frac{c_{l,n}}{\sqrt{2}}.
\end{eqnarray*}
The first series gives zero because
$(\xi_k\eta_l-\xi_l\eta_k)c_{l,n}=-(\xi_l\eta_k-\xi_k\eta_l)c_{k,n}$ 
for each pair of positive integers $(k,l)$ such that $k-l\in n\ZZ$.
Thus  $S_2(\kappa_n^l)=F_n$.
Likewise it follows that $S_2(\kappa_n^r)=-F_n$.
This in conjunction with (\ref{5nonex}) and (\ref{2nonex}) yields 
the claim of Lemma \ref{representation}.
\qed\end{proof}

Suppose that a subsequence of random variables $\{Z_{n_k}\colon\,
k\geq 1\}$ converges in probability to a random variable $Z$
as $k\to\infty$.
Since for each $k$ and any $0<M<\infty$,
$$
\Pr (\{|Z_{n_k}|>M\})\leq\Pr (\{|Z|>M/2\})+\Pr (\{|Z_{n_k}-Z|>M/2\}),
$$
it then follows that $\limsup_k\Pr (\{|Z_{n_k}|<M\})>0$ for all large
enough $M$.
Therefore to prove the theorem it is enough to prove that for each 
$M<\infty$,
\beq\label{nonex1}
\lim_{n\to\infty}\Pr (\{|Z_n|<M\})=0
\eeq
if $Z_n$ is either $S_{LC}(Y,X;\kappa_n)$ or $S_{RC}(Y,X;\kappa_n)$.
For such $Z_n$, we show that
\beq\label{nonex2}
EZ_n\geq C_1\ln n\qquad\mbox{and}\qquad
\sup_n\mbox{\rm Var}\, (Z_n)\leq C_2
\eeq
for some finite constants $C_1$ and $C_2$.
Then assuming this is true, for any $0<M, R<\infty$,
\begin{eqnarray}
\Pr (\{|Z_n|<M\})&\leq &\Pr (\{|EZ_n|<M+|Z_n-EZ_n|,\,\,
|Z_n-EZ_n|<R\})\nonumber\\
& &\quad +\Pr (\{|Z_n-EZ_n|\geq R\})\nonumber\\
&\leq &\Pr (\{|EZ_n|<M+R\})+R^{-2}\mbox{\rm Var}\, (Z_n).\nonumber
\end{eqnarray}
Taking $R$ large enough one can see that (\ref{nonex2})
yields (\ref{nonex1}).

We prove (\ref{nonex2}) only for $Z_n=S_{LC}(Y,X;\kappa_n)$,
$n\geq 1$, because a proof for the other case is the same.
By Lemma \ref{representation}, 
since $E(\xi_k\eta_l)=E(\eta_k\xi_l)=0$ for all $k, l$
and $E(\xi_k\xi_l)=E(\eta_k\eta_l)=0$ for $k\not=l$, 
the expectation of $Z_n$,
$$
EZ_n=EG_n+EF_n
=\frac{n}{2\pi^2}\sum_{(k,k)\in C_{-}(n)}\frac{1}{k^2}
\sin\Big (\frac{2\pi k}{n}\Big )
=\frac{n}{2\pi^2}\sum_{k=1}^{\infty}\frac{1}{k^2}
\sin\Big (\frac{2\pi k}{n}\Big ),
$$
because  $\sin (2\pi k/n)=0$ for $2k\in n\NN$.
For $0\leq x\leq\pi/2$, $\sin x\geq\frac{2}{\pi}\,x$.
Thus for $k\leq n/4$ we have $\sin (2\pi k/n)\geq 4k/n$.
For $n/4<k\leq n/2$, $\sin (2\pi k/n)\geq 0$,
and $\sum_{k\geq n/2}k^{-2}\leq 2/(n-2)$.
Thus for all $n\geq 4$
$$
EZ_n\geq\frac{2}{\pi^2}\Big (\sum_{1\leq k\leq n/4}\frac{1}{k}\Big )
-\frac{n}{2\pi^2}\frac{2}{n-2}\geq \frac{2}{\pi^2}(\ln n-1),
$$
and hence the first inequality in (\ref{nonex2}) holds.
For the variance of $Z_n$, we have $\sqrt{\mbox{Var}(Z_n)}\leq
4^{-1}\pi^{-2}\big (2\sqrt{V_1}+4\sqrt{V_2}+2\sqrt{V_3}+2\sqrt{V_4})$
for suitable $V_1,\dots,V_4$ to be bounded next for each $n\geq 1$.
Since $\sin^2x\leq \min\, (1,x^2)$ for all $x$, we have
$$
V_1=\mbox{Var}\Big\{n\sum_{k=1}^{\infty}\frac{\xi_k^2}{k^2}\sin
\Big (\frac{2\pi k}{n}\Big )\Big\}=n^2\sum_{k=1}^{\infty}
\frac{E(\xi_k^2-1)^2}{k^4}\sin^2\Big (\frac{2\pi k}{n}\Big )\leq
16\pi^2.
$$
Since $\sin (2\pi l/n)=0$ if $k+l\in n\NN$ and $k-l\in n\ZZ$,
we have
$$
V_2=E\Big\{n\sum_{(k,l)\in C_{-}(n),\,k>l}\frac{\xi_k\xi_l}{kl}
\sin\Big (\frac{2\pi l}{n}\Big )\Big\}^2
=n^2\sum_{m, l=1}^{\infty}\frac{\sin^2(2\pi l/n)}{(l+mn)^2l^2}
\leq n^2\sum_{m, l=1}^{\infty}\frac{1}{m^2n^2l^2}\leq 4.
$$
Denoting the integer part of an $r$ by $[r]$, likewise we get
$$
V_3=E\Big\{n\sum_{(k,l)\in C_{+}(n)}\frac{\xi_k\xi_l}{kl}\sin\Big (
\frac{2\pi l}{n}\Big )\Big\}^2=n^2\sum_{k=1}^{\infty}
\sum_{m=[k/n]+1}^{\infty}\frac{\sin^2(2\pi k/n)}
{k^2(nm-k)^2}\leq 12\pi^2+2.
$$
Since $1-\cos x=2\sin^2(x/2)\leq\min\, (2,x^2/2)$ for all $x$,
we have
$$
V_4=E\Big\{n\sum_{(k,l)\in C_{+}(n)}\frac{\xi_k\eta_l}{kl}(1-\cos\Big (
\frac{2\pi l}{n}\Big ))\Big\}^2=n^2\sum_{k=1}^{\infty}
\sum_{m=[k/n]+1}^{\infty}\frac{(1-\cos (2\pi k/n))^2}
{k^2(nm-k)^2}\leq 12\pi^4+2.
$$
Thus $\mbox{Var}(Z_n)<\pi^2$, and hence (\ref{nonex2}) holds
proving the Theorem \ref{nonexistence}.
\qed\end{proof}

\chapter{Extension of the class of semimartingales}\label{process}
\setcounter{thm}{0}

\vspace*{0.2truein}
\begin{quotation}{\footnotesize
Mathematics is endangered by a loss of unity and interaction

Richard Courant (1888-1972)}
\end{quotation}
\vspace*{0.2truein}

A stochastic process $X$ which is decomposable into a sum of
a local martingale $M$ and an adapted {\cadlag} stochastic process $\pv$ 
with locally bounded variation is called a semimartingale.
In this chapter we consider a stochastic process $X$ with the same
decomposition property except that $\pv$  is assumed to have locally 
bounded $p$-variation  for some $1\leq p<2$, so that $\pv$ may have almost
all sample functions with {\em unbounded} variation.
We call the extended  semimartingale $X$, a \nsem-semimartingale 
(see Definition \ref{p-semimartingale} below).
The main results are that an extended stochastic integral with respect to
$X$ and a continuous part of an extended quadratic variation of $X$
are {\em invariant} with respect to possibly different decompositions
$X-X_0=M+A$.
The new constructions in conjunction with the results of Chapter
\ref{variation} concerning the refinement Rimann-Stieltjes integral 
are then used to describe the unique solution to the linear extended
stochastic integral equation with respect to a \nsem-semimartingale.

\section{Stochastic processes and $p$-variation}

In this paper, a \emph{stochastic process} $X=\{X(t)\colon\,t\geq 0\}$
is a family of random variables $X(t)=X(t,\cdot)$ defined 
on a complete probability space  $(\Omega,{\cal F},\Pr)$.
Here a random variable is a real valued measurable function,
rather than an equivalence class of almost surely equal
measurable functions.
For each $\omega\in\Omega$, the function $X(\cdot)=X(\cdot,\omega)$
is called a {\em sample function} of $X$.
We say that two stochastic processes $X$ and $Y$ are {\em indistinguishable},
and write $X=Y$, if almost surely, $X(t)=Y(t)$ for each $t$.
Thus we identify two stochastic processes whose almost all sample functions
are the same.
This point of view implies that for two random variables $\xi$ and $\eta$
their equivalence means $\xi =\eta$ almost surely.
We say that two stochastic processes $X$ and $Y$ are {\em modifications}
of each other if for each $t$, $X(t)=Y(t)$ almost surely.

\paragraph*{Regulated and {\cadlag} stochastic processes.}
A stochastic process $X=\{X(t)\colon\,t\geq 0\}$ has regulated
sample functions if 
for almost all $\omega\in\Omega$, there exist the limits
$$
X(t-,\omega):=\lim_{u\uparrow t}X(u,\omega)\qquad\mbox{and}
\qquad X(s+,\omega):=\lim_{u\downarrow s}X(u,\omega)
$$
for $0\leq s<t<\infty$.
In this case we say that $X$ is a {\em regulated stochastic process}.
By Theorem 11.5 of Doob \cite[p.\ 361]{JLD}, if a separable stochastic
process $X$ is a martingale then $X$ is regulated.
In Probability Theory sample functions of stochastic processes
are often assumed to be regulated and right-continuous, that is,
{\cadlag}.
Let $X$ be a regulated stochastic process and let $\Omega_X$ be
the set of all $\omega\in\Omega$ such that $X(\cdot,\omega)$ is
regulated.
Due to completeness of the underlying probability space,
$\Omega_X\in {\cal F}$ and $\Pr (\Omega_X)=1$.
For each $t\in [0,\infty)$, let
$$
X_{+}(t):=X_{+}(t,\omega):=\left\{ \begin{array}{ll}
          X(t+,\omega) &\mbox{if $\omega\in\Omega_X$} \\
          X(t,\omega)  &\mbox{otherwise.} \end{array}
\right.
$$
Similarly define $X_{-}(t)$ for each $t\in (0,\infty)$ and
let $X_{-}(0):=X(0)$.
For each $t\in [0,\infty)$, $X_{+}(t)$ is a
random variable because it is the limit for $\omega\in\Omega_X$
of the random variables $X(r_n)$ for $r_n$ rational, $r_n\downarrow t$,
and equals $X(t)$ otherwise.
Likewise, $X(t-)$ is a random variable. 
Therefore $X_{+}=\{X_{+}(t)\colon\,t\geq 0\}$ and
$X_{-}=\{X_{-}(t)\colon\,t\geq 0\}$ are stochastic processes
on the same probability space as $X$.

Recall that $\Delta^{+}X(t):=X_{+}(t)-X(t)$ and 
$\Delta^{-}X(t):=X(t)-X_{-}(t)$.
For each $t\in [0,\infty )$, let
$$
\Omega_d(t):=\big\{\omega\in\Omega\colon\,\mbox{ either
$\Delta^{-}X(t)\not =0$ or $\Delta^{+}X(t)\not =0$}\big\}.
$$
For a regulated stochastic process $X$, 
a point $t\in [0,\infty)$ is called a point of {\em fixed
discontinuity} if $\Pr (\Omega_d(t))>0$.
If $t$ is not a point of fixed discontinuity then
\beq\label{cont}
\lim_{s\to t}X(s)=X(t)\qquad\mbox{almost surely}.
\eeq
Indeed, for $t\in (0,\infty)$ and 
$\omega\in\Omega_X\setminus\Omega_d(t)$,
we have
$$
\lim_{s\uparrow t}X(s,\omega)=X(t-,\omega)=X(t,\omega)
=X(t+,\omega)=\lim_{s\downarrow t}X(s,\omega).
$$
Since $\Pr (\Omega_X\setminus\Omega_d(t))=1$, (\ref{cont}) holds
when $t\in (0,\infty)$.
The same argument yields (\ref{cont}) when $t=0$.
If the stochastic process $X$ has no points of fixed discontinuity
then the three processes $X_{-}$, $X_{+}$ and $X$ are modifications
of each other, that is
$$
\Pr (\{X(t)=X(t-)\})=\Pr (\{X(t)=X(t+)\})=1
$$
for each $t\in [0,\infty)$.

A point $t\in [0,\infty)$ is a point of {\em stochastic
continuity} of $X$ if $X(s)\to X(t)$ in probability
as $s\to t$.

\begin{thm}\label{doob}
Let $X$ be a regulated stochastic process.
The set of points of fixed discontinuity of $X$ is at most countable
and coincides with the set of points of stochastic discontinuity
of $X$.
\end{thm}

\begin{proof}
Since the limit of a sequence convergent in probability is unique
almost surely, one can show that $X$ is not stochastically
continuous at some $t\in [0,\infty)$ if $t$ is a point
of fixed discontinuity. 
Also, $X$ is stochastically continuous at $t$ whenever
$t$ is not a point of fixed discontinuity.
This yields the second part of the claim.
The first part follows from Theorem 11.1, Ch. VII,
of Doob \cite{JLD}. 
\qed\end{proof}

By the preceding theorem, the set of fixed discontinuities of
any stochastic process with regulated sample functions is
at most countable.
However, almost every sample function of such a process
may have a non-fixed discontinuity, e.g.\
$X(t,\omega)=1_{\{t\geq \omega\}}$, where $\omega$ has a uniform distribution
in $[0,1]$. 

\paragraph*{The $p$-variation and its index.}
Suppose that almost all sample functions of a stochastic process $X$
have bounded $p$-variation on $[0,t]$ for some $0<p, t<\infty$. 
Since the set $\Xi [0,t]$ of all partitions of $[0,t]$ is
uncountable, the function 
$\omega\mapsto v_p(X(\cdot,\omega);[0,t])$ need not be measurable.
For example, let $\Pr$ be Lebesgue measure on $\Omega:=[0,1]$.
For a non-Lebesgue measurable set $A\subset [0,1]$, let 
$X=\{X(t)\colon\,0\leq t\leq 1\}$ be a stochastic process defined 
by $X(t,\omega)=0$ on $[0,1]$ except at the point $t=\omega\in A$ 
when $X(\omega,\omega)=1$.
Then $v_p(X;[0,1])=0$ if $\omega\not\in A$, and $=2$ if $\omega\in A$.
Notice that all sample functions of $X$ are regulated.

One possibility to overcome this difficulty is to understand
the supremum
$$
\sup\big\{s_p(X(\cdot,\omega);\kappa)\colon\,\kappa\in
\Xi [0,T]\big\},\quad \omega\in\Omega,
$$
as the lattice supremum in the space of equivalence classes of random variables
$L_0\!=\!L_0(\Omega,{\cal F},\Pr)$.
Recall that $L_0$ is a Dedekind complete Riesz space of countable
type (see \S IV.1 and \S VI.2 in Vulikh \cite{BZV}).
Therefore every infinite set $E$ in $L_0$ which is bounded above
has a supremum (in $L_0$), and there exists a finite or
countable subset $E'\subset E$ such that $\sup E'=\sup E$.

In Stochastic Analysis, traditionally a different route is used to deal with 
non-measurability problems.
Next we show that for a {\cadlag} stochastic process,
the $p$-variation over an interval 
is indistinguishable from the $p$-variation over a countable and 
everywhere dense set, which is always measurable.
Let $\lambda=\{\lambda_m\colon\,m\geq 1\}$ be a sequence of nested 
partitions $\lambda_m=\{t_i^m\colon\,i=0,\dots,n(m)\}$ of $[0,\infty)$
such that $\cup_m\lambda_m$ is dense in $[0,\infty)$, and let
$\Lambda [0,\infty)$ be the set of all such sequences.
For example, the sequence $\{i2^{-m}\colon\,i=0,\dots,m2^m\}$, 
$m\geq 1$, belongs to $\Lambda [0,\infty)$.
Let $f$ be a function on $[0,\infty)$, and let $0<p<\infty$. 
For each $m\geq 1$ and $t\geq 0$, let
$$
v_p(f;\lambda_m)(t):=\max\Big\{\sum_{j=1}^k\big |f(s_j\wedge t)
-f(s_{j-1}\wedge t)\big |^p\colon\,\{0,t_{n(m)}^m\}\subset\{s_j\colon\,
j=0,\dots,k\}\subset\lambda_m\Big\},
$$
which is the $p$-variation over the finite set $\{t_i^m\wedge t\colon\,
i=0,\dots,n(m)\}$, equal to the trace partition $\lambda_m\Cap [0,t]$
when $t\leq t_{n(m)}^m$ (see Notation \ref{function}.\ref{lmuv}).
Since $\lambda_m$, $m\geq 1$, are  nested partitions, the
sequence $v_p(f;\lambda_m)(t)$, $m\geq 1$, is nondecreasing
for each $t\geq 0$.
For $t\geq 0$, let
\beq\label{2loc-pv}
v_p(f)(t):=
v_p(f;\lambda)(t):=\sup_{m\geq 1}v_p(f;\lambda_m)(t)
=\lim_{m\to\infty}v_p(f;\lambda_m)(t).
\eeq
For a stochastic process $X$ and each $t\geq 0$,
$v_p(X)(t,\omega):=v_p(X(\cdot,\omega))(t)$ is possibly unbounded
but measurable function of $\omega\in\Omega$.

\begin{thm}\label{measurability}
Let $X=\{X(t)\colon\,t\geq 0\}$ be either a separable stochastic
process continuous in probability, or a {\cadlag} stochastic process,
let $\lambda\in\Lambda [0,\infty)$, and let $0<p<\infty$.
If for each integer $k\geq 1$,
\beq\label{loc-pv}
v_p(X;\lambda)(k)<\infty\qquad\mbox{almost surely,}
\eeq
then $\{v_p(X;[0,t])\colon\,t\geq 0\}$ is a stochastic process 
indistinguishable from $v_p(X)$.
Moreover, if in addition $p\geq 1$ then almost surely
for each $t>0$,
$$
v_p(X)(t)^{1/p}-v_p(X)(t-)^{1/p}\leq \big |\Delta^{-}X(t)\big |
\leq\big\{v_p(X)(t)-v_p(X)(t-)\big\}^{1/p},
$$
and for each $t\geq 0$,
$$
v_p(X)(t+)^{1/p}-v_p(X)(t)^{1/p}\leq \big |\Delta^{+}X(t)\big |
\leq\big\{v_p(X)(t+)-v_p(X)(t)\big\}^{1/p}.
$$
\end{thm}

\begin{proof}
Let $\lambda=\{\lambda_m\colon\,m\geq 1\}\in\Lambda [0,\infty)$.
First, suppose that $X$ is a separable stochastic process continuous
in probability.
Then by Theorem 2.2 of Doob \cite[Section II.2]{JLD}, $S:=\cup_m
\lambda_m$ is a separating set.
Thus there exists a $\Pr$-null set $\Omega_0$ such that if $A$ is a closed
subset of $\RR$ and $I$ is an open subset of $[0,\infty)$, then
$$
\big\{\omega\colon\,X(t,\omega)\in A,\,\,t\in S\cap I\big\}
\setminus\{\omega\colon\,X(t,\omega)\in A,\,\,t\in I\big\}
\subset\Omega_0.
$$
The following consequence of separability is what we need:
for each $\omega\not\in\Omega_0$ and any open interval $I\subset 
[0,\infty)$,
\beq\label{3loc-pv}
\overline{X(I,\omega)}=\overline{X(I\cap S,\omega)}.
\eeq
Here $\overline{X(I,\omega)}$ denotes the closure in $\RR$ of the set of 
values assumed by the sample function $X(\cdot,\omega)$ as $t\in I$,
and the other set has similar meaning.
In particular, for every $t\geq 0$ there exists a sequence
$\{s_k\colon\,k\geq 1\}\subset S$ such that $s_k\to t$ and
$X(s_k,\omega)\to X(t,\omega)$ for each $\omega\not\in\Omega_0$
(the sequence $\{s_k\}$ may depend on $\omega$).
Second, if $X$ is a {\cadlag} stochastic process then (\ref{3loc-pv})
holds for each {\cadlag} sample function $X(\cdot,\omega)$.

Let $\Omega_0$ be a $\Pr$-null set such that for each
$\omega\not\in\Omega_0$, (\ref{3loc-pv}) holds and 
$v_p(X;\lambda)(k,\omega)<\infty$ for each integer $k\geq 1$.
Let $0\leq r<t<\infty$.
If $r\in S$ then
\beq\label{4loc-pv}
v_p(X;\lambda)(r,\omega)\leq v_p(X,\lambda)(t,\omega)<\infty
\eeq
holds for each $\omega$ by definition (\ref{2loc-pv}).
If $r\not\in S$ then for each $\omega\not\in\Omega_0$ and
for each $m\geq 1$,
$$
v_p(X;\lambda_m)(r,\omega)=\lim_{s_k\to r,\,\,s_k\in S}
v_p(X;\lambda_m)(s_k,\omega)\leq v_p(X;\lambda)(t,\omega).
$$
Therefore (\ref{4loc-pv}) holds for each $\omega\not\in\Omega_0$
and any $0\leq r<t<\infty$.
Let $t>0$ and let $\kappa$ be a partition of $[0,t]$.
By (\ref{3loc-pv}), for each $\omega\not\in\Omega_0$, 
$s_p(X(\cdot,\omega);\kappa)$
can be approximated arbitrarily closely by sums 
$s_p(X(\cdot,\omega);\kappa_m)$ with $\{0,t\}\subset\kappa_m\subset
\lambda_m\Cap [0,t]$ as $m\to\infty$.
Therefore for each $\omega\not\in\Omega_0$,
$$
v_p(X(\cdot,\omega);[0,t])= v_p(X;\lambda)(t,\omega)
<\infty,
$$
proving the first part of the theorem.
The second part of the theorem follows from the inequalities
$$
v_p(f;[a,c])+v_p(f;[c,b])\leq v_p(f;[a,b])
\leq \big\{v_p(f;[a,c])^{1/p}+v_p(f;[c,b])^{1/p}\big \}^p
$$
valid for $a<c<b$ (Lemma 4.6 in \cite{DNb}),
and from the relations, 
$$
\lim_{z\uparrow x}v_p(f;[z,x])=\big |\Delta^{-}f(x)\big |^p
\quad\mbox{and}\quad
\lim_{z\downarrow y}v_p(f;[y,z])=\big |\Delta^{+}f(y)\big |^p
$$
valid for $a\leq y<x\leq b$ (Lemma 2.19 in \cite[Part II]{DNa}).
\qed\end{proof}

For a {\cadlag} stochastic process $X$ if (\ref{loc-pv}) holds
for some $\lambda\in\Lambda [0,\infty)$ then it also holds
for any other such $\lambda$, and $v_p(X)$ is a {\cadlag}
stochastic process.
It is clear from definition (\ref{2loc-pv}) that if $X$ is 
adapted then $v_p(X)$ is also adapted stochastic process 
whose almost all sample functions are nondecreasing.

\begin{defn}
{\rm Let $0<p<\infty$.
We say that a stochastic process $X$ is {\em locally of bounded
$p$-variation} if $X$ is adapted, {\cadlag} and (\ref{loc-pv}) holds
for some $\lambda\in\Lambda [0,\infty)$ and each integer $k\geq 1$.
Also, $v_p(X)$ is called the {\em $p$-variation process} for $X$.}
\end{defn}

For a function $f\colon\,[0,T]\mapsto\RR$, the $p$-variation index of $f$ 
is defined by
\beq\label{p-var-index}
\upsilon (f):=\upsilon (f;[0,T]):=\left\{
\begin{array}{ll}
  \inf\{p>0\colon\,v_p(f;[0,T])<\infty\} &\mbox{if the set is nonempty,}\\
  +\infty &\mbox{otherwise.}
\end{array} \right.
\eeq
For a stochastic process $X$ on a complete probability space 
$(\Omega,{\cal F},\Pr)$,
the $p$-variation index $\upsilon (X;[0,T])$ is defined provided
$\upsilon (X(\cdot,\omega);[0,T])$ is a constant for almost all
$\omega\in\Omega$.
We say that the \emph{local $p$-variation index $\pindex (X)=p$}
if with probability $1$, the $p$-variation index $\pindex (X;[0,T])=p$
for each $T>0$.

\paragraph*{Local martingales.}
We continue with recalling a standard framework of Stochastic
Analysis.
Let $(\Omega,{\cal F},\Pr)$ be a complete probability space,
and let $\FF=\{{\cal F}_t\colon\,0\leq t <\infty\}$ be a filtration
of sub-$\sigma$-algebras of ${\cal F}$
satisfying the usual hypotheses, that is, $\FF$ is right-continuous
and ${\cal F}_0$ contains all the $\Pr$-null sets of ${\cal F}$.
A stochastic process $X=\{X(t)\colon\,t\geq 0\}$ on 
$(\Omega,{\cal F},\Pr)$ is said to be {\em adapted} if $X(t)$ is 
${\cal F}_t$ measurable for each $t$.
An adapted, {\cadlag} stochastic process $X$ is called a {\em local
martingale} with respect to the filtration $\FF$ if there exists
a sequence of increasing $\FF$-stopping times, $\tau_n$, with
$\lim_{n\to\infty}\tau_n=+\infty$ a.s. such that $X^{\tau_n}1_{\{
\tau_n>0\}}=\{X(t\wedge\tau_n)1_{\{\tau_n>0\}}\colon\,t\geq 0\}$ 
is a uniformly integrable martingale for each $n$.
This definition of a local martingale is used in the book
by Protter \cite[p.\ 33]{PP}, which is used for references in this chapter.

The following result is due to L\'epingle \cite{DL}.
A different proof of this result was given by Pisier and Xu \cite{PandX}.

\begin{thm}\label{lepingle}
A local martingale is locally of bounded $p$-variation for each
$p>2$.
Moreover, for $p>2$ and $1\leq r<\infty$ there is a finite
constant $K_{p,r}$ such that for any martingale
$M=\{M(t)\colon\,0\leq t\leq T\}$, $T>0$, in $L^r(\Omega,{\cal F},
\Pr)$,
\beq\label{1lepingle}
E\big [v_p(M;[0,T])\big ]^{r/p}
\leq  K_{p,r}E\big [\sup_{0\leq t\leq T}|M(t)|\big ]^r.
\eeq
\end{thm}

\begin{proof}
Let $p>2$, $1\leq r<\infty$ and let $M=\{M(t)\colon\,0\leq t\leq T\}$
be a martingale in $L^r$.
Let $\lambda=\{\lambda_m\colon  m\geq 1\}$ be a nested sequence of 
partitions $\lambda_m=\{t_i^m\colon\,i=0,\dots,n(m)\}$ of $[0,T]$ such 
that $\cup_m\lambda_m$ is dense in $[0,T]$.
By Theorem 2.4 of Pisier and Xu \cite{PandX}, there is finite
constant $K_{p,r}$ such that for each $m\geq 1$,
$$
E\big [s_p(M;\lambda_m)(T)\big ]^{r/p}
\leq K_{p,r}E\big [\max_{0\leq i\leq n(m)}|M(t_{i}^m)|\big ]^r
\leq K_{p,r}E\big [\sup_{0\leq t\leq T}|M(t)|\big ]^r.
$$
Thus (\ref{1lepingle}) follows by Fatou lemma and Theorem
\ref{measurability}.
Suppose now that $M$ is a local martingale.
By the Fundamental Theorem of Local Martingales (Theorem
III.13 in \cite{PP}), $M=N+A$, where $N$ is locally bounded
martingale and $A$ is locally of bounded $1$-variation.
Thus an application of the first part of the proof yields
the conclusion of the theorem.
\qed\end{proof}

It is well known that sample functions of martingales which are 
continuous and nonconstant must have unbounded 1-variation
(see Lemma 3.2.1 of Fisk \cite{DLF}).
In fact they must have unbounded $p$-variation for every $p<2$.

\begin{thm}\label{pv-mart}
If for some $p<2$ a sample continuous local martingale is locally 
of bounded $p$-variation then almost surely its sample functions are 
constants.
\end{thm}

The proof of this statement is based on a stopping time
technique as follows.
We will say that a stochastic process $X$ is sample uniformly
continuous if for each $0<T<\infty$ there exists a $\Pr$-null set
$\Omega_0$ such that given $\epsilon >0$ one can choose $\delta>0$
such that $|X(t,\omega)-X(s,\omega)|<\epsilon$ whenever $|t-s|<\delta$,
$t, s\in [0,T]$ and $\omega\not\in\Omega_0$.
For a stopping time $\tau$, let $X^{\tau}:=\{X(t\wedge\tau)\colon\,
t\geq 0\}$.

\begin{lem}\label{fisk}
Let $X$ be a sample continuous stochastic process adapted to a 
filtration $\FF$.
For each $0<T<\infty$ there exists a sequence $\{\tau_k\colon\,
k\geq 1\}$ of increasing $\FF$ stopping times such that
$\Pr (\cup_{k=1}^{\infty}\{\tau_k=T\})=1$, and each $X^{\tau_k}$
is sample uniformly continuous and bounded by $k$.
\end{lem}

\begin{proof}
Let $0<T<\infty$.
Since almost every sample function of $X$ is uniformly continuous
on $[0,T]$, there exists a sequence $\{\delta_{kn}\colon\, k, n\geq 1\}$
of positive numbers such that, for each $k$, $\delta_{k1}>\delta_{k2}>
\dots$, and for each $k, n$,
\beq\label{1fisk}
\Pr\big (\{\sup\{|X(t)-X(s)|\colon\,|t-s|\leq\delta_{kn}\}\geq 1/k\}\big )
\leq 2^{-n-k}.
\eeq
Let $\tau_{kn}'$ be the least $t\in [0,T]$ such that
$$
\sup\big\{|X(s'')-X(s')|\colon\,|s''-s'|\leq\delta_{kn},\,\,
s'', s'\in [0,t]\big\}\geq 1/k.
$$
If there is no such $t$, let $\tau_{kn}':=T$.
Since $X$ is sample continuous each $\tau_{kn}'$ is an $\FF$ stopping time.
Let $\tau_k':=\inf_n\tau_{kn}'$.
Then for each $k$ and any real $r$, 
$\{\tau_k'<r\}=\cup_n\{\tau_{kn}'<r\}\in\FF_r$.
Since $\FF$ is right continuous, each $\tau_k'$ is an $\FF$ stopping
time.
Then each stopped process $X^{\tau_k'}$ is sample uniformly continuous.
Also, by (\ref{1fisk}), $\Pr (\{X^{\tau_k'}\not =X\})\leq 2^{-k}$.
It then follows from the Borell-Cantelli lemma that
$\Pr (\cup_k\{\tau_k'=T\})=1$.
Let $\tau_k''$ be the least $t\in [0,T]$ such that $|X(t)|>k$.
If no such $t$ exists, let $\tau_k'':=T$.
Then each $\tau_k''$ is an $\FF$ stopping time and $X^{\tau_k''}$
is bounded by $k$.
Clearly the sequence of $\FF$ stopping times 
$\tau_k:=\tau_k'\wedge \tau_k''$, $k\geq 1$, satisfies the
conclusion of the lemma.
\qed\end{proof}

Now we are ready to prove Theorem \ref{pv-mart}.

\medskip
\noindent
{\bf Proof of Theorem \ref{pv-mart}.}
Let a sample continuous local martingale $X$ be locally of bounded
$p$-variation for some $1\leq p<2$.
It is enough to prove that for any $0<T<\infty$, almost surely,
$X(t)=X(0)$ for each $0\leq t\leq T$.
Since $X$ is sample continuous, the $p$-variation process
$v_p(X)$ is also sample continuous by Theorem \ref{measurability}.
By Lemma \ref{fisk} applied to $X$ and $v_p(X)$, and since
$X$ is a local martingale, there exists a sequence $\{\tau_k\colon\,
k\geq 1\}$ of $\FF$ stopping times such that 
$\Pr (\cup_k\{\tau_k=T\})=1$, each $X^{\tau_k}$ is sample uniformly
continuous martingale bounded by $k$ and each $v_p(X)^{\tau_k}$
is bounded by $k$.
Notice that each $v_p(X)^{\tau_k}$ is the $p$-variation process
for $X^{\tau_k}$.
Let $\lambda=\{\lambda_m\colon\,m\geq 1\}$ be a sequence of nested
partitions $\lambda_m=\{t_i^m\colon\,i=0,\dots,n(m)\}$ of $[0,T]$
such that $\max_i(t_i^m-t_{i-1}^m)\to 0$ as $m\to\infty$.
Then $\cup_m\lambda_m$ is dense in $[0,T]$.
For each $k, m\geq 1$, let
$$
\epsilon_{k,m}:=\esssup_{\omega}\max_i\big |X^{\tau_k}(t_i^m)
-X^{\tau_k}(t_{i-1}^m)\big |.
$$
Let $t\in\cup_m\lambda_m$.
Then there exists an integer $m_t$ such that $t\in\lambda_m$ for each
$m\geq m_t$.
Since $X^{\tau_k}$ has orthogonal increments, for each $m\geq m_t$,
we have
$$
E\big [X^{\tau_k}(t)-X^{\tau_k}(0)\big ]^2
=Es_2\big (X^{\tau_k};\lambda_m\big )\leq E\big (\epsilon_{k,m}^{2-p}
v_p(X)^{\tau_k}(t)\big )\leq kE\big (\epsilon_{k,m}^{2-p}\big ).
$$
Since for each $k\geq 1$, $\epsilon_{k,m}\leq 2k$ and
$\epsilon_{k,m}\to 0$ almost surely as $m\to\infty$, by the
dominated convergence theorem, the left side of the preceding inequality 
is zero for each $k\geq 1$.
Therefore, almost surely, for each $k\geq 1$, $X^{\tau_k}(t)
=X^{\tau_k}(0)$ for each $t\in\cup_m\lambda_m$.
Since $X$ is sample continuous and since for almost every
$\omega\in\Omega$, there is an integer $k\geq 1$ such that
$X^{\tau_k}(t,\omega)=X(t,\omega)$ for $0\leq t\leq T$, it follows
that almost surely, $X(t)=X(0)$ for each $0\leq t\leq T$.
The proof of Theorem \ref{pv-mart} is complete.
\qed

\section{Stochastic integral and the Left Young integral}

We show here that values of the Left Young integral for pairs of
sample functions of stochastic processes and values of the corresponding
stochastic integral agree almost surely under conditions ensuring
the existence of both.
Let $Y$, $X$ be two regulated stochastic processes such that
for almost all $\omega\in\Omega$, the Left Young integral 
$$
(LY)\int_0^tY(\cdot,\omega)\,dX(\cdot,\omega)
$$ 
exists for all $t>0$.
For the rest of this chapter we consider an integration with respect
to a stochastic process $X$ which in addition is right-continuous.
In that case, by Definition \ref{variation}.\ref{LYandRY} of the 
Left Young integral, for $0\leq s<t<\infty$,
$$
(LY)\int_s^tY\,dX=(RRS)\int_s^tY_{-}^{(s)}\,dX
=(RRS)\int_s^tY_{-}^{(0)}\,X.
$$
The second equality holds also due to  the right-continuity of $X$.
We skip the superscript $(0)$ in what follows for notation
simplicity, and always assume that $Y(0-)=Y(0)$.
So that the Left Young integration  of $Y$ with respect to a {\cadlag} stochastic 
process $X$ reduces to the refinement Riemann-Stieltjes integration of the left 
modification $Y_{-}$ with respect to $X$. 
Then we define the {\it integral stochastic process} by
\beq\label{intsp}
\int Y_{-}\,dX:=\Big\{(RRS)\int_0^tY_{-}\,dX\colon\,t\geq 0\Big\}.
\eeq
Also for $1\leq p<\infty$, we say that a function $f$ on $[0,\infty)$ belongs to
$\dual ({\cal W}_p)^{{\rm loc}}$ if for each $0<T<\infty$,
$f$ belongs to $\dual ({\cal W}_p)[0,T]$ defined by (\ref{dual-Wp}).

\begin{prop}\label{MPS=Ito}
Let $1\leq p<\infty$, let $M$ be a local martingale of locally bounded
$p$-variation,
and let $H$ be an adapted stochastic process with almost all 
sample functions in {\rm $\dual ({\cal W}_p)^{{\rm loc}}$}.
Then the stochastic integral $H_{-}\cdot M$ and the integral
stochastic process $\smallint H_{-}\,dM$ are defined and
indistinguishable.
\end{prop}

\begin{proof}
The integral stochastic process $\smallint H_{-}\,dM$ exists by 
Theorem \ref{variation}.\ref{approximation}.
The stochastic integral $H_{-}\cdot M$ exists by Theorem III.17
in \cite[p.\ 106]{PP}.
We have to prove that for almost all $\omega\in\Omega$, 
the equality
\beq\label{1MPS=Ito}
(RRS)\int_0^tH(s-,\omega)\,dM(s,\omega)
=(H_{-}\cdot M)(t,\omega)
\eeq
holds for all $t\geq 0$.
There exists $\Omega_0\subset\Omega$ such that $\Pr (\Omega_0)=1$ and
for all $\omega\in\Omega_0$ the left and right sides of
(\ref{1MPS=Ito}) are regulated right-continuous functions of $t$.
Therefore it is enough to prove that (\ref{1MPS=Ito}) holds for
an arbitrary $t$ and for almost all $\omega$. 
Fix $t>0$.
To prove (\ref{1MPS=Ito}) we use 
term by term integration theorems for the $RRS$ integral and
the dominated convergence in probability for the stochastic integral.

For each integer $m\geq 1$, let $\lambda_m:=\{u_j^m:=tj2^{-m}\colon\,
j=0,\dots,2^m\}$ be a partition of $[0,t]$.
Let $H_m(s):=H(u_{j-1}^m)$ if $s\in [u_{j-1}^m,u_j^m)$ for some
$j=1,\dots,2^m$ and $H_m(t):=H(t)$.
We claim that for each $m\geq 1$, for each right-continuous
regulated sample function of $M$, and for each bounded sample function 
of $H$, we have
\beq\label{2MPS=Ito}
(RRS)\int_0^tH_m(s-)dM(s)=\sum_{j=1}^{2^m}H(u_{j-1}^m)\big [
M(u_j^m)-M(u_{j-1}^m)\big ]=S_{LC}(H,M;\lambda_m),
\eeq
where $H_m(0-)=H_m(0):=0$.
Indeed, let $m\geq 1$, and let $\kappa=\{([t_{i-1},t_i],s_i)\colon\,
i=1,\dots,n\}$ be a tagged partition of $[0,t]$ and a refinement of
$\lambda_m$.
Therefore for each $j=0,\dots,2^m$, there exists an index $i(j)\in
\{0,\dots,n\}$ such that $t_{i(j)}=u_j^m$.
Then the Riemann-Stieltjes sum based on $\kappa$ is given by
\begin{eqnarray*}
\lefteqn{S_{RS}((H_m)_{-},M;\kappa)=\sum_{i=1}^nH_m(s_i-)
\big [M(t_i)-M(t_{i-1})\big ]}\\[2mm]
&=&\sum_{j=1}^{2^m}\Big\{H_m(s_{i(j-1)+1}-)\big [M(t_{i(j-1)+1})
-M(u_{j-1}^m)\big ]+H(u_{j-1}^m)\big [M(u_j^m)-M(t_{i(j-1)+1})\big ]\Big\}.
\end{eqnarray*}
The right side of the preceding display approaches the right
side of (\ref{2MPS=Ito}) as $t_{i(j-1)+1}\downarrow u_{j-1}^m$
for each $j=1,\dots,2^m$ because sample functions of $M$ are
right-continuous and sample functions of $H$ are bounded.
Therefore given $\epsilon >0$, one can choose a partition
$\lambda$ of $[0,t]$ containing all points of $\lambda_m$
and points $v_j\in (u_{j-1}^m,u_j^m)$, $j=1,\dots,2^m$ sufficiently
close to $u_{j-1}^m$, so that the Riemann-Stieltjes sums
based on refinements of $\lambda$ are within $\epsilon$ to
the right side of (\ref{2MPS=Ito}).
Thus the $RRS$ integral in (\ref{2MPS=Ito}) exists and has
the stated value.

Next we prove that 
\beq\label{3MPS=Ito}
\lim_{m\to\infty}(RRS)\int_0^tH_m(s-)\,dM(s)
=(RRS)\int_0^tH(s-)\,dM(s).
\eeq
We have that almost every sample function of $M$ is right-continuous
and when restricted to $[0,t]$ is in ${\cal W}_p[0,t]$, while almost every
sample function of $H$ when restricted to $[0,t]$ is in 
$\dual ({\cal W}_p)[0,t]$.
First suppose that $p>1$.
In this case we use Theorem 5.6 on term by term integration of
L.\ C.\ Young \cite[p.\ 602]{LCY38b}.
Let $q<\infty$ be such that the sample function of $H$ is
in ${\cal W}_q[0,t]$ and $p^{-1}+q^{-1}>1$.
It is clear that the $q$-variation $v_q(H_m;[0,t])\leq
v_q(H;[0,t])$ for each $m\geq 1$.
The other assumptions of L.\ C.\ Young's theorem will follow once we show 
that $(H_m)_{-}$ converge to $H_{-}$ uniformly on the left at each
point of $(0,t]$, that is, for $s\in (0,t]$  and given $\epsilon >0$
there is an integer $M$ and a $\delta >0$ such that for all $m\geq M$
and all $x\in [s-\delta, s)$, $|H_m(x-)-H(x-)|<\epsilon$.
Let $s\in (0,t]$ and let $\epsilon >0$.
Choose a $\delta >0$ such that $|H(x)-H(s-)|<\epsilon /2$
for each $x\in [s-\delta, s)$.
Let $M$ be the minimal positive integer  such that $t2^{-M}\leq
\delta /2$.
Then for all $x\in [s-\delta/2,s)$ and all $m\geq M$,
$|H_m(x-)-H(x-)|\leq\epsilon /2+|H(s-)-H(u_{i-1}^m)|
<\epsilon$ for some $u_{i-1}^m\in [s-\delta,s)$.
Thus all assumptions of L.\ C.\ Young's theorem are satisfied,
and hence (\ref{3MPS=Ito}) holds in the case $p>1$.
Now suppose that $p=1$.
In this case we use the Osgood convergence theorem proved
for the Young-Stieltjes integrals by Hildebrandt 
\cite[Theorem 19.3.14]{THH63}.
Since $\|(H_m)_{-}\|_{\sup}\leq \|H\|_{\sup}<\infty$ and
$(H_m)_{-}$ converge pointwise to $H_{-}$, (\ref{3MPS=Ito})
holds by the Osgood theorem in the case $p=1$.
Therefore (\ref{3MPS=Ito}) holds.
On the other hand, by Theorem VIII.15 of Dellacherie and Meyer
\cite[p. 324]{DM82},
the limit
$$
\lim_{m\to\infty}S_{LC}(H,M;\lambda_m)=(Y_{-}\cdot M)(t)
$$
exists in probability.
This in conjunction with (\ref{2MPS=Ito}) and (\ref{3MPS=Ito})
yields that (\ref{1MPS=Ito}) holds for almost all $\omega\in\Omega$.
The proof of Proposition \ref{MPS=Ito} is complete.
\qed\end{proof}

\section{\nsem-semimartingales and the stochastic $\lsi$-integral}

In this section we define the class of \nsem-semimartingales 
and the stochastic $\lsi$-integral with respect to a \nsem-semimartingale.

\begin{defn}\label{p-semimartingale}
{\rm Let $1\leq p<2$, and let $X$ be an adapted {\cadlag} stochastic process.
We say that $X$ is a \emph{$p$-semimartingale} if there exist stochastic 
processes $M$ and $\pv$ such that 
\beq\label{1ps}
X-X(0)=M+\pv\qquad\mbox{almost surely,}
\eeq
$M(0)=\pv (0)=0$, $M$ is a local martingale and 
$\pv$ is a stochastic process of locally bounded $p$-variation.
The set of all such pairs $(M,\pv)$ is denoted by
${\cal D}_p(X)$, and ${\cal D}(X):=\cup\{{\cal D}_p(X)\colon\,
1\leq p<2\}$.
A stochastic process $X$ is called a {\it \nsem-semimartingale}
if $X$ is a $p$-semimartingale for some $1\leq p<2$, that is, if
the set ${\cal D}(X)$ is nonempty.}
\end{defn}

The 1-semimartingale is the same as semimartingale.
As in the case of semimartingales, decomposition (\ref{1ps})
may not be unique.
For example, the L\'evy decomposition
can be viewed in two different ways: either as a $1$-semimartingale
(see Theorem I.40 in \cite[p.\ 31]{PP}), or as a $p$-semimartingale
with the decomposition into a Brownian motion part and a pure
jump part.
However we prove that in most interesting cases a "natural" integral with respect
to a $p$-semimartingale and a continuous part of the quadratic variation
of a $p$-semimartingale are unique in spite of a possible non-uniqueness 
of decomposition (\ref{1ps}).

In the definition of a \nsem-semimartingale we require the second component
$A$ in decomposition (\ref{1ps}) to have locally finite $p$-variation for some
$1\leq p<2$.
Instead, one may define a class of stochastic processes with decomposition
(\ref{1ps}), where $A$ has almost all sample functions locally in
${\cal W}_2^{\ast}$.
Stochastic processes with such a property are reminiscent to Dirichlet stochastic
processes defined by F\"olmer \cite{HFb}, and further modified 
in different directions by other authors.

Let $X$ be a \nsem-semimartingale,
and let $H$ be an adapted regulated stochastic process.
Let $(M,\pv)\in {\cal D}(X)$.
The stochastic integral $H_{-}\cdot M=\{(H_{-}\cdot M)(t)\colon\,t\geq 0\}$ 
of the left-continuous modification
$H_{-}$ with respect to a local martingale $M$ is well-known.
We can and do assume that $H_{-}\cdot M$ is a {\cadlag} stochastic process.
It is natural to define a stochastic integral of $H$
with respect to a \nsem-semimartingale $X$ by
\beq\label{SI}
(H{\lsi}X)(t):=(H_{-}\cdot M)(t)+(RRS)\int_0^tH_{-}\,d\pv,\qquad t\geq 0,
\eeq
provided the integral stochastic process $\smallint H_{-}\,d\pv$
is defined (see (\ref{intsp})).
If $\pv$ is a stochastic process of locally bounded variation then
the refinement  Riemann-Stieltjes integral in (\ref{SI}) is defined and 
has the same value as the Lebesgue-Stieltjes integral.
If $\pv$ is a stochastic process of locally bounded $p$-variation
with $1<p<2$ then the refinement Riemann-Stieltjes integral in (\ref{SI}) 
is defined provided $H$ is a stochastic process of locally bounded 
$q$-variation with  $p^{-1}+q^{-1}>1$ by the L.\ C.\ Young Theorem on Stieltjes 
integrability (see (\ref{-1intr})).
In particular, the refinement  Riemann-Stieltjes integral in (\ref{SI}) is 
defined provided $H$ is a local martingale by Theorem \ref{lepingle}. 
However, since decomposition (\ref{1ps}) may not be unique
a value of sum (\ref{SI}) also need not be the same
for different decompositions of a \nsem-semimartingale $X$.
We define the stochastic $\lsi$-integral with respect to
a \nsem-semimartingale provided a sum of the two integrals
in (\ref{SI}) does not depend on decomposition (\ref{1ps}).

\begin{defn}\label{stochint}
{\rm Let $X$ be a \nsem-semimartingale, and let $H$ be an adapted regulated
stochastic process  such that the integral stochastic process 
$\smallint H_{-}\,d\pv$ is defined for each $(M,\pv)\in {\cal D}(X)$.
If the stochastic processes
\beq\label{1stochint}
\Big\{(H_{-}\cdot M)(t)+(RRS)\int_0^tH_{-}\,d\pv\colon\,t\geq 0
\Big\}
\eeq
are indistinguishable for different pairs $(M,\pv)\in {\cal D}(X)$
then define the \emph{stochastic $\lsi$-integral}
$H{\lsi}X=\{(H{\lsi} X)(t)\colon\,t\geq 0\}$ by (\ref{SI}).}
\end{defn}

The next statement follows from the bilinearity of the
stochastic integral and of the refinement Riemann-Stieltjes integral. 

\begin{prop}\label{bilinear}
Let $X$, $Y$ be \nsem-semimartingales and let $H$, $G$ be adapted
regulated stochastic processes.
Then
$$
(H+G){\lsi}(X+Y)=H{\lsi}X+H{\lsi}Y+G{\lsi}X+G{\lsi}Y
$$
provided all five stochastic $\lsi$-integrals are defined.
\end{prop}

For a stochastic process $Y=\{Y(t)\colon\,t\geq 0\}$,
the \emph{jump process} $\Delta^{-}Y=\{\Delta^{-}Y(t)\colon\,
t\geq 0\}$ is defined by $\Delta^{-}Y(t)=Y(t)-Y(t-)$ if $t>0$,
and $\Delta^{-}Y(0):=0$.

\begin{prop}\label{sijumps}
Let $X$ be a \nsem-semimartingale, and let $H$ be an adapted
regulated stochastic process such that the stochastic $\lsi$-integral
$H{\lsi}M$ is defined.
Then almost all sample functions of $H{\lsi}X$ are {\cadlag} and
$$
\big\{\Delta^{-}(H{\lsi}X)(t)\colon\,t\geq 0\big\}
=\big\{H_{-}(t)\Delta^{-}X(t)\colon\,t\geq 0\big\}.
$$
\end{prop}

\begin{proof}
By the same property for the stochastic integral with respect
to a local martingale (Theorem II.13 in \cite[p.\ 53]{PP}),
and by Proposition \ref{variation}.\ref{LYjumps}, we have
$$
\Delta^{-}(H{\lsi}X)=\Delta^{-}(H_{-}\cdot M)+\Delta^{-}(\smallint H_{-}
\,d\pv)=H_{-}\Delta^{-}M+H_{-}\Delta^{-}\pv=H_{-}\Delta^{-}X
$$
almost surely, and the proof is complete.
\qed\end{proof}

Next is a sufficient condition for the stochastic $\lsi$-integrability with 
respect to a \nsem-semimartingale. 

\begin{thm}\label{Estochint}
Let $X$ be a \nsem-semimartingale, and let $H$ be an adapted
stochastic process with the local $p$-variation index $\pindex (H)=2$.
Then the stochastic $\lsi$-integral $H{\lsi}X$ is defined, is a 
\nsem-semimartingale, and for each $t>0$ and $\lambda=\{\lambda_m
\colon\,m\geq 1\}\in\Lambda [0,t]$,
\beq\label{1Estochint}
\lim_{m\to\infty} S_{LC}(H,X;\lambda_m)
=(H{\lsi}X)(t)\qquad\mbox{in probability.}
\eeq
\end{thm}

\begin{proof}
Let $(M,\pv)\in {\cal D}(X)$.
Then $M$ is a local martingal, $\pv$ is an adapted {\cadlag} stochastic 
process of locally bounded $p$-variation for some $1\leq p<2$, 
$M(0)=\pv (0)=0$, and the decomposition (\ref{1ps}) holds.
By Theorem III.17 in \cite[p.\ 106]{PP}, the stochastic integral
$H_{-}\cdot M$ is defined and is a local martingale.
Since $1\leq p<2$ and the local $p$-variation index $\pindex (H)=2$,
there exists $2<q<\infty$ such that $p^{-1}+q^{-1}>1$
and $H$ has almost all sample functions in ${\cal W}_q^{\rm loc}$.
Thus by Theorem \ref{variation}.\ref{approximation}, there exists
the integral stochastic process $\smallint H_{-}\,d\pv$.
Hence  stochastic process (\ref{1stochint}) is defined for $(M,\pv)$.
We have to show that this process is unique up to indistinguishability
for different elements of ${\cal D}(X)$.
If the set ${\cal D}(X)$ has a single element then there is nothing to prove.
Suppose that there are two differenet elements in ${\cal D}(X)$,
that is, $(M_1,\pv_1)\in {\cal D}_{p_1}(X)$ and $(M_2,\pv_2)\in
{\cal D}_{p_2}(X)$ for some $p_1, p_2\in [1,2)$. 
Taking the maximum between $p_1$ and $p_2$ if they are different,
we can and do assume that $p_1=p_2=p$.
By the decomposability property of $X$, we have that stochastic
process $M_1-M_2=-(\pv_1-\pv_2)$ has almost all sample functions
in ${\cal W}_p^{\rm loc}$.
On the other hand, there exists $2<q<\infty$ such that $p^{-1}+q^{-1}>1$
and $H$ has almost all sample functions in ${\cal W}_q^{\rm loc}$,
so that the refinement Riemann-Stieltjes integral 
$$
(RRS)\int_0^tH_{-}(\cdot,\omega)\,d\big (M_1-M_2\big )(\cdot,\omega)
$$
exists for almost all $\omega\in\Omega$ and for all $t>0$.
Thus by Proposition \ref{MPS=Ito}, for almost all $\omega\in\Omega$
and for all $t>0$,
$$
(H_{-}\cdot M_1)(t,\omega)-(H_{-}\cdot M_2)(t,\omega)
=(RRS)\int_0^tH_{-}(\cdot,\omega)\,d(M_1-M_2)(\cdot,\omega)
$$
$$
=(RRS)\int_0^tH_{-}(\cdot,\omega)\,d\pv_2(\cdot,\omega)
-(RRS)\int_0^tH_{-}(\cdot,\omega)\,d\pv_1(\cdot,\omega).
$$
Therefore the stochastic processes (\ref{1stochint}) corresponding
to $(M_1,\pv_1)$ and $(M_2,\pv_2)$ are indistinguishable,
and so the stochastic $\lsi$-integral $H{\lsi}X$ is defined.

By Theorem \ref{variation}.\ref{approximation} and Proposition 
\ref{variation}.\ref{LYjumps}, the integral stochastic process 
$\smallint H_{-}\,d\pv$ is an adapted {\cadlag} stochastic process.
By Corollary \ref{variation}.\ref{indefinite}, $\smallint H_{-}\,d\pv$
has the same $p$-variation property as $\pv$.
Since $H_{-}\cdot M$ is a local martingale, it follows that 
the stochastic $\lsi$-integral $H{\lsi}X=H_{-}\cdot M+\smallint H_{-}\,d\pv$ 
is a \nsem-semimartingale with value $0$ at $0$.
The approximation in probability of the stochastic $\lsi$-integral by 
Left Cauchy sums follows from Theorem \ref{variation}.\ref{approximation} 
and from Theorem II.21 in \cite[p.\ 57]{PP}.
The proof of Theorem \ref{Estochint} is complete.
\qed\end{proof}

The preceding theorem shows that the property of being 
a \nsem-semimartingale is preserved by the stochastic $\lsi$-integration.
Also, it provides a class of stochastic $\lsi$-integrable processes
which is small as compared with the class of all predictable
stochastic processes integrable in the sense of the stochastic
integral with respect to a semimartingale.
However, the class of possible integrators allowed by Theorem 
\ref{Estochint} is larger.

In Stochastic Analysis the property proved next is called associativity
of the stochastic integral.

\begin{thm}\label{associativity}
Let $X$ be a \nsem-semimartingale, and let $G$, $H$ be adapted
stochastic processes with the local $p$-variation indeces 
$\pindex (G)=\pindex (H)=2$.
Then the stochastic $\lsi$-integral $H{\lsi}X$ is a \nsem-semimartingale and
\beq\label{1associativity}
G{\lsi}(H{\lsi}X)=(GH){\lsi}X.
\eeq
\end{thm}

\begin{proof}
The stochastic $\lsi$-integral $H{\lsi}X$ is defined and is a 
\nsem-semimartingale by Theorem \ref{Estochint}.
Let $(M,\pv)\in {\cal D}(X)$.
Then $(H_{-}\cdot M,\smallint H_{-}\,d\pv)\in {\cal D}(H{\lsi}X)$ 
by Theorem III.17 
in \cite[p.\ 106]{PP} and by Corollary \ref{variation}.\ref{indefinite}.
By the associativity of the stochastic integral with respect to
a local martingale (Theorem II.19 in \cite[p.\ 55]{PP}) and
by Proposition \ref{variation}.\ref{MPSsrule}, for any $0<t<\infty$,
the equalities
\begin{eqnarray*}
\big (G{\lsi}(H{\lsi}X)\big )(t) &=&
\big (G_{-}\cdot(H_{-}\cdot M)\big )(t)+(MPS)\int_0^tG_{-}\,d\Big (\int
H_{-}\,d\pv\Big )\\[2mm]
&=& \big ((GH)_{-}\cdot M\big )(t)+(MPS)\int_0^t(GH)_{-}\,d\pv
=\big ((GH){\lsi}X\big )(t)
\end{eqnarray*}
hold almost surely.
Since almost all sample functions of the stochastic integrals are {\cadlag},
the relation (\ref{1associativity}) follows.
\qed\end{proof}

\section{Quadratic covariation of \nsem-semimartingales}

Let $X$, $Y$ be two \nsem-semimartingales.
Since the local $p$-variation index of a \nsem-semimartingale
is equal to $2$ by Theorem \ref{lepingle}, the stochastic $\lsi$-integrals
$X{\lsi}Y$ and $Y{\lsi}X$ are defined by Theorem \ref{Estochint}.
Therefore the stochastic process $[X,Y]=\{[X,Y](t)\colon\,t\geq 0\}$ 
defined by
\beq\label{qcsp}
[X,Y](t):=X(t)Y(t)-X(0)Y(0)-\big (X{\lsi}Y\big )(t)-\big (Y{\lsi}
X\big )(t),\qquad t\geq 0,
\eeq
exists and is called the {\it quadratic covariation} of the
pair $X$, $Y$. 
In the case $X$ and $Y$ are semimartingales, definition (\ref{qcsp})
differ from  the definition of the quadratic covariation
in Section II.6 of \cite[p.\ 58]{PP} by the random variable
$X(0)Y(0)$.
Let $\bar X:=X-X(0)$ and $\bar Y:=Y-Y(0)$.
By (\ref{1Estochint}), we have $X{\lsi}\bar Y=X{\lsi}Y$ and
$Y{\lsi}\bar X=Y{\lsi}X$.
Therefore it follows that
\beq\label{3qcsp}
[X,Y]=\bar X\bar Y -\bar X{\lsi}\bar Y -\bar Y{\lsi}\bar X.
\eeq
The above notion of the quadratic covariation for \nsem-semimartingales 
is consistent with the quadratic $\lambda$-covariation for functions 
(Definition \ref{function}.\ref{8bracket}) in the following sense.

\begin{prop}\label{consistent}
Let $X$, $Y$ be \nsem-semimartingales such that for some 
$\lambda\in\Lambda [0,T]$, $0<T<\infty$,
the $2$-vector function $(X(\cdot,\omega),Y(\cdot,\omega))$ has
the quadratic $\lambda$-covariation for almost all $\omega\in\Omega$.
Then for almost all $\omega\in\Omega$,
\beq\label{1consistent}
[X,Y](t,\omega)=[X(\cdot,\omega),Y(\cdot,\omega)]_{\lambda}(t),
\qquad 0\leq t\leq T.
\eeq
\end{prop}

\begin{proof}
Let $\lambda=\{\lambda_m\colon\,m\geq 1\}\in\Lambda [0,T]$ be the sequence 
of partitions $\lambda_m=\{t_i^m\colon\,i=0,\dots,n(m)\}$
of $[0,T]$ satisfying the assumptions of the proposition. 
By (\ref{1Estochint}), we have that for each $0\leq t\leq T$,
\begin{eqnarray*}
[X,Y](t)&=&\lim_{m\to\infty}\sum_{i=1}^{n(m)}\big [X(t_i^m\wedge t)
-X(t_{i-1}^m\wedge t)\big ]\big [Y(t_i^m\wedge t)
-Y(t_{i-1}^m\wedge t)\big ]\\[2mm]
&=&\lim_{m\to\infty}C(X,Y;\lambda_m\Cap [0,t])
\qquad\mbox{in probability}.
\end{eqnarray*}
It then follows that there exists a countable dense subset $D$ of $[0,T]$,
and a subset $\Omega_0\subset\Omega$ such that $\Pr (\Omega_0)=1$
and for each $\omega\in\Omega_0$, the functions $[X,Y](\cdot,\omega)$
and $[X(\cdot,\omega),Y(\cdot,\omega)]_{\lambda}$ on $[0,T]$
are {\cadlag} and
$$
[X,Y](t,\omega)=[X(\cdot,\omega),Y(\cdot,\omega)]_{\lambda}(t),
\qquad t\in D.
$$
Therefore (\ref{1consistent}) holds.
\qed\end{proof}

Let $X$, $Y$ be \nsem-semimartingales, and let
$(M,\pv)\in {\cal D}(X)$, $(N,B)\in {\cal D}(Y)$.
By (\ref{3qcsp}) and by Proposition \ref{bilinear}, we have
the equality
$$
[X,Y]=[M+\pv,N+B]=[M,N]+[M,B]+[A,N]+[A,B].
$$
Since $M(0)=N(0)=0$, by definition (\ref{qcsp}), the equality
$$
[M,N]=MN-M_{-}\cdot N-N_{-}\cdot M
$$
holds and shows that $[M,N]$ agree with the quadratic
covariation for semimartingales used in Stochastic Analysis, and
hence inherits standard properties of the quadratic
covariation.
For example, almost all sample functions of $[M,N]$ are {\cadlag},
have locally bounded variation, and
$\Delta^{-}[M,N]=\Delta^{-}M\Delta^{-}N$.
The same is true for the quadratic covariations
$[M,B]$, $[A,N]$ and $[A,B]$.
This follows from the preceding proposition in conjunction
with Proposition \ref{function}.\ref{qcforpq}.
Also the last three quadratic covariation stochastic processes
have pure jump sample functions.
Since almost all sample functions of $[X,Y]$ are {\cadlag},
have locally bounded variation and $\Delta^{-}[X,Y]=\Delta^{-}X
\Delta^{-}Y$, almost all sample functions of the stochastic process
$$
[X,Y]^c:=\Big\{[X,Y](t)-\sum_{(0,t]}\Delta^{-}X\Delta^{-}Y
\colon\,t\geq 0\Big\}
$$ 
are continuous.
The following shows that the continuous part of the quadratic covariation
of $(2-\epsilon)$-semimartingales agree with the analogous Stochastic Analysis 
construction and does not depend on their decompositions
into local martingale and bounded $p$-variation parts.

\begin{prop}
Let $X$, $Y$ be \nsem-semimartingales.
Then for each $(M,\pv)\in {\cal D}(X)$ and $(N,B)\in {\cal D}(Y)$,
$[M,N]^c$ is indistinguishable from $[X,Y]^c$.
\end{prop}

\begin{proof}
Let $(M,\pv)\in {\cal D}(X)$ and $(N,B)\in {\cal D}(Y)$.
Since $M(0)=N(0)=\pv (0)=B(0)=0$, by (\ref{3qcsp}) and Proposition
\ref{bilinear}, we have for each $t>0$,
\begin{eqnarray*}
[X,Y](t)&=&(M+\pv )(t)(N+B)(t)-\big ((M+\pv){\lsi}
(N+B)\big )(t)-\big ((N+B){\lsi}(M+\pv )\big )(t)\\[2mm]
&=&[M,N](t)+\Big\{M(t)B(t)-\big (B_{-}\cdot M\big )(t)
-(RRS)\int_0^tM_{-}\,dB\Big\}\\[2mm]
& &+\Big\{\pv(t)N(t)-\big (\pv_{-}\cdot N\big )(t)
-(RRS)\int_0^tN_{-}\,d\pv\Big\}\\[2mm]
& &+\Big\{\pv(t)B(t)-(RRS)\int_0^t\pv_{-}\,dB
-(RRS)\int_0^tB_{-}\,d\pv\Big\}.
\end{eqnarray*}
Let $\pv=\pv (\cdot,\omega)$ and $B=B(\cdot,\omega)$ be {\cadlag} 
sample functions which have locally bounded 
$p$-variation for some $1\leq p<2$.
For any $t>0$ and any $\lambda\in\Lambda [0,t]$, by Propositions 
\ref{function}.\ref{qcforpq}, \ref{function}.\ref{property5}
and Theorem \ref{variation}.\ref{approximation}, there exists the
quadratic $\lambda$-covariation $[\pv,B]_{\lambda}$ on $[0,t]$ 
and the equality
$$
\pv(t)B(t)-(RRS)\int_0^t\pv_{-}\,dB
-(RRS)\int_0^tB_{-}\,d\pv=[\pv,B]_{\lambda}(t)
=\sum_{(0,t]}\Delta^{-}\pv\Delta^{-}B
$$
holds.
In addition to sample functions $\pv$ and $B$, let
$M=M(\cdot,\omega)$ and $N=N(\cdot,\omega)$ be {\cadlag} sample functions 
which have the local $p$-variation index equal to $2$.
For any $t>0$ and any $\lambda\in\Lambda [0,t]$, by the preceding argument, 
there exist the quadratic $\lambda$-covariations $[M,B]_{\lambda}$,
$[\pv,N]_{\lambda}$ on $[0,t]$, and the equalities
$$
M(t)B(t)-(RRS)\int_0^tB_{-}\,dM-(RRS)\int_0^tM_{-}\,dB
=[M,B]_{\lambda}(t)=\sum_{(0,t]}\Delta^{-}M\Delta^{-}B,
$$
$$
\pv(t)N(t)-(RRS)\int_0^t\pv_{-}\,dN-(RRS)\int_0^tN_{-}\,d\pv
=[\pv,N]_{\lambda}(t)=\sum_{(0,t]}\Delta^{-}\pv\Delta^{-}N
$$
hold.
By Proposition \ref{MPS=Ito}, the stochastic integrals 
$B_{-}\cdot M$ and $\pv_{-}\cdot N$ are indistinguishable
from the integral stochastic processes $\smallint B_{-}\,dM$
and $\smallint \pv_{-}\,dN$, respectively.
Then using the linearity of the unconditional sums, it follows that
the quadratic covariation $[X,Y]$ is
indistinguishable from the stochastic process
$$
\Big\{[M,N](t)-\sum_{(0,t]}\Delta^{-}M\Delta^{-}N
+\sum_{(0,t]}\Delta^{-}X\Delta^{-}Y\colon\,t\geq 0\Big\}.
$$ 
The proof of the proposition is complete.
\qed\end{proof}

\begin{prop}\label{1qcsp}
Let $X$, $Z$ be \nsem-semimartingales, let $H$ be an adapted {\cadlag} 
stochastic process with the local $p$-variation index $\pindex (H)=2$,
and let $Y:=C+H{\lsi}Z$ for some constant $C$.
Then
\beq\label{2qcsp}
[X,Y]=(RRS)\int H_{-}\,d[X,Z].
\eeq
\end{prop}

\begin{proof}
By the assumption there are local martingales $M, N$ and
adapted {\cadlag} stochastic processes $\pv, B$
having locally bounded $p$-variation for some $1\leq p<2$
such that $M(0)=N(0)=\pv (0)=B(0)=0$, $X-X(0)=M+\pv$ and
$Z-Z(0)=N+B$ almost surely.
Since $H_{-}\cdot M$ is a local martingale and $\smallint H_{-}\,d\pv$
is an adapted {\cadlag} stochastic process having
locally bounded $p$-variation, $Y$ is a \nsem-semimartingale
with $Y(0)=C$.
Therefore in (\ref{2qcsp}), the two quadratic covariations 
and the integral stochastic process are defined. 
Since almost all sample functions of the stochastic processes in 
(\ref{2qcsp}) are {\cadlag}, it is enough to prove that for any
$0<t<\infty$
$$
[X,Y](t)=(RRS)\int_0^t H_{-}\,d[X,Z]\qquad\mbox{almost surely.}
$$
To this aim fix an arbitrary $0<t<\infty$.
By Proposition \ref{consistent} applied to the two
pairs of \nsem-semimartingales $M$, $B$ and $\pv$, $Z$ in conjunction
with Proposition \ref{function}.\ref{qcforpq}, the equalities
$$
[M,B](t)=\sum_{(0,t]}\Delta^{-}M\Delta^{-}B,\qquad
[\pv,Z](t)=\sum_{(0,t]}\Delta^{-}\pv\Delta^{-}Z,
$$
hold almost surely.
By the preceding argument in conjunction with Propositions 
\ref{variation}.\ref{LYjumps} and \ref{sijumps}, 
the two pairs of \nsem-semimartingales $\pv$, $H{\lsi}Z$ and 
$M$, $\Phi:=\smallint H_{-}\,dB$,  satisfy the relations
$$
[M, \Phi ](t)=\sum_{(0,t]}\Delta^{-}M\Delta^{-}\Phi
=\sum_{(0,t]}H_{-}\Delta^{-}M\Delta^{-}B
=(RRS)\int_0^tH_{-}\,d[M,B],
$$
$$
[\pv,H{\lsi} Z](t)=\sum_{(0,t]}\Delta^{-}\pv\Delta^{-}(H{\lsi} Z)
=\sum_{(0,t]}H_{-}\Delta^{-}\pv\Delta^{-}Z
=(RRS)\int_0^tH_{-}\,d[\pv,Z]
$$
almost surely.
Since $M$, $N$ are local martingales such that $M(0)=N(0)=0$,
by Theorem II.29 in \cite[p.\ 68]{PP}, the equality 
$$
[M,H_{-}\cdot N](t)=(LS)\int_{[0,t]}H_{-}d[M,N]=(RRS)\int_0^tH_{-}d[M,N]
$$
holds almost surely, where $LS$ refers to the Lebesgue-Stieltjes
integration.
Since the quadratic covariation and the $RRS$ integral are bilinear,
by relation (\ref{3qcsp}), we have that equalities
\begin{eqnarray*}
[X,Y](t)&=&[M+\pv,H{\lsi} Z](t)
=[M,H_{-}\cdot N](t)+[M,\Phi ](t)
+[\pv,H{\lsi} Z](t)\\[2mm]
&=&(RRS)\int_0^tH_{-}d[M,N]+(RRS)\int_0^tH_{-}\,d[M,B]
+(RRS)\int_0^tH_{-}\,d[\pv,Z]\\[2mm]
&=&(RRS)\int_0^tH_{-}\,d[X,Z]
\end{eqnarray*}
hold almost surely.
The proof of Proposition \ref{1qcsp} is complete.
\qed\end{proof}

\section{The It\^o formula and the linear
stochastic $\lsi$-integral equation}

First we extend the It\^o formula to \nsem-semimartingales using
$\lsi$-integral, and second we show that the stochastic Dol\'eans
exponential is the unique \nsem-semimartingale satisfying the
linear $\lsi$-integral equation.

As before we say that a function $\phi\colon\,\RR^d\mapsto\RR$ is a 
$C^2$ class function if its all second order partial derivatives exist and 
are continuous.
Then we write $\phi_k':=\frac{\partial\phi}{\partial x_k}$
and $\phi_{kl}'':=\frac{\partial^2\phi}{\partial x_k\partial x_l}$
for $k,l=1,\dots,d$.

\begin{thm}\label{Ito}
Let $X=(X_1,\dots,X_d)$ be a $d$-tuple of \nsem-semimartingales, 
and let $\phi\colon\,\RR^d\mapsto\RR$ be a $C^2$ class function.
Then {\rm (1)} for $k=1,\dots,d$, the stochastic $\lsi$-integral 
$\big (\phi_k'{\circ}X){\lsi} X_k$ is defined, 
{\rm (2)} for $k,l=1,\dots,d$, the Riemann-Stieltjes integral stochastic process
$$
(RS)\int(\phi_{kl}''{\circ}X)_{-}\,d[X_k,X_l]^c:=
\Big\{(RS)\int_0^t(\phi_{kl}''{\circ}X)_{-}\,d[X_k,X_l]^c
\colon\,t\geq 0\Big\}
$$
is defined, 
{\rm (3)} the unconditional sum stochastic process
$$
\sum\Big\{\Delta^{-}(\phi{\circ}X)-\sum_{k=1}^d
(\phi_k'{\circ}X)_{-}\Delta^{-}X_k\Big\}
:=\Big\{\sum_{(0,t]}\big\{\Delta^{-}(\phi{\circ}X)-\sum_{k=1}^d
(\phi_k'{\circ}X)_{-}\Delta^{-}X_k\big\}\colon\,t\geq 0\Big\}
$$
is adapted {\cadlag} stochastic process with locally
bounded variation, and {\rm (4)}
the composition $\phi{\circ}X$ is a \nsem-semimartingale
and 
\begin{eqnarray}
\phi{\circ}X&=&(\phi{\circ}X)(0)+\sum_{k=1}^d(\phi_k'{\circ}X)
{\lsi}X_k+\frac{1}{2}\sum_{k,l=1}^d(RS)\int(\phi_{kl}''{\circ}X)_{-}\,
d[X_k,X_l]^c\label{2ps}\\[2mm]
& &+\sum\Big\{\Delta^{-}(\phi{\circ}X)-\sum_{k=1}^d
(\phi_k'{\circ}X)_{-}\Delta^{-}X_k\Big\}.\nonumber
\end{eqnarray}
\end{thm}

\begin{proof}
Composing with a Lipschitz smooth function does not change the 
$p$-variation property.
Thus the stochastic $\lsi$--integrals $(\phi_k'{\circ}X){\lsi}X_k$,
$k=1,\dots,d$, exist by Theorem \ref{Estochint}.
The Riemann-Stieltjes stochastic processes are defined and have
almost all sample functions continuous because each $[X_k,X_l]^c$
has almost all sample functions continuous with locally bounded
variation and each integrand is a regulated stochastic process.
Applying the mean value theorem and H\"older's inequality it follows
that the unconditional sum stochastic process is defined because
almost surely
$$
\sum_{(0,t]}\big\{\Delta^{-}X\big\}^2\leq [X](t)<\infty
$$
for each $t>0$.
By Theorem \ref{Estochint} again, the stochastic $\lsi$-integral is a
\nsem-semimartingale, so that the composition $\phi{\circ}X$ is
also a \nsem-semimartingale provided formula (\ref{2ps}) holds.
For notation simplicity the proof of (\ref{2ps}) is given in the 
case $d=1$.

The proof will be complete once we show that the equality 
\begin{eqnarray}\label{10ps}
(\phi{\circ}X)(t)&=&(\phi{\circ}X)(0)+\big ((\phi'{\circ}X)_{-}{\lsi}
X\big )(t)+\frac{1}{2}(RS)\int_0^t(\phi''{\circ}X)_{-}\,d[X]^c\\[2mm]
& &+\sum_{(0,t]}\Big\{\Delta^{-}(\phi{\circ}X)-(\phi'{\circ}X)_{-}
\Delta^{-}X\Big\}\nonumber
\end{eqnarray}
holds almost surely for an arbitrary $0<t<\infty$. 
Indeed, it then holds almost surely for countably many
$t$'s, and since all the stochastic processes are {\cadlag} the equality 
will hold almost surely for all $0<t<\infty$. 
Thus let $0<t<\infty$ be fixed, and let $(M,\pv)\in {\cal D}(X)$, that is,
$X=X(0)+M+\pv$, $M$ is a local martingale, $\pv$ is an adapted 
stochastic process with locally bounded $p$-variation,
and $M(0)=\pv (0)=0$.
For a positive integer $m$, let $\tau_m:=\inf\{s\in [0,t]\colon\,
|M(s)|>m$ or $|\pv (s)|>m\}$.
Then each $\tau_m$ is a stopping time such that
$\Pr (\cup_m\{\tau_m=t\})=1$.
Therefore (\ref{10ps}) holds provided it holds when
$M$ and $\pv$ are replaced by $M^{\tau_m}$ and $\pv^{\tau_m}$ respectively.
For notation simplicity again, let $|M(s)|\vee |\pv (s)|<K$ for some finite 
$K$ and all $s\in [0,t)$ almost surely.
Since $\phi''$ is uniformly continuous over $[-K,K]$,
we then are able to control increments $\phi''(X(s))-\phi''(X(r))$
almost surely with $s$ and $r$ having values in the domain of $X$ 
outside of an arbitrary small set.

Let $\epsilon \in (0,1)$.
Then there exists $\delta >0$ such that $|\phi ''(u)-\phi ''(v)|
<\epsilon$ whenever $|u-v|<\delta$ and $u, v\in [-K,K]$.
For almost all $\omega\in \Omega$ there exists a finite set
$\mu=\{s_j\colon\,j=1,\dots,k\}\subset (0,t]$ (depending on
$\omega$) such that sample functions $M=M(\cdot,\omega)$ and
$\pv=\pv (\cdot,\omega)$ satisfy the conditions:
\beq\label{7ps}
\Big |\sum_{\mu}\big\{\Delta^{-}(\phi{\circ}X)-(\phi'{\circ}X)_{-}
\Delta^{-}X\big\}
-\sum_{(0,t]}\big\{\Delta^{-}(\phi{\circ}X)-(\phi'{\circ}X)_{-}
\Delta^{-}X\big\}\Big |<\epsilon,
\eeq
\beq\label{8ps}
\Big |\sum_{\mu}\big\{(\phi''{\circ}X)_{-}(\Delta^{-}M)^2\big\}
-\sum_{(0,t]}\big\{(\phi''{\circ}X)_{-}(\Delta^{-}M)^2\big\}\Big |
<\epsilon,
\eeq
\beq\label{9ps}
\max_{1\leq j\leq k} Osc\,(X;[s_{j-1},s_j))<\delta
\quad\mbox{and}\quad
\sum_{j=1}^kv_2(\pv;[s_{j-1},s_j))<\epsilon,
\eeq
where $s_0:=0$.
The first condition in (\ref{9ps}) holds because almost all sample
functions of $X$ are {\cadlag}, and the second one follows from
Lemma \ref{variation}.\ref{lemma1}.
Let $\{\lambda_m\colon\,m\geq 1\}$ be a nested sequence of partitions
$\lambda_m=\{t_i^m\colon\,i=0,\dots,n(m)\}$ of $[0,t]$ such that
$\cup_m\lambda_m$ is dense in $[0,t]$, and
let $m_0=m_0(\omega)$ be the minimal integer such that
each $[s_{j-1},s_j)\cap\lambda_m\not =\emptyset$.
For each $m\geq m_0$, let $I_1(m):=\{i=1,\dots,n(m)\colon\,\mu\cap
(t_{i-1}^m,t_i^m]\not =\emptyset\}$ 
and $I_2(m):=\{1,\dots,n(m)\}\setminus I_1(m)$. 
By a telescoping sum and Taylor's theorem (\ref{3ps}),
for each $m\geq m_0$ we have
\begin{eqnarray}
\lefteqn{(\phi{\circ}X)(t)-(\phi{\circ}X)(0)
=\sum_{i=1}^{n(m)}\phi'(X(t_{i-1}^m))\Delta_i^mX
+\frac{1}{2}\sum_{i\in I_2(m)}\phi''(X(t_{i-1}^m))(\Delta_i^mM)^2}
\label{4ps}\\
&+&\sum_{i\in I_1(m)}\big\{\Delta_i^m(\phi{\circ}X)
-\phi'(X(t_{i-1}^m))\Delta_i^mX\big\}
+\,\frac{1}{2}\sum_{i\in I_2(m)}\phi''(y_i^m)
\big [2\Delta_i^mM\Delta_i^m\pv +(\Delta_i^m\pv)^2\big ]
\nonumber\\[2mm]
&+&\,\frac{1}{2}\sum_{i\in I_2(m)}\big [\phi''(y_i^m)
-\phi''(X(t_{i-1}^m))](\Delta_i^mM)^2=:\sum_{l=1}^5S_l(m),
\nonumber
\end{eqnarray}
where $\Delta_i^mf:=f(t_i^m)-f(t_{i-1}^m)$, $y_i^m:=X(t_{i-1}^m)
+\theta_i^m\Delta_i^mX\subset [-K,K]$ and $\theta_i^m\in (0,1)$
for $i=1,\dots,n(m)$.
By Theorem \ref{Estochint} and Theorem \ref{variation}.\ref{approximation}, 
we have
\begin{eqnarray*}
S_1(m)&=&\sum_{i=1}^{n(m)}\phi'(X(t_{i-1}^m))\Delta_i^mM
+\sum_{i=1}^{n(m)}\phi'(X(t_{i-1}^m))\Delta_i^m\pv\\
&\to&\big ((\phi'{\circ}X)_{-}\cdot M\big )(t)
+(RRS)\int_0^t(\phi'{\circ}X)_{-}\,d\pv
=\big ((\phi'{\circ}X){\lsi}X\big )(t)
\end{eqnarray*}
in probability as $m\to\infty$.
By Theorem II.30 in \cite[p.\ 69]{PP} and by definition of
the set $I_1(m)$, we have
\begin{eqnarray}
S_2(m)&=&\frac{1}{2}\sum_{i=1}^m\phi''(X(t_{i-1}^m))(\Delta_i^mM)^2
-\frac{1}{2}\sum_{i\in I_1(m)}\phi''(X(t_{i-1}^m))(\Delta_i^mM)^2
\nonumber\\
&\to&\frac{1}{2}(RRS)\int_0^t(\phi''{\circ}X)_{-}\,d[M]
-\frac{1}{2}\sum_{\mu}\big\{(\phi''{\circ}X)_{-}(\Delta^{-}M)^2
\big\}\label{5ps}
\end{eqnarray}
in probability as $m\to\infty$.
Since almost all sample functions of $X$ are {\cadlag},
by definition of $I_1(m)$ and because $\cup_m\lambda_m$ is dense in 
$[0,t]$, it follows that
\beq
S_3(m)\,\to\,\sum_{\mu}\big\{\Delta^{-}(\phi{\circ}X)-(\phi'{\circ}X)_{-}
\Delta^{-}X\big\}\label{6ps}
\eeq
almost surely as $m\to\infty$.
Notice that for each $i\in I_2(m)$, $[t_{i-1}^m,t_i^m]\subset
[s_{j-1},s_j)$ for some $j=1,\dots,k$.
Then by H\"older's inequality, and using the second condition in (\ref{9ps}),
we get the bound
\begin{eqnarray*}
|S_4(m)|&\leq&\frac{1}{2}\|\phi''\|_{\infty}
\Big\{2s_2(M;\lambda_m)^{1/2}\Big (\sum_{i\in I_2(m)}(\Delta_i^m\pv
)^2\Big )^{1/2}+\sum_{i\in I_2(m)}(\Delta_i^m\pv)^2\Big\}\\
&\leq &\frac{\sqrt{\epsilon}}{2}\|\phi''\|_{\infty}
\big\{2s_2(M;\lambda_m)^{1/2}+1\big\}.
\end{eqnarray*}
Finally using the first condition in (\ref{9ps}), we get the bound
$$
|S_5(m)|\leq\frac{1}{2}\max_{i\in I_2(m)}\big |\phi''(y_i^m)
-\phi''(X(t_{i-1}^m))\big |s_2(M;\lambda_m)
\leq \frac{\epsilon}{2} s_2(M;\lambda_m).
$$
Since the sums $S_1(m)$, $S_2(m)$ and $S_3(m)$ converge in probability
as $m\to\infty$, due to identity (\ref{4ps}), it follows that so does the sum
$S_4(m)+S_5(m)$.
Let $R(\epsilon)$ be the limit of $S_4(m)+S_5(m)$.
Notice that bounds of $S_4(m)$ and $S_5(m)$ also converge in
probability as $m\to\infty$.
Thus letting $\epsilon\downarrow 0$ it follows that $R(\epsilon)
\to 0$ in probability, while the sums over $\mu$ in (\ref{5ps})
and (\ref{6ps}) converge to the corresponding sums over $(0,t]$
by (\ref{7ps}) and (\ref{8ps}), respectively.
By Lemma \ref{function}.\ref{bv}, we have that
$$
(RRS)\int_0^t(\phi''{\circ}X)_{-}\,d[M]
=(RS)\int_0^t(\phi''{\circ}X)_{-}\,d[M]^c
+\sum_{(0,t]}\big\{(\phi''{\circ}X)_{-}(\Delta^{-}M)^2\big\}.
$$
Therefore equality (\ref{10ps}) holds.
The proof of Theorem \ref{Ito} is complete. 
\qed\end{proof}

By the preceding theorem and by Theorem \ref{Estochint} we have:

\begin{cor}
The class of all $(2-\epsilon)$-semimartingales is closed by taking
a composition with a $C^2$ class function.
\end{cor}

Now we turn to the linear stochastic $\lsi$-integral equation.
Let $X$ be a \nsem-semimartingale.
Since almost surely $\sum_{(0,t]}\{\Delta^{-}X\}^2\leq [X](t)
<\infty$ for each $t>0$, by Proposition \ref{function}.\ref{gamma}
and Lemma \ref{function}.\ref{bvofV}, the product
$$
V(X;t):=\prod_{(0,t]}(1+\Delta^{-}X)e^{-\Delta^{-}X},
\qquad t>0,
$$
converges absolutely for almost all sample functions of $X$.
The {\it stochastic Dol\'eans exponential} ${\cal E}(X)
=\{{\cal E}(X;t)\colon\,t\geq 0\}$ is then defined by
$$
{\cal E}(X;t):=\exp\big\{X(t)-X(0)-\frac{1}{2}[X]^c(t)\big\}
V(X;t),\qquad t\geq 0,
$$
where $V(X;0):=1$.
Suppose that for some $\lambda\in\Lambda [0,T]$, $0<T<\infty$,
almost all sample functions of $X$ have  the quadratic 
$\lambda$-variation on $[0,T]$.
Then by Proposition \ref{consistent}, for almost all $\omega\in
\Omega$, 
$$
[X]^c(t,\omega)=[X(\cdot,\omega)]_{\lambda}^c(t)
\quad\mbox{and}\quad
{\cal E}(X;t)(\omega)={\cal E}_{\lambda,0}(X(\cdot,\omega);t)
$$
for $0\leq t\leq T$.
That is the stochastic Dol\'eans exponential for $X$ agree
with the Dol\'eans exponential for sample functions of $X$
considered in Section \ref{doleansexp}.

For a semimartingale $X$, Dol\'eans-Dade \cite{DD} proved that
${\cal E}(X)$ is the unique semimartingale which satisfies
the linear stochastic integral equation
$$
Y=1+ Y_{-}\cdot X.
$$
Next this result is generalized to \nsem-semimartingales $X$.

\begin{thm}\label{1LSIE}
Let $X$ be a \nsem-semimartingale.
The stochastic Dol\'eans exponential ${\cal E}(X)$ is the unique 
\nsem-semimartingale which satisfies the linear stochastic
$\lsi$-integral  equation{\rm :} 
\beq\label{LSIE}
Y=1+ Y{\lsi}X.
\eeq
\end{thm}

\begin{proof}
To begin with we show that ${\cal E}(X)$ is a \nsem-semimartingale
and satisfies equation (\ref{LSIE}).
Let $\bar X=\{\bar X(t)\colon\,t\geq 0\}$ and $V=\{V(t)\colon\,t\geq 0\}$ 
be stochastic processes defined by 
\beq\label{2LSIE}
\bar X(t):=X(t)-X(0)-\frac{1}{2}[X]^c(t)
\eeq
and $V(t):=V(X;t)$ for $t\geq 0$.
By Lemma \ref{function}.\ref{bvofV}, $V=\{V(t)\colon\,t\geq 0\}$
is a {\cadlag} stochastic process with almost all pure jump sample
functions and locally of bounded variation.
Since $[X]^c$ has continuous sample functions locally of bounded variation,
$\bar X$ and $V$ are \nsem-semimartingales.
Let $\phi (u,v):=e^uv$ for $u,v\in\RR$ .
Then ${\cal E}(X)=\phi (\bar X,V)$, $(\bar X,V)$ is a $2$-tuple 
\nsem-semimartingale
and $\phi\colon\,\RR^2\mapsto\RR$ is a $C^2$ class function.
By Theorem \ref{Ito}, ${\cal E}(X)$ is a \nsem-semimartingale and 
\begin{eqnarray*}
{\cal E}(X) &=&1 +{\cal E}(X){\lsi}\bar X+
e^{\bar X}{\lsi}V
+\frac{1}{2}(RS)\int {\cal E}(X)_{-}\,d[X]^c\\[2mm]
& &+\sum\Big\{\Delta^{-}{\cal E}(X)-{\cal E}(X)_{-}
\Delta^{-}X-e^{\bar X_{-}}\Delta^{-}V\Big\}\\[2mm]
&=& 1+{\cal E}(X){\lsi}X+e^{\bar X}{\lsi}V
-\sum \Big\{e^{\bar X_{-}}\Delta^{-}V\Big\}.
\end{eqnarray*}
The last equality holds by Proposition \ref{bilinear} and because
$\Delta^{-}{\cal E}(X)={\cal E}(X)_{-}\Delta^{-}X$ by (\ref{2bvofV}).
Since almost all sample functions of $V$ are pure jump functions of 
locally bounded variation, the equality
$$
e^{\bar X}{\lsi}V
=\Big\{(RRS)\int_0^te^{\bar X_{-}}\,dV\colon\,t\geq 0\Big\}
=\Big\{\sum_{(0,t]}e^{\bar X_{-}}\Delta^{-}V\colon\,t\geq 0\Big\}
=\sum \Big\{e^{\bar X_{-}}\Delta^{-}V\Big\}
$$
holds almost surely.
Therefore ${\cal E}(X)$ is indistinguishable from $1+{\cal E}(X){\lsi}X$,
that is, the stochastic Dol\'eans exponential ${\cal E}(X)$
satisfies equation (\ref{LSIE}).

Now we prove that ${\cal E}(X)$ is the unique solution 
in the class of all \nsem-semimartingales.
Let $Y=\{Y(t)\colon\,t\geq 0\}$ be a \nsem-semimartingale
and a solution of (\ref{LSIE}). 
Let $V=\{V(t)\colon\,t\geq 0\}$ be a stochastic process defined by
$V(t):=Y(t)e^{-\bar X(t)}$, $t\geq 0$, where $\bar X(t)$ is defined
by (\ref{2LSIE}), and let $\psi (u,v):=e^{-u}v$ for $u,v\in\RR$.
Then $V=\psi (\bar X,Y)$, $(\bar X,Y)$ is a $2$-tuple 
\nsem-semimartingale, and $\psi\colon\,\RR^2
\mapsto\RR$ is a $C^2$ class function.
It is enough to prove that $V$ is indistinguishable from the
product stochastic process 
\beq\label{4LSIE}
\big\{V(X;t)\colon\,\,\, t\geq 0\big\}=
\Big\{\prod_{(0,t]}(1+\Delta^{-}X)e^{-\Delta^{-}X}
\colon\,\,\, t\geq 0\Big\}.
\eeq
To this aim we again apply Theorem \ref{Ito}, which in this case
yields the first equality
\begin{eqnarray}
V &=& 1-V{\lsi}\bar X+e^{-\bar X}{\lsi} Y
+\frac{1}{2}(RS)\int V_{-}d[X]^c
-(RS)\int e^{-\bar X_{-}}\,d[X,Y]^c\nonumber\\[2mm]
& &+\sum\Big\{\Delta^{-}V+V_{-}\Delta^{-}X
-e^{-\bar X_{-}}\Delta^{-}Y\Big\}\nonumber\\[2mm]
&=& 1-V{\lsi} X+\frac{1}{2}(RS)\int V_{-}\,d[X]^c +V{\lsi}X
+\frac{1}{2}(RS)\int V_{-}d[X]^c
-(RS)\int V_{-}\,d[X]^c\nonumber\\[2mm]
& &+\sum\Big\{\Delta^{-}V+V_{-}\Delta^{-}X
-V_{-}\Delta^{-}X\Big\}
= 1+\sum\Delta^{-}V.\label{3LSIE}
\end{eqnarray}
The second equality follows by Proposition \ref{bilinear}, Theorem
\ref{associativity} and Proposition \ref{1qcsp}.
By (\ref{3LSIE}), almost all sample functions of $V$ are 
right-continuous and pure jump functions of locally bounded variation.
Since $Y$ is a solution to the equation (\ref{LSIE}), we have
the equality
$$
\Delta^{-}V=e^{-\bar X_{-}-\Delta^{-}X}\big (Y_{-}+\Delta^{-}Y\big )
-e^{-\bar X_{-}}Y_{-}=V_{-}\Big\{\big (1+\Delta^{-}X\big )
e^{-\Delta^{-}X}-1\Big\}
$$
valid almost surely.
For $t>0$, let
$$
W(t):=\sum_{(0,t]}\Big\{\big (1+\Delta^{-}X\big )
e^{-\Delta^{-}X}-1\Big\}
$$
and $W(0):=0$.
Due to relation (\ref{7bvofV}), $W=\{W(t)\colon\,t\geq 0\}$ is a stochastic 
process with almost all {\cadlag} and pure jump sample functions 
locally of bounded variation.
This in conjunction with (\ref{3LSIE}), implies that $V$ is a
solution to the linear stochastic integral equation 
$$
V=1+(RRS)\int V_{-}\,dW=1+(LY)\int V\,dW.
$$
By Theorem \ref{function}.\ref{DN}, the solution is unique in the class of all
stochastic processes with regulated and right-continuous sample functions.
Thus $V$ is indistinguishable from the product stochastic process
(\ref{4LSIE}).
The proof of Theorem \ref{1LSIE} is complete.
\qed\end{proof}

\section{Concluding remarks}

The initial motivation for introducing the class of \nsem-semimartingales
was the wish to make a relation between calculi based on
the $p$-variation and martingale properties more transparent.
Due to time and space constrains we touch only most elementary and
basic aspects of the \nsem-semimartingale calculus.

In Stochastic Analysis, it is a tradition (with some exceptions)
to assume that stochastic processes are  {\cadlag}.
However it would be interesting to develop a calculus for stochastic
processes whose sample functions are regulated functions.
Such stochastic processes have different names in different books.
Gihman and Skorohod \cite[Chapter III, \S 4]{GS74} use the name
stochastic process without discontinuities of the second kind,
while Dellacherie and Meyer \cite[IV.20]{DM82} use the name: 
functions which are free of oscillatory discontinuities.
In these texts one can find conditions of existence of regulated 
stochastic processes.
Concerning this assumption Gihman and Skorohod \cite[p.\ 174]{GS74} notice 
that when dealing with regulated functions one does not distinguishes 
between two functions having at each point the same left and right limits (sic!).
Therefore it is natural to choose a certain convention concerning
the value of such a function at the discontinuity point.
Gihman and Skorohod \cite{GS74} concider the space of regulated functions
which are continuous either from the left or from the right at each point.
However from the point of view of applications it seems useful to
have a calculus which makes no such restrictions in advance.

\chapter{Stock price modelling in continuous time}\label{modelling}
\setcounter{thm}{0}

\vspace*{0.2truein}
\begin{quotation}{\footnotesize
In practice, both market practitioners and regulators are aware
that equilibrium is an illusion.
It is rare to find a field in which theory and practice are so
far apart, leaving ample room for alchemy and other forms of magic.

G.\ Soros. The Crisis of Global Capitalism [Open Society 
Endangered]. Public Affairs, New York, 1998, p.\ 41.} 
\end{quotation}
\vspace*{0.2truein}

In continuous time finance the uncertainty related to a future stock 
price behaviour is modelled by introducing a probability space with
a (semi)martingale $X=\{X_t\equiv X(t)\colon\,0\leq t\leq T\}$.
A stock price dynamics $P=\{P_t\equiv P(t)\colon\,0\leq t\leq T\}$
is then described by means of a solution of a stochastic differential equation
driven by $X$.
Usually, it is justifiable by the fact that the driving process (the return) 
satisfies a one of several variants  of the Random Walk Hypotheses.                         
In this chapter we continue developing a continuous time stock price
model undertaken in \cite{RN00} without using Stochastic Analysis 
constructions and suggest a different a motivation.

\section{Evolutionary asset pricing model}

\paragraph*{Model introduction.}
Prevailing stock market models are based on Probability Theory.
In financial applications, a probability space $(\Omega,{\cal F},\Pr)$
is considered as a set of possible "states of the world", 
where the probability $\Pr$ attach suitable weights on admissable events 
in the $\sigma$-algebra ${\cal F}$.
In continuous time stock price models, a set of possible "states of the world" 
can be expressed by a set of possible price trajectories in continuous time.
Under the Random Walk Hypotheses, typical distributions of price 
trajectories have support in a set of functions having unbounded variation,
and so we need a calculus applicable to such functions.
Stochastic Analysis provides several forms of calculi, such as the 
semimartingale theory or the Malliavin calculus, applicable to families
of functions with unbounded variation, which are sample functions
of suitable stochastic processes.
The fact that possible "states of the world",
that is, possible stock price trajectories, are indistinguishable by  
intrinsic constructions of Stochastic Analysis, such as
the It\^o integral, is an imperfection of this calculus.
Therefore it is natural to look for alternative approaches to
continuous time stock price models which allow a consideration of
a single price trajectory, a single trading strategy and so on
(see Willinger and Taqqu \cite{WandT89}  for a thorough discussion). 

An interpretation of theoretical constructions, such as a set of
possible "states of the world", is an important aspect of a model building.
In fact, an interpretation renders a formal theory into a model.
Here by a formal theory we mean an area of Mathematics, and not
an Economic Theory, and therefore a resulting model do not provide
economically meaningful explanations or predictions. 
Rather our model may serve as a motivation for further developments
both in an area of Mathematics and in an area of Financial Theory.
As it is usual in social sciences, a suggested interpretation includes
a subjective character of stock price behaviour.
Due to a human (scientific) activity, stock markets differ fundamentally
from usual subjects of natural sciences.
While the latter do not depend on existing theories, stock markets 
are very sensitive to various scientific interpretations.
Therefore it is common to incorporate a human factor into market models.
For example, in stochastic stock price models, a flow of $\sigma$-algebras 
of events is interpreted as a flow of information available for market traders.
As compared to stock price models motivated by the Random Walk 
Hypothesis, our interpretation concerns a tranformation relating
a price and its return, or equivalently, a form of a return.
Such an interpretation has no meaning in models based on
Stochastic Analysis due to the aforementiond indistinguishibility of
price trajectories. 

For simplicity we consider a market which consists of a single stock and an 
investor-specialist to be called a participant of the market.
The basic assumption of our model is that participant's thinking effects and 
is effected by a stock price behaviour.
In other words, we want to model a connection between a stock and market 
participants' expectations of the stock.
Due to this connection a stock price and participants'
expectations becomes equally (un)predictable.
To say informally, the model consists of an ordered pair $(M,A)$ 
of two sets of functions defined on $[0,T]$, and a one-to-one mapping $R$ 
from $M$ onto $A$, called a \emph{reflexivity mapping for $(M,A)$}.
Elements of the set $M$ are called \emph{price scenarios} 
$P=\{P(t)\colon\,0\leq t\leq T\}$ during the time period $[0,T]$,
and the stock is identified with the set $M$.
Participants' expectations are represented by the set $A$
of functions $X=\{X(t)\colon\,0\leq\ t\leq T\}$, each of which is called a 
\emph{return} of a price scenario $P\in M$ assigned by the reflexivity
mapping $R$.
Different elements of the pair $(M,A)$ may be
interpreted as different states of the world.

In continuous time finance it is customary to define a price evolution
by a unique (strong) solution $P=\{P_t\equiv P(t)\colon\,0\leq t\leq T\}$
of a stochastic differential equation 
\beq\label{SDE}
dP_t=\mu (t,P_t)\,dt+\sigma (t,P_t)\,dX_t,\qquad 0\leq t\leq T, 
\eeq
with the initial condition $P_0\!=\!x$,
driven by a semimartingale $X\!=\!\{X_t\!\equiv\! X(t)\colon\,0\leq t\leq T\}$.
Following the suggested construction, we take the sets $M$ and $A$ to be
supports of the distributions of $P$ and $X$, respectively.
The reflexivity mapping between $M$ and $A$ then could be
defined provided one can solve (\ref{SDE}) path by path.
Typically such a mapping may not be possible 
to express using intrinsic constructions of Stochastic Analysis.
Unique solutions of (\ref{SDE}) give rise to the solution mapping, 
sometimes called the It\^o mapping, defined on a class of driving
processes $X$.
But except for special cases, the It\^o mapping need not be given by 
a mapping defined on sample functions of the driving process $X$
because the It\^o integral is not defined for sample functions.
This aspect is important in econometric analysis of stock price
models, and sometimes makes simple exponential models preferable 
over the models given by solutions of (\ref{SDE}).
The pathwise aspect of the stochastic calculus
is also related to the problem of robustness
as discussed e.g.\ in Section 5.C of Bouleau and L\'epingle \cite{BandL94}.
By the above introductory remarks,  the reflexivity mapping, if exists, 
need not be continuous in a usual sense.
However roughly speaking our reflexivity mapping may be considered
as a substitute of the stochastic differential equation (\ref{SDE})
in continuous time models based on Stochastic Analysis.

In (\ref{SDE}) taking $X$ to be a standard Brownian motion $B$, 
$\sigma (t,u)\equiv u$, $\mu\equiv 0$ and the initial condition $x=1$, 
we get a simple case of the Black and Scholes model with the price evolution 
$P_B(t):=\exp\{B(t)-t/2\}$ for $0\leq t\leq T$.
Let $A$ be a set of continuous sample functions of $B$. 
Since the mapping $b\mapsto \{\exp\{b(t)-t/2\}\colon\,0\leq t\leq T\}$,
$b\in A$, is invertible on its range, say $M$, the set $(M,A)$ gives an 
example of a pair of sets with a reflexivity mapping.
If $A$ is a subset of the set of all right-continuous functions
on $[0,T]$ of bounded $p$-variation for some $0<p<2$ then by the
results of Dudley and Norvai\v sa \cite[Part II]{DNa}, some
natural reflexivity mappings can be extended to analytic mappings.

\paragraph*{Reflexivity mapping.}
More specifically, the reflexivity mapping will be modelled by means of a 
duality mapping between multiplicative and additive functions on
the simplex of an extended interval (see Section \ref{add&mult}).
Recall that a function $\pi=\{\pi (s,t)\colon\,(s,t)\in S\loc 0,T\roc \}$ 
is mutiplicative if $\pi(s,r)=\pi(s,t)\pi(t,r)$ for each $s, t, r\in\loc 0,T\roc$
such that $0\leq s\leq t\leq r\leq T$, 
and $\pi (t,t)=1$ for $t\in \loc 0,T\roc$.
Also, a function $\mu=\{\mu(s,t)\colon\,(s,t)\in S\loc 0,T\roc\}$ is
addative if $\mu(s,r)=\mu(s,t)+\mu(t,r)$ for each $s,t,r\in\loc 0,T\roc$
such that $0\leq s\leq t\leq r\leq T$,
and $\mu(t,t)=0$ for $t\in \loc 0,T\roc$.
For a given regulated function $f$ on $[0,T]$, $\mu (f)$ is additive
and upper continuous function on $S\loc 0,T\roc$ defined by
$\mu (f;s,t):=f(t)-f(s)$ for $(s,t)\in S\loc 0,T\roc$, and $\pi (f)$
is multiplicative and upper continuous function on $S\loc 0,T\roc$
defined by $\pi (f;s,t):=f(t)/f(s)$ for $(s,t)\in S\loc 0,T\roc$.

\begin{defn}
{\rm Let ${\cal M}={\cal M}[0,T]$ be a set of multiplicative functions, and let
${\cal A}={\cal A}[0,T]$ be a set of additive functions.
If there is a one-to-one mapping $\cal R$ from ${\cal M}$ onto ${\cal A}$,
then we say that $({\cal R},{\cal M},{\cal A})$ is a \emph{duality triple} over
$[0,T]$, and ${\cal R}$ is called a \emph{reflexivity mapping}.  }
\end{defn}

It will be seen from what follows that a reflexivity mapping for a given
${\cal M}$ and ${\cal A}$ may \emph{not} be unique.
We start with a simple example.
Let ${\cal M}[0,T]$ be the set of all multiplicative functions on 
$S\loc 0,T\roc$  which are positive, and let ${\cal A}[0,T]$ be the set 
of all real-valued additive functions on $S\loc 0,T\roc$.
Then taking the logarithm ${\cal R}_{\log}\pi (s,t):=\log \pi (s,t)$, $(s,t)\in
S\loc 0,T\roc$, induces the reflexivity mapping ${\cal R}_{\log}$ 
from ${\cal M}[0,T]$ onto ${\cal A}[0,T]$.
In the next subsection we describe a different reflexivity mapping
acting between subsets of ${\cal M}[0,T]$ and ${\cal A}[0,T]$.

In the cases considered below a natural domain of a reflexivity mapping
is quite a large set.
Therefore for modelling purposses, it is reasonable restrict a consideration 
to pairs of subsets defined as follows:

\begin{defn}\label{R-system}
{\rm Let $({\cal R},{\cal M},{\cal A})$ be a duality triple on $[0,T]$,
and let $M$, $A$ be two sets of regulated functions on $[0,T]$ such that
$f(0)=1$ for $f\in M$, $g(0)=0$ for $g\in A$,
$$
\pi (M):=\{\pi (f)\colon\,f\in M\}\subset {\cal M}\quad\mbox{and}
\quad\mu (A):=\{\mu (g)\colon\,g\in A\}\subset {\cal A}.
$$
If ${\cal R}(\pi (M))=\mu (A)$ then we say that the pair $(M,A)$ is an 
\emph{${\cal R}$-system} and the mapping $R\colon\,M\mapsto A$ 
defined by $R(f):=R_{{\cal R}\pi(f)}$, $f\in M$, is called the 
\emph{reflexivity mapping for $(M,A)$}, 
where $R_{\mu}(t):=\mu (0,t)$, $t\in [0,T]$ is the right 
distribution  function of $\mu$ restricted to $[0,T]$.}
\end{defn}

\begin{figure}
\begin{picture}(100,100)
\linethickness{0.2mm}

\put(40,10){$\pi (M)$}

\put(75,12){\vector(1,0){40}}

\put(90,20){${\cal R}$}

\put(130,10){$\mu (A)$}

\put(135,30){\vector(0,4){20}}

\put(90,70){$R$}

\put(45,60){$M$}

\put(75,62){\vector(1,0){40}}

\put(50,50){\vector(0,-4){20}}

\put(35,35){$\pi$}

\put(133,60){$A$}

\put(150,35){$R_{{\cal R}\pi}$}

\put(210,35){and}

\put(290,10){$\pi (M)$}

\put(365,12){\vector(-1,0){40}}

\put(340,20){${\cal R}^{-1}$}

\put(380,10){$\mu (A)$}

\put(385,50){\vector(0,-4){20}}

\put(340,70){$R^{-1}$}

\put(295,60){$M$}

\put(365,62){\vector(-1,0){40}}

\put(300,30){\vector(0,4){20}}

\put(260,35){$R_{{\cal R}^{-1}\mu}$}

\put(383,60){$A$}

\put(400,35){$\mu$}

\end{picture}
\caption{A reflexivity mapping for $(M,A)$ and its inverse.}\label{diagram}
\end{figure}

Let $({\cal R},{\cal M},{\cal A})$ be a duality triple on $[0,T]$,
and let $(M,A)$ be an $\cal R$-system.
By Corollaries \ref{variation}.\ref{projad} and \ref{variation}.\ref{projmult},
the reflexivity mapping $R$ for $(M,A)$ is one-to-one with the
inverse $R^{-1}\colon\,A\mapsto M$, and diagrams of both mappings are 
given by Figure \ref{diagram}.
For $(P,X)\in (M,A)$ and for each $0\leq s\leq t\leq T$, we have
\beq\label{assumption}
\frac{P(t)}{P(s)}=\pi (P;s,t)=({\cal R}^{-1}\mu (X))(s,t)
\quad\mbox{and}\quad
X(t)-X(s)=\mu (X;s,t)=({\cal R}\pi (P))(s,t).
\eeq
As it was said earlier a stock is identified with a set $M$
of possible price scenarious, and market participants'
expectations are represented by a set of returns $A$.
This model may be considered as a form of the ``Concept of Reflexivity'' 
introduced and advocated by Soros \cite[p.\ 42]{GS94}. 
The first relation in (\ref{assumption}) means that investors
make decisions in terms of the percentage changes in prices rather
than in terms of absolute prices.
In the context of finance this may be called the \emph{property of relative
price changes}.
This idea was first formulated and discussed in the literature by
Osborne \cite{MFMO59} and Samuelson \cite{PAS73} (see also Bernstein 
\cite[pp. 103-105]{PLB92}).
It differs in an important way from Bachelier's hypothesis
concerning differences of a stock price values.
The second relation in (\ref{assumption}) means that returns $X$
allow for aggregation over different time periods, the property 
preferable when one is
trying to test various assumptions made on the return $X$.
Let us call it as the \emph{aggregational property of returns}.
The two "properties" are so important for the present interpretation that
in a formal theory of a market they should be called as axioms.

\paragraph*{Evolutionary model.}
Here we describe a reflexivity mapping induced by the $\lambda$-generator
of a regulated evolution introduced in Section \ref{e-representation}.

\begin{defn}\label{evol-model}
{\rm Let $\lambda=\{\lambda_m\colon\,m\geq 1\}\in\Lambda [0,T]$.
A set ${\cal M}={\cal M}[0,T]$ of multiplicative interval functions on 
$S\loc 0,T\roc$ will be called a \emph{set of $\lambda$-evolutions},
and a mapping ${\cal R}_{\lambda}$ on ${\cal M}$ will be called
a \emph{$\lambda$-reflexivity mapping} if $(a)$ and $(b)$ hold,  where
\begin{enumerate}
\item[$(a)$] for each $\pi\in{\cal M}$ there is a function 
${\cal R}_{\lambda}(\pi):=\mu$ on $S\loc 0,T\roc$ such that for each 
$(s,t)\in S\loc 0,T\roc $,
\beq\label{Dol-exp}
\mu (s,t)=\lim_{m\to\infty} S(\pi;\lambda_m\Cap\loc s,t\roc ),
\eeq
where $S(\pi;\kappa):=\sum_{i=1}^n[\pi (t_{i-1},t_i)-1]$ for a partition 
$\kappa=\{t_i\colon\,i=0,\dots,n\}$;
\item[$(b)$] for each $\pi_1,\pi_2\in {\cal M}$, if
${\cal R}_{\lambda}(\pi_1)={\cal R}_{\lambda}(\pi_2)$ then
$\pi_1=\pi_2$.
\end{enumerate}}
\end{defn}

It is clear that a function $\mu$  arising in condition $(a)$ of the
preceding definition,  if it exists, is additive and unique, 
and so ${\cal R}_{\lambda}$ is a mapping on  $\cal M$ with
values in a class of additive interval functions on $[0,T]$.
By $(b)$, ${\cal R}_{\lambda}$ is a one-to-one mapping onto
the range  ${\cal A}:={\cal R}_{\lambda}({\cal M})$.
Therefore, $({\cal R}_{\lambda},{\cal M},{\cal A})$ is a duality triple 
over $[0,T]$ provided a set $\cal M$ of $\lambda$-evolutions 
and a $\lambda$-reflexivity mapping ${\cal R}_{\lambda}$ on 
$\cal M$ exist.

To construct a set of  $\lambda$-evolutions we use the set 
$L_{\lambda}[0,T]$ defined by (\ref{class-L}).
This and related set $E_{\lambda}[0,T]$ defined by (\ref{class-E}),
were used in Section \ref{e-representation} to solve the evolution 
representation problem for functions having the quadratic $\lambda$-variation.
For $\lambda\in\Lambda [0,T]$, let
\beq\label{class-M}
{\cal M}_{\lambda}={\cal M}_{\lambda}[0,T]:=\pi \big (L_{\lambda}[0,T]\big )
:=\big\{\pi (f)\colon\,f\in L_{\lambda}[0,T]\big\},
\eeq
\beq\label{class-A}
{\cal A}_{\lambda}={\cal A}_{\lambda}[0,T]:=\mu \big (E_{\lambda}[0,T]\big )
:=\big\{\mu (g)\colon\,g\in E_{\lambda}[0,T]\big\}.
 \eeq
For each $f\in L_{\lambda}[0,T]$, let 
\beq\label{l-reflex}
{\cal R}_{\lambda}\pi (f):=(LC)\smallint f^{-1}\,d_{\lambda}f,
\eeq
where the right side is defined by Proposition \ref{function}.\ref{Elogarithm}.
By Proposition \ref{function}.\ref{regul-evol2}, each 
$\pi (f)$, $f\in L_{\lambda}[0,T]$, is a regulated evolution on $[0,T]$,
and condition $(a)$ of Definition \ref{evol-model} holds provided $\pi (f)$ 
has the $\lambda$-generator of Definition \ref{function}.\ref{regul-evol3}. 
The latter holds by Theorem \ref{function}.\ref{generator1}, and so
condition $(a)$ of Definition \ref{evol-model} holds for the set 
${\cal M}_{\lambda}$.
While condition $(b)$ of Definition \ref{evol-model} holds for
${\cal M}_{\lambda}$ by Theorem \ref{function}.\ref{ratio}.
Moreover
by Propositions \ref{function}.\ref{Elogarithm}  and \ref{function}.\ref{property1},
${\cal R}_{\lambda}({\cal M}_{\lambda})\subset {\cal A}_{\lambda}$, and
by Theorem \ref{function}.\ref{Ulogarithm}, ${\cal R}_{\lambda}
({\cal M}_{\lambda})\supset {\cal A}_{\lambda}$.
We thus proved the following:

\begin{cor}
For $\lambda\in\Lambda [0,T]$, ${\cal M}_{\lambda}$ is a 
$\lambda$-evolutionary set, ${\cal R}_{\lambda}$ is a 
$\lambda$-reflexivity mapping from ${\cal M}_{\lambda}$
onto ${\cal A}_{\lambda}$, and
$({\cal R}_{\lambda},{\cal M}_{\lambda},{\cal A}_{\lambda})$ is
a duality triple.
\end{cor}

In our financial interpretation, ${\cal M}_{\lambda}$ is the set of price scenarios,
${\cal A}_{\lambda}$ is the set of returns  of price scenarios assigned
by the reflexivity mapping ${\cal R}_{\lambda}$.
Therefore given $\lambda\in\Lambda [0,T]$, the duality triple 
$({\cal R}_{\lambda},{\cal M}_{\lambda},{\cal A}_{\lambda})$
is called the \emph{$\lambda$-evolutionary asset pricing model}.
If $\mu\in{\cal A}_{\lambda}$ then $\mu={\cal R}_{\lambda}\pi (f)$
for some $f\in L_{\lambda}[0,T]$, and hence by Theorem 
\ref{function}.\ref{ratio}, the inverse to ${\cal R}_{\lambda}$ is given by
\beq\label{l-inverse}
{\cal R}_{\lambda}^{-1}\mu=\prodi (1+d_{\lambda}({\cal L}_{\lambda}f))
=\pi (f),
\eeq
where ${\cal L}_{\lambda}f(t)=\mu (0,t)=(LC)\smallint_0^tf^{-1}\,d_{\lambda}f$
for $t\in [0,T]$.
Next is a special case of Definition \ref{R-system}.

\begin{defn}\label{E-system}
{\rm Let $({\cal R}_{\lambda},{\cal M}_{\lambda},{\cal A}_{\lambda})$
be the $\lambda$-evolutionary asset pricing model for some 
$\lambda\in\Lambda [0,T]$,
and let $P$, $X$ be two sets of regulated functions on $[0,T]$ such that
$f(0)=1$ for $f\in P$, $g(0)=0$ for $g\in X$,
$P\subset L_{\lambda}[0,T]$ and $X\subset E_{\lambda}[0,T]$.
The pair $(P,X)$  is called a \emph{$\lambda$-evolutionary system}
if it is an ${\cal R}_{\lambda}$-system, that is, if 
${\cal R}_{\lambda}(\pi (P))=\mu (X)$. 
Moreover, the reflexivity mapping for $(P,X)$, say $R_{\lambda}$, 
is called the \emph{$\lambda$-reflexivity mapping for $(P,X)$}.
The set $P$ is called a \emph{stock of the $\lambda$-evolutionary
model}, or simply a stock, and each element of $P$ is called
a \emph{stock price scenario}.}
\end{defn}

\begin{rem}\label{convention}
{\rm For notation convenience $P$ may denote a set (= a stock) and
an element of a set (= a stock price scenario).
If there is a possibility of a confusion we also may use a notation
for a stock price scenario different from a stock.
The same rule applies to a return $X$.
Typical examples of a $\lambda$-evolutionary system are pairs of
stochastic process, and so our rule conforms to the tradition in
Stochastic Analysis to denote a stochastic process and a sample function
by the same letter.}
\end{rem}

Let $(P,X)$ be a $\lambda$-evolutionary system for some $\lambda\in
\Lambda [0,T]$.
By (\ref{l-reflex}), the $\lambda$-reflexivity mapping for $(P,X)$
is ${\cal L}_{\lambda}$, that is, the mapping from 
$P$ onto $X$ induced by the right distribution  function of the Left Cauchy 
$\lambda$-integral.
By (\ref{l-inverse}),  the inverse of $R_{\lambda}={\cal L}_{\lambda}$
is the right distribution function ${\cal P}_{\lambda}$ of the
product $\lambda$-integral (see Figure \ref{diagram}). 
By Theorem \ref{function}.\ref{ExistLprod}, the inverse $R_{\lambda}^{-1}
={\cal P}_{\lambda}$ is also equal to the mapping induced by the forward 
Dol\'eans exponential ${\cal E}_{\lambda}(g):={\cal E}_{\lambda,0}(g)$ 
defined by (\ref{Dol1}) with $[a,b]=[0,T]$, for $g\in E_{\lambda}[0,T]$.

\begin{prop}\label{E-stock}
For $\lambda\in\Lambda [0,T]$ and $X\subset E_{\lambda} [0,T]$,
$({\cal E}_{\lambda}(X),X)$ is a $\lambda$-evolutionary system.
Also for each $g\in X$, $f:={\cal E}_{\lambda}(g)\in {\cal E}_{\lambda}(X)$ 
and $0\leq t\leq T$,
\beq\label{1E-stock}
f(t)=1+(LC)\int_0^tf\,d_{\lambda}g,
\eeq
\beq\label{2E-stock}
[f]_{\lambda}(t)=(LY)\int_0^tf^2\,d[g]_{\lambda}\qquad\mbox{and}\qquad
[f,g]_{\lambda}(t)=(LY)\int_0^tf\,d[g]_{\lambda}.
\eeq
\end{prop}

\begin{proof}
Let $\lambda\in\Lambda [0,T]$ and $X\subset E_{\lambda} [0,T]$.
For $g\in X$, by Theorem \ref{function}.\ref{qvofEa}, ${\cal E}_{\lambda}(g)$
has the quadratic $\lambda$-variation.
Also since $1+\Delta^{-}g\wedge\Delta^{+}g\gg 0$, by Lemma 
\ref{function}.\ref{gamma} and Corollary \ref{variation}.\ref{projmult},
${\cal E}_{\lambda}(g)\gg 0$.
Since $g\in X$ is arbitrary, ${\cal E}_{\lambda}(X)\subset L_{\lambda}[0,T]$.
Again let $g\in X$.
By Theorem \ref{variation}.\ref{interv} and Corollary \ref{function}.\ref{beta}, 
$\pi ({\cal E}_{\lambda}(g))=\beta_{\lambda}(g)$, and so
by Theorem  \ref{function}.\ref{ExistLprod}, $\prodi (1+d_{\lambda}g)
=\pi({\cal E}_{\lambda} (g))$.
Therefore by Theorem \ref{function}.\ref{Ulogarithm},
$$
{\cal R}_{\lambda}\pi ({\cal E}_{\lambda}(g))
=(LC)\int ({\cal P}_{\lambda}g)^{-1}\,d_{\lambda}({\cal P}_{\lambda}g)
=\mu(g).
$$
Since $g\in X$ is arbitrary, ${\cal R}_{\lambda}\pi ({\cal E}_{\lambda}(X))
=\mu (X)$, that is, $({\cal E}_{\lambda}(X),X)$ is a $\lambda$-evolutionary
system.
Finally (\ref{1E-stock}) is a consequence of Theorem \ref{function}.\ref{LIE},
and  (\ref{2E-stock}) is a consequence of Theorem \ref{function}.\ref{qvofEa}.
\qed\end{proof}

The following is a $\lambda$-evolutionary version of the  
Black and Scholes model.

\begin{exmp}\label{P_B}
{\rm Let $B$ be a standard Brownian motion on a complete probability space 
$(\Omega, {\cal F},\Pr)$.
Given $\lambda\in\Lambda [0,T]$, almost all sample functions
$b:=B(\cdot,\omega)$, $\omega\in\Omega$, of $B$ are continuous 
and have the quadratic  $\lambda$-variation.
Thus denoting by $B_{\lambda}$ the set of all such functions $b$, and denoting 
by $P_{B_{\lambda}}$ the set of all functions $\exp\{b(t)-t/2\}$, $0\leq t\leq T$, 
we have that $(P_{B_{\lambda}},B_{\lambda})$ is a $\lambda$-evolutionary 
system.
To see this notice that $P_{B_{\lambda}}=\{{\cal E}_{\lambda}(b)\colon\,
b\in B_{\lambda}\}$ and apply the preceding proposition.
Letting $N(\lambda):=\{\omega\in\Omega\colon\,B(\cdot,\omega)\not\in
B_{\lambda}\}$, we have $\cup\{N(\lambda)\colon\,\lambda\in\Lambda
[0,T]\}=\Omega$, $\Pr$-almost surely.}
\end{exmp}

\section{Almost sure approximation by a discrete time model}

Here we prove that the $\lambda$-evolutionary version of the Black
and Scholes model of Example \ref{P_B} is a limit of discrete time
binomial models.
It is boring to say once again that the point here is an almost sure convergence
because in probability variant of such approximation is well-known.

\paragraph*{Returns in discrete time.}
We start by showing that the two continuous time price dynamics 
described by the reflexivity mappings ${\cal R}_{\lambda}$ and
${\cal R}_{\log}$, accord to two discrete time models based on
different notions of returns.
Consider a discrete time model of a stock price $\{P(t)\colon\,
t=0,1,\dots,T\}$ with the {\em simple net return} $R_{net}(t)$
defined by
$$ 
R_{net}(t):=\frac{P(t)}{P(t-1)}-1,\qquad t=1,\dots,T.
$$
Also we assume that $P(0)=1$ and $P(t)>0$ for $t=1,\dots,T$.
As it is usual for discrete time models, their returns are defined at the 
right endpoints $t$ of unitary time lags $[t-1,t]$ since such returns always
represent a change on the stock price between dates $t-1$ and $t$.
To imbed a discrete time model into a continuous time model one has to
redefine the discrete time returns as follows: 
let $X(0):=0$ and $X(t)-X(t-1):=R_{net}(t)$ for  $t=1,\dots,T$.
Adding up the increments of $X$, we get the relations
\beq\label{1return}
X(t)=\sum_{s=1}^tR_{net}(s)=\sum_{s=1}^t
\frac{P(s)-P(s-1)}{P(s-1)},\qquad t=1,\dots,T.
\eeq
Let $\lambda=\{\lambda_m\colon\,m\geq 1\}\in\Lambda [0,T]$ 
with $\lambda_m=\{t_i^m\colon\,i=0,\dots,n(m)\}$.
Suppose that a price scenario $P$ has its positive values defined on the interval
$[0,T]$.
Then for each $m\geq 1$, $P$ restricted to $\lambda_m$ gives the
discrete time model with the simple net return $X_m$ such that
$$
X_m(t_i^m)-X_m(t_{i-1}^m)=\frac{P(t_i^m)}{P(t_{i-1}^m)}-1,
\qquad i=1,\dots,n(m).
$$
Relation (\ref{Dol-exp}) with
$\pi (t_i^m,t_{i-1}^m)=P(t_i^m)/P(t_{i-1}^m)$ for each $i$,
means that $\mu (t,s)$ can be approximated arbitrary 
closely by $X_m(t)-X_m(s)$ as $m\to\infty$. 
Thus the function $\mu$ satisfying (\ref{Dol-exp}), if exists, extends 
the simple net return $X$ to a continuous time setting.

The following relations in the discrete setting illustrate the duality relations 
(\ref{assumption}) used to define a reflexivity mapping ${\cal R}$.
Relations (\ref{1return}) and
\beq\label{2return}
P(t)=\prod_{s=1}^t\big [1+X(s)-X(s-1)\big ],\qquad t=1,\dots,T,
\eeq
provide a one-to-one correspondence between 
a positive price $P$ and a simple net return $X$ having jumps 
bigger than $-1$.
Instead of (\ref{2return}) one can use the equivalent relation
$$
P(t)=1+\sum_{s=1}^tP(s-1)\big [X(s)-X(s-1)\big ].
$$
For each $t=1,\dots,T$, let ${\cal L}(P)(t)$ be defined by the
right side of (\ref{1return}), and let ${\cal E}(X)(t)$ be
defined by the right side of (\ref{2return}).
Then for $t=1,\dots,T$, we have
$$
{\cal L}({\cal E}(X))(t)=\sum_{s=1}^t\Big [\frac{{\cal E}(X)(s)}
{{\cal E}(X)(s-1)}-1\Big ]=X(t)
$$
and
$$
{\cal E}({\cal L}(P))(t)=\prod_{s=1}^t\Big (1+{\cal L}(P)(s)
-{\cal L}(P)(s-1)\Big )=P(t).
$$
This duality is used in discrete time stock price models by specifying simple
net returns rather than prices (see \cite[Section 3.2]{SRP}).

In econometric analysis instead of the simple
net return, one uses the {\em log-return} $R_{log}$, defined by
$$
R_{log}(t):=\log\Big ( \frac{P(t)}{P(t-1)}\Big ),\qquad t=1,\dots,T.
$$
The log return is preferable to the simple net return 
in data analysis because of the aggregational (additivity)  property:
$$
R_{log}(t)+R_{log}(t-1)+\dots
+R_{log}(t-s+1)=\log\Big ( \frac{P(t)}{P(t-s)}\Big )
$$
valid for any $s,t\in\{0,1,\dots,T\}$, $s<t$.
Clearly such a property does not hold for the simple net return 
$R_{net}$.
In contrast, we see that if (\ref{Dol-exp}) holds then 
the continuous time analog of the simple net 
return possess the aggregational property.
In sum, a continuous time stock price model based on the reflexivity
mapping induced by the logarithm is a limit of discrete time models
based on the log-returns, while discrete time models based on
the simple net returns give rise to a continuous time models
based on the $\lambda$-reflexivity mapping to be discussed in more detail
in what follows.

\paragraph*{Approximation by the binomial model.}
The binomial model is a simple but yet very important discrete time model 
of price evolution of a single stock.
We show that the binomial price evolution approximates the 
$\lambda$-evolutionary system $(P_{B_{\lambda}},B_{\lambda})$ defined 
in Example \ref{P_B}.
To recall a simplest variant of the binomial model,
let $\epsilon_1,\epsilon_2,\dots$ be a sequence of independent
identically distributed random variables on $(\Omega,{\cal F},\Pr)$ 
such that $\Pr(\{\epsilon_1=1\})=1-\Pr (\{\epsilon_1=0\})=1/2$.
For each $t=1,2,\dots$ and for any $a\in (0,1)$, let
$$
P(t):=P(0)(1+a)^{N_t}(1-a)^{t-N_t},
$$
where $N_t:=\epsilon_1+\dots +\epsilon_t$.
According to this time evolution,  called the {\em binomial model}, 
the price either goes up by the factor $1+a$, or it goes down by the factor $1-a$.
The probability of an up move or down move is the same, and the moves over
time are independent of each other.
The binomial model is a discrete time analogue of the Black and Scholes
continuous time model.
To justify this one usually refers to a weak convergence of corresponding
returns.
Next we show that this analogy is much stronger in the sense that the 
approximation of price scenarios holds with probability $1$.

To this aim we use the random walk approximation of a standard
Brownian motion originated in Knight \cite{FBK62} and further developed in  
Section 6.2 of  R\'ev\'esz \cite{PR90} and in Szabados \cite{TS96}.
Let $\{r_1,r_2,\dots\}$ be a sequence of independent identically
distributed random variables on $(\Omega,{\cal F},\Pr)$ with
the distribution $\Pr(\{r_1=1\})=1-\Pr (\{r_1=-1\})=1/2$.
Letting $S(0):=0$ and $S(k):=r_1+\dots +r_k$ for each $k\geq 1$,
we have the {\em random walk} process $S=\{S(k)\colon\,k=0,1,\dots\}$
on the lattice of all non-negative integers.
On a probability space carrying a countable array of independent random walks
one can construct a sequence of ''twisted'' and ''shrinked" random
walks giving rise to a standard Brownian motion.
More specificaly the following holds:

\begin{thm}\label{RW}
On a rich enough probability space $(\Omega,{\cal F},\Pr)$ one can
construct stochastic processes $B$ and $\{W_m\colon\,m\geq 1\}$
such that
\begin{enumerate}
\item[$(a)$] for a random walk process $S$ and for each $m\geq 1$,
$$
\big\{W_m(k2^{-2m})\colon\,k=0,1,\dots\big\}
\stackrel{\rm f.d.d.}{=}\big\{2^{-m}S(k)\colon\,k=0,1,\dots\big\};
$$
\item[$(b)$] for each $m\geq 1$, there is a random partition 
$\{T_{m+1}(k)\colon\,k\geq 1\}$ of $[0,\infty)$ such that for each
$k\geq 1$, $W_{m+1}(T_{m+1}(k))=W_m(k2^{-2m})${\rm ;}
\item[$(c)$] $B$ is a standard Brownian motion, and for each $T>0$,
$\|B-\widehat{W}_m\|_{[0,T],\sup}\to 0$ almost surely as 
$m\to \infty$, where $\widehat{W}_m$ is a process indexed by $t\in[0,\infty)$
 obtained from $W_m$ by linear interpolation{\rm ;}
\item[$(d)$] with probability $1$, for any $T'>T>0$  and for all but 
finitely many $m$ there exists a random partition $\{\tau_m(k)\colon\,
0\leq k\leq T2^{2m}\}$ of $[0,T')$ such that $B(\tau_m(k))=W_m(k2^{-2m})$
for $0\leq k\leq T2^{2m}$, where $\max_k|\tau_m(k)-\tau_m(k-1)|\to 0$
and  $\max_k|\tau_m(k)-k2^{-2m}|\to 0$ both as $m\to\infty$.
\end{enumerate}
\end{thm}

The construction of $B$ and $\{W_m\colon\,m\geq 1\}$ satisfying
$(a)$, $(b)$ and $(c)$ is given in \cite[Section 6.2]{PR90} and
in \cite{TS96}.
Statement $(d)$ follows from Lemma 7(b) in \cite{TS96}.

Let $m\geq 1$, and let the (discrete time) simple net return $X_m$
be defined by 
$$
X_m(t):=W_m(t2^{-2m})\qquad\mbox{for $t=0,1,\dots ,T2^{2m}.$}
$$
Then by $(a)$ of the preceding theorm, $X_m(0)=0$ almost surely and
$\{X_m(t)-X_m(t-1)\colon\,t=1,\dots,T2^{2m}\}$ are independent identically 
distributed random variables on $(\Omega,{\cal F},\Pr)$ 
such that $\Pr (\{X_m(1)=2^{-m}\})=1-\Pr (\{X_m(1)=-2^{-m}\})=1/2$.
For each $k=1,\dots,T2^{2m}$, let $\epsilon_k:=1$ if 
$X_m(k)-X_m(k-1)=2^{-m}$ and let $\epsilon_k:=0$ if 
$X_m(k)-X_m(k-1)=-2^{-m}$.
Then the discrete time evolutionary
price $P_m$ is given by $P_m(0):=1$ and for each $t=1,\dots, T2^{2m}$,
$$
P_m(t)=\prod_{k=1}^{t}\big [1+W_m(k2^{-2m})-W_m((k-1)2^{-2m})\big ]
=\big (1+2^{-m}\big )^{N_t}(1-2^{-m})^{t-N_t},
$$
where $N_t:=\epsilon_1+\dots +\epsilon_t$.
Therefore for the probability space $(\Omega,{\cal F},\Pr)$ of
Theorem \ref{RW}, $\{(P_m(\cdot,\omega),X_m(\cdot,\omega))
\colon\,\omega\in\Omega\}$ are discrete time binomial models
for each integer $m\geq 1$.
We show next that on the same probability space we can construct
the $\lambda$-evolutionary system $(P_{B_{\lambda}},B_{\lambda})$ 
of Example \ref{P_B} for $\lambda=\{\lambda_m\colon\,m\geq 1\}$ with
$\lambda_m=\{k2^{-2m}\colon\,k=0,\dots,T2^{2m}\}$.

\begin{thm}\label{approx}
Let $B$ and $W_m$, $m\geq 1$, be stochastic processes
on the probability space $(\Omega,{\cal F},\Pr)$ satisfying conditions
$(a) - (d)$ of Theorem {\rm \ref{RW}}, and let $0<T<\infty$.
Then with probability $1$, for each $0\leq t\leq T$,
$$
\lim_{m\to\infty}\prod_{k=1}^{[t2^{2m}]}\big [1+W_m(k2^{-2m})
-W_m((k-1)2^{-2m})\big ]=\exp\big\{B(t)-t/2\big\}.
$$
\end{thm}

\begin{proof}
With probability $1$, for all but finitely many $m$, let
$\{\tau_m(k)\colon\,0\leq k\leq T2^{2m}\}$ be random partitions of
$[0,T')$ with $T'>T$ satisfying condition $(d)$ of Theorem \ref{RW}.
Then for any $0<t\leq T$, and for all but finitely many $m$,
$$
\prod_{k=1}^{[t2^{2m}]}\big [1+W_m(k2^{-2m})-W_m((k-1)2^{-2m})\big ]
=\prod_{k=1}^{[t2^{2m}]}\big [1+B(\tau_m(k))-B(\tau_m(k-1))\big ].
$$
Since a Brownian motion is sample continuous, it is enough to prove
that for any $0<t\leq T$, with probability $1$,
$$
\lim_{m\to\infty}\prod_{k=1}^{T2^{2m}}\big [1+B(\tau_m(k)\wedge t)
-B(\tau_m(k-1)\wedge t)\big ]=\exp\{B(t)-t/2\}.
$$
Let $0<t\leq T$.
By Proposition 4.10 in \cite[Part II]{DNa},
this is so once we show that with probability $1$,
$$
\lim_{m\to\infty}\max_k|B(\tau_m(k))-B(\tau_m(k-1))|=0
$$
and
$$
\lim_{m\to\infty}\sum_{k=1}^{T2^{2m}}\big [B(\tau_m(k)\wedge t)
-B(\tau_m(k-1)\wedge t)\big ]^2=t.
$$
Indeed, for each $m$, let $n=n_{m,t}$ be the integer such that 
$(n-1)2^{-2m}< t\leq n2^{-2m}$.
By statements $(d)$ and $(a)$ of Theorem \ref{RW}, we have for all
but finitely many $m$,
$$
|B(\tau_m(k))-B(\tau_m(k-1))|=|W_m(k2^{-2m})-W_m((k-1)2^{-2m})|
=2^{-m}
$$
for $1\leq k\leq T2^{2m}$, and
$$
\sum_{k=1}^{T2^{2m}}\big [B(\tau_m(k)\wedge t)
-B(\tau_m(k-1)\wedge t)\big ]^2=(n-1)2^{-2m}+[B(t)-B(\tau_m(n-1))]^2.
$$ 
Since $\tau_m(n-1)-t=[\tau_m(n-1)-(n-1)2^{-2m}]+[(n-1)2^{-2m}-t]\to 0$
as $m\to\infty$, the right side of the last display tends to $t$, 
proving the above two limiting conditions.
The proof of Theorem \ref{approx} is complete. 
\qed\end{proof}

\section{Option pricing and hedging}\label{option}

\paragraph*{Trading strategies.}
For $\lambda\in\Lambda [0,T]$ and  $k=0,\dots,n$, let $(P_k,X_k)$
be a $\lambda$-evolutionary system (see Definition \ref{E-system}).
Then $P=(P_0,\dots,P_n)$ is called a portfolio of $n+1$  
stocks of the $\lambda$-evolutionary model during a time period $[0,T$.

\begin{defn}\label{trade-str}
{\rm Let $P=(P_0,\dots,P_n)$ be a portfolio of $n+1$ stocks of the 
$\lambda$-evolutionary model for some $\lambda\in\Lambda [0,T]$,
and let$\gamma_k$, $k=0,\dots,n$, be real-valued functions on $[0,T]$
A vector function $\gamma=(\gamma_0,\dots,\gamma_n)$ is called
a \emph{$\lambda$-trading strategy for $P$} if for each $k=0,\dots,n$,
and a price scenario $P_k$, the Left Cauchy $\lambda$-integral 
$(LC)\smallint \gamma_k\,d_{\lambda}P_k$ is defined on $[0,T)$,
$\Delta^{-}P_k(T)=0$ whenever $\gamma_k(T-)$ is undefined,
and there exists the limit 
\beq\label{1trade-str}
(LC)\int_0^{T-}\gamma_k\,d_{\lambda}P_k:=
\lim_{u\uparrow T}(LC)\int_0^u\gamma_k\,d_{\lambda}P_k.
\eeq
Also, each $\gamma_k$ is then called a \emph{$\lambda$-trading strategy
for $P_k$}, and for  we let
\beq\label{2trade-str}
(LC)\int_0^T\gamma_k\,d_{\lambda}P_k:=\left\{ \begin{array}{ll}
\gamma_k(T-)\Delta^{-}P_k(T)
+(LC)\int_0^{T-}\gamma_k\,d_{\lambda}P_k
&\mbox{if $\gamma_k(T-)$ exists,}\\
(LC)\int_0^{T-}\gamma_k\,d_{\lambda}P_k&\mbox{otherwise.}
\end{array} \right. 
\eeq}
\end{defn}

It will be seen below (Proposition \ref{NK1}) that there is a price scenario
$P_k$ and a $\lambda$-trading strategy $\gamma_k$ for $P_k$ such that
$\gamma_k(T-)$ does not exist but the limit (\ref{1trade-str}) exists.
By Definition \ref{function}.\ref{LCint}, if the Left Cauchy $\lambda$-integral
$(LC)\smallint \gamma_k\,d_{\lambda}P_k$ is defined on $[0,T]$ then
its value on $[0,T]$ agree with (\ref{2trade-str}).

Following our convention spelled out by Remark \ref{convention}, 
the indefinite integrals for (\ref{1trade-str}) and (\ref{2trade-str}) may be
viewed as families of functions corresponding to different strategies
of a stock.
Also a $\lambda$-trading strategy may denote a set of different strategies as well
as a single strategy corresponding to a price scenario.

\begin{defn}
{\rm For $\lambda\in\Lambda [0,T]$, let $P=(P_0,\dots,P_n)$ be a portfolio
of $n+1$  stocks of the $\lambda$-evolutionary model during a time period 
$[0,T]$,
and let $\gamma=(\gamma_0,\dots,\gamma_n)$ be a $\lambda$-trading
strategy for $P$.
\begin{enumerate}
\item[$(a)$] For each price scenario $P$ define
the functions $V=V(P)$ and $G=G(P)$ defined by
$$
V(P;t):=\sum_{k=0}^n\gamma_k(t)P_k(t)\qquad\mbox{and}\qquad
G(P;t):=\sum_{k=0}^n(LC)\int_0^t\gamma_k\,d_{\lambda}P_k
$$
are called the \emph{portfolio value function} and the \emph{portfolio
gain function}, respectively.
\item[$(b)$] The $\lambda$-trading strategy $\gamma$ is called
\emph{self-financing} if $V(P;t)=V(P;0)+G(P;t)$ for each 
price scenario $P$ and $t\in [0,T]$.
\end{enumerate}}
\end{defn}

\paragraph*{Black-Scholes formula.}
Let $\lambda\in\Lambda [0,T]$, $1\leq p<2$ and $\sigma >0$.
Let $CQ_{\lambda,\sigma}[0,T]$ be the class of all continuous functions 
$g$ on $[0,T]$ having the quadratic $\lambda$-variation
with the bracket function $[g]_{\lambda}(t)=\sigma^2 t$, $0\leq t\leq T$,
and let $C\CLW_p[0,T]$ be the class all continuous functions 
on $[0,T]$ having bounded $p$-variation.
Consider the class of functions
\beq\label{BandS1}
X_{\lambda,p,\sigma}[0,T]:=
\Big\{w:=g+h\colon\,g\in CQ_{\lambda,\sigma}[0,T],\quad
h\in C\CLW_p[0,T],\quad g(0)=h(0)=0\Big\}.
\eeq
By Corollary \ref{function}.\ref{qvclassQ}, each function $w\in X$ is
continuous and has the quadratic $\lambda$-variation with the bracket function
$[w]_{\lambda}(t)=\sigma^2 t$, $0\leq t\leq T$, that is, 
$w\in CQ_{\lambda,\sigma}[0,T]$.
Let $X$ be a subset of $X_{\lambda,p,\sigma}[0,T]$.
Thus $X\subset E_{\lambda}[0,T]$, and so by Proposition \ref{E-stock},
$({\cal E}_{\lambda}(X),X)$ is a $\lambda$-evolutionary system.
By Theorem \ref{function}.\ref{LIE}, for each $w\in X$, 
$P:={\cal E}_{\lambda}(w)$ and $0\leq t\leq T$,
$$
P(t)=1+(L)\int_0^tP\,dw=1+(RS)\int_0^tP\,dh+(LC)\int_0^tP\,d_{\lambda}g.
$$
For example, $X$ could be the set $B_{\lambda}$ of all sample functions of a 
standard Brownian  motion which are continuous and have the quadratic 
$\lambda$-variation (cf.\ Example \ref{P_B}), and so 
$B_{\lambda}\subset X_{\lambda,1,1}$.
For another example, let $B_H$ be a fractional Brownian motion with
the Hurst index $1/2<H<1$ defined on the same complete probability
space $(\Omega,{\cal F},\Pr)$ as a standard Brownian motion $B$.
Then there is $\Omega_0\subset\Omega$ such that $\Pr (\Omega_0)=1$
and $X:=\{B(\cdot,\omega)+B_H(\cdot,\omega)\colon\,\omega\in
\Omega_0\}\subset X_{\lambda,p,1}$ for any $1/H<p<2$.

Consider a portfolio $(Q,P)$ with $P={\cal E}_{\lambda}(X)
=\{{\cal E}_{\lambda}(\omega)\colon\,\omega\in X\}$ for some 
$X\subset X_{\lambda,p,\sigma}[0,T]$,  and a set $Q$ consisting of a single 
price scenario $Q(t):=\exp\{rt\}$, $0\leq t\leq T$, depending on $r>0$.
The set $Q$ may be considered as a bond, while $P$ is a risky asset.
Given $K>0$ and $r>0$, we show that there is a self-financing 
$\lambda$-trading strategy $(\alpha,\beta)$ for the portfolio 
$(Q,P)$ with the portfolio value on a date $T$:
$$
V(T)=\alpha (T)Q(T)+\beta (T) P(T)=\max\big\{0,P(T)-K\}
$$
for any price scenario $(Q,P)$.
To define the desired $\lambda$-trading strategy $(\alpha,\beta)$,
for $t\in [0,T)$ and $x>0$, let
$$
d_1(t,x)=\big [\ln (x/K)+(r+\sigma^2/2)(T-t)]/(\sigma\sqrt{T-t})
\quad\mbox{and}\quad
d_2(t,x)=d_1(t,x)-\sigma\sqrt{T-t}.
$$
Denoting the standard normal distribution function by $N$, let
$\phi\equiv \phi_{\sigma,K}$ be the function on $[0,T]\times [0,+\infty)$
defined by
\beq\label{2spm}
\phi(t,x):=\left\{\begin{array}{ll}
xN (d_1(t,x))-Ke^{-r(T-t)}N (d_2(t,x)) 
&\mbox{if $(t,x)\in [0,T)\times (0,+\infty)$,}\\
\max\{0,x-K\} &\mbox{if $(t,x)\in\{T\}\times (0,+\infty)$,}\\
0 &\mbox{if $(t,x)\in [0,T]\times\{0\}$}.
\end{array}\right.
\eeq
Notice that $\phi (T,x)=\lim_{t\uparrow T}\phi (t,x)$ for each $x>0$,
and $\phi (t,0)=\lim_{x\downarrow 0}\phi (t,x)$ for each $t\in [0,T]$.
For a price scenario $P$, let $\Phi P(t):=\phi (t,P(t))$,  $0\leq t\leq T$.
We use notation (\ref{Fdt}) and (\ref{Fdx}) for partial derivatives of $\phi$. 
Then for a price scenario $(Q,P)$ and $0\leq s\leq T$, let
$$
\beta (s):=\left\{\begin{array}{ll}
\Phi_x'P(s)=N (d_1(s,P(s))) &       \mbox{if  $0\leq s <T$,}\\
0 &\mbox{if $s=T$,}
\end{array}\right.
$$
and 
$$
\alpha (s):=\big\{\Phi P(s)-P(s)\beta (s)\big\}/Q(s).
$$
Therefore the portfolio value function $V$ for a price scenario
$(Q,P)$ and $s\in [0,T]$ is given by
\beq\label{BandS2}
V(s)=\alpha (s)Q(s)+\beta (s)P(s)=\Phi P(s).
\eeq
The values of $\alpha$ and $\beta$ at $s=T$ do not affect the
self-financing condition to be proved next, and they can be chosen 
arbitrarily so that (\ref{BandS2}) holds for $s=T$.
An arbitrage argument assure that the portfolio value function $V$
gives a fair price of a European call option on the stock 
$P={\cal E}_{\lambda}(X)$.
This price is called the \emph{Black-Scholes formula}.

\begin{thm}\label{BandS}
For $\lambda\in\Lambda [0,T]$, $1\leq p<2$ and $\sigma >0$,
let  $X$ be a subset of $X_{\lambda,p,\sigma}[0,T]$, 
and let $P={\cal E}_{\lambda}(X)$.
Then $(\alpha,\beta)$  is a self-financing $\lambda$-trading strategy
for the portfolio $(Q,P)$.
\end{thm}

\begin{proof}
Let $(Q,P)$ be a price scenario, that is,
$Q(t)=\exp\{rt\}$, $0\leq t\leq T$, and $P={\cal E}_{\lambda}(w)$ 
for some $w\in X$.
Therefore $w=g+h$ for some $g\in CQ_{\lambda,\sigma}[0,T]$ and
$h\in CW_p[0,T]$.
Thus by Corollary \ref{function}.\ref{qvclassQ}, $w$ has the
quadratic $\lambda$-variation with the bracket function
$[w]_{\lambda}(t)=[g]_{\lambda}^c(t)=\sigma^2 t$, $0\leq t\leq T$.
By Proposition \ref{E-stock}, $P={\cal E}_{\lambda}(w)$
has the quadratic $\lambda$-variation with the bracket function
\beq\label{1BandS}
[P]_{\lambda}(t)=(LY)\int_0^t{\cal E}_{\lambda}(w)^2\,d[w]_{\lambda}
=\sigma^2 \int_0^tP^2(s)\,ds,\qquad 0\leq t\leq T.
\eeq
The function $\alpha$ is bounded on $[0,T]$ and continuous on $[0,T)$.
The function $Q$ is continuous and has bounded variation on $[0,T]$.
Thus by Theorem II.13.16 in Hildebrandt \cite{THH63}, $\alpha$
is Riemann-Stieltjes integrable on $[0,T]$ with respect to $Q$, and
for $0\leq t\leq T$, 
\beq\label{2BandS}
(RS)\int_0^t\alpha\,dQ=r\int_0^t\alpha (s)Q(s)\,ds=-r\int_0^t\big [
P(s)\Phi_x'P(s)-\Phi P(s)\big ]\,ds.
\eeq
To prove that $\beta$ is a $\lambda$-trading strategy for 
${\cal E}_{\lambda}(w)$ and $(\alpha,\beta)$ is self-financing 
we use a chain rule for the composition $\Phi P$.
By definition (\ref{2spm}), $\phi$ is continuous on $[0,T]\times (0,\infty)$,
it is a $C^2$ class function on $[0,T)\times (0,\infty)$, and  
\beq\label{4BandS}
\Phi_t'P(s)=\left\{ \begin{array}{ll}
-\frac{\sigma P(s)}{2\sqrt{T-s}}N'(d_1(s,P(s))-rKe^{-r(T-s)}N(d_2(s,P(s))
&\mbox{if $s\in [0,T)$,}\\ 
0 &\mbox{if $s=T$,} 
\end{array} \right. 
\eeq
\beq\label{5BandS}
\Phi_{xx}''P(s)=\left\{ \begin{array}{ll}
\frac{1}{\sigma P(s)\sqrt{T-s}}N'(d_1(s,P(s)) &\mbox{if $s\in [0,T)$,}\\ 
0 &\mbox{if $s=T$.} 
\end{array} \right. 
\eeq
Notice that $\phi $ satisfies the partial differential equation
\beq\label{3spm}
\left\{\begin{array}{ll}
\frac{\partial u}{\partial t}+Lu=0&\mbox{on $[0,T)\times (0,+\infty)$},\\
u(T,x)=\max\{0,x-K\} &\mbox{for $x\in (0,+\infty)$},\\
u(t,0)=0 &\mbox{for $t\in [0,T]$,}
\end{array}\right.
\eeq
where
$$
Lu(t,x):=\frac{\sigma^2}{2}x^2
\frac{\partial^2u}{\partial x^2}(t,x)+rx\frac{\partial u}{\partial x}
(t,x)-ru(t,x).
$$
By the chain rule of Theorem \ref{function}.\ref{mdchrule}(1)
(cf. (\ref{2chrule-HK})), the Left Cauchy $\lambda$-integral
$(LC)\smallint \Phi_x'P\,d_{\lambda}P$ is defined on $[0,T)$,
and by (\ref{BandS2}), for each $0\leq u < T$, 
\begin{eqnarray}
V(P;u)-V(P;0)&=&(R)\int_0^u\Phi_t'P+(LC)\int_0^u\Phi_x'P\,d_{\lambda}P
+\frac{1}{2}(RS)\int_0^u\Phi_{xx}''P\,d[P]_{\lambda}\nonumber\\
\mbox{by (\ref{1BandS}) and (\ref{2BandS})}\quad
& =&(RS)\int_0^u\alpha\,dQ+(LC)\int_0^u\beta \,d_{\lambda}P
+(R)\int_0^u\Big [\frac{\partial \phi}{\partial t}+L\phi\Big ](s,P(s))\,ds
\nonumber\\ \mbox{by (\ref{3spm})}\quad
&=&(RS)\int_0^u\alpha\,dQ+(LC)\int_0^u\beta\,d_{\lambda}P=G(P;u).
\label{3BandS}\end{eqnarray}
The last equality holds because the Left Cauchy $\lambda$-integral extends 
the Riemann-Stieltjes integral.

To prove that the left and right sides of (\ref{3BandS}) are equal for
$u=T$ we use the chain rule of Proposition \ref{function}.\ref{chrule-HK}.
Letting $\chi (s):=1/\sqrt{T-s}$ for $s\in [0,T)$ and $\chi (T):=0$,
there is a finite constant $C$ such that for all $0\leq s\leq T$,
$$
\big |\Phi_t'P(s)\big |\leq C\chi (s)\qquad\mbox{and}\qquad
\big |\Phi_{xx}''P(s)\big |\leq C\chi (s).
$$
By the Hake theorem (Section 7.3 in \cite{RMM})
and by the substitution theorem for the Henstock-Kurzweil integral
(Section 7.2 in \cite{RMM}),
$\chi $ is Henstock-Kurzweil integrable on $[0,T]$ with respect to
the identity function and with respect to the bracket function $[P]_{\lambda}$.
By the dominated convergence theorem for the Henstock-Kurzweil
integral (Section 7.8 in \cite{RMM}), it then follows that
the Henstock-Kurzweil integrals $(HK)\smallint_0^T\Phi_t'P$ 
and $(HK)\smallint_0^T\Phi_{xx}''P\,d[P]_{\lambda}$ are defined.
By Proposition \ref{function}.\ref{chrule-HK}, we then have that
$$
(LC)\int_0^{T-}\Phi_x'P\,d_{\lambda}P
=V(P;T)-V(P;0)-(HK)\int_0^T\Phi_t'P-\frac{1}{2}(HK)\int_0^T
\Phi_{xx}''P\,d[P]_{\lambda}.
$$
Since $\Delta^{-}P(T)=0$, by (\ref{2trade-str}), $(\alpha,\beta)$
is a $\lambda$-trading strategy for $(Q,P)$ and 
$V(P;T)=V(P;0)+G(P;T)$.
The proof is complete.
\qed\end{proof}

The use of the Henstock-Kurzweil integrals in the preceding proof is
essential because the functions $\Phi_t'P$ and $\Phi_{xx}''P$ defined
by (\ref{4BandS}) and (\ref{5BandS}), respectively, can be unbounded
near $T$.
Before showing this, following Bick and Willinger \cite{BandW94},
we notice that
$$
\lim_{t\uparrow T}d_1(t,P(t))=\lim_{t\uparrow T}\big\{\ln \big (P(t)/K
\big )/\sigma\sqrt{T-t}\big\}=\left\{ \begin{array}{ll}
-\infty &\mbox{if $P(T)<K$,} \\
+\infty &\mbox{if $P(T)>K$.}
\end{array} \right. 
$$
Thus in the case $P(T)\not =K$, $\Phi_t'P(T-)=\Phi_{xx}''P(T-)=0$.
However the case $P(T)=K$ is different.
Bick and Willinger (see Remark following Proposition 1 in 
\cite[p.\ 366]{BandW94}) conjectured that in this case $(\alpha,\beta)$
may not be trading strategies in their sense. 
In contrary the preceding theorem shows that $(\alpha,\beta)$
is always a self-financing $\lambda$-trading strategy.
The following shows that $\Phi_t'P$ and $\Phi_{xx}''P$ can be unbounded
and gives a positive answer to the second part of
the problem of Bick and Willinger \cite{BandW94}.

\begin{prop}\label{NK1}
Let $\sigma=1$ and $K=\sqrt{e}$.
There exist $\lambda\in\Lambda [0,1]$ and
a price scenario $P={\cal E}_{\lambda}(w)$ with 
$w\in X_{\lambda,1,\sigma}[0,1]$, $w(1)=1$, such that the limit 
$\lim_{t\uparrow 1}\beta(t)$ does not exist, and the functions
$\Phi_t'P$, $\Phi_{xx}''P$ are unbounded.
\end{prop}

\begin{proof}
To construct the desirable return $w$ we choose a self-affine function of Example
\ref{function}.\ref{NK}.
By Theorem 1 of K\^ono  \cite{NK88}, a continuous self-affine function 
$w$ on $[0,1]$ with the scale parameter $H\in (0,1)$ 
to base $r\geq 4$ satisfying $w(0)=0$ and $w(1)=1$ is completely determined
by the sequence $x(k)\in\{1,-1\}$, $k=0,\dots,r-1$, 
and such that $\sum_{k=0}^{r-1}x(k)=r^H$.
Then for $0\leq t\leq 1$,
$$
w(t)=\sum_{n=1}^{\infty}y_{n-1}(t)s(\delta_n)r^{-nH},
$$
where $t=\sum_{n=1}^{\infty}\delta_nr^{-n}$ with $\delta_n\in\{0,1,
\dots,r-1\}$, $s(0)=0$, $s(j)=\sum_{i=0}^{j-1}x(i)$, $j=1,\dots,r$,
$y_0(t)=1$ and 
\beq\label{1NK1}
y_n(t)=\prod_{k=1}^nx(\delta_k).
\eeq
For $m=1,2,\dots$, 
\beq\label{2NK1}
0\leq h =\sum_{n=m+1}^{\infty}\delta_nr^{-n}\leq r^{-mH}
\quad\mbox{and}\quad t=\sum_{n=1}^m\delta_nr^{-n}=ir^{-m},
\eeq
letting $\epsilon_{i,m}:=y_m(ir^{-m})$, we have
\begin{eqnarray*}
w(ir^{-m}+h)-w(ir^{-m})&=&\sum_{n=m+1}^{\infty}y_{n-1}(ir^{-m}+h)
s(\delta_n)r^{-nH}\\[2mm]
&=&y_m(ir^{-m})r^{-mH}\sum_{n=1}^{\infty}y_{n-1}(r^mh)
s(\delta_{m+n})r^{-nH}\\[2mm]
&=&\epsilon_{i,m}r^{-mH}w(r^mh).
\end{eqnarray*}
Now let $w$ be such that $H=1/2$, $r=4$ and $\{x(k)\colon\,k=0,\dots,3\}
=\{1,1,1,-1\}$.
For $h=r^{-m}$ with $m=1,\dots$ and $t=1-h$ in (\ref{2NK1}),
all $\delta_n=r-1=3$ and $i=r^m-1=3^m-1$.
Thus by (\ref{1NK1}), we have for $m=1,2,\dots$,
\beq\label{3KN1}
w(1)-w(1-4^{-m})=(-1)^m4^{-m/2}w(1).
\eeq
Letting $\lambda_m:=\{ir^{-m}\colon\,i=0,\dots,r^m\}$,
we have $\lambda=\{\lambda_m\colon\,m\geq 1\}\in\Lambda [0,1]$,
and so $w\in X_{\lambda,1,\sigma}$.
For the price scenario $P(t)=\exp\{w(t)-t/2\}$, $0\leq t\leq 1$,
we then have
$$
\lim_{t\uparrow 1}d_1(t,P(t))
=\lim_{t\uparrow 1}\frac{\ln (P(t)/K)}{\sqrt{1-t}}
=\lim_{h\downarrow 0}\frac{w(1-h)-w(1)}{\sqrt{h}}.
$$ 
However by (\ref{3KN1}), the right side does not exist, proving the
proposition.
\qed\end{proof}

\section{Returns and the $p$-variation: examples}
 
Let $0<p<2$ and let  $X$ be a stochastic process on a complete probability space 
$(\Omega,{\cal F},\Pr)$ such that almost all its sample functions
have bounded $p$-variation  and are right-continuous.
By Proposition \ref{function}.\ref{qcforpq} (see also Example 
\ref{function}.\ref{qv-for-W2}),  given $\lambda\in\Lambda [0,T]$,
there is $\Omega_0\in{\cal F}$ of full probability such that
$X(\cdot,\omega)$ has the quadratic $\lambda$-variation for each
$\omega\in\Omega_0$.
Let $\tau$ be the first (random) moment in $(0,T]$ such that 
$X(\tau)-X(\tau-)\leq -1$, and let $X^{\tau-}:=X(\cdot\wedge \tau-)$.
Letting $X:=\{X^{\tau-}(\cdot,\omega)\colon\,\omega\in\Omega_0\}$ and 
$P_X:=\{{\cal E}_{\lambda}(X^{\tau-}(\cdot,\omega))
\colon\,\omega\in\Omega_0\}$, by Proposition \ref{E-system},
$(P_X,X)$ is a $\lambda$-evolutionary system, that is,
in the $\lambda$-evolutionary asset pricing model $P_X$ is a set of
price scenarios of a stock and $X$ is a set of returns.
So the quadratic $\lambda$-variation and the $p$-variation property 
are the main indicators of a stochastic process for its applicability in the
$\lambda$-evolutionary asset pricing model.   

The $p$-variation property has been investigated for several classes
of stochastic processes.
However, the results are scattered over many different journals.
Below we provide some of these results related to financial models.
More information can be found in the annotated reference list
on $p$-variation in \cite[Part IV]{DNa}.

 \paragraph*{Fractional Brownian motion.}
A fractional Brownian motion $B_H=\{B_H(t)\colon\,t\geq 0\}$ 
with index $H\in (0,1)$ is a mean zero Gaussian process with 
the covariance function
\beq\label{fbm1}
E\big\{B_H(t)B_H(s)\big\}=\frac 12 \big\{
t^{2H}+s^{2H}-|t-s|^{2H}\big\}\quad\mbox{for $t, s\geq 0$}
\eeq
and $B_H(0)=0$ almost surely.
Since the right side of (\ref{fbm1}) is equal to
$t\wedge s$ for $H=1/2$, $B_H$ is a Brownian motion
in this case.
Many important properties of a fractional Brownian motion can be found in
Kahane \cite[Sect. 18]{J-PK85}.
In particular, from the a.s. bound $w_{B_H}(t)=O(t^H\sqrt{\log (1/t)})$ 
as $t\downarrow 0$ of the modulus of continuity on $[0,T]$ 
it follows that $v_p(B_H;[0,T])<\infty$ a.s. whenever $p>1/H$. 
The next statement provides a more precise description of the 
$p$-variation of a fractional Brownian motion.
Let $s_{\psi}(f;\kappa):=\sum_{i=1}^n\psi (|f(x_i)-f(x_{i-1})|)$
for $\kappa=\{x_i\colon\,i=0,\dots,n\}$ and let 
$\psi_H(u):=u^{1/H}/(2LLu)^{1/(2H)}$ for $u>0$, where $LLu:=
\log_e|\log_e u|$, and $\psi_H(0):=0$.

\begin{prop}
For $0<T<\infty$, almost all sample functions of a fractional 
Brownian motion  $B_H$ with $H\in (0,1)$ satisfy the relation
\beq\label{fbm2}
\lim_{\delta\downarrow 0}\,\sup\big\{s_{\psi_H}(B_H;\kappa )
\colon\,\mbox{$\kappa\in\Xi [0,T]$ and the mesh $|\kappa|\leq\delta$}
\big\}=T.
\eeq
\end{prop}

\begin{rem}
{\rm For the case $H=1/2$, (\ref{fbm2}) was proved by Taylor \cite[Theorem 1]{SJT}.
The general case of (\ref{fbm2}) with $\leq$ instead of $=$ follows from
Theorem 3 of Kawada and K\^ono \cite{KK} because almost all sample functions
of $B_H$ are continuous.
The converse inequality does not follow from Theorem 4 of Kawada and K\^ono 
\cite{KK} because their condition (v) fails to hold for $\sigma^2(u)=u^{2H}$, 
$u>0$.} 
\end{rem}

\begin{proof}
By the preceding remark, it is enough to prove (\ref{fbm2}) with $\geq$
instead of $=$.
Let $\phi_H(u):=u^H\sqrt{2LLu}$ for $u>0$.
Oodaira \cite{HO72} proved that $\eta_u:=\{B_H(tu)/\phi_H(u)
\colon\,t\in [0,1]\}$ satisfies the functional law of the iterated logarithm
as $u\downarrow 0$.
Taking a continuous functional $\Phi\colon C[0,1]\to\RR$ 
such that $\Phi (f):=f(1)$ for each $f\in C[0,1]$, by the
continuous mapping theorem, it then follows that, with probability 1,
\beq\label{fbm3}
\limsup_{u\downarrow 0}\frac {|B_H(u)|}{\phi_H(u)}=\sup\big\{
f(1)\colon\,f\in K_H\big\}\geq 1,
\eeq
where $K_H$ is the unit ball in the Reproducing Kernel Hilbert Space
corresponding to the covariance of $B_H$.
Denoting by $\phi_H^{-1}$, the inverse function of $\phi_H$, we
have that $\psi_H(u)/\phi_H^{-1}(u)\to 1$ as $u\downarrow 0$.
Thus, by (\ref{fbm3}), we have with probability 1
$$
\limsup_{u\downarrow 0}\frac {\psi_H(|B_H(u)|)}{u}\geq 1.
$$
This in conjunction with Fubini's theorem yields that
the Lebesgue measure
$$
|\big\{t\in [0,T]\colon\,\limsup_{u\downarrow 0}
\frac {\psi_H(|B_H(t+u)-B_H(t)|)}{u}\geq 1\big\}|=T.
$$
Now the proof can be completed using the Vitali covering theorem
in the same way as for the case $H=1/2$ in Taylor \cite[Theorem 1]{SJT}.
\qed\end{proof}

Statistical analysis of a financial data under the fractional Brownian 
motion hypthesis is discussed in Beran \cite{JB94}.

\paragraph*{Ornstein-Uhlenbeck process.}
A mean zero Gaussian stochastic process $u=\{u(t)\colon\,-\infty<t<+\infty\}$
with the covariance
$$
E\{u(t)u(s)\}=\exp\{-\beta |t-s|\}\qquad\mbox{for $t,s\in
(-\infty,+\infty)$ and $\beta >0$}
$$
is called the OU (Ornstein-Uhlenbeck) process.
It is the unique solution of the Langevin equation driven by a
Brownian motion.
Therefore $u$ is also called the OU velocity process
to distinguish it from the corresponding OU position
process $x=\{x(t)\colon\,t\geq 0\}$ defined by $x(t)-x(0)=
\int_0^tu(s)\,ds$.
The latter stochastic process models the $x$-coordinate of the free particle
at the time $t$.
The OU position process $x$ is a mean zero Gaussian process
with variance
$$
E[x(t)-x(s)]^2=2\beta^{-2}\big [e^{-\beta |t-s|}-1
+\beta |t-s|\big ]\quad\mbox{for $t, s\geq 0$}.
$$
When restricted to a bounded interval, almost all
sample functions of $x$ have bounded variation.
This follows from the result of Ibragimov \cite[Theorem 6]{IAI73}:

\begin{thm}
If $X=\{X(t)\colon\,0\leq t\leq T\}$ is a separable measurable stochastic 
process such that $E|X(t)-X(s)|\leq K|t-s|$ for $0\leq s,t\leq T$,
then with probability $1$,
sample functions of $X$ are of bounded variation.
\end{thm}

Goldenberg \cite{DHG86} argued that future prices in markets with frictions 
should be modelled by stochastic processes with smooths sample functions
such as of an OU position process.
Ross \cite[Sect. 6.3.4]{SMR83} indicated a different financial application 
of an integrated Brownian motion having sample functions
of bounded variation.
Notice that the $p$-variation of an OU 
velocity process is the same as of a Brownian motion. 
The OU velocity process was scrutinized by Borell
\cite{ChB90} as the return process on the basis of Danish stocks.

\paragraph*{Homogeneous L\'evy process.}
Let $X=\{X(t)\colon\,t\geq 0\}$ be a separable, continuous 
in probability stochastic process with independent increments.
Then $X$ is called the \emph{L\'evy process} if almost all its
sample functions are regulated, right continuous and $X(0)=0$
almost surely.
A L\'evy process $X$ is called \emph{homogeneous} if a distribution of 
$X(t+s)-X(t)$ with $t, s\geq 0$ does not depend on $t$.
Given a real number $a$, a positive number $b$,
and a $\sigma$-finite measure $L_X$ on $\RR\setminus\{0\}$ such 
that  $L_X(\{x\in\RR\colon\,|x|>\delta\})<\infty$ for $\delta >0$ and
$\int_{\RR\setminus \{0\}} (1\wedge x^2)\,L_X(dx)<\infty$, let
\beq\label{Levy1}
\Phi (u):=iau-bu^2+\int_{\RR\setminus\{0\}}\big (e^{iux}-1+
iuh(x)\big )L_X(dx)
\eeq
for $u\in\RR$, where $h$ is a bounded Borel function.
The function $\Phi$ is called the \emph{characteristic exponent} and 
the measure $L_X$ is called the \emph{L\'evy measure}.
Then the characteristic function of a homogeneous L\'evy process $X$ is
given by
$
E\exp\{ iuX(t)\}=\exp\big\{t\Phi (u)\big\}
$
for each $t\geq 0$ and $u\in\RR$.

It is well-known that sample functions of a homogeneous L\'evy
process $X$ with the characteristic exponent (\ref{Levy1}) have bounded 
variation if and only if $b\equiv 0$ and
\beq\label{Levy2}
\int_{\RR\setminus \{0\}} (1\wedge |x|)\,L_X(dx)<\infty
\eeq
(cf. e.g. Theorem 3 on p. 279 of Gikhman and Skorokhod \cite{GS69}).
The following result of Bretagnolle \cite[Th\'eor\`eme III b]{JLB72}
is less well-known:

\begin{thm}\label{Bretagnolle}
Let $1<p<2$ and let $X=\{X(t)\colon\,t\geq 0\}$ be a mean zero homogeneous 
L\'evy process with the characteristic exponent {\rm (\ref{Levy1})} such that
$b\equiv 0$.
Then $v_p(X;[0,1])<\infty$ with probability $1$ if and only if
$$
\int_{\RR\setminus \{0\}} (1\wedge |x|^p)\,L_X(dx)<\infty.
$$  
\end{thm}

Less precise but more general result is a characterization of
the $p$-variation index (see (\ref{p-var-index}) for the definition)
of sample functions of a L\'evy process by the Blumenthal-Gettor index.
Let $X$ be a homogeneous L\'evy process.
The \emph{Blumenthal-Getoor index} of $X$ is defined by
$$
\beta (X):=\inf\big\{\alpha >0\colon\,\int_{\RR\setminus\{0\}}
\big (1\wedge |x|^{\alpha}\big )\,L_X(dx)<\infty\big\},
$$
where $L_X$ is the L\'evy measure of $X$.
Notice that $0<\beta (X)\leq 2$.
If $X$ has no Gaussian part, then for any $0<T<\infty$, the $p$-variation index
$\upsilon (X;[0,T])=\beta (X)$ almost surely.
This follows from Theorems 4.1 and 4.2 of Blumenthal and Getoor \cite{BG},
and from Theorem 2 of Monroe \cite{MI}.

\paragraph*{Hyperbolic L\'evy motion.}
A homogeneous L\'evy process $X$ is called a hyperbolic L\'evy motion
if in its characteristic exponent (\ref{Levy1}) $b\equiv 0$ and
the L\'evy measure $L_X(dx)=g(x;\alpha,\delta)\,dx$, where
$$
g(x;\alpha,\delta)=\frac {1}{\pi^2|x|}\int_0^{\infty}\frac
{\exp\{-|x|\sqrt {2y+(\alpha )^2}\}}{y(J_1^2(\delta\sqrt{2y})
+Y_1^2(\delta\sqrt {2y}))}\,dy+\frac {\exp \{-|x|\}}{|x|}
$$
for $\alpha >0$ and $\delta\geq 0$.
Here $J_1$ and $Y_1$ are the Bessel functions of the first and
second kind, respectively.
Eberlein and Keller \cite{EandK} used a hyperbolic L\'evy motion to model 
stock price returns and showed (ibid, p. 295) that 
$g(x;\alpha,\delta )\sim x^{-2}$  as $x\to 0$.
Therefore, by Theorem \ref{Bretagnolle}, almost all sample 
functions of a hyperbolic L\'evy motion have bounded $p$-variation on $[0,1]$
for each $p>1$, that is, the $p$-variation index $\upsilon (X;[0,1])=1$.

\paragraph*{Normal inverse Gaussian L\'evy process.}
A homogeneous L\'evy process $X$ is called a normal inverse 
Gaussian L\'evy process if in its characteristic exponent (\ref{Levy1})
$b\equiv 0$,
$a=(2\alpha\delta/\pi)\smallint_0^1 \mbox{sinh} (\beta x)K_1(\alpha x)\,dx$
and the L\'evy measure $L_X(dx)=f(x;\alpha,\beta,\delta)\,dx$, where
$$
f(x;\alpha,\beta,\delta)=\frac {\alpha\delta}{\pi}\frac {e^{-\beta x}}
{|x|}K_1(\alpha |x|)=\frac {\delta}{\pi\sqrt {2}}\frac {e^{\beta x}}
{|x|}\int_0^{\infty}y^{-1/2}e^{-|x|\sqrt{2y+\alpha ^2}}dy
$$
for $\alpha >0$, $0\leq |\beta |<\alpha$, $\mu\in\RR$ and $\delta
>0$.
Here $K_1$ denotes the modified Bessel function of third order and
index 1.
Barndorff-Nielsen \cite{OEB-N97} showed that
$f(x;\alpha,\beta,\delta)\sim (\delta/\pi )x^{-2}$ as $x\to 0$.
Therefore, by Theorem \ref{Bretagnolle}, almost all sample functions
of a normal inverse Gaussian L\'evy process have bounded $p$-variation on 
$[0,1]$ for each  $p>1$, that is, the $p$-variation index $\upsilon (X;[0,1])=1$.

\paragraph*{The V. G. process.} This process (variance gamma) is defined
to be the composition $N(t):=B(G(t))$, $t\geq 0$, of a Brownian motion $B$
and the process $G$ which has independent gamma increments with mean
$t$ and variance $\nu t$.
Then $N$ is a homogeneous L\'evy process with $a=b\equiv 0$ and 
with the L\'evy measure
$L_X(dx)=(\nu |x|)^{-1}\exp\{-|x|\sqrt {2}/\sigma\sqrt\nu\}\,dx$.
Almost all sample functions of $N$ have bounded variation
because (\ref{Levy2}) holds.
Madan and Seneta \cite{MandS} introduced the V. G. process to replace
a Brownian motion in models for stock price returns.

\paragraph*{$\alpha$-stable L\'evy motion.}
A homogeneous L\'evy process $X$ is called an $\alpha$-stable L\'evy motion
of index $\alpha$ if its characteristic exponent (\ref{Levy1}) is given by 
$b\equiv 0$ and the L\'evy measure $L_X=L_{\alpha,r,q}$, where  
$L_{\alpha,r,q}(dx):=rx^{-1-\alpha}dx$ if $x>0$ and $L_{\alpha,r,q}(dx)
:=q(-x)^{-1-\alpha}dx$ if $x<0$                    
for $\alpha\in (0,2)$ and $r, q\geq 0$ with $r+q>0$.
If $\alpha <1$ then (\ref{Levy2}) holds and 
$\smallint h\,dL_X<\infty$.
In this case, it is said that an $\alpha$-stable L\'evy motion 
has no drift if $a+\smallint h\,dL_X=0$.
The following result is due to Fristedt and Taylor \cite[Theorem 2]{FT}.

\begin{thm}
Let $X_{\alpha}=\{X_{\alpha}(t)\colon\,t\in [0,1]\}$ be an 
$\alpha$-stable L\'evy motion with $\alpha\in (0,2)$ having no drift 
if $\alpha <1$ and with $r=q$ if $\alpha =1$.
For an increasing function $\psi\colon\,[0,\infty)\to [0,\infty )$,
with probability $1$
\beq\label{Levy3}
\lim_{|\kappa |\to 0}s_{\psi}(X_{\alpha};\kappa)\equiv
\lim_{|\kappa |\to 0}\sum_{i=1}^n\psi (|X_{\alpha}(t_i)-X_{\alpha}
(t_{i-1})|)=\sum_{(0,1]}\psi (|\Delta^{-}X_{\alpha}|).
\eeq
\end{thm}

To apply this result it is useful to recall that the right
side of (\ref{Levy3}) is finite almost surely if and only if
$\int_{-1}^1\psi (u)u^{-1-\alpha }\,du<\infty$.
Xu \cite{QHX} established necessary and/or sufficient conditions
for the boundedness of the $p$-variation of
a symmetric $\alpha$-stable processes with possibly dependent
increments.
More information about $p$-variation of stable prcesses can be
found in Fristedt \cite{BF74}.

Multivariate stable distributions are often used as an alternative 
to the Gaussian hypothesis in financial models.
Usually statistical analysis of financial data under the
stable distribution hypothesis is based on the log-return.
Therefore the conclusions of such analysis concern the exponential
model.
On the other hand, most theoretical work assuming the non-Gaussian
stable distribution is based on the simple net return.
This discord have been discussed by Elton et al. \cite{EGK}.  
Financial applications of stable distributions are reviewed by 
McCulloch \cite{JHMC}.
According to Campbell et al. \cite[p. 18]{CLMK}: ``Although stable 
distributions were popular in the 1960's and early 1970's, they are less
commonly used today. They have fallen out of favor partly because they
make theoretical modelling so difficult ...''.

\section{Estimating the $p$-variation index}

In the evolutionary asset pricing model a price scenario is represented by
a function having the quadratic $\lambda$-variation with some 
$\lambda\in\Lambda [0,T]$.
Therefore it is important to have a method to test this property
when a function on $[0,T]$ is given by a finite set of values.
This looks like a task comparable with the familiar problem of testing
the hypothesis of finitiness of a second moment of a random variable
having a finite sample.
Instead, it seems simpler to start with a statistical analysis of whether
a historical data with high probability correspond to a sample function
of a stochastic process having bounded $p$-variation for some $p<2$. 
This would be very usefull because a robustness of the evolutionary asset
pricing model changes drastically if a price scenario 
follow a function having bounded $p$-variation for some $0<p<2$.
Recall that the $p$-variation index of a function is the infimum of all $p$th such 
that the $p$-variation is finite (cf.\ (\ref{p-var-index})).
In this section we discuss two methods of estimating the $p$-variation index 
of a function from a given finite set of its values.

\paragraph*{Gladyshev class and the Orey index.}
A method of estimation of the $p$-variation index for a class of stochastic 
processes was suggested in \cite{NS99}.
It is based on the result of Gladyshev \cite{EGG} concerning 
a rescalled quadratic variation property of a Gaussian stochastic process $X$ 
with stationary increments.
Let $X=\{X(t)\colon\,t\geq 0\}$ be a Gaussian stochastic
process with stationary increments and continuous in quadratic mean.
Let $\sigma =\sigma_X$ be the incremental variance of $X$ given by
$$
\sigma_X(u)^2:=E\big [(X(t+u)-X(t))^2\big ], \qquad t, u\geq 0.
$$
Following Orey (1970), let
$$
\gamma_{\ast}:=\inf\big\{\gamma >0\colon\,u^{\gamma}/\sigma (u)\to 0,
\,\mbox{as}\,u\downarrow 0\big\}
=\limsup_{u\downarrow 0}\Big (\frac{\log \sigma (u)
}{\log u}\Big )
$$
and
$$
\gamma^{\ast}:=\sup\big\{\gamma >0\colon\,u^{\gamma}/\sigma (u)\to 
+\infty,\,\mbox{as}\,u\downarrow 0\big\}
=\liminf_{u\downarrow 0}\Big (\frac{\log \sigma (u)
}{\log u}\Big ).
$$
We always have that $0\leq\gamma^{\ast}\leq\gamma_{\ast}\leq \infty$.
If $\gamma_{\ast}=\gamma^{\ast}$ then we say that $X$ has the
{\em Orey index} $\gamma_X:=\gamma_{\ast}=\gamma^{\ast}$.
If $X$ is a fractional Brownian motion $B_H$ with the Hurst index $H$,
then its Orey index $\gamma_X=H$. 
If $X$ has the Orey index $\gamma_X$, then its $p$-variation index is 
$1/\gamma_X$, so that an estimation of the $p$-variation index 
reduces to an estimation of the Orey index.

For each integer $m\geq 1$, let $\lambda_m=\{i2^{-m}\colon\,
i=0,\dots,2^m\}$ be a nested sequence of dyadic partitions of $[0,1]$.
Gladyshev \cite{EGG} proved under mild conditions on $X$ that if
the Orey index $\gamma_X$ exists then almost surely
$$
\lim_{m\to\infty}\frac{\log\sqrt{s_2(X;\lambda_m)/2^m}}{\log (1/2^m)}
=\gamma_X.
$$
Motivated by Gladyshev's result,
the following definition was suggested in \cite{NS99}.

\begin{defn}
{\rm Let $\{N_m\colon\,m\geq 1\}$ be a sequence of integers increasing to
infinity and let $\lambda=\{\lambda_m\colon\,m\geq 1\}$ be a sequence of
partitions $\lambda_m=\lambda (T/N_m)=\{iT/N_m\colon\,i=0,\dots,N_m\}$ 
of $[0,T]$.
We say that a stochastic process $X=\{X(t)\colon\,t\geq 0\}$ belongs to 
the \emph{Gladyshev class}, or \emph{class ${\cal G}^{\lambda}$} 
whenever the limit
\beq\label{formula}
G_{\lambda}(X):=\lim_{m\to\infty}\frac{\log\sqrt{s_2(X;\lambda_m)/N_m}}
{\log (1/N_m)}
\eeq
exists for almost all sample functions of $X$.}
\end{defn}

We do not know how large Gladyshev's class is.
If $X$ is a Gaussian stochastic process, then some information can be derived 
from known results about the \emph{strong limit theorem}:  for some sequence
$\{N_m\colon\,m\geq 1\}$ of integers increasing to infinity and
for some constant $C_1$, there exists the limit
\beq\label{slth}
\lim_{m\to\infty}\frac{s_2(X;\lambda_m)}{N_m\sigma_X(1/N_m)^2}=C_1^2,
\qquad\mbox{almost surely}.
\eeq
Indeed for each $m\geq 1$ and any positive constant $C_1$, we have the identity
$$
\log\sqrt{\frac{s_2(X;\lambda_m)}{N_m}}=\log\sigma_X(1/N_m)
+\log C_1+\frac{1}{2}\log\left (1+
\frac{s_2(X;\lambda_m)}{C_1^2N_m\sigma_X(1/N_m)^2}-1\right ).
$$
Then assuming that the Orey index $\ind_X$ exists,
$X$ belongs to the Gladyshev class ${\cal G}^{\lambda}$ and 
$\ind_X=G_{\lambda}(X)$, almost surely.
This relation was used in \cite{NS99} to estimate the Orey index from a 
sample function of a
stochastic process $X$ given its values at finitely many points.

The rest contains a brief review of known results related to
the strong limit theorem (\ref{slth}) for Gaussian processes.
In addition to this theorem, there are several results in a literature
concerning a restricted form of $p$-variation of stochastic processes.
Corollary \ref{particular} below relates such results with the strong
limit theorem used in the present paper.
At the end of this section we formulate a part of Gladyshev's result
concerning Gaussian stochastic processes with stationary increments.

For a Brownian motion $B$, L\'evy \cite{PL40}, 
and independently Cameron and Martin \cite{CM},  proved that the limit
\beq\label{levy}
\lim_{m\to\infty}s_2(B;\lambda (T/2^m))=T
\eeq
exists almost surely.
Therefore the limit (\ref{formula}) with $X=B$ and $N_m=2^m$ also exists 
almost surely and $G_{\lambda}(B)=1/2$.
Each Gaussian process $X$ satisfying the strong limit theorem of Baxter 
\cite{GB} also belongs to Gladyshev's class and $G_{\lambda}(X)=1/2$.
Gladyshev \cite{EGG} extended Baxter's result to a large class of Gaussian
stochastic processes.
In the case a Gaussian process $X$ has stationary increments and has the Orey
index $\ind_X\in (0,1)$, the limit 
\beq\label{glad}
\lim_{m\to\infty}\frac{s_2(X;\lambda (T/2^m))}{Es_2(X;\lambda (T/2^m))}=T
\eeq
exists almost surely provided the spectral density of $X$ has a power 
behaviour (see Theorem \ref{stationary} below for a precise statement).
This yields (\ref{slth}) for $N_m=2^m$.
Yaglom \cite{AMY} discussed several applications of this type
strong limit theorems.
He also suggested that relation (\ref{glad}), and its further
modifications, should hold for a variety of stochastic processes.  

L\'evy \cite{PL40} proved relation (\ref{levy})  for a sequence
of dyadic partitions $\{kT2^{-n}\colon\,k=0,\dots,2^n\}$ of the time interval
$[0,T]$.
He also noticed that the same holds for any {\em nested} sequence
$\{\kappa_n\colon\,n\geq 1\}$ of partitions $\kappa_n$ such that 
$\cup_n\kappa_n$ is dense in $[0,T]$.
Dudley \cite{RMD73} proved that the sequence of dyadic partitions
in (\ref{levy}) can be replaced by an {\em arbitrarily} sequence of partitions
$\kappa_n=\{0=t_0^n<t_1^n<\cdots <t_{N_n}^n=T\}$ as long as the mesh 
$|\kappa_n|:=\max_i(t_i^n-t_{i-1}^n)$ tends to zero faster
than $1/\log n$ as $n\to\infty$.
The rate $o(1/\log n)$ cannot be improved as was shown 
by Fern\'andez de la Vega \cite{WFV}.
Gin\'e and Klein \cite{GK} extended the results of Baxter \cite{GB}
and Gladyshev \cite{EGG} to arbitrary sequences $\{\kappa_n\colon\,
n\geq 1\}$ of partitions subject to a certain rate of convergence to zero of 
$|\kappa_n|$.
For partitions $\kappa_n=\{kT/N_n\colon\,k=0,\dots,N_n\}$, $n=1,2,\dots$, 
with equally spaced points we have $|\kappa_n|=1/N_n$.
It follows from Gin\'e and Klein \cite{GK} that (\ref{glad}) holds
provided $N_n^{-\ind}=o(1/\log n)$ if $\ind_X=\ind\in (0,1/2)$,
$N_n^{-1}\log N_n=o(1/\log n)$ if $\ind_X=1/2$ and
$N_n^{-1}=o(1/\log n)$ if $\ind_X\in (1/2,1)$.
Recently Shao \cite{QMS} proved that (\ref{glad}) does hold
for Gaussian stochastic processes $X$ with stationary increments
such that $\sigma_X$ is either convex or concave, and $N_n$
satisfies conditions similar to the ones of Gin\'e and Klein \cite{GK}.
These conditions on $\sigma_X$ in Shao's result are quite general
and they do not seem to imply the existence of the Orey index $\ind_X$.

\begin{thm}\label{general}
Let $X=\{X(t)\colon\,t\geq 0\}$ be a mean zero Gaussian stochastic 
process with stationary increments and continuous in quadratic mean.
Suppose that for some $\ind\in (0,1)$, the Orey index $\ind_X=\ind$ and 
there exists the limit
\beq\label{0index}
\lim_{m\to\infty}\frac{T/N_m}{\sigma_X(T/N_m)^{1/\ind}}
\,s_{1/\ind}(X;\lambda (T/N_m)) = c\qquad\mbox{a.s.}
\eeq
for some sequence of integers $\{N_m\colon\,m\geq 1\}$ 
such that $N_m\to\infty$ as $m\to\infty$ and for some positive constants
$c$, $T$.
Then $X$ belongs to the Gladyshev class ${\cal G}^{\lambda}$ and
$G_{\lambda}(X)=\ind$ almost surely.
\end{thm}

\begin{proof}
For each integer $m\geq 1$, let $h_m:=T/N_m$. 
First we prove that
\beq\label{1index}
\liminf_{m\to\infty}\frac{\log\sqrt{s_2(X;\lambda_m)/N_m}}
{\log (1/N_m)}=
\frac{1}{2}+\liminf_{m\to\infty}\frac{\log s_2(X;\lambda_m)}{2\log h_m}
\leq\ind\qquad\mbox{a.s.}.
\eeq
For strictly positive functions $\phi$ and $\psi$ on $(0,1]$ such that
$\lim_{h\downarrow 0}\phi (h)=\lim_{h\downarrow 0} \psi (h)=0$ the relation
\beq\label{2index}
\liminf_{h\downarrow 0}\frac{\log\,\phi(h)}{\log\,\psi(h)}
=\sup\big\{\alpha>0\colon\,\mbox{$\psi(h)^{\alpha}/\phi(h)\to +\infty$
as $h\downarrow 0$}\big\},
\eeq
holds (see e.g.\ Annex A.4 in \cite{CT}).
Therefore relation (\ref{1index}) holds if and only if
\beq\label{3index}
\sup\Big\{\alpha >0\colon\,\mbox{$h_m^{\alpha-1}s_2(X;\lambda_m)
\to +\infty$
as $m\to\infty$}\Big\}\geq 2-2\ind\qquad\mbox{a.s.},
\eeq
provided $s_2(X;\lambda_m)/h_m\to +\infty$ almost surely.
Let $0<\ind''<\ind_X=\ind<\ind'\leq 1$.
Since almost all sample functions of $X$ satisfy H\"older's condition
of order $\ind''$ (see e.g.\ Section 9.4 in \cite{CL}), 
there is a positive random variable $K$ such that
\beq\label{4index}
\big|X(ih_m)-X((i-1)h_m)\big |\leq Kh_m^{\ind''}
\eeq
for all $i=1,\dots,N_m$ and $m\geq 1$.
Denoting $\epsilon :=\ind''(1/\ind-1/\ind')>0$, we then have
$$
s_{1/\ind'}(X;\lambda_m)\geq K^{-\epsilon/\ind''}\left [
\frac{\sigma(h_m)}{h_m^{\ind(1+\epsilon)}}\right ]^{1/\ind}
\frac{h_m}{\sigma(h_m)^{1/\ind}}\,s_{1/\ind}(X;\lambda_m).
$$
Since $\ind (1+\epsilon )>\ind=\ind_{\ast}$, the quantity 
in square brackets above tends to infinity as $m\to\infty$.
By assumption (\ref{0index}), it then follows that, for each 
$\ind'\in (\ind,1]$, 
$$
\lim_{n\to\infty}s_{1/\ind'}(X;\lambda_m)=+\infty\qquad\mbox{a.s.}.
$$
To prove (\ref{3index}),  first suppose that $1/2\leq\ind<1$.
For any $\ind'\in (\ind,1]$, by H\"older's inequality, we have
that
$$
s_{1/\ind'}(X;\lambda_m)^{2\ind'}\leq N_m^{2\ind'-1}s_2(X;\lambda_m)^2
$$
for each $m\geq 1$.
Thus for each $\ind'\in (\ind,1]$, 
\beq\label{5index}
\lim_{m\to\infty}h_m^{1-2\ind'}s_2(X;\lambda_m)=+\infty
\qquad\mbox{a.s.}.
\eeq
This yields (\ref{3index}) in the case $1/2\leq\ind<1$.
Now suppose $0<\ind<1/2$.
For $\epsilon\in (0,1/2-\ind)$, 
let $\ind':=\ind+\epsilon$ and $\ind'':=\ind-\epsilon$.
Since almost all sample functions of $X$ satisfy H\"older's
condition (\ref{4index}), we have
$$
s_{1/\ind'}(X;\lambda_m)\leq K^{1/\ind'-2}h_m^{\ind''(1/\ind'-2)}
s_2(X;\lambda_m)
$$
for each $m\geq 1$.
Since $\ind''(1/\ind'-2)=1-2\ind-2\epsilon [1/(\ind+\epsilon)-1]
<1-2\ind-\epsilon$, it follows that for an arbitrarily small
$\epsilon >0$,
\beq\label{6index}
\lim_{m\to\infty}h_m^{1-2\ind-\epsilon}s_2(X;\lambda_m)=+\infty
\qquad\mbox{a.s.}.
\eeq
This yields (\ref{3index}) in the case $0<\ind< 1/2$.
In both cases, by (\ref{5index}) and (\ref{6index}), we also have
$s_2(X;\lambda_m)/h_m\to\infty$ almost surely.
Therefore (\ref{3index}) implies (\ref{1index}).

To prove that $X\in {\cal G}^{\lambda}$ and $G_{\lambda}(X)=\ind$
almost surely, by (\ref{1index}), it is enough to show that
\beq\label{8index}
\limsup_{m\to\infty}\frac{\log\sqrt{s_2(X;\lambda_m)/N_m}}
{\log (1/N_m)}=
\frac{1}{2}+\limsup_{m\to\infty}\frac{\log s_2(X;\lambda_m)}
{2\log h_m}\geq\ind\qquad\mbox{a.s.}.
\eeq
Under the same conditions as (\ref{2index}) we have the relation:
\beq\label{9index}
\limsup_{h\downarrow 0}\frac{\log\,\phi(h)}{\log\,\psi(h)}
=\inf\big\{\alpha>0\colon\,\mbox{$\psi(h)^{\alpha}/\phi(h)\to 0$
 as $h\downarrow 0$}\big\}.
\eeq
Since $s_2(X;\lambda_m)/h_m\to\infty$ as $m\to\infty$ almost surely,
(\ref{8index}) holds if and only if
\beq\label{10index}
\inf\big\{\alpha>0\colon\,\mbox{$h_m^{\alpha-1}s_2(X;\lambda_m)\to 0$
as $m\to\infty$}\big\}\leq 2-2\ind\qquad\mbox{a.s.}.
\eeq
Let $0<\ind''<\ind$.
Since almost all sample functions of $X$ satisfy H\"older's
condition (\ref{4index}), we have
$$
s_{1/\ind''}(X,\lambda_m)\leq K^{\epsilon/\ind''}\left [
\frac{\sigma (h_m)}{h_m^{\ind (1-\epsilon)}}\right ]^{1/\ind}
\frac{h_m}{\sigma (h_m)^{1/\ind}}s_{1/\ind}(X;\lambda_m)
$$ 
for each $m\geq 1$, where $\epsilon :=1-\ind''/\ind>0$.
Since $\ind (1-\epsilon)<\ind=\ind^{\ast}$, the quantity in
square brackets tends to zero as $m\to\infty$.
By assumption (\ref{0index}), it then follows that
for each $\ind''\in (0,\ind)$,
$$
\lim_{m\to\infty}s_{1/\ind''}(X;\lambda_m)=0\qquad\mbox{a.s.}.
$$ 
To prove (\ref{10index}), first suppose that $0<\ind\leq 1/2$.
For any $\ind''\in (0,\ind)$, by H\"older's inequality, we have that
$$
s_2(X;\lambda_m)\leq N_m^{1-2\ind''}s_{1/\ind''}(X;\lambda_m)^{2\ind''}
$$
for each $m\geq 1$.
Thus for each $\ind''\in (0,\ind)$,
\beq\label{11index}
\lim_{m\to\infty}h_m^{1-2\ind''}s_2(X;\lambda_m)=0\qquad\mbox{a.s.}.
\eeq
This yields (\ref{10index}) in the case $0<\ind\leq 1/2$.
Now suppose that $1/2<\ind<1$.
Let $\ind''\in (1/2,\ind)$.
Since almost all sample functions of $X$ satisfy H\"older's
condition (\ref{4index}), we have
$$
s_2(X;\lambda_m)\leq K^{2-1/\ind''}h_m^{2\ind''-1}
s_{1/\ind''}(X;\lambda_m)
$$
for each $m\geq 1$.
Thus (\ref{11index}) holds for each $\ind''\in (1/2,\ind)$.
This yields (\ref{10index}) in the case $1/2<\ind<1$.
Since (\ref{10index}) implies  (\ref{8index}),
the proof of Theorem \ref{general} is complete.
\qed\end{proof}

\begin{cor}\label{particular}
Let $X=\{X(t)\colon\,t\geq 0\}$ be a mean zero Gaussian stochastic 
process with stationary increments and continuous in quadratic mean.
Suppose that for some $\ind\in (0,1)$,
\begin{enumerate}
\item[$(a)$] the incremental variance $\sigma_X$ of $X$
satisfies the condition $\sigma_X(h)=C_1h^{\ind}(1+\epsilon (h))$
where $C_1$ is a positive constant and $\epsilon (h)\to 0$
as $h\downarrow 0${\rm ;}
\item[$(b)$] there exists the limit
\beq\label{12index}
\lim_{m\to\infty}s_{1/\ind}(X;\lambda (T/N_m))=C_2\qquad
\mbox{a.s.}
\eeq
for some sequence of integers $\{N_m\colon\,m\geq 1\}$ such that $N_m\to\infty$
as $m\to\infty$ and some positive constants $C_2$, $T$.
\end{enumerate}
Then the Orey index $\ind_X=\ind$,
$X$ belongs to the Gladyshev class ${\cal G}^{\lambda}$
and $G_{\lambda}(X)=\ind$ almost surely, where 
$\lambda=\{\lambda (T/N_m)\colon\,m\geq 1\}$.
\end{cor}

Conditions for (\ref{12index}) to hold are given by
Marcus and Rosen \cite{MR} in the case $0<\ind\leq 1/2$,
and by Shao \cite{QMS} in the case $0<\ind <1$.
Let $X$ satisfy the assumptions of Corollary \ref{particular}
except possibly (\ref{12index}) and let $0<T<\infty$.
Suppose that $\sigma_X^2$ is nondecreasing and concave
on $[0,T]$ if $0<\ind\leq 1/2$, and that $\sigma_X^2$ is
nondecreasing and convex on $[0,T+\epsilon]$ for some
$\epsilon >0$ if $1/2<\ind <1$.
Also, suppose that $1/N_m=o(1/(\log m)^a)$ as $m\to\infty$,
where $a:=1/\{2[\ind\wedge (1-\ind)]\}$.
Then by Corollaries 1.1 and 1.2 of Shao \cite{QMS}, (\ref{12index})
holds with $C_2=TC_1^{1/\ind}E|\eta|^{1/\ind}$, where
$\eta$ is a standard normal random variable 
and $C_1$ is the constant from $(a)$.

Further conditions on existence of the Orey index
will be given in terms of a spectral density
of a Gaussian stochastic process with stationary increments.
A nondecreasing function $F$ on $\RR\setminus\{0\}$ is called
a L\'evy-Khinchin function if
$$
\int_{\RR\setminus\{0\}}\min \{1,u^2\}\,dF(u)<\infty.
$$
Let $X$ be a real-valued mean zero Gaussian stochastic process with
stationary increments, continuous in probability and $X(0)=0$
almost surely.
Then there is a unique Le\'vy-Khinchin function $F_X$ on 
$\RR\setminus\{0\}$ such that
$$
E\big [X(t)X(s)\big ]=\int_{\RR\setminus\{0\}}\big (e^{it\lambda}-1
\big )\big (e^{-is\lambda}-1\big )\,dF_X(\lambda).
$$
If $F_X$ is absolutely continuous, the derivative $f_X=F_X'$
is called the spectral density of $X$.

The spectral density of a fractional Brownian motion $B_H$ with the Hurst 
exponent $0<H<1$ is given by the power function 
$f_H(\lambda)=c_H|\lambda|^{-1-2H}$, 
where $c_H=(2\pi)^{-1}\sin(\pi H)\Gamma (1+2H)$.
This follows from the relation:
$$
\sigma_{B_H}(u)^2=\int_{-\infty}^{\infty}\big |e^{iu\lambda}-1\big |^2
f_H(\lambda)\,d\lambda=8c_Hu^{2H}\int_0^{\infty}\Big (\sin
\frac{\lambda}{2}\Big )^2\frac{d\lambda}{\lambda^{1+2H}}=u^{2H}.
$$
So that the Orey index $\ind_{B_H}$ exists and equals $H$.

A spectral density which differ from the spectral density $f_H$
by a slowly varying function gives rise to a Gaussian stochastic process
with the same Orey index $H$.
Consider a quasi-monotone slowly varying (near infinity) function $\ell$ defined
in \cite[Section 2.7]{BGT} and let $f_{\ind,\ell}(\lambda):=\ell (|\lambda|)
|\lambda|^{-1-2\ind}$ for some $0<\ind<1$.
By Abelian Theorem 4.1.5 from \cite{BGT}, we have
\beq\label{3spectral}
\int_{-\infty}^{\infty}\big |e^{iu\lambda}-1\big |^2f_{\ind,\ell}
(\lambda)\,d\lambda
=8u^{2\ind}\int_0^{\infty}\ell\Big (\frac{\lambda}{u}\Big )\Big (
\sin\frac{\lambda}{2}\Big )^2\frac{d\lambda}{\lambda^{1+2\ind}}
\sim\ell\Big (\frac{1}{u}\Big )u^{2\ind}c_{\ind}^{-1}
\eeq
as $u\downarrow 0$, where $\phi(u)\sim\psi(u)$ as $u\downarrow 0$
means that $\lim_{u\downarrow 0}\phi(u)/\psi(u)=1$.
Let $X$ be a Gaussian stochastic process with spectral density
$f_{\ind,\ell}$.
By elementary properties of slowly varying functions (e.g.\
Proposition 1.3.6 in \cite{BGT}), it follows that the Orey index $\ind_X$
exists and equals $\ind$ for any quasi-monotone slowly varying
function $\ell$.

The following statement contains a part of Theorem 3 of
Gladyshev \cite{EGG}.

\begin{thm}\label{stationary}
Let $X$ be a mean zero Gaussian stochastic process with stationary
increments having the spectral density $f_X$ such that, for some
$0<\ind<1$,
$$
f_X(\lambda)=|\lambda|^{-1-2\ind}\big [\ell(|\lambda|)+b(\lambda)\big ],
$$
where $\ell$ is a quasi-monotone slowly varying function,
$b$ is zero on a neighbourhood of zero, bounded elsewhere
and $b(\lambda)=O(\lambda^{-2})$ as $\lambda \to\infty$.
Then the Orey index $\ind_X=\ind$.
If in addition $\ell (|\lambda|)=c$ for some $c>0$, then
$$
\frac{1}{2}-\lim_{m\to\infty}\frac{\log s_2(X,\lambda (2^{-m}))}{2m}=\ind
\qquad\mbox{a.s.}.
$$
\end{thm}

\begin{proof}
Suppose that $b$ is zero on $[-1,1]$.
Then for some constant $C$ depending on $b$, and for each $u>0$, we  have
$$
\int_{-\infty}^{\infty}\Big (\sin\frac{u\lambda}{2}\Big )^2
\frac{b(\lambda)\,d\lambda}{|\lambda|^{1+2\ind}}
\leq Cu^{2+2\ind}\int_u^{\infty}\Big (\sin\frac{\lambda}{2}
\Big )^2\frac{d\lambda}{\lambda^{3+2\ind}}
$$
\beq\label{6spectral}
\leq Cu^{2\ind}\int_u^{u^{1/(2+2\ind)}}\Big (\sin\frac{\lambda}{2}\Big )^2
\frac{d\lambda}{\lambda^{1+2\gamma}}+Cu^{1+2\ind}/(2+2\ind).
\eeq
Since the right side of (\ref{6spectral}) is of order $o(u^{2\ind})$ as 
$u\downarrow 0$, by (\ref{3spectral}), it follows that 
$\sigma_X(u)\sim\sqrt{c_{\ind}^{-1}\ell(1/u)}u^{\ind}$ as $u\downarrow 0$.
Therefore the Orey index $\ind_X$ exists and equals $\ind$,
proving the first part of the theorem.
The second part follows from Theorem 3 of Gladyshev \cite{EGG}
with $\alpha=1+2\gamma$.
\qed\end{proof}

\appendix
\chapter{Convergence of directed functions}\label{convergence}
\setcounter{thm}{0}


In this paper, the Moore-Smith limit theory as developed by McShane 
\cite{EJM} (see also \cite{MB}) is used to define
the convergence with respect to refinements of partitions.
Namely, McShane's notion of the direction is used to define several Stieltjes type 
integrals, the unconditional convergence of sums and products, and 
the local $p$-variation.
For readers convenience we state the main facts about directed
functions in this appendix.
A connection with the better known theories of limits based on nets 
and filters is discussed in the last section of \cite{EJM}.


\paragraph{Directed functions.}
We start with the notion of a directed function and its convergence
which are basic to McShane's limit theory 
(see p.\ 10 in \cite{EJM}, or p.\ 33 in \cite{MB}).

\begin{defn}\label{direction}
{\rm Let $D$ be a set.
\begin{enumerate}
\item[(a)] A nonempty family $\ugotik$ of nonempty subsets of $D$
is called a {\em direction in} $D$ if it is
directed downward by inclusion, that is, if $A_1$ and $A_2$ belong to 
$\ugotik$, then there is a set $A_3$ in $\ugotik$ such that 
$A_3\subset A_1\cap A_2$.
\item[(b)] Let $\ugotik$ be a direction in $D$.
If $P(x)$ is a property valid for $x\in D$ then we say that
``{\em ultimately} $P(x)$'' (or ``$P(x)$ is ultimately true'')
provided that there is a set $A$ in $\ugotik$ such that
$P(x)$ holds for each $x\in A$.
\item[(c)]
Let $\ugotik$ be a direction in $D$.
If a real-valued function $f$ is ultimately defined on $D$,
that is ultimately ``$f(x)$, $x\in D$, is a real number'',
then the ordered pair $(f,\ugotik)$ is called a {\em directed function}.
A directed function $(f,\ugotik)$ converges to a real number $r$,
in symbols
$$
\lim_{\ugotik}f=\lim_{x,\ugotik}f(x)=r,
$$
if for each neighborhood $V$ of $r$, $f$ is ultimately in $V$. 
\end{enumerate}}
\end{defn}

The convergence of directed functions satisfies all natural properties
of limits including the uniqueness of the limit 
(p.\ 12 in \cite{EJM}).
Next is the Cauchy test for the convergence of directed functions 
(see p.\ 37 in \cite{MB} for the proof).

\begin{thm}\label{Cauchy}
Let $f$ be a real-valued function defined on a set $D$, and let
$\ugotik$ be a direction on $D$.
The directed function $(f,\ugotik )$ has a limit if and only if
for each $\epsilon >0$ there exists $A\in\ugotik$ such that
$|f(x)-f(y)|<\epsilon$ for each $x, y\in A$.
\end{thm}

\paragraph{Order convergence.}
Let $(f,\ugotik)$ be a real-valued directed function.
Define the {\em upper} and {\em lower limits} of $(f,\ugotik)$,
respectively, by
\beq\label{1orderconv}
\limsup_{\ugotik}f:=\mathop{\inf}_{A\in\ugotik}\sup_{x\in A}f(x)
\quad\mbox{and}\quad
\liminf_{\ugotik}f:=\sup_{A\in\ugotik}\inf_{x\in A}f(x).
\eeq

\begin{lem}
For a directed function $(f,\ugotik)$, $\liminf_{\ugotik}f
\leq \limsup_{\ugotik}f$.
\end{lem}

A real-valued directed function $(f,\ugotik)$ is {\em order convergent to 
$b\in\RR$} if
$$
\liminf_{\ugotik}f=\limsup_{\ugotik}f=b
$$
The following theorem is a special case of Theorem 11.4
of McShane and Botts \cite[p.\ 55]{MB} proved for functions
with values in the extended set of real numbers.

\begin{thm}\label{orderconv}
Let $(f,\ugotik)$ be a real-valued directed function and let
$b\in\RR$.
Then $(f,\ugotik)$ is order convergent to $b$ if and only if
it converges to $b$.
\end{thm}

We use this theorem in Section \ref{Wiener} to characterize the Wiener 
class of functions having the local $p$-variation, 
and at the end of Section \ref{extended}
to relate several definitions of extended 
Riemann-Stieltjes integrals.

\paragraph{Unconditionally convergent sums.}
Suppose that $f$ is a real-valued function defined on an interval $J$.
Let ${\cal F}(J)$ be the family of all {\em finite} sets of points
of $J$.
For each $\sigma\in{\cal F}(J)$, let 
$$
S(f;\sigma):=\sum_{x\in\sigma}f(x).
$$
The {\em sum-function} $S(f)=\{S(f;\sigma)\colon\,\sigma\in {\cal F}(J)\}$ 
is then defined on ${\cal F}(J)$.
For each finite subset $\sigma$ of $J$, let $A(\sigma)$ be the set of all
finite subsets $\sigma'$ of $J$ which contain $\sigma$.
Then $A(\sigma)$ is a nonempty subset in the domain of the sum-function
$S$.
Let $\fgotik$ be the family of all sets $A(\sigma)$ for 
$\sigma\in{\cal F}(J)$, that is
\beq\label{Fgotik}
\fgotik=\fgotik (J):=\big\{A(\sigma)\colon\,\sigma\in
{\cal F}(J)\big\}.
\eeq
Let $A(\sigma_1)$ and $A(\sigma_2)$ be any two members of $\fgotik$.
Then $A(\sigma_1\cup\sigma_2)$ belongs to $\fgotik$, and every finite set 
$\sigma$ of $J$ which contains $\sigma_1\cup\sigma_2$ also contains both
$\sigma_1$ and $\sigma_2$, that is, $\sigma$ belongs to $A(\sigma_1)\cap
A(\sigma_2)$.
Hence $\fgotik$ is a direction in ${\cal F}(J)$.

\begin{defn}\label{unconditional}
{\rm Let $f$ be a real-valued function defined on an interval $J$.
We say that the {\em non-ordered sum $\sum_Jf$ converges unconditionally}
if the directed function $(S(f),\fgotik )$ has a limit and let
$$
\sum_Jf:=\lim_{\sigma,\fgotik}S(f;\sigma).
$$ 
Also we say that the {\em non-ordered sum $\sum_Jf$ converges absolutely} 
if there exists $K$ such that
$\sum_{x\in\sigma}|f (x)|\leq K$ for each $\sigma\in {\cal F}(J)$.}
\end{defn}

\begin{thm}\label{boundedness}
For a real-valued function $f$ on an interval $J$, the non-ordered sum 
$\sum_Jf$ converges unconditionally if and only if it converges absolutely.
\end{thm}

{\proof}
First, suppose that $f$ is non-negative on $J$.
In this case, the directed function $(S(f),\fgotik)$ converges if and 
only if the set $\{S(f;\sigma)\colon\,\sigma\in {\cal F}(J)\}$ has an upper 
bound, and if both statements hold then the limit is the least upper bound.
Now, suppose that $f$ may have arbitrary real values.
Let $f^{+}$ and $f^{-}$ be positive and negative parts of $f$,
respectively, so that $f=f^{+}-f^{-}$ and $|f|=f^{+}+f^{-}$.  
If $\{S(|f|;\sigma)\colon\,\sigma\in {\cal F}(J)\}$ is bounded
then also the sets  $\{S(f^{+};\sigma)\colon\,\sigma\in {\cal F}(J)\}$
and  $\{S(f^{-};\sigma)\colon\,\sigma\in {\cal F}(J)\}$ are bounded.
By the first part, the directed functions $(S(f^{+}),\fgotik)$
and $(S(f^{-}),\fgotik)$ converge.
Hence by linearity of the limit, the sum $\sum_Jf$ converges unconditionally.
For the converse implication we use the Cauchy test (Theorem \ref{Cauchy}
above) to show that if the sum $\sum_Jf$ converges unconditionally
then so do the sums $\sum_Jf^{+}$ and $\sum_Jf^{-}$.
Indeed, for $k=1,2$ and for any finite sets $\sigma_0\subset\sigma_k$,
let 
$$
\sigma_0\subset\sigma_k':=\sigma_0\cup\{x\in\sigma_k\setminus\sigma_0
\colon\,f(x)>0\}\subset\sigma_k.
$$
Then $S(f^{+};\sigma_1)-S(f^{+};\sigma_2)=S(f;\sigma_1')
-S(f;\sigma_2')$.
A similar relation holds between $f$ and $f^{-}$.
This implies the claim.
By the above argument it follows that the set $\{S(|f|;\sigma)\colon\,
\sigma\in {\cal F}(J)\}$ is bounded.
The proof of Theorem \ref{boundedness} is complete.
\qed

For a function $f$ with real values,
by the preceding theorem, the sum $\sum_Jf$ converges unconditionally
if and only if it converges absolutely.
This is still true for functions with values in
any finite-dimensional Banach space.
In an infinite-dimensional Banach space, absolute convergence implies
unconditional convergence but unconditional convergence never
implies absolute convergence (Dvoretzky and Rogers \cite{DR}).

\begin{thm}\label{sum-add}
Let $J$, $J_1$, $J_2$ be intervals such that $J=J_1\cup J_2$
and $J_1\cap J_2=\emptyset$.
Then the non-ordered sum $\sum_Jf$ converges unconditionally if and only if
both $\sum_{J_1}f$ and $\sum_{J_2}f$ converge unconditionally,
and then
\beq\label{sumadditive}
\sum_Jf=\sum_{J_1}f+\sum_{J_2}f.
\eeq
\end{thm}

{\proof}
Suppose that the non-ordered sum $\sum_Jf$ converges unconditionally.
By the preceding theorem, the set $\{S(|f|;\sigma)\colon\,
\sigma\in {\cal F}(J)\}$ is bounded.
Therefore the same sets with ${\cal F}(J)$ replaced by
${\cal F}(J_1)$ or ${\cal F}(J_2)$ also are bounded,
and hence the sums on the right side of (\ref{sumadditive})
converge unconditionally.
The converse implication follows similarly.
The conclusion of the theorem now follows by linearity of the limit
from equality (\ref{sumadditive}) with
$J$, $J_1$, $J_2$ replaced by finite sets $\sigma\subset J$,
$\sigma_1\subset J_1$, $\sigma_2\subset J_2$ such that
$\sigma=\sigma_1\cup\sigma_2$.
\qed

\paragraph{Unconditionally convergent products.}
Here we extend the concept of unconditional convergence to 
non-ordered products.
Again let $f$ be a real-valued function defined on an interval $J$,
and let ${\cal F}(J)$ be the family of all finite sets of points
of $J$.
For each $\sigma\in{\cal F}(J)$, let
$$
P(f;\sigma):=\prod_{x\in\sigma}f(x).
$$
The {\em product-function} $P(f)=\{P(f;\sigma)\colon\,\sigma\in {\cal F}(J)\}$ 
is then defined on ${\cal F}(J)$.
Let $\fgotik$ be the direction on ${\cal F}(J)$ defined by
(\ref{Fgotik}).

\begin{defn}\label{conv-prod}
{\rm Let $f$ be a real-valued function defined on an interval $J$.
We say that the {\em non-ordered product $\prod_Jf$ converges 
unconditionally} if the directed function $(P(f),\fgotik )$ has a 
{\em non-zero} limit, and then let
$$
\prod_Jf:=\lim_{\sigma,\fgotik}P(f;\sigma).
$$ }
\end{defn}

If $f$ is positive everywhere on $J$ then $\log f(x)$ is defined
for each $x\in J$. 
Also, if in this case the non-ordered product $\prod_Jf$ converges
then the directed function $(S(\log f),\fgotik)$ has a limit 
$\sum_J\log f$ such that
\beq\label{1log-prod}
\log\Big ( \prod_Jf\Big )=\sum_J\log f.
\eeq
We show below that the unconditional convergence of a non-ordered
sum $\sum_Jf$ is equivalent to the unconditional convergence of
a non-ordered product $\prod_Jf$.
Thus the unconditional convergence of non-ordered products
extends the absolute convergence of products of sequences.   

\begin{prop}\label{log-prod}
For a real-valued function $f$ on $J$, if the non-ordered product 
$\prod_Jf$ converges unconditionally then for each $\epsilon >0$, 
ultimately $|f-1|<\epsilon$.
\end{prop}

\begin{proof}
Let $P=\prod_Jf$, and let $\epsilon\in (0,|P|)$.
Then there is $\sigma_0\in {\cal F}(J)$ such that
$|P(f;\sigma)-P|<\epsilon$ for each $\sigma\supset\sigma_0$.
For $x\in J\setminus\sigma_0$, we then have
$$
|f(x)-1|=\frac{|P(f;\sigma_0\cup\{x\})-P(f;\sigma_0)|}{|P(f;
\sigma_0)|}<\frac{2\epsilon}{|P|-\epsilon}.
$$
Since $\epsilon$ is arbitrary, the conclusion follows.
\qed\end{proof}

It follows from the preceding proposition that $f$ is equal $1$ 
everywhere on $J$ except on a countable set, 
and $f$ is positive on $J$ except on a finite set.
In particular, the sum-function $S(\log f)$ is ultimately
defined on ${\cal F}(J)$ with respect to the direction $\fgotik$,
whenever the non-ordered product $\prod_Jf$ converges
unconditionally.
Then the same arguments used to show (\ref{1log-prod}) for a positive
$f$, yield the following statement:

\begin{thm}\label{2log-prod}
Let $f$ be a real-valued function defined on an interval $J$.
The non-ordered product $\prod_Jf$ converges unconditionally
if and only if the non-ordered sum $\sum_J\log f$ converges 
unconditionally.
\end{thm}

The following criterion for the unconditional convergence of 
non-ordered products is used in Lemma \ref{function}.\ref{gamma}.

\begin{thm}
Let $\unc$ be a real-valued function on $J$ such that $\unc (x)\not =-1$
for each $x\in J$.
Then the non-ordered product $\prod_J(1+\unc)$ converges unconditionally
if and only if the non-ordered sum $\sum_J\unc$ converges unconditionally.
\end{thm}

\begin{proof}
Suppose that the non-ordered product $\prod_J(1+\unc)$ converges 
unconditionally.
By Proposition  \ref{log-prod}, for each $\epsilon >0$,
ultimately $|\unc |<\epsilon$.
Therefore ultimately $|1-\unc^{-1}\log (1+\unc)|<\epsilon$.
By Theorem  \ref{2log-prod}, it then follows that the
directed function $(S(\unc),\fgotik)$ has a limit
proving the unconditional convergence of the non-ordered sum
$\sum_J\unc$.
The converse implication follows similarly  because
a limit of the directed function $(P(1+\unc),\fgotik)$ is non-zero.
\qed\end{proof}

\addcontentsline{toc}{chapter}{\protect\numberline{}{Bibliography}}

\end{document}